\documentclass[12pt, final, a4paper]{template}

\DeclareRobustCommand{\subtitle}[2]{\vspace{.5\baselineskip}\\#1: #2}

\begin{document}

\setcounter{page}{1}
\normalem

\date{}
\title{Gluing invariants of Donaldson--Thomas type\subtitle{Part I}{the Darboux stack}}

\authorbenjamin
\authorjulian
\authormarco

\begin{abstract}
  Let $X$ be a $(-1)$\nobreakdashes-shifted symplectic derived Deligne--Mumford stack.
  In this paper we introduce the Darboux stack of $X$, parametrizing local presentations of $X$ as a derived critical locus of a function $f$ on a smooth formal scheme $U$.
  Local invariants such as the Milnor number $\mu_f$, the perverse sheaf of vanishing cycles $\PJoyce_{U,f}$ and the category of matrix factorizations $\MF(U,f)$ are naturally defined on the Darboux stack, without ambiguity.
  The stack of non-degenerate flat quadratic bundles acts on the Darboux stack and our main theorem is the contractibility of the quotient stack when taking a further homotopy quotient identifying isotopic automorphisms.
  As a corollary we recover the gluing results for vanishing cycles by Brav--Bussi--Dupont--Joyce--Szendr\H oi.

  In a sequel to this paper we apply this general mechanism to glue the locally defined categories of matrix factorizations $\MF(U,f)$ under the prescription of 
  additional orientation data, thus answering positively conjectures by Kontsevich--Soibelman and Toda in motivic Donaldson--Thomas theory.
\end{abstract}

\maketitle
\renewcommand{\subtitle}[2]{ -- #1}

\tableofcontents

\section{Introduction}
\label{Introduction}

\subsection{Main result}
\label{sectionmainresult}
Throughout this article, we will be working over the field of complex numbers $\basefield$.
Motivated by applications to Donaldson--Thomas invariants (see \cref{motivations}), this paper studies the gluing of local invariants attached to a $(-1)$\nobreakdashes-shifted symplectic derived Deligne--Mumford stack $X$ \cite{1111.3209}. The simplest example of such $X$ is the derived critical locus $\dCrit(U,f)$ of a function $f$ on a smooth scheme $U$ (see \cref{symmetryderivedcriticalocusexplicit}) and a Darboux lemma proved in \cite{MR3904157} shows that in fact, Zariski locally, any $(-1)$\nobreakdashes-shifted symplectic $X$ is symplectically equivalent to some $\dCrit(U,f)$. This paper starts with the observation that such local presentations are not unique:

\begin{example}
  \label{ambiguitypresentationdcrit}
  Given the two pairs  $(\affineline{\basefield}, x^2)$ and $(\affine^2_\basefield, x^2+y^2)$, we have a canonical symplectic identification $\dCrit(\affineline{\basefield}, x^2)=\dCrit(\affine^2_\basefield, x^2+y^2)=\Spec(\basefield)$, where $\Spec(\basefield)$ is equipped with the zero $(-1)$\nobreakdashes-shifted symplectic form.  More generally, we can always reproduce this phenomena by adding a quadratic variable. For instance, $\dCrit(\affineline{\basefield}, x^3)=\dCrit(\affine^2_{\basefield}, x^3+y^2)$.
\end{example}
Our goal is to investigate the moduli of such local Darboux parametrizations of a $(-1)$\nobreakdashes-shifted symplectic derived Deligne--Mumford stack $X$ in a way that resolves all ambiguities.

To start with, part of the moduli data should concern an \emph{exact structure} for the symplectic form.
Indeed, if $X=\dCrit(U,f)$ then the symplectic structure  has a canonical {exact structure} (see \cref{symmetryderivedcriticalocusexplicit}).
More generally, a result of Deligne implies that every $(-1)$\nobreakdashes-shifted symplectic derived stack admits an {exact structure} -- see \cite[Cor.\,5.3]{MR3285853} and our refinement in \cref{canonicalexactstructure}.
For this reason, we will consider $(-1)$\nobreakdashes-shifted symplectic derived Deligne--Mumford stack $X$ \emph{together with a choice of exact structure}.

The naive attempt to define a moduli functor would be to assign to every affine Zariski open $S=\Spec(A)$ in $X$, the  $\infty$\nobreakdashes-groupoid classifying \emph{exact symplectic identifications} of $S$\footnote{equipped with the exact symplectic structure restricted from $X$.} with a derived critical locus $\dCrit(U,f)$ for $U$ a smooth scheme and $f$ a function on $U$. Unfortunately, this naive attempt is not even functorial: if $S'\subseteq S$ is an open Zariski subset of $S$, there is a priori no canonical way to produce a Zariski open $U'$ of $U$ with a function $f'$ whose derived critical locus is $S'$. In order to overcome this difficulty we observe  that $\dCrit(U,f)$ with its exact $(-1)$\nobreakdashes-symplectic structure only depends on the smooth formal scheme given by the formal completion $\formalU$ of the derived critical locus $\dCrit(U,f)$ inside $U$ together with the Taylor expansion of $f$ (see \cref{dCritonlydependsformal}). Therefore, we can replace \emph{algebraic} smooth schemes $U$, by \emph{formal} smooth schemes $\formalU$ in our Darboux local models. The fact that étale covers can be uniquely extended along infinitesimal thickenings solves the functoriality problem.

Following this discussion, we can introduce the central object of this work, namely, the Darboux moduli \ifunctor (where $\lambda$ denotes the $(-1)$\nobreakdashes-shifted exact symplectic form of $X$)
\[
  \Darbstack^{\lambda}_X \colon  X_\et^\op \to \inftygpd
\]
defined on the small étale site of $X$ and sending $S\to X$ to the $\infty$\nobreakdashes-groupoid $\Darbstack^{\lambda}_X(S)$ of \emph{exact symplectic} identifications $S\simeq \dCrit(\formalU,f)$ where $\formalU$ is a smooth formal scheme which is a nilpotent thickening of $S$. See \cref{definitiondarbouxstackasstacks}\footnote{As in many geometric situations (for instance, with the orientation sheaf on a manifold) $\Darbstack_X^\lambda$ might not admit global sections but the Darboux's lemma of \cite{MR3904157} guarantees the existence of local sections.}.

To address \cref{ambiguitypresentationdcrit} describing quadratic bundles as a source of ambiguities, we consider  the stack $\stackQuadnabla_X$ classifying étale locally trivial quadratic bundles with compatible flat connection $(Q,q,\nabla)$ (see \cref{sectionStackQuad}). Generalizing \cref{ambiguitypresentationdcrit}, the stack  $\stackQuadnabla_X$ acts on $\Darbstack^{\lambda}_X$ and we focus our attention on the (hypercomplete) quotient stack
\[
  \quot{\Darbstack_X^{\lambda}}{\stackQuadnabla_X.}
\]
Unfortunately, not all redundancies emerge from the action of quadratic bundles. The following example identifies a second source given by automorphisms:

\begin{example}
  \label{ambiguitymorphisms}
  Let $\formalU=\Spf(\basefield[\![x,y]\!])$ be the formal completion of $0$ in $\affine^2$ and $f(x,y) \coloneqq x^3+y^4$. Let $X=\dCrit(\formalU, f) = \Spec\left( {\basefield[\![x,y]\!]}/{(x^2,y^3)} \right)$ and consider the automorphism $\phi \colon \formalU\to \formalU$ obtained by sending $x\mapsto x+y^4$ and $y\mapsto y.h(x,y)$, where $h=\sqrt[4]{1-3x^2- 3xy^4- y^8}\in \basefield[\![x,y]\!]$. This automorphism is the identity on $X$ and in fact preserves the shifted symplectic structure. It provides a non-trivial section $[\phi]$ of the automorphism sheaf
  \[
    \homotopysheaf_1\left(\quot{\Darbstack^{\lambda}_X}{\stackQuadnabla_X}, (\formalU, f)\right)
  \]
  at the point given by $(\formalU, f)$. This shows that the quotient stack is not contractible.
\end{example}

The next remark explains the strategy to solve the non-contractibility of the quotient $\Darbstack^{\lambda}_X/\stackQuadnabla_X$:

\begin{observation}
  \label{isotopiessolveambiguity}
  In \cref{ambiguitymorphisms}, it is easy to exhibit a $1$\nobreakdashes-parameter family (isotopy) $\phi_t$ interpolating between $\phi$ and the identity map, given by $x \mapsto x+ty^4$ and $y \mapsto y \sqrt[4]{1-3tx^2- 3t^2xy^4- t^3y^8}$.
\end{observation}

\cref{isotopiessolveambiguity} suggests that to solve the ambiguities posed by \cref{ambiguitymorphisms}, we need, at the very least, to identify $\affineline{}$\nobreakdashes-isotopic families.
We achieve this by constructing further quotient stacks $\DA$  and  $\QA$, keeping the same objects but identifying $\affineline{}$\nobreakdashes-isotopic families of isomorphisms (see \cref{Liouvilleuptoisotopy}). Finally, we can state our main  result:

\begin{theoremnonumber}[Contractibility -- see \cref{contractibilitytheoremdetailedversion}]
  \label{theoremcontractibility}
  Let $X$ be a Deligne--Mumford derived stack equipped with a $(-1)$\nobreakdashes-shifted exact symplectic form $\lambda$. Then the canonical morphism to the final object of the small (hypercomplete) étale $\infty$\nobreakdashes-topos
  \[
   \quotA \to \constantsheaf \ast_X
  \]
  is an equivalence of (hypercomplete) $\infty$\nobreakdashes-stacks on $X_\et$. In other words, the quotient stack is contractible.
\end{theoremnonumber}

\begin{remark}
  To extend the above theorem to the case of derived Artin stacks requires some work.
  Local models are now given by Lagrangian correspondences to a derived critical locus.
  It is straightforward to see that the added data amounts to coisotropic structures on smooth covers.
  It is linear in nature and will hence cause no additional difficulties.
  We thus strongly believe that extending \cref{theoremcontractibility} to the Artin case will be a formal (albeit not short) argument.
  We will not explore this question in this paper.
\end{remark}

\subsection{Application: Gluing local invariants of singularities}\label{sectiongluinglocalinvariants}
We now discuss the applications of \cref{theoremcontractibility}. The Darboux lemma of \cite{MR3904157} allows us to attach local invariants to a $(-1)$\nobreakdashes-shifted symplectic derived Deligne--Mumford stack $X$:
\begin{example}\label{Milnornumber}
  To a smooth scheme $U$ with a function $f$ and a closed point $x$ in the classical critical locus of $f$, we can assign the Milnor number $\mu_{f,x}\in \integers$ given by
  \[(-1)^{\dim U}(1-\chi(F_x))\]
  where $F_x$ is the Milnor fiber at $x$. See \cite{milnor1968singular} and \cite[\S 1.2]{MR2600874}.
\end{example}

\begin{example}\label{perversesheaf}
  The theory of perverse sheaves \cite{MR751966, Dimca2004, MR3353002} allows us to categorify the construction of the Milnor number: for every smooth scheme $U$ with a function $f$ we can assign a perverse sheaf $P_{U,f}$ on the classical critical locus of $f$, computing the cohomology of vanishing cycles for the function $f$.
  Finally, we recover the Milnor number $\mu_f$ of \cref{Milnornumber} as the Euler characteristic of the stalk of the perverse sheaf $P_{U,f}$ at $x$
  \[
    \mu_{f,x}= \chi(P_{U,f,x}).
  \]
\end{example}

\begin{example}\label{exampleMF}
  In \cite{MR570778}, Eisenbud introduced the notion of matrix factorizations, which was then extensively studied, in particular by Orlov \cite{MR2101296, MR2735755}.
  They form a $2$\nobreakdashes-periodic (dg)-category $\MF(U,f)$ which is a categorical invariant associated to the singularities of function $f$ on a smooth scheme $U$.
  Thanks to several contributions \cite{MR2824483,MR3121870,1212.2859, MR3877165} this invariant provides a categorification of the cohomology of vanishing cycles of \cref{perversesheaf}, \ie $\HP(\MF(U,f))$ recovers the $2$\nobreakdashes-periodized vanishing cohomology.
\end{example}

Some invariants, such as the Milnor number of \cref{Milnornumber} are insensitive to the ambiguity posed in \cref{ambiguitypresentationdcrit} and therefore, can be globalized along $X$:

\begin{theorem}[Behrend \cite{MR2600874}]\label{behrendfunctiontheorem}
  Let $X$ be a proper $(-1)$\nobreakdashes-shifted symplectic derived Deligne--Mumford stack. Then the locally defined Milnor numbers $\mu_{f,x}$ are well-defined independently of the local presentations of $X$ as a derived critical locus and the assignment $x\mapsto \nu(x) \coloneqq  \mu_{f,x}$ defines a constructible function $\nu \colon X \to \integers$ such that
  \[
    \int_{[X]^\virtual} 1 = \chi(X, \nu)
  \]
  where in the r.h.s we have the Euler characteristic of $X$ weighted by $\nu$ and the l.h.s is the volume of the virtual fundamental class of \cite{MR1437495}.
\end{theorem}

This paper concerns precisely those invariants that are \emph{sensitive} to the ambiguity posed by \cref{ambiguitypresentationdcrit}:

\begin{observation}\label{ambiguityforsheaves}
  The ambiguity exhibited in \cref{ambiguitypresentationdcrit} extends to the sheaves $P_{U,f}$. For instance, the sheaves $P_{\affineline{\basefield}, x^2}$ and $P_{\affine^2_{\basefield}, x^2+y^2}$ are isomorphic but not canonically so: the identification requires the choice of an orientation of the Milnor fiber, which is a sphere.
\end{observation}

A key result by Brav-Bussi-Dupont-Joyce-Szendr\H oi  (BBDJS) establishes a condition under which the ambiguity in \cref{ambiguityforsheaves} can be solved and the globalization of the locally defined sheaves $P_{U,f}$, is possible:

\begin{theorem}[\cite{MR3353002}]\label{joycegluing}
  Let $X$ be a $(-1)$\nobreakdashes-shifted symplectic derived Deligne--Mumford stack. Assume that there exists a line bundle $L$ together with an equivalence $\phi \colon L \otimes L \simeq \det(\tangent_X)$  (called  \emph{orientation data}). Then:
  \begin{enumerate}
    \item The locally defined perverse sheaves of vanishing cycles $P_{U,f}$ glue to a globally defined perverse sheaf $\PJoyce\in \Perv(X)$.
    \item The Euler characteristic of the perverse sheaf $\PJoyce$ recovers Behrend's function of \cref{behrendfunctiontheorem}: $\chi(\PJoyce)=\nu$.
  \end{enumerate}
\end{theorem}

\begin{observations}{}
  \item \label{squarerootbundle}
  The obstruction to the existence of such orientation data is a class $\beta \in \cohomology^2(X, \quot{\integers}{2})$.
  \item The proof of \cref{joycegluing} is obtained via gluing by hand along local presentations.
  It is however impossible to generalize such an approach to more complicated invariants, and in particular to invariants of higher categorical nature such as matrix factorizations.
\end{observations}

In \cite{MR2851153,MR2681792} Kontsevich-Soibelman conjectured the possibility of gluing more refined invariants of singularities such as the motives of vanishing cycles \cite{ayoub1, ayoub2} as objects in the Morel-Voevodsky motivic homotopy theory of schemes \cite{voevodsky-morel}. More recently, Toda \cite{Toda2019,Toda2023} conjectured the possibility of gluing  the 2-periodic dg-categories of matrix factorizations of \cref{exampleMF}.
Our main result  \cref{theoremcontractibility} together with the observation that all the local invariants listed in \cref{Milnornumber,perversesheaf,exampleMF} only depend on the formal completions (see \cref{sectionalgebraictoformulainvariants}), provides  a mechanism to answer both Kontsevich-Soibelman and Toda's conjectures in full generality, upon the prescription of the appropriate orientation data. More precisely:

\begin{itemize}
	\item \cref{theoremcontractibility} implies that any local invariant satisfying simultaneously
	\begin{itemize}
		\item a Thom-Sebastiani formula for quadratic bundles, 
		\item isotopy-invariance,
		\end{itemize}
	 can be glued up to a \emph{twist} determined by the action of quadratic bundles.
\item	To \emph{untwist} this gluing one needs further \emph{orientation data} which depend on the nature of the local invariant:
	\begin{itemize}
		\item In \cref{examplemilnornatural} we use \cref{theoremcontractibility} to recover Behrend's function of \cref{behrendfunctiontheorem}.
		\item In \cref{naturaltransformationperversesection}  we use  \cref{theoremcontractibility} to recover BBDJS's perverse sheaf of \cref{joycegluing} and the notion of orientation data therein.
	\end{itemize}
	\end{itemize}

In the follow-up article \cite{hennionholsteinrobaloII}, we use \cref{theoremcontractibility} to study the more complicated example of matrix factorizations, thus
partially answering conjectures of Kontsevich-Soibelman and Toda:

\begin{theoremnonumber}[\cite{hennionholsteinrobaloII}]\label{theoremgluingMF}
  Let $X$ be a $(-1)$\nobreakdashes-shifted symplectic derived Deligne--Mumford stack.
  Then the locally defined $2$\nobreakdashes-periodic dg-categories of matrix factorizations $\MF(U,f)$
  of \cref{exampleMF} can be glued as a crystal of 2-periodic dg-categories "up to isotopies", provided the following orientation data:
  \begin{itemize}
    \item A continuous function $X \to \quot{\integers}{2}$ controlling the parity of the dimension of the local Darboux models -- corresponding to a trivialization of the zero class $\beta_1 = 0 \in \cohomology^1(X, \quot{\integers}{2})$;
    \item A square root of the canonical bundle $\canonicalbundle_X$ -- corresponding to a trivialization of the class $\beta_2 = \beta \in \cohomology^2(X, \quot{\integers}{2})$ of \cref{squarerootbundle};
    \item A sort of shifted spin structure on $X$ -- corresponding to the trivialization of a class $\beta_3 \in \cohomology^3(X, \quot{\integers}{2})$.
  \end{itemize}
  The space of such choices forms a gerbe over $X$, over which such a gluing always canonically exists.
\end{theoremnonumber}

\begin{warning}
The appearance of "up to isotopies"  in \cref{theoremgluingMF} is necessary as we do not yet know if $\MF$ itself satisfies the necessary local $\affineline{}$-invariance property.
\end{warning}

\subsection{Motivation:  Donaldson--Thomas invariants}\label{motivations}
To conclude this introduction let us  explain how  \cref{theoremcontractibility} fits in the framework of Donaldson--Thomas theory and $(-1)$\nobreakdashes-shifted symplectic geometry.

Following ideas of Kaluza-Klein (see \cite{Bourguignon1989} for a mathematical discussion), string theory \cite{CANDELAS198546} suggests that spacetime is locally of the form $\bbR^4\times Y$, with  $Y$ a Calabi-Yau threefold. For the purposes of this paper, this means $Y$ is a smooth algebraic variety of dimension 3 over $\complex$, together with a trivialization of the canonical bundle $\omega_Y\simeq \structuresheaf_Y$, such as the Fermat quintic.
The trajectories of string-like particles through space-time, define complex algebraic curves $C$ in $Y$ (see \cite[\S 25.4]{Zwiebach2009})  and an important part of modern enumerative algebraic geometry concerns the counting of such curves. There are several counting approaches.
The Gromov-Witten  approach counts \emph{parametrized curves} via intersection theory in the  moduli space $\stablemaps_{g,n}(Y,\beta)$ of stable maps $f \colon C\to Y$ passing through $n$ fixed points $y_1,.., y_n$ in $Y$ and with fixed degree $f_*[C]=\beta\in \cohomology_2(X,\integers)$.
See \cite{MR1291244, MR1363062, MR1412436, MR1431140}. The fact that this counting is well-defined is a consequence of the deformation theory of stable maps being concentrated in degrees $[-1,0]$.  A second counting strategy is to look at the \emph{Hilbert scheme}, $\Hilbert_2(Y)$ of \emph{embedded curves} in $Y$.
Contrary to the GW approach, the deformation theory of an embedded curve is no longer concentrated in degrees $[-1,0]$ unless the embedding is a local complete intersection. This paper deals with a third approach, the Donaldson--Thomas (DT) counting introduced in \cite{MR1818182}, which replaces embedded curves by their associated ideal sheaves $I_C\subseteq \structuresheaf_Y$ (and more generally, coherent sheaves). In this situation we are interested in the moduli space $\MCoh(Y)^\sigma$ of coherent sheaves on $Y$ with a prescribed stability condition.
As observed by R. Thomas \cite{MR1818182}, the deformation theory of sheaves is better-behaved: not only it is concentrated in degrees $[-1,0]$ but furthermore, a combination of Serre duality and the CY-condition forces a duality on the deformation theory of $E\in \MCoh(Y)^\sigma$:
\begin{equation}\label{thomasformula}
  \{1^{st} \text{ order def. of } E\} \simeq \left( \{\text{Obstructions to def. of } E\}\right)^{\vee}.
\end{equation}
Therefore, by  \cite{MR1437495},  $\MCoh(Y)^\sigma$ admits a virtual fundamental class and we define $\DT$\nobreakdashes-invariants by
\begin{equation}\label{definitionDTinvariants}
  \DT(Y) \coloneqq \int_{[\MCoh(Y)^\sigma]^{\virtual}} 1.
\end{equation}
DT invariants are explicitly related to GW-invariants. See \cite{Maulik2006,Bridgeland2011,Toda2010,Pandharipande2014} and more recently \cite{pardon2024universally}.

Our paper deals with the geometrical interpretation of the symmetry latent in \cref{thomasformula}, which Behrend \cite{MR2600874} incorporated in the notion of \emph{symmetric perfect obstruction theory} enhancing the notion of perfect obstruction theory  of \cite{MR1437495}.
In \cite{1111.3209}, Pantev-To\"en-Vaquié-Vezzosi discovered that the extra symmetry \cref{thomasformula}, and more generally, Behrend's notion of symmetric perfect obstruction theory, are really shadows of a $(-1)$\nobreakdashes-shifted symplectic structure on the derived enhancement $\RMCoh(Y)^\sigma$ of $\MCoh(Y)^\sigma$.
It goes as follows: since $Y$ is smooth, $\RMCoh(Y)^\sigma$ can be found as an open substack of the derived mapping stack $\derivedMap(Y, \Perf)$ where $\Perf$ is the derived moduli stack of perfect complexes of \cite{toen-vaquie}.
By \cite{1111.3209}, $\Perf$ has a $2$\nobreakdashes-shifted symplectic structure given by the universal trace map and therefore the derived mapping stack $\derivedMap(Y, \Perf)$ inherits a $(2-3)=(-1)$\nobreakdashes-shifted symplectic form, obtained by integrating the $2$\nobreakdashes-shifted symplectic form on $\Perf$. Therefore, as an open substack, $\RMCoh(Y)^\sigma$  inherits the $(-1)$\nobreakdashes-shifted symplectic form.
It is explained in \cite[Section 3.2]{toen-vaquie} how to recover the symmetric perfect obstruction theory from a $(-1)$\nobreakdashes-shifted symplectic structure.

This discussion places $\DT$\nobreakdashes-theory within the framework of $(-1)$\nobreakdashes-shifted symplectic geometry:  \cref{behrendfunctiontheorem} combined with \cref{definitionDTinvariants} tells us Behrend's function $\nu$ computes $\DT$\nobreakdashes-invariants and \cref{joycegluing}-(ii) exhibits Joyce's sheaf $\PJoyce$ as a \emph{categorification} of $\DT$\nobreakdashes-invariants, in the sense that $\chi(\PJoyce)=\nu$ gives back the $\DT$\nobreakdashes-counting.
Finally, our main result \cref{theoremcontractibility} provides a general mechanism to glue refined invariants of Donaldson--Thomas type, such as the categories of matrix factorization in  \cite{hennionholsteinrobaloII}.

\subsection{Outline}

After reviewing the theory of shifted symplectic structures (including isotopic and Lagrangian fibrations, distributions, existence of exact forms) in \cref{remindersshiftedsymplectic}
we turn to constructing the Darboux stack in \cref{SectionDarboux}.
In \cref{dCritandformalcompletions} we introduce the derived critical locus of a function $f$ on a smooth formal scheme $\formalU$.
We then view derived critical loci in terms of {Liouville data} in \cref{sectioncharacterizationsformaldarboux}, establishing that given a $(-1)$-shifted exact symplectic derived affine scheme $S$, the data of an exact symplectic equivalence $S \simeq \dCrit(\formalU, f)$ is tantamount to a closed immersion $i \colon S \hookrightarrow \formalU$ together with a structure of Lagrangian fibration on $i$.
In \cref{sectionconstructionformaldarboux} we then construct the moduli functor classifying Liouville data, and define the {Darboux stack} $\Darbstack^{\lambda}_X$ of a derived Deligne-Mumford stack with a $(-1)$-shifted exact symplectic form $\lambda$ in \cref{definitiondarbouxstackasstacks}.
We then compare our presentation with the more familiar language of Landau-Ginzburg pairs in \cref{sectionLGpairs}.

In \cref{section-Quadraticbundles} we turn to quadratic bundles and their action on the Darboux stack.
In \cref{sectionStackQuad} we construct the stack of flat quadratic bundles and then show in \cref{quadandliouville,section-Action} that it acts the Darboux stack $\Darbstack^{\lambda}_X$.
This section culminates with \cref{thm-transitivity}, proving this action is transitive.

We then turn to the proof of \cref{theoremcontractibility} in \cref{sectionisotopies}.
After reminders and preliminary results on classical $\affineline{}$-localisation in \cref{subsectionclassicalA1localization}, we explain how to impose $\affineline{}$-invariance on mapping spaces only in \cref{subsectionisotopies}. 
This lets us consider in \cref{darbouxstackuptoisitopysection} the $\affineline{}$-isotopic quotients $\QA$ and $\DA$ whose quotient is our main focus.
In \cref{contractibilitytheoremdetailedversion} we give a more precise formulation of \cref{theoremcontractibility} and reduce its proof to \cref{prop:A1-automorphisms}, which we prove in \cref{sectionproofcontractibility}.

We consider applications of \cref{theoremcontractibility} in \cref{sectionapplications}. We introduce the notion of critical invariants in \cref{sectioncriticalinvariants}, which covers Milnor number, vanishing cycles and categories of matrix factorizations -- see \cref{examplesofalgebraiccriticalinvariants}.
For discrete and 1-categorical critical invariants we show such invariants only depend on formal completions in \cref{sectionalgebraictoformulainvariants}.
In \cref{examplemilnornatural} we use \cref{theoremcontractibility} to  recover Behrend's function of \cref{behrendfunctiontheorem}.
In \cref{naturaltransformationperversesection}  we use  \cref{theoremcontractibility} to recover BBDJS's perverse sheaf and mixed Hodge module of \cref{joycegluing} and the notion of orientation data therein.

We conclude with some thoughts on universal orientation data in \cref{sectionuniversalorientation}, including an intriguing relation to semi-topological Grothendieck--Witt spectra.

\subsection{Acknowledgments}
This work has been in development for nearly four years, and during this period, we had the opportunity to discuss our results with several people and gain valuable insights. We would like to express our gratitude to Bertrand To\"en, who initially suggested us to explore the quotient of the stack of foliations, and to Mauro Porta, who was involved in the early phases of this project.
 Additionally, we extend our thanks to several people who provided us with highly encouraging feedback at different stages, and answered to some of our naive questions: Dario Beraldo, Chris Brav, Damien Calaque, Jean Fasel, Julien Grivaux, Stéphane Guillermou, Dominic Joyce, Mikhail Kapranov, Young-Hoon Kiem, Gérard Laumon, Tony Pantev, Hyeonjun Park, Massimo Pippi, Nick Rozenblyum, Pavel Safronov, Olivier Schiffmann, Vivek Shende, Hiro Lee Tanaka,  Richard Thomas and Gabriele Vezzosi.
\footnote{The first and third authors were supported by the grant ANR-17-CE40-0014.
	
	The second author acknowledges support by the Deutsche Forschungsgemeinschaft (DFG, German Research Foundation) through EXC 2121 ``Quantum Universe'' -- project number 390833306 -- and SFB 1624 ``Higher structures, moduli spaces and integrability'' -- project number 506632645. 
}

\subsection{Conventions and notations}
\label{notations}

Throughout the paper we work exclusively over the base field $\basefield$.

\begin{notations}{Unless mentioned otherwise, we use the notations of \cite{lurie-htt, lurie-ha} for higher categories. In particular:} \label{notationhighercategories}
  \item We denote by $\inftygpd$ the $(\infty,1)$\nobreakdashes-category of $\infty$\nobreakdashes-groupoids;
  \item We denote by $\inftycats$ the $(\infty,1)$\nobreakdashes-category of $(\infty,1)$\nobreakdashes-categories\footnote{Set theoretic issues will not play an important role and will thus stay implicit throughout.}, just called \icategories from now on;
  \item We write  $\calC^\maximalgroupoid$ for the maximal $\infty$\nobreakdashes-groupoid of an \icategory $\calC$;
  \item We use $\dgMod_\basefield$ for the underlying \icategory of chain complexes of $\basefield$\nobreakdashes-modules, obtained by inverting quasi-isomorphisms.
  However, and contrary to \cite{lurie-ha} we opt for the \emph{cohomological} convention.
\end{notations}

\begin{notations}[Derived Geometry]{We assume the reader is familiar with derived algebraic geometry at the level of the introductory surveys \cite[\S 1-4]{MR3285853} and \cite[\S 1-4]{Pantev2021}. For conventions:}\label{notationderivedgeometry}
  \item Derived affine rings are modeled by $\cdga$, the \icategory of cohomologically negatively graded commutative differential graded algebras of finite presentation over $\basefield$. For $A\in \cdga$ we denote by $\Spec(A)$ the corresponding derived affine scheme. We denote by $\dAff_\basefield$ the \icategory of derived affine schemes over $\basefield$, equivalent to $(\cdga)^\op$ via the functor $\Spec$.

  \item Let $\dSt_\basefield \subseteq \PSh(\dAff_\basefield, \inftygpd)$ denote the \icategory of derived stacks over $\basefield$, with étale hyperdescent.

  \item \label{descriptionsofslicecategoriesofstacks}For any $X \in \dSt_\basefield$, we denote by $\dSt_X$ the slice \icategory of stacks over $X$. Throughout this paper we will use other descriptions of this category, namely:
  \begin{enumerate}
    \item $\dSt_X \simeq \Sh(\dAff_X, \inftygpd)$ where $\dAff_X \subset \dSt_X$ denotes the full subcategory spanned by derived affine schemes equipped with the induced Grothendieck topology.
    This equivalence is a consequence of the equivalence $\PSh(\dAff_X)\simeq \PSh(\dAff)_{X}$ of \cite[\href{https://kerodon.net/tag/04BV}{04BV}]{kerodon} together with the fact the forgetful functor $\dAff_{X} \to \dAff_\basefield$ preserves covers and therefore the right Kan extension preserves sheaves \cite[Prop.\,3.12.11]{barwick2020exodromy}.
    \item  $\dSt_X \simeq \Sh(\dSt_X, \inftygpd) \simeq \Fun^{\mathcal{R}}( \dSt_X^\op, \inftygpd)$ where on the $\dSt_X$ is equipped with its canonical topology, and the RHS denote the category of limit preserving functors. See \eg \cite[Prop.\,1.1.12]{lurie-structuredspaces}.
  \end{enumerate}

  \item For two derived stacks $X$ and $Y$ we denote by $\stackMap(X,Y)$ the derived mapping stack in $\dSt_\basefield$.

  \item  For $X\in \dSt_\basefield$ we denote by $\classicaltruncation(X)$ its classical truncation.

  \item For $X\in \dSt_\basefield$ we denote by $\DQCoh(X)$ its derived \icategory of quasi-coherent sheaves \cite[\S 4.1]{MR3285853} and $\Perf(X)\subseteq \DQCoh(X)$ the full subcategory of perfect complexes. If $X=\Spec(A), \DQCoh(X)=\dgMod_A$.

  \item\label{defcotangent} For a derived stack $X$ we denote by $\cotangent_X\in \DQCoh(X)$ the cotangent complex of $X$ relative to $\basefield$ \cite[\S  3.1 and  \S4.1]{MR3285853} whenever it exists. This is the case for derived Artin stacks \cite[Cor.\,2.2.3.3.]{MR2394633}. Furthermore, as our stacks are all locally of finite presentation, the cotangent complex $\cotangent_{X}$ is perfect with dual $\tangent_{X}$ the \emph{tangent complex} of $X$ -- see \cite[Prop.\,2.2.2.4]{MR2394633} or \cite[Thm.\,7.4.3.18]{lurie-ha}. For a map of derived stacks $X\to Y$ we denote by $\cotangent_{X/Y}\in\DQCoh(X)$ its \emph{relative cotangent complex} if it exists.
\end{notations}

\begin{reminders}[Formal stacks and formal completions]{%
    We consider formal derived stacks and formal completions using the formalism of \cite[\S 2.1]{MR3653319}.}
  \label{formalstacks}
  \item For $A \in \cdga$, we denote by $\nilpotents_A \subseteq \pi_0(A)$ the ideal of nilpotent elements, $A_\red \coloneqq \pi_0(A)/\nilpotents_A$ the quotient, and  $i \colon \CAlg_\basefield^{\red} \hookrightarrow \cdga$ the inclusion of the full subcategory spanned by reduced discrete commutative algebras.
  The assignment $A \mapsto A_\red$ defines a left adjoint to the inclusion $i$.
  \item\label{reducedisadjoint} We denote by $\dSt_\basefield^\red$ the \icategory of stacks on the site $\CAlg_\basefield^{\red}$ with the étale topology.
  Composition with $i$ induces a restriction functor $i^\ast \colon \dSt_\basefield\to \dSt_\basefield^{\red}$ that admits both a left adjoint $i_!$ and a right adjoint $i_\ast$, both fully faithful.
  \item\label{remindersonderhamstack} For a derived stack $X$, we denote by $X_{\deRham} \coloneqq i_\ast i^\ast(X)$ its de Rham stack (introduced by Simpson). Its functor of points is given by $X_\deRham(A) \coloneqq X(A_\red)$.
  By adjunction, we get a functorial morphism $X \to X_\deRham$.
  The de Rham stack plays a crucial role, for its ties to flat connections and $\calD$\nobreakdashes-modules and for its relation to formal completions:
  \begin{remindersinner}
    \item\label{sheavesoverderhamareDmodules} When $X$ is a scheme, $\DQCoh(X_\derham)$ is equivalent to the \icategory of complexes of left D-modules \cite[Part I, Chapter 4]{GR-II}.\item\label{derhamandformalcompletion} The quotient map $X \to X_\deRham$  plays the role of the universal formal completion: given a map $f \colon Z=\Spec(A)\to X$, the formal completion $\widehat{Z}$ of $X$ along $Z$ coincides with the derived fiber product in $\dSt_\basefield$
    \[
      \begin{tikzcd}
        \widehat{Z} \tikzcart \ar{d} \ar{r} &
        X \ar{d} \\
        Z_{\deRham} \ar{r}[swap]{f_{\deRham}} &
        X_{\deRham}
      \end{tikzcd}
    \]
    (\cf \cite[Prop.\,2.1.8]{MR3653319}).
    Since $X_{\deRham}$ has no deformation theory, its cotangent complex vanishes.
    In particular, the canonical map $\cotangent_{X/k} \to \cotangent_{X/X_\deRham}$ is an equivalence and the map $\widehat{Z}\to X$ is formally étale.
  \end{remindersinner}
  \item As in \cite[\S 2.1]{MR3653319}, throughout this paper, formal stacks will be seen as a full subcategory of the \icategory of derived stacks $\dSt_\basefield$.
  In particular, smooth formal schemes $\FSchsm$  form a full subcategory of derived stacks.
  See \cite[Thm.\,1.4.2, Lem.\,1.4.5]{zbMATH06760099} and \cite[\S I. 1.4]{zbMATH06829730} for the presentation as inductive limits of schemes.
  \item For a derived stack $X\in \dSt_\basefield$ we set $X_\red \coloneqq i_!\circ i^\ast(X)$ its reduced substack.
  By adjunction it comes with a canonical map $X_\red \to X$.
  When $X=\Spec(A)$ is an affine derived scheme, we have $X_\red=\Spec(A_\red)$. It enjoys the following properties:
  \begin{remindersinner}
    \item By construction, the assignment $X\mapsto X_\red$ defines a left adjoint to the functor $(-)_{\deRham} \colon \dSt_\basefield\to \dSt_\basefield$ (see \eg \cite[2.1.3 and 2.1.4]{MR3653319}).
    \item\label{equivalencederhamvsreduced} There are canonical equivalences $(X_\red)_\deRham\iso (X)_\deRham$ and $(X)_\red\iso (X_\deRham)_\red$.
    \item\label{equivalencedreducediffequivalencederham}
    A morphism $X \to Y$ induces an equivalence $S_\red \to F_\red$ if and only if it induces an equivalence $S_\derham\to F_\derham$.
    \item \label{reducedoflimits} For any diagram $\alpha \mapsto X_\alpha$, the canonical morphism is an equivalence
    \begin{equation}\label{extractingreducedsubstack}
      (\lim_\alpha X_\alpha)_{\red}\simeq (\lim_\alpha (X_\alpha)_{\red})_{\red}.
    \end{equation}
  \end{remindersinner}
  \item\label{reducedeqplusetaleisiso} Consider a morphism of derived stacks $f \colon X \to Y$ (locally of finite presentation). If both $X$ and $Y$ admit a deformation theory, if $f_\red$ is an equivalence and if $f$ is formally étale, then $f$ is an equivalence (see \cite[I-1 8.3.2]{GR-II}).
  \item\label{etaleformalcompletionisnoop} The above specializes to: if $f \colon X \to Y$ is formally étale, then $X \to \widehat X$ is an equivalence: \ie we have a Cartesian square
  \[
    \begin{tikzcd}
      X \tikzcart \ar{d} \ar{r} &
      Y \ar{d} \\
      X_{\deRham} \ar{r}[swap]{f_{\deRham}} &
      Y_{\deRham}.
    \end{tikzcd}
  \]
\end{reminders}

\begin{notations}[Étale site]{}\label{notationetale}
  \item When $X$ is a derived stack, we denote by $X_\et$ the small étale site of $X$, with objects given by maps $S \to X$ which are étale (meaning, representable, locally of finite presentation and the relative cotangent complex $\cotangent_{S/X}$ vanishes).
  \item For a presentable \icategory $\calC$, we denote by
  \[
    \Sh(X_\et,\calC)\subseteq \PSh(X_\et, \calC) \coloneqq \Fun(X_\et^\op, \calC)
  \]
  the full subcategory spanned by hypercomplete sheaves.

  \item We denote by $X_\et^\daff \subseteq X_\et$ the full subcategory spanned by étale maps $S \to X$ with $S$ an \emph{affine} derived scheme. It follows \cite[VII Prop.\,3.1]{MR0354653} that if $X$ admits an affine étale atlas, the inclusion $X_\et^\daff \subseteq X_\et$ induces an equivalence on sheaves
  \[
    \Sh(X_\et^\daff, \calC)\iso \Sh(X_\et,\calC).
  \]
  For this reason we will mainly work with $X_\et^\daff$ throughout the paper.
\end{notations}

\begin{reminders}{}\label{remindersonequivalencesofsites}
  \item \label{equivalenceetalesites}
  For any derived stack $X$, the closed immersions $X_\red \to \classicaltruncation(X) \to X$ induce equivalences
  \[
    X_\et \simeq (\classicaltruncation(X))_\et \simeq (X_\red)_\et.
  \]
  Indeed, it follows from \cite[\href{https://stacks.math.columbia.edu/tag/04DZ}{Thm.\,04DZ}]{stacks-project} that the étale site is invariant under nilpotent extensions of schemes.
  It follows from \cite[Cor.\,2.2.2.9]{MR2394633} or \cite[Thm.\,7.5.0.6]{lurie-ha} that the étale site is impervious to derived structure as well.
  \item\label{identificationmorphismofsites} For an étale schematic morphism $u \colon X \to Y$ between derived stacks, the pullback functor $u^\ast \colon Y_\et \to X_\et$ admits a left adjoint $u_\sharp$ sending an étale map $S \to X$ to the composition $S \to X \to Y$. This defines a morphism of sites and therefore at the level of sheaves, the induced pullback functor
  \[
    u^\ast \colon \Sh(Y_\et,\calC)\to \Sh(X_\et,\calC)
  \]
  as in \cite[Ex.\,17.3]{MR2182076}, identifies with $u^\ast(F)\simeq F\circ u_\sharp$.
  \item \label{enoughpoints} By \cite[VIII  Thm.\,3.5 and Thm.\,7.9]{MR0354652}\footnote{See also \cite[\href{https://stacks.math.columbia.edu/tag/04JU}{04JU}]{stacks-project}} and (a), when $X$ is a derived Deligne-Mumford stack, the classical topos $\Sh(X_\et, \Sets)$ has enough points, coinciding with the collection of geometric points of $X$.
\end{reminders}

\section{Reminders on shifted symplectic geometry}
\label{remindersshiftedsymplectic}

In this section we review the tools related to shifted symplectic geometry used in the paper. The reader can consult  \cite[\S 7-8]{Pantev2021} and \cite[\S 5]{MR3285853} for a survey and  \cite{1111.3209, MR3940792, MR3653319} for the technical details.

\subsection{Closed forms and shifted symplectic structures}
\label{sectionmixedgradedstuff}

The main technical obstacle in defining derived symplectic forms lies in the definition of closed forms. When $X$ is a smooth classical affine scheme of dimension $n$ the situation is easy:
$p$\nobreakdashes-forms are defined as global sections of the wedge power  $\bigwedge^p\Omega_X^1$ of  $\Omega_X^1$  the vector bundle of Kähler differentials. To define \emph{closed $p$\nobreakdashes-forms} we consider the \emph{naively truncated} de Rham complex of vector bundles on $X$
\begin{equation}\label{naivetruncatedderhamcomplex}
  \Omega_X^{\geq p} \coloneqq \left[
    \begin{tikzcd}[column sep=scriptsize]
      0 \ar[r] & \cdots \ar[r] & 0 \ar[r] & \underbrace{\Omega_X^p}_{\mathclap{\text{coh. deg. } p}} \ar{r}{\dderham} & \underbrace{\Omega_X^{p+1}}_{\mathclap{\text{coh. deg.  }p+1}}\ar{r}{\dderham} & \cdots \ar{r}{\dderham} & \underbrace{\Omega_X^n}_{\mathclap{\text{coh. deg. }n}}
    \end{tikzcd}\hspace{1em} \right]
\end{equation}
and recover the usual definition of closed $p$\nobreakdashes-forms on $X$ as cohomology classes $\cohomology^0(X,\Omega_X^{\geq p}[p])$.

For an affine derived scheme $S=\Spec(A)$ we must replace Kähler differentials by the cotangent complex $\cotangent_S$ (\cref{defcotangent}). To access the naively truncated de Rham complex, we consider first the \emph{de Rham algebra} $\DR(X)$ which assembles the different derived wedge products $(\bigwedge^p \cotangent_X)[p]$ as graded pieces of a graded complex,
\begin{equation}
  \label{derhamalgebra}
  \DR(X) \coloneqq  \Sym_{A}(\cotangent_A[1])= \bigoplus_{p\geq 0} \left(\bigwedge^p_A \cotangent_A \right)[p]
\end{equation}
with $\cotangent_A[1]$ in weight $1$, and implements the de Rham differential as an action of an extra operator $\epsilon$ that shifts the cohomological degree by $-1$ and the grading by $1$,
\begin{equation}\label{derhamalgebramixed}
  \epsilon=\dderham \colon  \left(\bigwedge^p_A \cotangent_A\right)[p]\to \left(\bigwedge_A^{p+1} \cotangent_A\right)[p+1][-1].
\end{equation}
Here, we see each $\bigwedge^p_A \cotangent_A \in \DQCoh(S)=\dgMod_A$ as an object in $\dgMod_\basefield$ via the forgetful functor.

Finally, the \emph{totalization} -- $|\DR(X)|$ -- with respect to the internal differentials of each graded piece, together with the extra differential defined by $\dderham$, gives the correct construction of $\Omega_X^{\geq 0}$ in the non-smooth setting. Let us now quickly survey the structures underlying this discussion, summarizing the results of \cite[\S 1.1 and \S 1.2]{1111.3209}:

\begin{reminders}[Graded modules and mixed graded modules]{}
  \item \label{remindergraded} Denote by $\gradedmodules_\basefield$ the \icategory of $\integers$\nobreakdashes-graded chain complexes. Objects in $\gradedmodules_\basefield$ are families $E= (E_n)_{n \in \integers}$ with $E_n \in \dgMod_\basefield$. By construction,  $\gradedmodules_\basefield$ comes with a canonical direct sum functor $\bigoplus \colon  \gradedmodules_\basefield \to \dgMod_\basefield$ sending $E= (E_n)_{n \in \integers} \to \bigoplus_n E_n$ and we will sometimes identify $E$ with the associated direct sum. Also, $\gradedmodules_\basefield$ comes with operators $E\mapsto E(p)$ that shifts the weights by $p$, namely, by sending $E= (E_n)_{n \in \integers}\mapsto E(p) \coloneqq  (E_{n+p})_{n\in \integers}$. To conclude, we observe that $\gradedmodules_\basefield$ carries a symmetric monoidal structure given on objects by \[
    E\derivedtensor{\basefield} F= \left(\bigoplus_{n+m=\ell} E_n\derivedtensor{\basefield} F_m\right)_{\ell \in \integers}
  \]
  \item \label{remindermixedgraded}
  Denote by $\basefield[\epsilon]_{\gr}$ the  graded  strictly associative dg-algebra freely generated by an element $\epsilon$ in cohomological degree $-1$ and weight $1$, and \emph{strictly} verifying $\epsilon^2=0$.
  A \emph{mixed graded module} is a left $\basefield[\epsilon]_{\gr}$\nobreakdashes-module in $\gradedmodules_\basefield$.
  Informally, this amounts to the data of an object $E=(E_n)_{n\in \integers}$  together with an extra operator
  \begin{equation}\label{conventionmixedgradeddirection}
    \epsilon \colon  E(1)[1]\to E
  \end{equation}
  with $\epsilon \circ \epsilon =0$. We denote by
  $\mixedgradedmodules_\basefield$ the \icategory of mixed graded modules. The algebra $\basefield[\epsilon]_{\gr}$ carries a strictly commutative graded Hopf structure
  $\basefield[\epsilon]_\gr\to \basefield[\epsilon]_\gr\otimes_{\gr}\basefield[\epsilon]_\gr
  $ determined by $
    \epsilon\mapsto \epsilon\otimes 1 + 1\otimes \epsilon
  $.
  This induces a symmetric monoidal on  $\mixedgradedmodules_\basefield$ and we defined \emph{mixed graded algebras} $\mixedgradedalgebras_\basefield$ as commutative algebra objects in mixed graded modules.
\end{reminders}

\begin{remark}
  In this paper we use the definition of mixed graded modules of \cite{1111.3209} where the action of $\epsilon$ \emph{decreases} the cohomological degree by one \cref{conventionmixedgradeddirection}. This is \emph{opposite} to the conventions used in \cite{MR3653319}. Thanks to the red-shift equivalence, both approaches are equivalent and in particular, yield compatible filtered Tate totalizations (in the sense of \cref{reminderrealizations} below). See \cite[Rem.\,1.1.3]{MR3653319} and \cite[p.\,9]{Toen2023b}.
\end{remark}
\begin{reminder}\label{reminderderhamalgebra}[de Rham algebra]
  The \emph{de Rham} algebra $\DR$ of \cref{derhamalgebra} with the mixed graded structure of \cref{derhamalgebramixed} is implemented by an \ifunctor
  \[
    \DR \colon \cdga \to \mixedgradedalgebras_\basefield.
  \]
\end{reminder}

\begin{reminder}\label{reminderrealizations}[Totalizations]
  Let $\basefield(p)$ be the mixed graded module given by a copy of $\basefield$ in weight $p$ and zero otherwise, equipped with the zero action of $\epsilon$.  Let $E$ be a mixed graded complex. We define the \emph{weight $p$ totalization} $|E|^{\geq p}$ as the image of $E$ under the  limit preserving lax monoidal\footnote{See \cite[\S 1.1]{Toen2020a}} \ifunctor given by
  \begin{equation}\label{formulafortruncatedderhamasRHOM}
  \begin{aligned}
    |-|^{\geq p} \colon \mixedgradedmodules_\basefield &\to \dgMod_\basefield
    \\
    E &\mapsto |E|^{\geq p} \coloneqq \RHom_{\basefield[\epsilon]_\gr}(\basefield(p), E).
  \end{aligned}
  \end{equation}
  Explicitly,
  \[
    |E|^{\geq p} = \left(\prod_{n\geq p} E_n[-2(n-p)], d+\epsilon\right)
  \]
  with differential given by $d+\epsilon$, where $\epsilon$ is the mixed graded differential and $d$ is the internal differential.  The assembly of all totalizations $|-|^{\geq p }$ provides a functor with values in graded modules
  \[
    |-|^\gr \colon \mixedgradedmodules_\basefield \to \dgMod_\basefield^\gr
    \hspace{2em}
    E \mapsto (|E|^{\geq p})_{p \in \integers}.
  \]
  The composition with the map of mixed graded modules
  \begin{equation}\label{canonicalmapmixgradedmodules}
  \basefield[\epsilon]_\gr(p) \to \basefield(p)
  \hspace{1em} \text{sending} \hspace{1em}
  \epsilon \mapsto 0
  \end{equation}
  defines a natural transformation $\eta$ between $|-|^\gr$ and the forgetful functor that discards the $\basefield[\epsilon]_\gr$\nobreakdashes-action $\mixedgradedmodules_\basefield \to \dgMod_\basefield^\gr$.
  On each weight $\eta$ is given by
  \begin{equation}\label{naturaltransformationunderlyingform}
    \eta_E^p \colon |E|^\gr_p \coloneqq  |E|^{\geq p} \to E_p.
  \end{equation}
\end{reminder}

\begin{example}\label{closedformssmoothcase}
  Let $S=\Spec(A)$ where $A$ is a smooth $\basefield$\nobreakdashes-algebra. Then we have
  \[
    \Omega_A^{\geq p}\simeq |\DR(A)|^{\geq p}[-2p]
  \]
  where $\Omega_A^{\geq p}$ denotes the complex in \cref{naivetruncatedderhamcomplex}.
\end{example}

\begin{reminder}\label{constructiontruncatedderhamcomplex}[Derived complete truncated de Rham complexes]
  Let $S=\Spec(A)$ be a derived affine scheme. The \emph{derived complete $p$\nobreakdashes-truncated de Rham complex}, denoted $\derhamcomplex_S^{\geq p}$, is defined as
  \[
    \derhamcomplex_S^{\geq p} \coloneqq  |\DR(A)|^{\geq p}[-2p].
  \]
  This recovers Illusie's derived de Rham complex. See \cite[Construction 4.1]{1207.6193}. The complex of \emph{closed} $p$\nobreakdashes-forms is $\derhamcomplex_S^{\geq p}[p]$. More generally,  the complex of \emph{$n$\nobreakdashes-shifted closed $p$\nobreakdashes-forms} is given by $\derhamcomplex_S^{\geq p}[p+n]$. In particular, we can present the data of a $n$\nobreakdashes-shifted closed $p$\nobreakdashes-form as a map of mixed graded modules
  \[
    \basefield[p-n](p)\to \DR(A)
  \]
  where $\basefield[p-n](p)$ is the mixed graded module given by $\basefield[p-n]$ in weight $p$ with trivial action of $\epsilon$.  The natural transformation \cref{naturaltransformationunderlyingform} defines a map in $\dgMod_\basefield$
  \begin{equation}\label{underlying2form}
    \derhamcomplex_S^{\geq p}[p+n]=|\DR(A)|^{\geq p}[n-p]\to \left(\bigwedge^p_A \cotangent_A\right)[p][n-p]=\left(\bigwedge^p_A \cotangent_A\right)[n]
  \end{equation}
  which forgets the data rendering the $p$\nobreakdashes-form closed.
\end{reminder}

Finally, we can recall the definition of $(-1)$\nobreakdashes-shifted symplectic form in the affine case:

\begin{definition}\label{definitionsymplectic}[Shifted Symplectic structures]
  Let $S=\Spec(A)$ be a derived affine scheme.
  A \emph{$(-1)$\nobreakdashes-shifted closed $2$\nobreakdashes-form} on $S$ is a map in $\mixedgradedmodules_\basefield$
  \[
    \omega \colon \basefield[2-(-1)](2)=\basefield[3](2) \to \DR(S).
  \]
  Equivalently, it corresponds to a map in $\dgMod_\basefield$
  \[
    \omega \colon \basefield\to |\DR(S)|^{\geq 2}[-3] = \derhamcomplex^{\geq 2}_S[1].
  \]
  The \emph{underlying} $2$\nobreakdashes-form is given by the composition in
  \begin{equation}\label{underlying2form2}
    \basefield\to \derhamcomplex^{\geq 2}_S[1]\to (\cotangent_S \wedge \cotangent_S )[-1]
  \end{equation}
  and via the extension of scalars adjunction $\DQCoh(S)\to \dgMod_\basefield$, \cref{underlying2form2} corresponds to a map in $\DQCoh(S)=\dgMod_A$
  \begin{equation}
    \label{underlying2form22}
    \structuresheaf_S \to (\cotangent_S\wedge \cotangent_S )[-1]\rlap{.}
  \end{equation}
  We say that $\omega$ is \emph{symplectic} if furthermore, the $\structuresheaf_S$\nobreakdashes-linear pairing associated to the composition with \cref{underlying2form22} is a non-degenerate pairing, \ie induces an equivalence $\omega \colon \tangent_S \simeq \cotangent_S[-1]$.
\end{definition}

\begin{reminder}
  \label{basicelementaction}
  \cite[p.8]{2007.09251} shows that $\RHom_{\basefield[\epsilon]_\gr}(\basefield(0), \basefield(1)[2])\simeq \basefield$. Through this equivalence, the element $1\in \basefield$ defines a map of mixed graded complexes $\basefield(0) \to \basefield(1)[2]$ which is the boundary map of the cofiber of \cref{canonicalmapmixgradedmodules} in mixed graded modules:
  \begin{equation}\label{basiccofiber}
    \begin{tikzcd}
      \basefield[\epsilon]_\gr \ar{d} \ar[r] \tikzcart & \basefield(0) \ar{d}{1} \\ 0 \ar[r] & \basefield(1)[2]\rlap{.}
    \end{tikzcd}
  \end{equation}
  The iteration of the map $1$ in \cref{basiccofiber} defines a sequence of maps in $\mixedgradedmodules_\basefield$
  \[
    \begin{tikzcd}[column sep=scriptsize]
      \cdots \ar[r] & \basefield(-2)[-4] \ar[r] & \basefield(-1)[-2] \ar[r] & \basefield(0) \ar[r] & \basefield(1)[2] \ar[r] & \basefield(2)[4] \ar[r] & \cdots
    \end{tikzcd}
  \]
  so that, given $E\in \mixedgradedmodules_\basefield$, and the formula in \cref{formulafortruncatedderhamasRHOM}, this produces a sequence
  \begin{equation}\label{sequencefiltered}
    \begin{tikzcd}
      \cdots \ar{r} & {} |E|^{\geq 2}[-4] \ar{r} & {} |E|^{\geq 1}[-2] \ar{r} & {} |E|^{\geq 0}\rlap{.}
    \end{tikzcd}
  \end{equation}
  In particular, for any derived affine scheme $S$ and for any $p \geq 0$, the composition of the maps in \cref{sequencefiltered} defines a map
  \begin{equation}\label{inclusionnaivelytruncationmap}
    \derhamcomplex_S^{\geq p} \to \derhamcomplex_S^{\geq 0}
  \end{equation}
  corresponding to the inclusion of the naively truncated de Rham complex.  Moreover, the cofiber sequence \cref{basiccofiber} induces for any derived affine scheme $S=\Spec(A)$,  cofiber sequences in $\Mod_\basefield$
  \begin{equation}\label{gradedpiecestruncatedderhamfiltration}
    \begin{tikzcd}
      \derhamcomplex_S^{\geq p+1 } \ar{r} \ar{d} \tikzcocart & \derhamcomplex_S^{\geq p} \ar{d} \\ 0 \ar{r} & \left(\bigwedge^p_S \cotangent_S\right)[-p]\rlap{.}
    \end{tikzcd}
  \end{equation}
\end{reminder}

\begin{definitions}[{\cite{1111.3209}}]{}
  \item \label{stackofclosedforms} The \ifunctors $\stackforms^{p, \closed}(- , n)$ and  $\stackforms^{p}(- , n)$ sending an affine derived scheme $S$ to the Dold-Kan constructions
  \[
    \stackforms^{p, \closed}(- , n) \colon  S\mapsto \DK\left(\derhamcomplex_S^{\geq p}[p+n]\right)
    \hspace{1em}\text{and}\hspace{.5em}
    \stackforms^{p}(- , n) \colon S\mapsto \DK\left(\left(\bigwedge^p_S \cotangent_S\right)[n]\right)
  \]
  satisfy étale hyperdescent and therefore define derived stacks classifying respectively, closed $p$\nobreakdashes-forms and $p$\nobreakdashes-forms (see \cite[Prop.\,1.11]{1111.3209}).
  \item For $X$ a derived Artin stack, we define $n$\nobreakdashes-shifted symplectic forms as maps of derived stacks $X\to \stackforms^{p, \closed}(-,n)$ satisfying the non-degeneracy condition of \cref{definitionsymplectic}.
\end{definitions}

\subsection{Exact symplectic structures}
\label{sectionfixingexactstructureequalsconstant}

\begin{notation}\label{truncatedabovederhamcomplex}
  Let $S$ be a derived affine scheme. We denote $\derhamcomplex_S^{\leq p-1}$ the homotopy cofiber of \cref{inclusionnaivelytruncationmap} in $\dgMod_\basefield$
  \begin{equation}\label{cofibersequencesfortruncations}
    \begin{tikzcd}
      \derhamcomplex_S^{\geq p} \ar{r} \ar{d} \tikzcocart & \ar{d} \derhamcomplex_S^{\geq 0} \\
      0 \ar{r} & \derhamcomplex_S^{\leq p-1}\rlap{.}
    \end{tikzcd}
  \end{equation}
\end{notation}

\begin{definitions}{}
  \item \label{stackofexactpforms} We define the stack $\stackforms^{p,\exact}(-,n)$ of $n$\nobreakdashes-shifted exact $p$\nobreakdashes-forms:
  \begin{multline*}
    \stackforms^{p,\exact}(-,n) \coloneqq \fiber\left(\stackforms^{p,\closed}(-,n) \to \stackforms^{0,\closed}(-,n+p)\right) \colon \\S \mapsto \DK\left(\derhamcomplex^{\leq p-1}_S[p+n-1]\right).
  \end{multline*}
  \item A \emph{$p$\nobreakdashes-shifted exact symplectic derived stack} is a derived Artin stack $X$ equipped with an exact $2$\nobreakdashes-form $\lambda \colon X \to \stackforms^{2,\exact}(-, p)$ such that the underlying (closed) $2$\nobreakdashes-form is non-degenerate (thus symplectic).
  \item \label{definitionexact} Let $(X,\omega_X)$ be a $p$\nobreakdashes-shifted symplectic derived Artin stack. An \emph{exact structure} on $X$ is a lift of $\omega_X \colon X \to \stackforms^{2,\closed}(-, p)$ to $\stackforms^{2,\exact}(-, p)$.
\end{definitions}

\begin{example}\label{symmetryderivedcriticalocusexplicit}
  The derived critical locus $\dCrit(U,f)$ of a function $f$ on a smooth affine scheme $U$ comes equipped with a $(-1)$\nobreakdashes-shifted exact symplectic form.
  Recall that $\dCrit(U,f)$ is defined as the derived intersection
  \begin{equation}\label{fiberproductdCrit}
    \begin{tikzcd}
      \dCrit(U,f) \ar{r}{i} \ar{d}[swap]{j} \tikzcart & U \ar{d}{0}\\
      U \ar{r}{df} & \cotangentstack U\rlap{.}
    \end{tikzcd}
  \end{equation}
  Notice that $i$ and $j$ are homotopic, as the composition with the canonical projection $\pi \colon \cotangentstack U \to U$ implies $i=\pi\circ df \circ i \sim \pi\circ 0\circ j=j$.
  In particular, the $2$\nobreakdashes-cell rendering the commutativity of \cref{fiberproductdCrit} amounts to a null-homotopy $\phi$ of the composition
  \begin{equation}\label{criticallocusasnullhomotopy}
    \begin{tikzcd} \basefield\ar{r}{f} & \structuresheaf_U \ar{r}{(\dderham)_U} & i^\ast \cotangent_U. \end{tikzcd}
  \end{equation}
  In explicit terms, the derived ring of functions $\structuresheaf_{\dCrit(U,f)}$ can be computed using the explicit Koszul resolution for the fiber product \cref{fiberproductdCrit} and we obtain
  \begin{equation}
    \label{explicitmodulefunctionsdCrit}
    \structuresheaf_{\dCrit(U,f)}\simeq (\Sym(\tangent_U[-1]), df)=[\cdots \to \Lambda^2\tangent_\formalU \to \tangent_\formalU]
  \end{equation}
  with differential given by contraction with $df$. See \cite[\S 7.3]{Pantev2021}.

  Now, both the zero section and the section given by the $1$\nobreakdashes-form $df$ are Lagrangian for the standard symplectic form on $\cotangentstack U$: the zero section, tautologically, and $df$ because $d(df)=0$.
  By \cite[Cor.\,2.10]{1111.3209}, being a derived intersection of Lagrangians, $\dCrit(U,f)$ inherits a $(-1)$\nobreakdashes-shifted symplectic form.
  Moreover, this $(-1)$\nobreakdashes-symplectic form has a canonical exact structure (\cref{definitionexact}) resulting from the combination of the Liouville form on $\cotangentstack U$, with the fact both $0$ and $df$ are exact $1$\nobreakdashes-forms on $U$.

  For the sake of completion we include here the construction of this exact structure: let $\basefield \to \Omega_{\cotangentstack U}^{\geq 2}[2]$ be the map classifying the symplectic form in $\cotangentstack U$, according to \cref{definitionsymplectic}. The Liouville form $\lambda$ is given by a map
  \[
    \lambda \colon \basefield[-2]\to \Omega_{\cotangentstack U}^{\leq 1}[-1]
  \]
  which we can fit in a commutative diagram
  \[
    \begin{tikzcd}[row sep=scriptsize]
      && \Omega_{U}^{\leq 1}[-1] \ar{dr}{i^\ast} &
      \\
      \basefield[-2] \ar{r}{\lambda} \ar[bend left=1.8em]{urr}[swap]{df} \ar[bend right=1.8em]{drr}{0} &
      \Omega_{\cotangentstack U}^{\leq 1}[-1] \ar{ur}{(df)^\ast} \ar{dr}[swap]{0^\ast} &&
      \Omega_{\dCrit(U,f)}^{\leq 1}[-1].
      \\
      && \Omega_{U}^{\leq 1}[-1] \ar{ur}[swap]{j^\ast}&
    \end{tikzcd}
  \]
  But using \cref{closedformssmoothcase}, since $U$ is smooth, $\Omega_{U}^{\leq 1}[-1]$ is the 2-term complex
  \[
    \Omega_{U}^{\leq 1}[-1]=[\underbrace{\structuresheaf_U}_{\mathclap{\text{deg. }1}}\to \underbrace{\Omega^1_U}_{\mathclap{\text{deg. }2}} ]
  \]
  and $df$ factors through the inclusion of the component in degree 2,
  \[
    \basefield[-2]\to \Omega_U^1[-2]\to  \Omega_{U}^{\leq 1}[-1]
  \]
  But since $df$ is \emph{exact}, its primitive $f$ provides a factorization through the homotopy fiber
  \[
    \begin{tikzcd}
      &\structuresheaf_U[-2] \ar{d}{\dderham} \ar{r} \tikzcart & 0 \ar{d} \\
      \basefield[-2] \ar{r}{df} \ar{ur}{f} & \Omega^1_U[-2] \ar[r] & \Omega_{U}^{\leq 1}[-1]\rlap{.}
    \end{tikzcd}
  \]
  The above diagram becomes
  \[
    \begin{tikzcd}[row sep=scriptsize]
      && \Omega_{U}^{\leq 1}[-1] \ar{dr}{i^\ast} &
      \\
      \basefield[-2] \ar{r}{\lambda} \ar[bend left=1.8em]{urr}[swap]{0} \ar[bend right=1.8em]{drr}{0} &
      \Omega_{\cotangentstack U}^{\leq 1}[-1] \ar{ur}{(df)^\ast} \ar{dr}[swap]{0^\ast} &&
      \Omega_{\dCrit(U,f)}^{\leq 1}[-1]
      \\
      && \Omega_{U}^{\leq 1}[-1] \ar{ur}[swap]{j^\ast}&
    \end{tikzcd}
  \]
  and the induced delooping map
  \begin{equation}\label{explicitexactstructureondCrit}
    \lambda \colon \basefield[-2]\to \Omega_{\dCrit(U,f)}^{\leq 1}[-2]
  \end{equation}
  defines the exact structure on $\dCrit(U,f)$.

  To conclude this example we provide an explicit description for the underlying non-degenerate 2-form: the derived fiber product \cref{fiberproductdCrit} induces a pullback square of tangent complexes
  \begin{equation}\label{lagrangianconditiondCritHessian}
    \begin{tikzcd}
      \tangent_{\dCrit(U,f)} \tikzcart \ar{r}{\tangentmap{i}}\ar{d}{\tangentmap{i}} & \tikzcart i^\ast \tangent_U\ar{d}{\tangentmap{0}}\ar{r} & 0 \ar{d} \\
      i^\ast\tangent_U\ar{r}{\tangentmap{df}} & i^\ast 0^\ast \tangent_{\cotangentstack U} \simeq i^\ast\tangent_U\oplus i^\ast \cotangent_U\ar{r} & i^\ast\cotangent_U
    \end{tikzcd}
  \end{equation}
  and the bottom horizontal  composition is given by the Hessian of $f$, $\hessian(f)$, at each closed point.
  It follows that the induced map $\tangent_{\dCrit(U,f)/U}\to i^\ast \cotangent_U[-1]$ is an equivalence and $\tangent_{\dCrit(U,f)}$ is given by the two-term complex
  \begin{equation}\label{explicitmodelcotangetcrit22}
    \begin{tikzcd}[column sep=large]
      i^\ast \tangent_{U} \ar{r}{\hessian(f)} & i^\ast \cotangent_{U}.
    \end{tikzcd}
  \end{equation}
  Finally, thanks to the symmetry of the Hessian, the complexes $\tangent_{\dCrit(f)}$ and $\cotangent_{\dCrit(f)}$ differ only by a shift of $-1$, \ie
  \[
    \tangent_{\dCrit(f)}\simeq \cotangent_{\dCrit(f)}[-1].
  \]
  This equivalence is precisely the underlying 2-form of the \emph{$(-1)$\nobreakdashes-shifted exact symplectic structure} obtained from the (exact) Lagrangian intersection.
\end{example}

\begin{proposition}\label{exactformssplitting}
  Consider the stack (sheaf, really) $\stackforms^{0,\closed}(-,0) \colon S \mapsto \cohomology_\deRham^0(S)$.
  The fiber-cofiber sequence $\derhamcomplex^{\leq 1} \to \derhamcomplex^{\geq 2}[1] \to \derhamcomplex^{\geq 0}[1]$ induces, after rotation, a fiber sequence
  \[
    \stackforms^{0,\closed}(-,0) \to \stackforms^{2,\exact}(-,-1) \to \stackforms^{2,\closed}(-,-1).
  \]
  This sequence is also a cofiber sequence and admits a canonical retract
  \begin{align*}
    \stackforms^{2,\exact}(-,-1) & \to \stackforms^{0,\closed}(-,0) \\
    \lambda                      & \mapsto c_\lambda.
  \end{align*}
  In particular, $\stackforms^{2,\exact}(-,-1) \simeq \stackforms^{2,\closed}(-,-1) \oplus \stackforms^{0,\closed}(-,0)$.
\end{proposition}
\begin{remark}
  Note that the short exact sequence from the proposition corresponds under Dold-Kan and HKR to the fiber sequence of periodic cyclic, cyclic and negative cyclic chains on $A = \structuresheaf_S$:
  \[
    \operatorname{PC}(A)[-2] \to \operatorname{CC}(A)[-4] \to \operatorname{NC}(A)[-3].
  \]
\end{remark}
\begin{proof}[Proof of \cref{exactformssplitting}]
  \newcommand{\lc}{\mathrm{lc}}
  By \cite[Prop.\,5.6(a)]{MR3904157}, the morphism
  \[
    \stackforms^{2,\exact}(-,-1) \to \stackforms^{2,\closed}(-,-1)
  \]
  is an epimorphism, so the fiber sequence at hand is indeed a cofiber sequence as well.
  The functorial projection $\derhamcomplex_S^{\leq 1} \to \derhamcomplex_S^{\leq 0} = \globalsections(S,\structuresheaf_S)$ induces a morphism
  \[
    \pi \colon \stackforms^{2,\exact}(-,-1) \to \stackforms^{2,\exact}(-,-1)_\deRham \to \affineline{\deRham}.
  \]
  For $S$ affine (derived) scheme, any element of $\stackforms^{2,\exact}(S,-1)_\deRham \coloneqq \DK(\derhamcomplex_{S_\red}^{\leq 1})$ amounts to a function $f \colon S_\red \to \affineline{}$ and a nullhomotopy $df \sim 0 \in \cotangent_{S_\red}$.
  The function $f$ is in particular locally constant on the smooth locus of $S_\red$, and thus on all of $S_\red$ by continuity (and \cite[\href{https://stacks.math.columbia.edu/tag/056V}{056V}]{stacks-project}).
  The morphism $\pi$ thus factors (uniquely) through the subsheaf $\affineline{\deRham,\lc}$ of $\affineline{\deRham}$ whose $S$\nobreakdashes-points are \emph{locally constant} functions on $S_\red$.

  Since de Rham cohomology is invariant under nilpotent extensions, we moreover have $\affineline{\deRham,\lc} \simeq \stackforms^{0,\closed}(-,0)_\deRham \simeq \stackforms^{0,\closed}(-,0) \colon S \mapsto \cohomology_\deRham^0(S)$.
  The induced morphism
  \[
    \stackforms^{2,\exact}(-,-1) \to \affineline{\deRham,\lc} \simeq \stackforms^{0,\closed}(-,0)
  \]
  is then a retract of $\stackforms^{0,\closed}(-,0) \to \stackforms^{2,\exact}(-,-1)$.
\end{proof}

\begin{example}
  As explained in \cref{symmetryderivedcriticalocusexplicit}, for any $\LG$\nobreakdashes-pair $f \colon U \to \affineline{}$, the derived critical locus comes with a canonical $(-1)$\nobreakdashes-shifted exact symplectic form $\lambda$.
  The associated locally constant function $c_\lambda$ is nothing but the restriction of $f$ to the reduced critical locus $\Crit(f)_\red$.
  It therefore corresponds to the critical values of the given function.
\end{example}

\begin{definition}\label{locallyconstantfunctionunderlyingexact}
  Let $X$ be a derived Artin stack. For any $(-1)$\nobreakdashes-shifted exact $2$\nobreakdashes-form $\lambda$ on $X$, we call the \emph{critical value function} of $\lambda$ its image $c_\lambda \in \cohomology_\deRham^0(X)$ via the retract of \cref{exactformssplitting}.
\end{definition}

\begin{corollary}\label{canonicalexactstructure}
  Any $(-1)$\nobreakdashes-shifted symplectic derived Artin stack admits a unique exact structure with vanishing critical value function.
\end{corollary}

\subsection{Relative de Rham algebras and completions}
\label{sectionrelativederham}

We now recall the construction of the \emph{relative} de Rham algebra. This notion is required to define isotropic and Lagrangian fibrations and to state \cref{Bhatt}.

\begin{construction}\label{constructionrelativederham}
  Following \cite[1.3.12, 1.3.17]{MR3653319}, \cref{reminderderhamalgebra} extends to the relative setting: for any morphism of affine derived schemes $S \to S'$ we have a \emph{relative} de Rham algebra given by the mixed graded algebra (see \cref{remindermixedgraded})
  \begin{equation}\label{definitionrelativederham}
    \DR(S/S') \coloneqq  \Sym_{S}(\cotangent_{S/S'}[1]) \in \mixedgradedalgebras_\basefield
  \end{equation}
  where $\cotangent_{S/S'}[1]$ is in weight 1, equipped with the mixed graded structure given by the relative de Rham differential.
  By \cite[Prop.\,1.3.16]{MR3653319} $\DR(S/S')$ is in fact a $\structuresheaf_{S'}$\nobreakdashes-mixed graded algebra, \ie the relative de Rham differential is $\structuresheaf_{S'}$\nobreakdashes-linear.  Following the discussion \cite[before Lem.\,1.3.18]{MR3653319}, \cref{definitionrelativederham} is functorial on morphisms of affine derived schemes.
  As it satisfies étale descent in both source and target, it can be extended to morphisms of derived stacks (see \cite[2.4.2]{MR3653319})
  \[
    \DR \colon \Fun(\Delta^1, \dSt_\basefield)^\op \to \mixedgradedalgebras_\basefield.
  \]
  If moreover $X \to Y$ admits a relative cotangent complex, we have an equivalence of graded complexes
  \[
    \DR(X/Y) \simeq \bigoplus_{p\geq 0} \Gamma(X, \Sym_{\structuresheaf_X}(\cotangent_{X/Y}[1])).
  \]
  Further, for any $X \to Y \to Z$, the commutative squares
  \[
    \begin{tikzcd}[row sep=small]
      X \ar{r} \ar{d} & Y \ar{d} \\
      Y \ar{r} & Y
    \end{tikzcd}
    \hspace{1cm}
    \begin{tikzcd}[row sep=small]
      X \ar{r} \ar{d} & X \ar{d} \\
      Y \ar{r} & Z
    \end{tikzcd}
  \]
  induce morphisms
  \begin{align}\label{functionstorelativeforms}
    \mathllap{\globalsections(Y, \structuresheaf_Y) \simeq {}} \DR(Y/Y) &\to \DR(X/Y) \rlap{ and}
    \\ \label{formstorelativeforms}
    \DR(X/Z) &\to \DR(X/Y).
  \end{align}
\end{construction}

\begin{notation}\label{relativeclosedpforms}
  Let $X \to Y$ be a morphism of derived stacks. Similarly to \cref{constructiontruncatedderhamcomplex,truncatedabovederhamcomplex}, we define the complex of \emph{relative closed} $p$\nobreakdashes-forms
  \[
    \derhamcomplex_{X/Y}^{\geq p} \coloneqq |\DR(X/Y)|^{\geq p}[-2p]
  \]
  and the complex of \emph{relative} exact forms
  \[
    \derhamcomplex_{X/Y}^{\leq p} \coloneqq \cofiber\left(\derhamcomplex_{X/Y}^{\geq p+1} \to \derhamcomplex_{X/Y}^{\geq 0}\right).
  \]
  By construction, they come with functorial morphisms summed up in the following diagram, where $X \to Y \to Z \in \dSt_\basefield$
  \begin{equation}
    \label{functorialitiesderham}
    \begin{tikzcd}[row sep=small]
      & \derhamcomplex_{X/Y}^{\geq p+1} \ar{d} & \derhamcomplex_{X/Z}^{\geq p+1} \ar{l}{\cref{formstorelativeforms}} \ar{d}
      \\
      \globalsections(Y, \structuresheaf_Y) \ar{r}{\cref{functionstorelativeforms}} & \derhamcomplex_{X/Y}^{\geq 0} \ar{d} & \derhamcomplex_{X/Z}^{\geq 0} \ar{d} \ar{l}{\cref{formstorelativeforms}}
      \\
      & \derhamcomplex_{X/Y}^{\leq p} & \derhamcomplex_{X/Z}^{\leq p}\rlap{.} \ar{l}{\cref{formstorelativeforms}}
    \end{tikzcd}
  \end{equation}
\end{notation}

\begin{definitions}{}
  \item\label{stackofrelativeclosedforms} We define the stack of $n$\nobreakdashes-shifted relative closed $p$\nobreakdashes-forms, denote by $\stackforms^{p, \closed}_\relative(- , n)$, as the stack $\Fun(\Delta^1, \dSt_\basefield)^\op \to \inftygpd$ sending a morphism of stacks $X \to Z$ to the Dold-Kan constructions
  \[
    \stackforms^{p, \closed}_\relative(X/Z, n) \coloneqq \DK(\derhamcomplex_{X/Z}^{\geq p}[p+n]).
  \]
  \item\label{stackofrelativeexactpforms} We define the stack $\stackforms^{p,\exact}_\relative(-,n)$ of $n$\nobreakdashes-shifted relative exact $p$\nobreakdashes-forms:
  \[
    \stackforms^{p,\exact}_\relative(-,n) \coloneqq \fiber\left(\stackforms^{p,\closed}_\relative(-,n) \to \stackforms^{0,\closed}_\relative(-,n+p)\right).
  \]
\end{definitions}

Finally, we can state Bhatt's theorem, which will play a crucial role in our
description of derived critical loci:
\begin{theorem}\cite[Prop.\,4.16, Ex.\,4.18]{1207.6193}\label{Bhatt}
  Let $i \colon S\hookrightarrow \formalU$ be a closed immersion with $S$ a derived affine scheme, $\formalU$ a smooth affine formal scheme and assume that $i$ is an isomorphism on the reduced subschemes. Then the inclusion of relatively constant functions $\cref{functionstorelativeforms}\colon\structuresheaf_\formalU \to |\DR(S/\formalU)|$ is an equivalence in $\cdga$
  \[
    \structuresheaf_\formalU \simeq |\DR(S/\formalU)| = \derhamcomplex_{S/\formalU}^{\geq 0}.
  \]
\end{theorem}

\subsection{Lagrangian fibrations and  distributions}
\label{sectionLagrangianDistributions}
In this section, we introduce the key notion of Lagrangian fibrations. It relies on the notion of (Lagrangian) distribution that we will discuss first.

\subsubsection{Distributions}
\begin{definition}
  \label{underlyingisotropicdistribution}
  Let $X$ be a derived Artin stack locally of finite presentation. A \emph{distribution} on $X$ is the data of a perfect complex $\rmE \in \Perf(X)$ together with a map of perfect complexes $\ell \colon \rmE\to \tangent_X$. The \emph{normal complex} $\rmN_\ell$ of $\ell$ is its cofiber in $\Perf(X)$, so we have by definition a fiber-cofiber sequence
  \begin{equation}\label{normalofdistribution}
    \begin{tikzcd} \rmE \ar{r}{\ell} & \tangent_X \ar{r} & \rmN_\ell. \end{tikzcd}
  \end{equation}
  A morphism of distributions is a commutative triangle in $\Perf(X)$
  \[
    \begin{tikzcd}[row sep={between origins,1.1em}]
      E' \ar{dd}\ar{rd}{\ell'} & \\
      &\tangent_X.    \\
      E \ar{ur}[swap]{\ell}&
    \end{tikzcd}
  \]
  Distributions on $X$ form an \icategory formally defined as $\Dist(X) \coloneqq \Perf(X)_{/\tangent_X}$.
\end{definition}

\begin{example}
  \label{fiberdistribution}
  Let $f \colon X \to Y$ be a morphism of derived Artin stacks of finite presentation.
  Then the relative tangent $\ell_{X/Y} \colon \tangent_{X/Y} \to \tangent_X$ defines a distribution on $X$ with normal complex $f^\ast \tangent_Y$.
  We call it the \emph{fiberwise distribution}.
  Given a morphism of stacks under $X$
  \[ \begin{tikzcd}[row sep={between origins,1.5em}]
      & Y_1\ar{dd}{\rho} \\
      X \ar{ur}{f_1} \ar{rd}[swap]{f_2} & \\
      &  Y_2
    \end{tikzcd}
  \]
  we obtain an obvious morphism of distributions $\tangent_{X/\rho}\colon(\tangent_{X/Y_1},\ell_1)\to (\tangent_{X/Y_2}, \ell_2)$
  \begin{equation}\label{diagrammapinducedatlevelofdistributionsbyuniversalproperty}
    \begin{tikzcd}
      \arrow[rounded corners, "\ell_{1}", to path={|- ([yshift=1em]\tikztotarget.north)[near end]\tikztonodes -- (\tikztotarget)}]{rr} \tangent_{X/Y_1} \tikzcart\arrow[r, dashed, "\tangent_{X/\rho}"] \ar{d} & \ar{d} \tangent_{X/Y_2} \tikzcart\ar{r}{\ell_2} & \tangent_X \ar{d}
    \\
      0\ar{r} & i^\ast \tangent_{Y_1/Y_2} \ar{r} \ar{d} \tikzcart & \rmN_{\ell_1} \mathrlap{{} \coloneqq i^\ast \tangent_{Y_1}} \ar{d}{i^\ast \tangentmap{\rho}}
    \\
      & 0 \ar{r} & \rmN_{\ell_2} \mathrlap{{} \coloneqq j^\ast \tangent_{Y_2}}
    \end{tikzcd}
  \end{equation}
  where $\tangentmap{\rho}$ is the derivative of $\rho$.
\end{example}

\begin{definitions}{Let $\ell \colon \rmE\to \tangent_X$ be a distribution.}
  \item \label{definitionsmoothdistribution} We say that it is \emph{smooth} if $\rmE$ is a perfect complex in Tor-amplitude $-1$, \ie there exists a finite rank vector bundle $\rmM$ over $X$ together with a quasi-isomorphism $\rmM^\vee[-1]\simeq \rmE$.
  \item \label{definitiontransverselysmooth} We say that it is \emph{transversely smooth} if the normal complex $\rmN_\ell$ is a vector bundle. We will then refer to $\rmN_\ell$ as the normal bundle of the distribution.
\end{definitions}

\begin{example}
  \label{distributiondCrit}
  Let $U$ be a smooth affine $\basefield$\nobreakdashes-scheme with a function $f \colon U\to \mathbb{A}^1$ and $i \colon \dCrit(U,f)\hookrightarrow U$ the inclusion of the derived critical locus. Then the pullback square \cref{lagrangianconditiondCritHessian} provides an identification $\tangent_{\dCrit(U,f)/U}\simeq i^\ast \cotangent_U[-1]$ and by definition $\rmN_\ell=i^\ast \tangent_U$.
  Since $U$ is smooth the fiberwise distribution is simultaneously smooth and transversely smooth.
\end{example}

\subsubsection{Isotropic and Lagrangian distributions}
\begin{definition}\label{isotropicdistribution}
  Let $\ell \colon \rmE\to \tangent_S$ be a distribution on a $(-1)$\nobreakdashes-symplectic derived affine scheme $S$. Let $\ell^\vee \colon \cotangent_X\to \rmE^\vee$ be its $S$\nobreakdashes-linear dual map.
  An \emph{isotropic structure} on $\ell \colon \rmE\to \tangent_S$ is the data of a null-homotopy $\eta$ for the composition
  \begin{equation}\label{restrictionofunderlyingformtofoliation}
    \begin{tikzcd}[column sep=3em]
      \structuresheaf_S \ar{r}{\cref{underlying2form22}} & (\cotangent_S \wedge \cotangent_S )[-1] \ar{r}{\ell^\vee \wedge \ell^\vee} & (\rmE^\vee \wedge \rmE^\vee)[-1]=\Sym^2_S(\rmE^\vee[-1])[1] \rlap{.}
    \end{tikzcd}
  \end{equation}
  This is equivalent to the data of a \underline{self-dual} 2-cell $\eta$ rendering the commutativity of the square
  \begin{equation}\label[square]{square-lagdist}
    \begin{tikzcd}[column sep=small]
      \rmE \ar{rr}{\ell} \ar{d} && \tangent_S \ar{r}{\sim}[swap]{\omega_S} & \cotangent_S[-1] \ar{d}{\ell^\dual[-1]} \\
      0 \ar[phantom]{urrr}[description,sloped]{\sim}[above]{\scriptstyle \eta} \ar{rrr} &&& \rmE^\vee[-1]\rlap{.}
    \end{tikzcd}
  \end{equation}
  This in turn induces a morphism $\alpha_\eta \colon \rmN_\ell \to \rmE^\vee[-1]$ together with a commutative square
  \[
    \begin{tikzcd}
      \tangent_S \ar{d}\ar{r}{\sim}[swap]{\omega_S} &\cotangent_S[-1] \ar{d}{\ell^\vee[-1]}\\
      \rmN_\ell \ar{r}{\alpha_\eta}&\rmE^\vee[-1]\rlap{.}
    \end{tikzcd}
  \]
\end{definition}

\begin{definition}\label{categoryisotropicdistributionsconstruction}
  For any $(-1)$\nobreakdashes-symplectic derived affine scheme $S$, we denote by $\IsotDist(S)$ the \icategory of isotropic distributions on $S$. It is formally defined as the fiber product
  \[
    \begin{tikzcd}[column sep=small]
      \IsotDist(S) \tikzcart \ar{r}{} \ar{d} &  \Dist(S)\ar{d}  \\
      \left(\Perf(S)\right)^\op  \ar{r} & \left(\undercat{\structuresheaf_S}{\Perf(S)}\right)^\op
    \end{tikzcd}
  \]
  where $\Perf(S)\to \undercat{\structuresheaf_S}{\Perf(S)}$ is the pointing by the zero map, and the vertical map sends a distribution $\rmE \to \tangent_S$ to the map \cref{restrictionofunderlyingformtofoliation} induced by the symplectic structure.
\end{definition}

\begin{definition}\label{lagrangians}
  Let $(\ell \colon \rmE\to \tangent_S, \eta)$ be an isotropic distribution on $S$ as in \cref{isotropicdistribution}. We say that the isotropic distribution is \emph{Lagrangian} if the square \cref{square-lagdist} is coCartesian.
  We denote by $\LagDist_S \subseteq \IsotDist_S$ the full sub-\icategory of Lagrangian distributions.
\end{definition}

\begin{observation}\label{forlagrangianssmoothequalstransverselysmooth}
  It follows from the coCartesian property of the square \cref{square-lagdist} that a Lagrangian distribution is smooth if and only if it is transversely smooth.
\end{observation}

\begin{construction}\label{const-quadfromlagdist}\label{const:quadfromlagdist-lagrangian}
  Let $(\rmE,\ell,\eta)$ be an isotropic distribution over $S=\Spec(A)$.
  We will construct a natural symmetric $2$\nobreakdashes-form on the normal complex $\rmN_\ell$. Consider the fiber sequence $\begin{tikzcd}[cramped] \rmN^\dual_\ell \ar{r} & \cotangent_S \ar{r}[description]{\ell^\dual} & \rmE^\dual \end{tikzcd}$
  and denote by $\partial^\dual_\ell \colon \rmE^\dual \to \rmN^\dual_\ell[1]$ its boundary operator.
  The composite morphism
  \[
    \begin{tikzcd}
      (\cotangent_S \wedge \cotangent_S)[-1] \ar{r}{\ell^\dual \wedge \ell^\dual}
      & (\rmE^\dual \wedge \rmE^\dual)[-1] \ar{r}{\partial_\ell^\dual \wedge \partial_\ell^\dual}
      & (\rmN^\dual_\ell[1] \wedge \rmN^\dual_\ell[1])[-1] \simeq \Sym^2_S(\rmN^\dual_\ell)[1]
    \end{tikzcd}
  \]
  is then nullhomotopic. Using the isotropic structure $\eta$, we get a commutative diagram
  \[
    \begin{tikzcd}
      \structuresheaf_S \ar{r}{\omega_S} \ar{d} & (\cotangent_S \wedge \cotangent_S)[-1] \ar{r} \ar{d}{\ell^\dual \wedge \ell^\dual} & 0 \ar{d}
      \\
      0 \ar{r} \ar[phantom]{ur}[description,sloped]{\sim}[above]{\scriptstyle \eta} & (\rmE^\dual \wedge \rmE^\dual)[-1] \ar{r}{\partial_\ell^\dual \wedge \partial_\ell^\dual} & \Sym^2_S(\rmN^\dual_\ell)[1].
    \end{tikzcd}
  \]
  We get a loop $\structuresheaf_S \to \Sym^2_S(\rmN^\dual_\ell)[1][-1] \simeq \Sym^2_S(\rmN^\dual_\ell)$ providing the announced symmetric $2$\nobreakdashes-form on the normal bundle
  \[
    \rmB_\eta \colon \Sym^2_S(\rmN_\ell) \to \structuresheaf_S
  \]
  with induced self-dual pairing $\rmH_\eta$ given by the diagonal composition
  \[
    \begin{tikzcd}
      \ar{dr}{\rmH_\eta} \rmN_\ell \ar{r}{\alpha_\eta} \ar{d}{\partial_\ell} &\rmE^\dual[-1]\ar{d}{\partial_\ell^\vee}\\
      \rmE[1] \ar{r}{\alpha_\eta^\dual} & \rmN_\ell^\vee
    \end{tikzcd}
  \]
  Finally, assuming further that the isotropic structure is Lagrangian, the above square identifies $\rmH_\eta$ with the boundary operator $\partial_\ell$.
\end{construction}

\subsubsection{Lagrangian fibrations}
We now discuss the notions of Lagrangian fibrations used in this work.
\begin{definition}\label{definitionisotropic}
  Let $S=\Spec(A)$ an affine derived scheme equipped with a $(-1)$\nobreakdashes-shifted symplectic structure.
  Let $i \colon S\to F$ with $F$ a derived stack with cotangent complex.
  A structure of \emph{isotropic fibration} on $i$ is the data of a null-homotopy the map in $\dgMod_\basefield$
  \begin{equation}
    \label{underlying2form2isotropic}
    \basefield\to \derhamcomplex^{\geq 2}_S[1]\to \derhamcomplex^{\geq 2}_{S/F}[1].
  \end{equation}
\end{definition}

\begin{example}\label{distributionassociatedtofibration}
  Let $S$ a $(-1)$\nobreakdashes-symplectic derived affine scheme and let $i \colon S\to F$ be an isotropic fibration as in \cref{definitionisotropic}.
  Then the fiberwise distribution of \cref{fiberdistribution} inherits a natural isotropic structure induced by the naturality of the maps \cref{underlying2form}:
  \[
    \begin{tikzcd}
      \basefield \ar{r}{\omega} & \ar{r}{\cref{underlying2form2isotropic}} \derhamcomplex^{\geq 2}_S[1] \ar{d}{\cref{underlying2form}} & \derhamcomplex^{\geq 2}_{S/F}[1] \ar{d}{\cref{underlying2form}}
      \\
      & (\bigwedge^2 \cotangent_S)[-1] \ar{r}{\ell^\vee\wedge_S\ell^\vee} & (\bigwedge^2 \cotangent_{S/F})[-1].
    \end{tikzcd}
  \]
\end{example}

\begin{definition}\label{definition-oflagrangianfibration}
  Let $i \colon S\to F$ be an isotropic fibration as in \cref{definitionisotropic}.
  We say that it is a \emph{Lagrangian fibration} if the associated fiberwise isotropic distribution (of \cref{distributionassociatedtofibration}), is a Lagrangian distribution.
  That is to say if the induced map $\alpha_\eta \colon i^\ast\tangent_F\to\cotangent_{S/F}[-1]$ is an equivalence.
\end{definition}

The most important example of Lagrangian fibration in this paper is the following:

\begin{proposition}\cite[Thm.\,3.5, Rem.\,3.12]{Grataloup2020}\label{lagrangianfibrationdCrit}
  Let $U$ and $f$ be as in \cref{distributiondCrit}. Then the closed immersion $i \colon \dCrit(U,f)\hookrightarrow U$ has a canonical structure of Lagrangian fibration.
\end{proposition}

\begin{observation}\label{bilinearformisthehessian}
  \Cref{const:quadfromlagdist-lagrangian} shows that the normal bundle of a Lagrangian distribution has an associated symmetric bilinear pairing $\rmH_\eta$. In the case of \cref{lagrangianfibrationdCrit}, $\rmH_\eta$ is the Hessian matrix and the Lagrangian condition is the fact the square \cref{lagrangianconditiondCritHessian} is coCartesian.
\end{observation}

We will study this example in more detail in \cref{sectioncharacterizationsformaldarboux} below. To conclude this section with summarize the functoriality properties for Lagrangians distributions:

\begin{construction}\label{stackisotropicandlagrangiandistributions}
  Let $X$ be a $(-1)$\nobreakdashes-shifted symplectic derived Deligne--Mumford stack.
  Since for any étale map $u \colon S\to X$ we have $u^\ast\tangent_X\simeq \tangent_S$, the pullback of distributions is well-defined. Moreover, the condition of being Lagrangian is stable under such pullbacks. In this case, the assignments of distributions (\cref{underlyingisotropicdistribution}), isotropic distributions (\cref{categoryisotropicdistributionsconstruction}) and Lagrangians distributions (\cref{lagrangians}), are encoded by \ifunctors
  \begin{equation}\label{stacksoflagrangiansdistributionsformula}
    \stackDist_X, \stackLagDist_X, \stackIsotDist_X \colon  (X_\et^\daff)^\op \to \inftycats.
  \end{equation}
  Moreover, these satisfy hyperdescent in $X_\et^\daff$ due to descent for perfect complexes.
\end{construction}

\section{The stack of Darboux charts in \texorpdfstring{$(-1)$}{(-1)}-symplectic geometry}
\label{SectionDarboux}

In this section we construct the Darboux moduli functor discussed in \cref{sectionmainresult}.

\subsection{Derived critical loci and formal completions}\label{dCritandformalcompletions}
We start with the observation that the derived critical locus with the exact symplectic form of \cref{symmetryderivedcriticalocusexplicit}, only depends on the formal completion -- \cref{dCritonlydependsformal} below.

\begin{construction}  \label{dCritofformalscheme}
  Let $\formalU$ be a smooth (formally smooth and topologically of finite presentation) formal scheme over $\basefield$ with a function $f \colon \formalU \to \affineline{}$.
  Its cotangent stack $\cotangentstack \formalU$ admits an exact $0$\nobreakdashes-shifted symplectic structure, classified by a morphism of derived stacks
  \[
    \cotangentstack \formalU \to \stackforms^{2,\closed}(-,0).
  \]
  The de Rham differential $df$ determines a section $df \colon \formalU\to \cotangentstack \formalU$ and we can form the derived intersection in derived stacks
  \begin{equation}\label{fiberproductdCritformal}
    \begin{tikzcd}
      \dCrit(\formalU,f) \ar{r}{i} \ar{d}{i} \tikzcart & \formalU \ar{d}{df}\\
      \formalU \ar{r}{0} & \cotangentstack \formalU\rlap{.}
    \end{tikzcd}
  \end{equation}
  As in \cref{symmetryderivedcriticalocusexplicit}, the derived stack $\dCrit(\formalU,f)$ inherits a $(-1)$\nobreakdashes-shifted symplectic structure because both $0$ and $df$ are Lagrangians in $\cotangentstack \formalU$.
  Moreover, the fact that $df$ is \emph{exact} and the $0$\nobreakdashes-shifted form in $\cotangentstack \formalU$ is exact due to the Liouville form, implies once again that $\dCrit(\formalU,f)$ has a canonical exact structure.
\end{construction}

\begin{construction}
  \label{morphismofcriticalcharts}
  Let $(\formalU,f)$ and  $(\formalV,g)$ be smooth formal schemes over $\basefield$ equipped with functions $f \colon \formalU \to \affineline{}$ and $g \colon \formalV \to \affineline{}$.
  Let $u \colon \formalU \to \formalV$ be a map such that $g\circ u=f$.
  In this case we extract a correspondence between the associated derived critical loci:
  \begin{equation}
    \label{correspondancecriticallocus}
    \begin{tikzcd}
      &
      \mathllap{\dCrit(g)|_\formalU \coloneqq {}} \dCrit(g)\times_{\formalV} \formalU \ar{dl}[swap]{p} \ar{dr}{q} \\
      \dCrit(f) &
      &
      \dCrit(g)
    \end{tikzcd}
  \end{equation}
  obtained from the commutative diagram
  \begin{equation}
    \label{constructionoflagrangiancorrespondance}
    \begin{tikzcd}[column sep=small]
      &
      \dCrit(f) \ar{dl} \ar{dd} &
      &
      \dCrit(g)|_\formalU \ar[dashed]{ll}[swap]{p} \ar{rr}{q} \ar{dl} \ar{dd} &
      &
      \dCrit(g) \ar{dl} \ar{dd} \\
      \formalU \ar{dd}{df} &
      &
      \formalU \ar[equals, crossing over]{ll} \ar[crossing over]{rr}[pos=0.4]{u} &
      &
      \formalV \\
      &
      \formalU \ar{dl} &
      &
      \formalU \ar[equals]{ll} \ar{rr} \ar{dl} &
      &
      \formalV \ar{dl} \\
      \cotangentstack \formalU &
      &
      \cotangentstack \formalV \times_{\formalV} \formalU \ar{ll}{\mathrm{D}_u} \tikzcart \ar{d} \ar{rr} \ar[from=uu, crossing over, near start, "dg|_\formalU"] &
      &
      \cotangentstack \formalV \ar{d} \ar[from=uu, crossing over, near start, "dg"] \\
      &
      &
      \formalU \ar{rr}{u} &
      {} &
      \formalV\rlap{.}
    \end{tikzcd}
  \end{equation}
  The correspondence \cref{correspondancecriticallocus} is Lagrangian in the sense of \cite[Def.\,1.7]{MR3940792}: if $\omega_{f}$ denotes the symplectic form on $\dCrit(f)$ and $\omega_g$ on $\dCrit(g)$ then there is a homotopy between
  \begin{equation}
    \label{canonicallagrangiancorrespondance}
    p^\ast(\omega_f)\sim q^\ast(\omega_g)\rlap{.}
  \end{equation}
  In fact, not only the symplectic structures are homotopic but the exact structures also: this is a direct consequence of the fact that the correspondence of cotangent bundles
  \begin{equation}
    \label{lagrangiancorrespondancecotangentbundle}
    \begin{tikzcd}
      \cotangentstack \formalU &
      \displaystyle \cotangentstack \formalV \times_{\formalV} \formalU \ar{l} \ar{r} &
      \cotangentstack \formalV
    \end{tikzcd}
  \end{equation}
  is exact Lagrangian with respect to the standard symplectic forms, since the  pullback of the two Liouville forms to $\cotangentstack \formalV \times_{\formalV} \formalU$ agree. See \cite[Ex.\,2.4.9(a)]{2108.02473} for more details.
\end{construction}

\begin{construction}
  \label{comparingformalandnonformal}
  Let $U$ be smooth $\basefield$\nobreakdashes-scheme with a function $f\colon U \to \affine^1$. Denote by $u \colon \formalU \to U$ the formal completion of $U$ along $\dCrit(U,f)$ and $\widehat{f} \coloneqq f \circ u$.  We apply \cref{morphismofcriticalcharts} to the morphism $u \colon (\formalU,\widehat{f}) \to (U,f)$. Since $u$ is formally étale, the left morphism in the exact Lagrangian correspondence \cref{lagrangiancorrespondancecotangentbundle} is an isomorphism.
  Therefore, so is the morphism $p$ in the induced exact Lagrangian correspondence \cref{correspondancecriticallocus}
  \begin{equation}
    \label{correspondancecriticallocusformalnonformal}
    \begin{tikzcd}
      & \dCrit(\formalU, \widehat{f}) & \dCrit(U,f)\times_{U} \formalU \ar{l}[swap]{p}{\sim} \ar{r}{q} & \dCrit(U,f)\rlap{.}
    \end{tikzcd}
  \end{equation}
  Finally, as in top right face of diagram \cref{constructionoflagrangiancorrespondance}, we extract a pullback square of derived stacks
  \begin{equation}\label{pullbacketaledCritformal}
    \begin{tikzcd}
      \dCrit(\formalU, \widehat{f}) \ar{r}{q} \ar{d}  \tikzcart & \dCrit(U,f) \ar{d}\\
      \formalU \ar{r}{u}& U
    \end{tikzcd}
  \end{equation}
  and since $u$ is étale, so is $q$.
\end{construction}

\begin{proposition}
  \label{dCritonlydependsformal}  Let $U$ be a smooth $\basefield$\nobreakdashes-scheme with a function $f \colon U\to \affine^1$ and let $\formalU$ be the formal neighborhood of $\dCrit(U,f)$ in $U$ equipped with the restriction of $f$ to $\formalU$. Then the map $q$ of \cref{pullbacketaledCritformal}
  \[
    \begin{tikzcd}
      \dCrit(\formalU, \widehat{f}) \ar{r}{q} & \dCrit(U,f)
    \end{tikzcd}
  \]
  is an equivalence of $(-1)$\nobreakdashes-shifted exact symplectic derived schemes
\end{proposition}

\begin{proof}
  By definition, both $\formalU$ and $\dCrit(U,f)$ have the same underlying reduced scheme. The Cartesian square \cref{pullbacketaledCritformal} implies that so does $\dCrit(\formalU, \widehat{f})$ -- see \cref{extractingreducedsubstack}.
  Since $q$ is étale and an isomorphism on reduced schemes, the result follows from \cref{reducedeqplusetaleisiso}, using \cref{correspondancecriticallocusformalnonformal} for the exact symplectic compatibility.
\end{proof}

\begin{example}
  We have exact symplectic equivalences
  \[
    \Spec(\basefield[x]/(3x^2))=\dCrit(\affine^1, x^3)\simeq \dCrit(\widehat{\affine}^1_0, x^3).
  \]
\end{example}

\subsection{Derived critical loci as Lagrangian fibrations}
\label{sectioncharacterizationsformaldarboux}

We have seen in \cref{symmetryderivedcriticalocusexplicit} that $\dCrit(U,f)$ admits an \emph{exact structure}. In this section we use it to reformulate the structure of Lagrangian fibration of \cref{lagrangianfibrationdCrit} and establish the fundamental result -- \cref{exactequivalence} -- showing that for a $(-1)$\nobreakdashes-shifted exact symplectic derived affine scheme $S$, the data of a Lagrangian fibration is equivalent to the data presenting $S$ as a derived critical locus.

\begin{observation}\label{constructionLiouville}
  Let $S$ be a $(-1)$\nobreakdashes-shifted symplectic affine derived scheme and let $i \colon S \hookrightarrow \formalU$ be closed immersion, with $\formalU$ a smooth affine formal scheme of finite presentation and $i$ an isomorphism on the underlying reduced schemes, \ie $i_\red\colon S_\red\simeq \formalU_\red$.
  Then the data of an exact structure provides by definition a lifting of $\omega_S$:
  \begin{equation}\label{isotropicinexactcase1}
    \begin{tikzcd}[row sep={between origins,1.8em}]
      & \derhamcomplex^{\leq 1}_S \ar[r]\ar[dd] & \derhamcomplex^{\leq 1}_{S/\formalU} \ar[dd]
      \\
      \basefield \ar{ur}[near end]{\lambda} \ar{dr}[swap, near end]{\omega_S} \tikzhomotopy[pos=.65]{r}{\alpha} & {}
      \\
      & \derhamcomplex^{\geq 2}_S[1] \ar{r}{\cref{functorialitiesderham}} & \derhamcomplex^{\geq 2}_{S/\formalU}[1]\rlap{.}
    \end{tikzcd}
  \end{equation}
\end{observation}

\begin{observation}\label{constructionLiouvilleisotropic}
  Assume $S \to \formalU$ as in \cref{constructionLiouville}. Rotating the cofiber sequence \cref{cofibersequencesfortruncations} for the \emph{relative} de Rham complex of $S$ over $\formalU$ and using Bhatt's equivalence (see \cref{Bhatt}), we find a cofiber sequence
  \begin{equation}
    \label{isotropicfibersequenceforrelative}
    \structuresheaf_\formalU \simeq \derhamcomplex^{\geq 0}_{S/\formalU}\to \derhamcomplex^{\leq 1}_{S/\formalU}\to \derhamcomplex^{\geq 2}_{S/\formalU}[1]\rlap{.}
  \end{equation}
  Combining the diagrams \cref{isotropicinexactcase1} and \cref{isotropicfibersequenceforrelative}, we find that in the presence of an exact structure $\alpha$, the data of an isotropic fibration $\eta$ on the morphism $i \colon S\hookrightarrow \formalU$ amounts to a lift $f \colon \basefield \to \structuresheaf_\formalU$ and a 2-cell $\delta$ rendering the commutativity of the diagram
  \begin{equation}\label{isotropicinexactcase2}
    \begin{tikzcd}[row sep={between origins,1.8em}]
      && \structuresheaf_\formalU \ar{dd} \tikzhomotopy{ddl}{\delta}
      \\ \\
      & \derhamcomplex^{\leq 1}_S \ar[r]\ar[dd] & \derhamcomplex^{\leq 1}_{S/\formalU} \ar[dd]
      \\
      \basefield \ar{ur}[near end]{\lambda} \ar{dr}[swap,near end]{\omega_S} \ar[bend left]{uuurr}{f} \tikzhomotopy[pos=.65]{r}{\alpha} & {}
      \\
      & \derhamcomplex^{\geq 2}_S[1] \ar{r}{\cref{functorialitiesderham}} & \derhamcomplex^{\geq 2}_{S/\formalU}[1]\rlap{.}
    \end{tikzcd}
  \end{equation}
\end{observation}

\begin{definition}\label{Liouvillecomplex}
  Let $S$ be a derived affine scheme, $\formalU$ a smooth affine formal scheme and $i \colon S \hookrightarrow \formalU$ a closed immersion such that $i_\red \colon S \to \formalU_\red$ is an isomorphism.
  We define the \emph{Liouville complex} of $S\to \formalU$, $\Liouv(S/\formalU)$ as the homotopy pullback in $\dgMod_\basefield$
  \begin{equation}\label{isotropicinexactcase33}
    \begin{tikzcd}
      \Liouv(S/\formalU) \ar[r] \ar{d} \tikzcart & \structuresheaf_\formalU \ar{d}
      \\
      \derhamcomplex^{\leq 1}_S \ar[r] & \derhamcomplex^{\leq 1}_{S/\formalU}\rlap{.}
    \end{tikzcd}
  \end{equation}
\end{definition}

\begin{observation}\label{Liouvilleisotropic}
  It follows from \cref{constructionLiouvilleisotropic} and \cref{Liouvillecomplex} that having fixed an exact structure $\alpha$ on $S$, the data of an isotropic fibration on $i \colon S\hookrightarrow \formalU$ amounts to a lift $\kappa$:
  \begin{equation}\label{isotropicstructureviaLiouville}
    \begin{tikzcd}[row sep={between origins,3.5em}]
      & \ar{d} \Liouv(S/\formalU) \\
      \basefield \ar{ru}[name=K]{\kappa} \ar{r}[description]{\lambda} \ar{dr}[swap, pos=0.6, name=Omega]{\omega_S} & \ar{d} \derhamcomplex^{\leq 1}_S \tikzhomotopy{K}{\delta} \tikzhomotopy{Omega}{\alpha} \\
      & \derhamcomplex_S^{\geq 2}[1]\rlap{.}
    \end{tikzcd}
  \end{equation}
\end{observation}

The Liouville complex has a more explicit description, justifying to some extent its name:

\begin{lemma}\label{constructionexplicitLiouvillestructure}\label{dataLiouvillexplicitformula}
  Let $S$ be a derived affine scheme, $\formalU$ a smooth affine formal scheme and $i \colon S \hookrightarrow \formalU$ a closed immersion such that $i_\red \colon S \to \formalU_\red$ is an isomorphism.
  The complex $\Liouv(S/\formalU)$ is naturally equivalent to the total space of the bicomplex
  \begin{equation}\label{isotropicinexactcase4}
    \Liouv(S/\formalU) \simeq \left[\begin{tikzcd}[cramped]
        \structuresheaf_\formalU \ar{r}{\dderham} & i^* \Omega_\formalU
      \end{tikzcd}\right].
  \end{equation}
  In particular, a map $\kappa \colon \basefield \to \Liouv(S/\formalU)$ amounts to the data of a pair $(f, \phi)$, where $f \in \structuresheaf_\formalU$ and $\phi \colon \dderham f \sim 0$ is a nullhomotopy of $\dderham f$ in $i^* \Omega_\formalU$.
\end{lemma}
\begin{proof}
  By construction and by \cref{Bhatt}, the square \cref{isotropicinexactcase33} is equivalent to
  \[
    \begin{tikzcd}
      \Liouv(S/\formalU) \ar{r} \ar{d} \tikzcart & \structuresheaf_\formalU \ar{d} \\
      \left[ \structuresheaf_S \to \cotangent_S \right] \ar{r} & \left[ \structuresheaf_S \to \cotangent_{S/\formalU} \right].
    \end{tikzcd}
  \]
  The result follows since the fiber of $\cotangent_S \to \cotangent_{S/\formalU}$ is by definition $i^* \Omega_\formalU$.
\end{proof}

\begin{corollary}\label{functionassociatedtoisotropicfibration}
  Assume further that $S$ is equipped with an exact $(-1)$\nobreakdashes-symplectic structure $\lambda \in \derhamcomplex_S^{\leq 1}$.
  The structure of an isotropic fibration on $i \colon S \hookrightarrow \formalU$ amounts to a function $f \in \structuresheaf_\formalU$, a nullhomotopy $\phi \colon \dderham f \sim 0$ in $i^* \Omega_\formalU$ and a homotopy
  \[
    \delta \colon (f,\phi) \sim \lambda \in \derhamcomplex_S^{\leq 1}.
  \]
\end{corollary}

\begin{observation}
  \label{LiouvillestructureondCrit}
  A particular case of \cref{constructionexplicitLiouvillestructure} is when $S=\dCrit(\formalU,f)$ equals the derived critical locus of a function $f$ on a smooth formal scheme $\formalU$. In this case, the equivalence \cref{isotropicinexactcase4} provides a more explicit construction for the exact structure of \cref{explicitexactstructureondCrit}, namely the morphism $\kappa$ in \cref{isotropicstructureviaLiouville} corresponds precisely to the data of the original function $f$ and the homotopy $\phi \colon \dderham(f)\sim 0$ in $i^\ast \cotangent_\formalU$ as in \cref{criticallocusasnullhomotopy}. Also in this case, the 2-cell $\alpha \circ \delta$ of \cref{isotropicstructureviaLiouville} rendering $\dderham(\phi) \sim \omega_S$ is the one of \cite[5.7-(a)]{MR3904157} obtained using explicit resolutions.
\end{observation}

\begin{proposition}\label{exactequivalence}
  Let $S$ be an exact $(-1)$\nobreakdashes-symplectic derived affine scheme.
  Let $i \colon S \hookrightarrow \formalU$ be an isotropic fibration.
  Denote by $f \colon \formalU \to \affineline{\basefield}$ the function induced by \cref{functionassociatedtoisotropicfibration}. Then
  \begin{enumerate}
    \item The morphism $i$ factors canonically through $\dCrit(\formalU, f)$ and the induced morphism $\widebar{\imath} \colon S \to \dCrit(\formalU,f)$ is compatible with the symplectic and exact structures.
    \item If the isotropic fibration is moreover Lagrangian, then $\widebar{\imath}$ is an equivalence.
  \end{enumerate}
\end{proposition}
\begin{proof}
  By \cref{functionassociatedtoisotropicfibration}, we have a nullhomotopy $\phi$ of $\dderham f$ in $i^* \Omega_\formalU$. It provides exactly the announced factorization $\widebar{\imath} \colon S\to \dCrit(\formalU,f)$ by definition of the critical locus.

  We need to show that $\widebar{\imath}$ is compatible with the symplectic and exact structures.
  Since those structures can be recovered from the section of the Liouville complex, both for $S$ and for $\dCrit(\formalU,f)$, it is enough to compare those sections. Consider the diagram:
  \[
    \begin{tikzcd}[row sep=tiny]
      & \Liouv(\dCrit(\formalU,f)/\formalU) \ar{dd} \ar{r}[swap]{\sim}{\cref{isotropicinexactcase4}} & \left[\structuresheaf_\formalU \to j^* \Omega_\formalU \right] \ar{dd}
      \ar[start anchor=center, end anchor=center, phantom]{ddl}{\scriptstyle \mathrm{(\tau)}}
      \\
      \basefield \ar{ur} \ar{dr} \ar[phantom]{r}[near end]{\scriptstyle \mathrm{(\sigma)}} & {}
      \\
      & \Liouv(S/\formalU) \ar{r}{\sim}[swap]{\cref{isotropicinexactcase4}} & \left[\structuresheaf_\formalU \to i^* \Omega_\formalU \right]
    \end{tikzcd}
  \]
  where $j \colon \dCrit(\formalU,f) \to \formalU$ denotes the structural embedding.
  The square $\mathrm{(\tau)}$ tautologically commutes, and so does the outer triangle. Thus $\mathrm{(\sigma)}$ commutes, proving our claim.

  If we further assume the isotropic structure to be Lagrangian,  we have $\cotangent_{S/\formalU}\simeq i^\ast\tangent_\formalU[1]$. Therefore, the fact that $\widebar{\imath}$ commutes with the maps to $\formalU$, together with the formulas for the cotangent complex of $\dCrit(\formalU,f)$ in \cref{explicitmodelcotangetcrit22}, imply that $\widebar{\imath}$ is étale. Since by construction $\widebar{\imath}$ is also an isomorphism of truncations, it is an equivalence of derived schemes (see \cref{reducedeqplusetaleisiso}).
\end{proof}

\begin{summary}\label{descriptionlocaldarboux}
  In conclusion, we have established that for any $(-1)$\nobreakdashes-shifted exact symplectic derived affine scheme $S$, there is a bijection between the following type of data:
  \begin{enumerate}
    \item The data of a smooth affine formal scheme $\formalU$ with a function $f$ and an \emph{exact symplectic equivalence} $S \simeq \dCrit(\formalU,f)$;
    \item The data of a closed immersion  $i \colon S \hookrightarrow \formalU$ with $\formalU$ a smooth affine formal scheme and inducing an isomorphism of reduced sub-schemes, together with a structure of Lagrangian fibration on $i$.
  \end{enumerate}
  Indeed, the map from (i) to (ii) is given by \cref{lagrangianfibrationdCrit}, observing that the arguments in \cite[Thm.\,3.5, Rem.\,3.12]{Grataloup2020} work mutatis mutandis for a formal scheme $\formalU$. The fact it is a bijection follows from \cref{exactequivalence}.
\end{summary}

\subsection{Moduli functor of Darboux data}\label{sectionconstructionformaldarboux}
We now construct the moduli functor classifying Liouville data as in \cref{Liouvillecomplex}. As a first step we introduce a stack classifying factorizations of morphisms between derived stacks:

\begin{construction}\label{stackofstacks}
  Consider the evaluation map
  \[
    \ev_{02} \colon \Fun(\Delta^2, \dSt_\basefield) \to \Fun(\Delta^1, \dSt_\basefield).
  \]
  The existence of fiber products guarantees that $\ev_{02}$ is a Cartesian fibration. We  denote its straightening by
  \[
    \stackfactor \colon \Fun(\Delta^1, \dSt_\basefield)^\op\to \inftycats
  \]
  informally sending $X\to Z$ to the \icategory of \emph{factorizations} $\dSt_{X/./Z}$ and  a morphism $u \colon (X' \to Z') \to (X \to Z)$ in $\Fun(\Delta^1, \dSt_\basefield)$ to the pullback functor  $u^\ast \colon \stackfactor(X \to Z)\to \stackfactor(X' \to Z')$ defined on objects by sending the factorization $X \to Y \to Z$ to $X'\to Y \times_Z Z'\to Z'$. Following (\cite[6.1.3.9-(3), 6.1.3.10]{lurie-htt}), $\stackfactor$ is a stack.
\end{construction}

\begin{construction}\label{Liouvillefunctor4mars}
  The functorialities of \cref{functorialitiesderham} induce \ifunctors
  \[
    \begin{tikzcd}[row sep=0pt, column sep=small]
      \Liouv \colon \Fun(\Delta^2, \dSt_\basefield) \ar{r} & \Fun(\Lambda^2_2, \dgMod_\basefield)^\op \ar{r}{\lim} & \dgMod_\basefield^\op
      \\
      \phantom{\Liouv \colon {}} (X \to Y \to Z) \ar[mapsto]{r} & \left(\globalsections(Y, \structuresheaf_Y) \to \derhamcomplex_{X/Y}^{\leq 1} \from \derhamcomplex_{X/Z}^{\leq 1}\right) \ar[mapsto]{r} & \Liouv(X/Y/Z).
    \end{tikzcd}
  \]
  and
  \[
    \begin{tikzcd}[row sep=0pt, column sep=small]
      \isotropic\colon \Fun(\Delta^2, \dSt_\basefield) \ar{r} & \Fun(\Lambda^2_2, \dgMod_\basefield)^\op \ar{r}{\lim} & \dgMod_\basefield^\op\\
      \phantom{\isotropic \colon {}} (X \to Y \to Z) \ar[mapsto]{r} & \left(0 \to \derhamcomplex_{X/Y}^{\geq 2}[1] \from \derhamcomplex_{X/Z}^{\geq 2}[1]\right) \ar[mapsto]{r} & \isotropic(X/Y/Z).
    \end{tikzcd}
  \]
  We note in particular that the Liouville complex $\Liouv(S/\formalU)$ that we introduced in \cref{Liouvillecomplex} is nothing but $\Liouv(S/\formalU/S_\derham)$.

  Since both functors  $\Liouv$ and $\isotropic$ have values in modules, their associated Dold-Kan \igroupoids carry abelian group structures. We denote the associated abelian group stacks\footnote{These are indeed stacks, being obtained as fiber product of stacks for the étale topology.}:
  \begin{align*}
    \Fun(\Delta^2,\dSt_\basefield)^\op               & \to \inftygpd,                               \\\vspace{0.5cm}
    \stackofliouvillesections \colon (X \to Y \to Z) & \mapsto \DK\left( \Liouv(X/Y/Z) \right),     \\\vspace{0.5cm}
    \stackofisotropicsections \colon
    (X \to Y \to Z)                                  & \mapsto \DK\left( \isotropic(X/Y/Z) \right).
  \end{align*}
  Using the projections $\Liouv(X/Y/Z) \to \derhamcomplex_{X/Z}^{\leq 1}$
  and $\isotropic(X/Y/Z) \to \derhamcomplex_{X/Z}^{\geq 2}[1]$ together with the commutativity of \cref{functorialitiesderham} and \cref{stackofclosedforms,stackofexactpforms}, we obtain a commutative square of symmetric monoidal natural transformations of abelian group stacks
  \begin{equation}\label{squareofabeliangroupstacks}
    \begin{tikzcd}[column sep={between origins, 9.5em}]
      \stackofliouvillesections \ar{r} \ar{d} & \stackofisotropicsections \ar{d}
      \\
      \stackforms^{2,\exact}_\relative(\ev_{02}(-),-1) \ar{r} & \stackforms^{2,\closed}_\relative(\ev_{02}(-), -1)
    \end{tikzcd}
  \end{equation}
  where
  \[
    \stackforms^{2,\exact}_\relative(-,-1), \stackforms^{2,\closed}_\relative(- , -1) \colon \Fun(\Delta^1, \dSt_\basefield)^\op\to \inftygpd
  \]
  are, respectively, the stacks of \underline{relative} exact and symplectic $(-1)$\nobreakdashes-shifted forms.
\end{construction}

\begin{observation}\label{squarenotCartesianingeneral}
  In general, the square \cref{squareofabeliangroupstacks} is not Cartesian. Indeed, the fiber of the map $\stackofliouvillesections \to \stackforms^{2,\exact}_{\relative}(\ev_{02}(-),-1)$ is given by the assignment $(X\to Y \to Z)\mapsto \DK\left(\fiber\left(\globalsections(Y, \structuresheaf_Y) \to \derhamcomplex_{X/Y}^{\leq 1}\right)\right)$ while the fiber of $\stackofisotropicsections \to \stackforms^{2,\closed}_{\relative}(\ev_{02}(-), -1)$ is given by $(X\to Y \to Z) \mapsto \DK\left( \derhamcomplex_{X/Y}^{\geq 2}\right)$. However, thanks to
  Bhatt's equivalence (see \cref{Bhatt}), the canonical comparison map
  \[
    \fiber\left(\globalsections(Y, \structuresheaf_Y) \to \derhamcomplex_{X/Y}^{\leq 1}\right)\to \derhamcomplex_{X/Y}^{\geq 2}
  \]
  becomes an equivalence when restricted to the full subcategory of $\Fun(\Delta^2, \dSt_\basefield)$ spanned by morphisms $i \colon S \hookrightarrow \formalU$ (over some $Z$) with $S$ an affine derived scheme, $\formalU$ a smooth affine formal scheme and $i$ an isomorphism of reduced subschemes.
\end{observation}

We now construct the moduli functor of Liouville data:

\begin{construction}\label{constructionLiouville4mars}
  Let $\stackforms\colon\Fun(\Delta^2, \dSt_\basefield)^\op \to \inftygpd$ be a stack, classified by a Cartesian fibration $p_\stackforms\colon \dSt^{\Delta^2,\stackforms} \coloneqq \int \stackforms \to \Fun(\Delta^2, \dSt_\basefield)$.
  Denote by $q_\stackforms$ the functor
  \[
    q_\stackforms \colon
    \begin{tikzcd}[cramped]
      \dSt^{\Delta^2,\stackforms} \ar{r}{p_\stackforms} & \Fun(\Delta^2, \dSt_\basefield) \ar{r}{\ev_{02}} & \Fun(\Delta^1, \dSt_\basefield).
    \end{tikzcd}
  \]
  Since both $p_\stackforms$ and $\ev_{02}$ (\cref{stackofstacks}) are Cartesian fibrations, so is $q_\stackforms$ \cite[2.4.2.3]{lurie-htt}. We denote by
  \[
    \stackfactor^{\stackforms} \colon \Fun(\Delta^1, \dSt_\basefield)^\op \to \inftycats
  \]
  the resulting \ifunctor.
  It is a stack as both $\stackforms$ and $\stackfactor$ are, and it comes naturally equipped with a natural transformation $\stackfactor^{\stackforms}\to \stackfactor$ that forgets the information of $\stackforms$. We apply this construction to the functors $\stackofliouvillesections$ and $\stackofisotropicsections$ of \cref{Liouvillefunctor4mars} and  obtain a commutative square of stacks over $\Fun(\Delta^1, \dSt_\basefield)$
  \begin{equation}\label{squareofabeliangroupstacksGrothendieck2}
    \begin{tikzcd}
      \mathllap{\stackfactor^{\Liouv} \coloneqq {}}\stackfactor^{\stackforms^{\Liouv}} \ar{r} \ar{d} & \stackfactor^{\stackforms^\isotropic} \mathrlap{{}\eqqcolon \stackfactor^{\isotropic}} \ar{d}\\
      \stackforms^{2,\exact}_\relative(-,-1) \ar{r} & \stackforms^{2,\closed}_\relative(-, -1).
    \end{tikzcd}
  \end{equation}
\end{construction}

\begin{remark}
  Unpacking \cref{constructionLiouville4mars}, for any $X \to Z \in \Fun(\Delta^1,\dSt_\basefield)$, the objects of the \icategory $\stackfactor^{\Liouv}(X \to Z)$ are factorizations $X \to Y \to Z$ equipped with a section of $\Liouv(X/Y/Z)$.
\end{remark}

\begin{construction}\label{restrictionLiouvillefunctoralongfunctor}
  Let $\calD \colon \Fun(\Delta^1, \dSt_\basefield)^\op \to \inftycats$ be a stack for the étale topology.
  For any stack $Y$ and any functor $F \colon \dSt_Y \to \Fun(\Delta^1, \dSt_\basefield)$, we consider $\calD_{F}$ the restriction of $\calD$ along $F$.
  Whenever $F$ preserves étale covers and fiber products, the restriction $\calD_{F} \colon \dSt_Y^\op \to \inftycats$ remains a stack for the étale topology.
\end{construction}

\begin{example}\label{stackofstacksnofactorizations}
  Let $e \colon \dSt_\basefield\to \Fun(\Delta^1, \dSt_\basefield)$ be the functor sending $X\mapsto (\varnothing \to X)$.
  The restriction  $\stackfactor_e$ of \cref{restrictionLiouvillefunctoralongfunctor} is the stack of derived stacks of \cite[6.1.3.9-(3), 6.1.3.10]{lurie-htt} $\thestackofstacks \colon \dSt_\basefield^\op\to \inftycats$, sending $X$ to the \icategory $\dSt_X$ of derived stacks over $X$.
\end{example}

\begin{notation}\label{notationfunctorpderham}
  Let $X \in \dSt_\basefield$.
  We denote by $p_X \colon \dSt_{X_\derham} \to \Fun(\Delta^1, \dSt_\basefield)$ given by $S \mapsto ( S \times_{X_\derham} X \to S)$.
\end{notation}

\begin{construction}\label{constructionliouviellefunctor4mars2}
  Let $X \in \dSt_\basefield$. The functor $p_X$ preserves étale maps and fiber products.
  \cref{restrictionLiouvillefunctoralongfunctor} with $F=p_X$ applied to the functors \cref{squareofabeliangroupstacksGrothendieck2}, yields a commutative square of stacks $\dSt_{X_\derham}^\op \to \inftycats$
  \[
    \begin{tikzcd}[row sep=scriptsize]
      \stackfactor^{\Liouv}_{p_X} \ar{r} \ar{d} & \stackfactor^{\isotropic}_{p_X} \ar{d} \\
      \stackforms^{2,\exact}_\relative(p_X(-),-1) \ar{r} & \stackforms^{2,\closed}_\relative(p_X(-),-1)\rlap{.}
    \end{tikzcd}
  \]
\end{construction}

\begin{definition}\label{fixedexactformstack}
  Let $X$ be a derived stack equipped with an exact 2-form $\lambda$.
  The form $\lambda$ induces a global section of $\stackforms^{2,\exact}_\relative(p_X(-),-1)$, using the observation $\stackforms^{2,\exact}(X,-1) \simeq \stackforms^{2,\exact}_\relative(X/X_\deRham,-1)$ (see \cref{derhamandformalcompletion}). We define $\stackfactor^{\Liouv, \lambda}_{p_X}$ as the fiber product (in categorical stacks over $\dSt_{X_\derham}$)
  \begin{equation}\label{definitionLiouvillestackalphafiberproductconstruction}
    \begin{tikzcd}[row sep=scriptsize]
      \stackfactor^{\Liouv, \lambda}_{p_X} \ar{r} \ar{d} \tikzcart & \stackfactor^{\Liouv}_{p_X} \ar{d} \\
      \{\lambda\} \ar{r} & \stackforms^{2,\exact}(p_X(-),-1)\rlap{.}
    \end{tikzcd}
  \end{equation}
\end{definition}

Finally, as discussed in \cref{squarenotCartesianingeneral}  we restrict our attention to factorizations for which the square \cref{squareofabeliangroupstacks} is a pullback:
\begin{definitions}{Let $X$ be a derived Deligne--Mumford stack.}
  \item \label{notationforrestrictiontoXet} Denote by $(-)^X_\derham$ the functor $X_\et^\daff \to \dSt_{X_\derham}$ mapping $S$ to $S_\derham$.
  \item \label{Liouvillesubfunctor} Let $\lambda$ be an exact 2-form on $X$.
  For any $S \to X \in X_\et^\daff$, we denote by $\Liouvstack^\lambda_X(S)$ the full subcategory of $\stackfactor^{\Liouv, \lambda}_{p_X}(S_\derham)$ spanned by those objects $S
  \to Y \to S_\derham$ where $Y$ is an affine smooth formal scheme $\formalU$ topologically of finite presentation over $\basefield$, and $S \to \formalU$ is a closed immersion which induces an isomorphism on the reduced subschemes $S_\red \simeq \formalU_\red$.
\end{definitions}

\begin{proposition}\label{Liouvillesubfunctorprop}
  In the situation of \cref{Liouvillesubfunctor}, $S \mapsto \Liouvstack^\lambda_X(S)$ defines a full substack of $\stackfactor^{\Liouv,\lambda}_{p_X} \circ (-)_\derham^X$ on the small affine étale site $X_\et^\daff$.
\end{proposition}
\begin{proof}
  We first proves it defines a subfunctor. Fix an étale schematic morphism $u \colon S \to S'$ over $X$.
  Let $S' \hookrightarrow \formalU' \to S'_\deRham$ with $S'\hookrightarrow \formalU'$ a closed immersion which induces an isomorphism of reduced schemes and $\formalU'$ a smooth affine formal scheme topologically of finite presentation over $\basefield$.
  The functoriality $u^\ast$ combined with the fact $u$ is étale gives us two pullback squares (using \cref{etaleformalcompletionisnoop})
  \[
    \begin{tikzcd}
      S \ar[hook]{r}{i} \ar{d}[swap]{u} \tikzcart & u^\ast(\formalU') \ar{r} \ar{d}[swap]{p} \tikzcart & S_\deRham \ar{d} \\
      S' \ar[hook]{r} & \formalU' \ar{r} & S'_\deRham.
    \end{tikzcd}
  \]
  It follows that $p$ is étale, $i$ is a closed immersion and $u^\ast(\formalU')$ is a formal thickening of $S$.
  Moreover, since $u$ is schematic, $u^\ast(\formalU')$ is a also smooth formal scheme.
  Finally, by \cref{extractingreducedsubstack} $i$ is an isomorphism on reduced schemes and $\Liouvstack^\lambda_X$ indeed defines a subfunctor.
  The stack property follows as the condition of being a smooth formal scheme or a closed immersion are étale-local properties.
\end{proof}

Finally, we introduce the moduli of Darboux data of a $(-1)$\nobreakdashes-shifted exact symplectic derived Deligne--Mumford stack.

\begin{definition}\label{definitiondarbouxstackasstacks}
  Let $X$ be a derived Deligne--Mumford stack equipped with a $(-1)$\nobreakdashes-shifted exact symplectic form $\lambda$.
  The \emph{Darboux stack} of $X$, denoted by $\Darbstack_X^\lambda$, is the full substack\footnote{This is indeed a substack, as the non-degeneracy condition is stable under étale maps, and étale-local.} on $X_\et^\daff$
  \[
    \Darbstack^{\lambda}_X \subset \Liouvstack^{\lambda,\simeq}_X \colon \left(X_\et^\daff\right)^\op \to \inftygpd
  \]
  of the maximal subgroupoid stack of $\Liouvstack_X^{\lambda}$ whose $S$\nobreakdashes-points are factorizations $S\to \formalU \to S_\derham$ together with section $s \in \Liouv(S/\formalU)$ such that the induced isotropic fibration via the projection in \cref{constructionliouviellefunctor4mars2} is Lagrangian in the sense of \cref{definition-oflagrangianfibration}.
\end{definition}

\begin{observation}
  If $Y \to X$ is an étale mophism, then $\Darbstack_Y^{\lambda_{|Y}}$ is naturally equivalent to the restriction of $\Darbstack_X^\lambda$ along $Y_\et^\daff \to X_\et^\daff$.
 \end{observation}

\begin{remark}
   In view of \cref{Bhatt}, Darboux charts could be equivalently described as Lagrangian foliations: \ie foliations, in the sense of \cite{Toen2020a}, equipped with a non-degenerate isotropic structure (similarly to \cref{definitionisotropic}). More precisely, the data of symplectic identification $S\simeq \dCrit(\formalU,f)$ produces a Lagrangian foliation by the fibers of the Lagrangian fibration $S\simeq \dCrit(\formalU,f)\to \formalU$ of  \cref{lagrangianfibrationdCrit}. Equivalently, every Lagrangian foliation $\calF$ on  $S$ gives us a Darboux chart given by the formal quotient $\formalU \coloneqq S/\calF$. Although we will not need this in this paper, the approach via Lagrangian foliations was our first definition of the Darboux stack inspired by T. Pantev's talk \cite{toen-pantev-video}. We still believe it will become relevant when defining Darboux stacks for $n$\nobreakdashes-shifted symplectic structures in general. We intend to investigate the analogous of \cref{theoremcontractibility} for $n=-2$ in a future work.
\end{remark}

\subsection{Formal \texorpdfstring{$\LG$}{LG}-pairs}\label{sectionLGpairs}
To conclude this section we relate Darboux data to the more familiar language of (formal) Landau--Ginzburg pairs ($\LG$\nobreakdashes-pairs). We recall the terminology:

\begin{definitions}[{$\LG$\nobreakdashes-pairs}]{}\label{definitionsofalgebraicandformalLGpairs}
  \item An \emph{algebraic} $\LG$\nobreakdashes-pair over $\basefield$ is a pair $(U,f)$ consisting of a smooth $\basefield$\nobreakdashes-scheme $U$ with a function $f \colon U \to \affineline{\basefield}$. A morphism of $\LG$\nobreakdashes-pairs $(U,f) \to (V,g)$ is a morphism of schemes $u \colon U \to V$ such that $g \circ u = f$.
  \item A \emph{formal} $\LG$\nobreakdashes-pair consists of a pair $(\formalU, f)$ where $\formalU$ is a smooth formal scheme and $f \colon \formalU \to \affineline{\basefield}$.
  Morphisms are defined as in the algebraic case.
\end{definitions}

\begin{construction}\label{LGpairsasGrothendieckconstruction}
  Consider the (limit-preserving) \ifunctor of functions
  \[
    \LG \coloneqq \Map(-, \affineline{\basefield}) \colon \dSt_\basefield^\op \to \inftygpd.
  \]
  It comes with a canonical equivalence $\dSt_{\affineline{\basefield}}\simeq \int \LG$ where the right hand side is the Cartesian Grothendieck construction of the \ifunctor $\LG$.
  Finally, the canonical projection $\Liouv(X/Y/Z) \to \globalsections(Y, \structuresheaf_Y)$ in \cref{Liouvillefunctor4mars} induces a symmetric monoidal natural transformation of abelian group stacks on $\Fun(\Delta^2, \dSt_\basefield)$
  \begin{equation}\label{squareofabeliangroupstackstoLG}
    \begin{tikzcd}
      \stackofliouvillesections \ar{r}& \LG \circ \, \ev_1(-)
    \end{tikzcd}
  \end{equation}
  where $\ev_1\colon\Fun(\Delta^2,\dSt_\basefield) \to \dSt_\basefield$ is the evaluation at $1$ in $\Delta^2$.
  Applying \cref{constructionLiouville4mars}, we get an \ifunctor
  \[
    \stackfactor^{\LG} \coloneqq \stackfactor^{\LG \circ\, \ev_1} \colon (X \to Z) \mapsto \left\{ X \to Y \to Z \times \affineline{\basefield} \right\}.
  \]
  We now fix a derived stack $X$. Restricting the above functor along $p_X$ (\cf \cref{notationfunctorpderham}), the map \cref{squareofabeliangroupstackstoLG} induces a morphism of monoidal stacks $\dSt_{X_\derham}^\op \to \monoidalcats$
  \begin{equation}\label{naturaltransformationfromLiouvilletoLG}
    \begin{tikzcd}
      \stackfactor^{\Liouv}_{p_X} \ar{r} & \stackfactor^{\LG}_{p_X}\rlap{.}
    \end{tikzcd}
  \end{equation}
  We denote by $\stackLGder{X}{}$ the full substack
  \[
    \stackLGder{X}{} \subset \stackfactor^{\LG}_{p_X} \circ (-)_\derham^X \colon (X^\daff_\et)^\op \to \inftycats
  \]
  (\cf \cref{notationforrestrictiontoXet}) spanned by factorizations $S \to Y \to S_\derham$ with a function $f \colon Y \to \affineline{\basefield}$ such that:
  \begin{enumerate}
    \item The stack $Y$ is a smooth affine formal scheme over $\basefield$;
    \item The map $S \to Y$ is a closed immersion and a reduced equivalence.
    \item The form $df_{|{\classicaltruncation(S)}}$ vanishes and the induced map $\classicaltruncation(S) \hookrightarrow \Crit(Y,f)$ is an isomorphism.
  \end{enumerate}
  The same argument as in \cref{Liouvillesubfunctorprop} shows that this indeed defines a substack.
  Assume further that $X$ is Deligne--Mumford.
  For any $(-1)$\nobreakdashes-shifted exact symplectic structure $\lambda$, the natural transformation \cref{naturaltransformationfromLiouvilletoLG} restricts to a morphism of stacks over $X^\daff_\et$
  \begin{equation}\label{naturaltransformationfromLiouvilletoLG2}
    \Darbstack_X^\lambda \to \stackLGder{X}{}.
  \end{equation}
\end{construction}

\begin{observation}
  Unpacking \cref{LGpairsasGrothendieckconstruction}, the functor \cref{naturaltransformationfromLiouvilletoLG2} sends a factorization $S\hookrightarrow \formalU\to S_\derham$ with Liouville section $s\in \Liouv(S/\formalU/S_\derham)$ to the formal scheme $\formalU$ equipped with the function given by the projection $\Liouv(S/\formalU/S_\derham)\to \globalsections(\formalU, \structuresheaf_\formalU)$. A morphism of Darboux charts is sent to the underlying morphism of formal $\LG$\nobreakdashes-pairs under $S$
  \[
    \begin{tikzcd}[row sep=0]
      & \formalU\ar{dd}{\rho} \ar{dr}{f} \\
      S \ar{ur}{i} \ar{dr}[swap]{j} & {}& \affineline{\basefield}\rlap{.} \\
      & \formalV \ar{ur}[swap]{g}
    \end{tikzcd}
  \]
\end{observation}

\begin{observation}\label{remarkcriticalvalues}
  The role of the exact structure can be better explained at this point.
  Recall from \cref{locallyconstantfunctionunderlyingexact} that any $(-1)$\nobreakdashes-shifted exact form $\lambda$ induces a locally constant critical value function $c_\lambda$.
  Fix $S$ étale over $X$ and fix an $S$\nobreakdashes-point of $\Darbstack_X^\lambda$. Using \cref{naturaltransformationfromLiouvilletoLG2}, we get an $\LG$\nobreakdashes-pair
  \[
    S \hookrightarrow \formalU \to^f \affineline{}.
  \]
  A rapid diagram chase then shows that $f_{|S_\red} = c_\lambda$.
  It follows that the choice of an exact structure determines the critical values of the $\LG$\nobreakdashes-pairs allowed as Darboux charts.
  \Cref{canonicalexactstructure} then simply means that we can enforce Darboux charts to have only $0$ as critical value.
\end{observation}

\section{Quadratic bundles and their action on the Darboux stack}
\label{section-Quadraticbundles}

As shown in \cref{ambiguitypresentationdcrit}, the space of local presentations of a $(-1)$\nobreakdashes-shifted symplectic derived Deligne--Mumford stack as a derived critical locus has an internal symmetry by the action of quadratic bundles. In this section we introduce the convenient moduli of quadratic bundles acting on Darboux charts (\cref{sectionStackQuad}), construct the action (\cref{section-Action}) and show that it is transitive -- 
\cref{thm-transitivity}.

\subsection{The stack of quadratic bundles with flat connection.}
\label{sectionStackQuad}

\begin{notation}\label{trivialquadraticbundle}
  Let $X$ be a derived stack. We denote by $\trivialquadratic{X}$ the trivial rank $n$ vector bundle on $X$ equipped with the standard quadratic form $q_{\mathrm{std}}=x_1^2+ \cdots + x_n^2$ and the trivial flat connection.
\end{notation}

\begin{definition}
  \label{quadgroupoid}
  Let $S$ an affine derived scheme.
  We consider $\Quad^\nabla(S)$ the 1-groupoid classifying \emph{non-degenerate} quadratic bundles on $S$ with an étale locally trivial flat connection.
  More precisely, $\Quad^\nabla(S)$ classifies triples $ (Q,q, \nabla)$ where  $(Q, q)$ is a non-degenerate quadratic vector bundle on $S$ with $\nabla$ a compatible flat connection along $S$, and such that, there exists an étale cover $S'\to S$ such that the restriction of $ (Q,q, \nabla)$ is isomorphic to $\trivialquadratic{S'}$.
\end{definition}

\begin{observation}\label{definitionofQuadviaderhamstack}
  In practice it will be more convenient to reformulate \cref{quadgroupoid} using the de Rham stack, identifying the groupoid $\Quad^\nabla(S)$ with the full subgroupoid of $\Quad(S_\derham)$ classifying étale locally trivial connections.
  In particular, every object $(Q,q,\nabla)\in \Quad^\nabla(S)$ is obtained via pullback from a quadratic bundle $\formalQ$ on $S_\derham$, \ie, we have a Cartesian diagram
  \[
    \begin{tikzcd}
      Q \ar[r] \ar{d} \tikzcart & \formalQ \ar{d} \\
      S \ar[r] & S_{\derham}\rlap{.}
    \end{tikzcd}
  \]
\end{observation}

\begin{definition}\label{definitionBO}
  We denote by $\B\orthogonal$ the classical $1$\nobreakdashes-stack in groupoids
  \[
    \B\orthogonal = \coprod_n \B(\orthogonal_n),
  \]
  equipped with its direct sum monoidal structure induced by the direct sum of matrices
  \[
    \oplus \colon \orthogonal_p \times \orthogonal_q \to \orthogonal_{p+q}.
  \]
\end{definition}

\begin{definition}\label{stackquaddef}
  Let $X$ a derived Deligne--Mumford stack. We consider the
  stack on $X_\et^\daff$ classifying quadratic bundles on $X$ with an étale locally trivial flat connection, \ie
  \begin{align*}
    \stackQuadnabla_X \colon (X_\et^\daff)^\op & \to \inftygpd
    \\
    [S \to_\et X]                              & \mapsto \stackQuadnabla_X(S) \coloneqq \Quad^\nabla(S)
  \end{align*}
  as in \cref{quadgroupoid}. By definition,  this is the substack of
  \begin{align*}
    (X_\et^\daff)^\op & \to \inftygpd                          \\
    S                 & \mapsto \Map(S_\deRham, \B\orthogonal)
  \end{align*}
  spanned by locally trivial sections. The monoid structure $\oplus$ on $\B\orthogonal$ induces a monoidal structure on the stack $\stackQuadnabla_X$.
\end{definition}

\begin{observation}\label{constantquadstack}
  Consider $\Quad_\basefield = \B\orthogonal(\complex) = \coprod_n \B(\orthogonal_n(\complex))$ the 1-groupoid of non-degenerate quadratic forms over $\complex$. Let $p \colon X\to \Spec(\complex)$ be the canonical projection. Then the stack $\stackQuadnabla_X$ coincides with $p^{-1} \Quad_\complex$, where $p^{-1} \colon  \Sh_{\et}(\Spec(\complex))\to \Sh_{\et}(X)$ is the sheafification of the naive inverse image of presheaves $p^{\dagger}$. Indeed, notice first that $p^{\dagger}\Quad_\complex$ is the constant functor with value $\Quad_\complex$.
  Therefore, we can consider the map of prestacks
  \[
    (p^{\dagger}\Quad_\complex)(S) = \Quad_\complex \to \stackQuadnabla_X(S)
  \]
  sending $(\complex^{\oplus_n},q_{\standard})$ to $Q=\trivialquadratic{S}$.
  It induces an isomorphism on $\pi_1$ as part of the Riemann--Hilbert correspondence applied to the classical truncation of $S$.
  The defining condition of $\stackQuadnabla_X$ implies that after sheafification, the above morphism becomes an isomorphism on the $\pi_0$\nobreakdashes-sheaves showing that the induced map
  \[
    p^{-1}\Quad_\complex \to \stackQuadnabla_X
  \]
  is an isomorphism of stacks.
\end{observation}

\begin{observation}\label{trivialsheafonpizero}
  Notice that by \cref{stackquaddef}, $\homotopysheaf_0\left(\stackQuadnabla_{X}\right) \simeq \constantsheaf \integers^{\geq 0}_X$ as sheaves (here $\homotopysheaf_0$ denotes the étale hypersheafification of the $\pi_0$\nobreakdashes-presheaf) corresponding to the rank.
\end{observation}

\subsection{Quadratic bundles and Liouville complexes}\label{quadandliouville}
The goal of this section is to construct the action of the stack of non-degenerate flat quadratic bundles on Darboux charts. We start with a symmetric monoidal enhancement of the discussion in \cref{sectionconstructionformaldarboux}:

\begin{construction}\label{relativesymmetricmonoidalstructure1}
  The Cartesian fibration
  $\ev_{02} \colon \Fun(\Delta^2, \dSt_\basefield) \to \Fun(\Delta^1, \dSt_\basefield)$ of \cref{stackofstacks} parametrizes a family of symmetric monoidal \icategories \footnote{\cite[\S 2.4.1]{lurie-ha}} over $\Fun(\Delta^1,\dSt_\basefield)$, with the relative Cartesian monoidal structure.
  On objects over $X \to Z$, the tensor product of $X\to Y\to Z$ with $X\to Y'\to Z$ is the fiber product $X \to Y \times_Z Y'\to Z$.  Since this fiberwise  monoidal structure is given by pullbacks, \cite[2.4.1.9]{lurie-ha}  upgrades the stack $\stackfactor$ of \cref{stackofstacks} to an \ifunctor $\stackfactor^{\times}\colon\Fun(\Delta^1, \dSt_\basefield)^\op \to \monoidalcats$, which we can also read as an \ifunctor
  \[
    \Fun(\Delta^1, \dSt_\basefield)^\op\times \pFSets \to \inftycats
  \]
  satisfying the Segal-monoid condition of \cite[2.4.2.1]{lurie-ha}. Finally, we denote by
  \[
    p_{\stackfactor^{\times}} \colon \int_\cocart \stackfactor^{\times} \to \Fun(\Delta^1,\dSt_\basefield)^\op\times \pFSets
  \]
  the associated \emph{coCartesian} fibration.

\end{construction}

\begin{construction}\label{constructionLiouville4marsmonoidal}
  The abelian group structure on the Liouville complexes  extends the \ifunctor $\stackofliouvillesections \colon \Fun(\Delta^2, \dSt_\basefield)^\op\to \inftygpd$
  of \cref{Liouvillefunctor4mars} to an \ifunctor
  \[
    \Fun(\Delta^2, \dSt_\basefield)^\op\to \CAlg(\inftygpd)
  \]
  which we can write as a map
  \[
    \begin{tikzcd}[row sep=tiny, column sep=tiny]
      \Fun(\Delta^2,\dSt_\basefield)^\op\ar{rr}{\stackofliouvillesectionsplus} \ar{dr}[swap]{\ev_{02}^\op} & & \CAlg(\inftygpd) \times \Fun(\Delta^1, \dSt_\basefield)^\op\ar{dl}\\
      & \Fun(\Delta^1, \dSt_\basefield)^\op &
    \end{tikzcd}
  \]
  We now claim that the functor $\stackofliouvillesectionsplus$ is $\Fun(\Delta^1,\dSt_\basefield)^\op$\nobreakdashes-fiberwise lax symmetric monoidal with respect to the $\Fun(\Delta^1, \dSt_\basefield)^\op$\nobreakdashes-relative symmetric monoidal structure of \cref{relativesymmetricmonoidalstructure1} on the left hand side and the trivial $\Fun(\Delta^1, \dSt_\basefield)^\op$\nobreakdashes-relative symmetric monoidal structure on the right hand side given by the Cartesian symmetric monoidal structure on $\CAlg(\inftygpd)$.
  The relative lax monoidal structure over $(X\to Z)\in \Fun(\Delta^1, \dSt_\basefield)$ given by the canonical sum morphism
  \[
    \textstyle \Liouv(X/Y/Z) \oplus \Liouv(X/Y'/Z) \lra^{\boxplus} \Liouv\left(X/Y \times_Z Y'/Z\right).
  \]
  promotes $\stackofliouvillesectionsplus$ to a map between \emph{coCartesian} fibrations
  \[
    \begin{tikzcd}[row sep=tiny, column sep=tiny]
      \displaystyle \smash{\int_{\cocart}} \stackfactor^{\times} \ar{dr}\ar{rr}{\stackofliouvillesectionspluslax}&& \ar{dl} \Fun(\Delta^1,\dSt_\basefield)^\op \times \CAlg(\inftygpd)^{\times}
    \\
      &\Fun(\Delta^1,\dSt_\basefield)^\op \times \pFSets &
    \end{tikzcd}
  \]
  where $\CAlg(\inftygpd)^{\times} \to \pFSets$ is the Cartesian symmetric monoidal structure associated to $\CAlg(\inftygpd)$ and finite products therein \cite[2.4.1.5]{lurie-ha}.
  It comes with a functor $\pi \colon \CAlg(\inftygpd)^{\times} \to \CAlg(\inftygpd)$. We consider the composition
  \begin{equation}\label{compositionlaxforget}
    \begin{tikzcd}[row sep=tiny, column sep=small]
      \displaystyle \smash{\int_{\cocart}} \stackfactor^{\times} \ar{rrr}{\stackofliouvillesectionspluslax} &&& \Fun(\Delta^1,\dSt_\basefield)^\op \times \CAlg(\inftygpd)^{\times} \ar{r}{\pi} & \CAlg(\inftygpd) \ar{r} & \inftygpd
    \end{tikzcd}
  \end{equation}
  where the last functor extracts the underlying \igroupoid. Informally, it sends an object $(X\to Y_i\to Z)_{i\in \{1,\cdots n\}}$ lying over $(X\to Z, \langle n \rangle)\in \Fun(\Delta^1, \dSt_\basefield)^\op\times \pFSets$ to the product $\prod_{i\in \{1,\cdots ,n\}} \Liouv(X/Y_i/Z)\in \inftygpd$. Next, consider the associated \emph{coCartesian} fibration.
  \[
    \begin{tikzcd}
      \mathcal{V}_{\stackofliouvillesections}\ar{r} & \displaystyle \smash{\int_{\cocart}} \stackfactor^{\times}
    \end{tikzcd}
  \]
  together with the composition
  \[
    \begin{tikzcd}
      q^+_{\Liouv} \colon \mathcal{V}_{\stackofliouvillesections} \ar{r} & \displaystyle \smash{\int_{\cocart}} \stackfactor^{\times} \ar{r}{p_{\stackfactor^{\times}}} & \Fun(\Delta^1,\dSt_\basefield)^\op \times \pFSets
    \end{tikzcd}
  \]
  which, being a composition of coCartesian fibrations, is also a coCartesian fibration. Finally, we consider the  straightening of $q^+_{\Liouv}$:
  \[
    \Fun(\Delta^1, \dSt_\basefield)^\op\times \pFSets \to \inftycats
  \]
  and observe that by construction it also satisfies the Segal-monoid conditions of \cite[2.4.2.1]{lurie-ha}  relatively to  $\Fun(\Delta^1, \dSt_\basefield)^\op$. Indeed, by construction, the fiber of $p_{\stackfactor^{\times}}$ over $(X\to Z, \langle n \rangle)$ is the \icategory $\prod_{i=1}^n(\dSt_\basefield)_{X/./Z}$ and the restriction of the functor \cref{compositionlaxforget} to this fiber, factors, by construction, through the $n$\nobreakdashes-fold product
  \[
    \begin{tikzcd}[column sep=large]
      \displaystyle \prod_{i=1}^{n} (\dSt_\basefield)_{X/./Z} \ar{r}{\prod_{i=1}^{n}\stackofliouvillesections} & \displaystyle \prod_{i=1}^{n} \inftygpd \ar{r}{\times}&\inftygpd
    \end{tikzcd}
  \]
  so that its Grothendieck construction is the $n$\nobreakdashes-fold product of the Grothendieck constructions of each component. We obtain an \ifunctor
  \[
    \stackfactor^{\Liouv, \boxplus} \colon \Fun(\Delta^1,\dSt_\basefield)^\op \to  \monoidalcats
  \]
  together with a symmetric monoidal natural transformation $\stackfactor^{\Liouv, \boxplus}\to \stackfactor^{\times}$, providing a symmetric monoidal enhancement of \cref{constructionLiouville4mars}.
\end{construction}

\begin{observations}{}
  \item The monoidal structure in $  \stackfactor^{\Liouv, \boxplus}(X\to Z)$ maps $X \to Y \to Z$ and $X \to Y' \to Z$ to $X \to Y \times_Z Y' \to Z$, with Liouville section obtained using the canonical sum morphism
  \[
    \textstyle \Liouv(X/Y/Z) \oplus \Liouv(X/Y'/Z) \lra^{\boxplus} \Liouv\left(X/Y \times_Z Y'/Z\right).
  \]
  \item\label{liouvtoLGmonoidal} Consider the functor $\LG\circ\, \ev_{1}$ of \cref{LGpairsasGrothendieckconstruction} with its canonical abelian group structure and the lax monoidal structure induced by
  \[
    \globalsections(Y, \structuresheaf_Y)\oplus \globalsections(Y', \structuresheaf_{Y'})\to \globalsections\left(Y\fiberproduct{Z} Y', \structuresheaf_{Y\fiberproduct{Z} Y'}\right).
  \]
  We can run \cref{constructionLiouville4marsmonoidal} replacing $\stackforms^{\Liouv}$ by $\LG\circ\, \ev_{1}$ and obtain a symmetric monoidal stack $\stackfactor^{\LG, \boxplus} \colon \Fun(\Delta^1,\dSt_\basefield)^\op \to \inftycats^{\otimes}$, together with a symmetric monoidal upgrade of the natural transformation \cref{naturaltransformationfromLiouvilletoLG}
  \begin{equation}\label{laxnaturaltransformationfromLiouvilletoLG}
    \stackfactor^{\Liouv, \boxplus} \to \stackfactor^{\LG, \boxplus}.
  \end{equation}
\end{observations}

We are now finally ready to address the relation between quadratic bundles and Liouville structures:

\begin{proposition}\label{orthogonaltoliouv}
  Let $X \in \dSt_\basefield$.
  Taking the formal completion of the zero section, defines a symmetric monoidal morphism of stacks $\dSt_{X_\derham}^\op \to \monoidalcats$
  \[
    \B\orthogonal^{\oplus}_{X_\derham} \to \stackfactor^{\Liouv,\boxplus}_{p_X}
  \]
  such that the composite
  \[
    \B\orthogonal^{\oplus}_{X_\derham} \to \stackfactor^{\Liouv,\boxplus}_{p_X} \to \stackforms^{2,\exact}_\relative(p_X(-),-1)
  \]
  is nullhomotopic (where $p_X$ is as in \cref{notationfunctorpderham}).
\end{proposition}
\begin{proof}
  Write $\myuline{\FSch}^{\Liouv}$ for the full substack of $\stackfactor^{\Liouv}_{p_X}$ whose $S$\nobreakdashes-points are factorizations $S \to^s Y \to^\pi S$  such that the morphism $\pi$ is smooth and represented by a formal affine scheme.
  \begin{lemma*}
    The classical truncation $\classicaltruncation\myuline{\FSch}^{\Liouv}$ of $\myuline{\FSch}^{\Liouv}$ is a $1$\nobreakdashes-stack.
    \begin{proof}
      Let $S \in \Aff_\basefield$ and pick an $S$\nobreakdashes-point of $\myuline{\FSch}^{\Liouv}$.
      It corresponds to a factorization $S \to^s Y \to^\pi S$ (plus a Liouville section).
      Denote by $\rmN^\dual \simeq s^* \Omega_{Y/S}^1$ the conormal bundle of $s$.
      We can compute
      \[
        E \coloneqq \fiber\left(\Liouv(S/Y/S) \to \derhamcomplex_{S/S}^{\leq 1} \simeq \Gamma(S, \structuresheaf_S)\right) \simeq \fiber\left(\Gamma(Y,\structuresheaf_Y) \to \derhamcomplex_{S/Y}^{\leq 1}\right)
      \]
      as a (split) extension $0 \to \rmN^\dual \to E \to \ker(\Gamma(Y,\structuresheaf_Y) \to \Gamma(S,\structuresheaf_S)) \to 0$.
      It is in particular concentrated in degree $0$ and the result follows.
    \end{proof}
  \end{lemma*}
  We can now return to the proof of \cref{orthogonaltoliouv}.
  The inclusion $\classicaltruncation\myuline{\FSch}^{\Liouv} \to \classicaltruncation\stackfactor^{\Liouv}_{p_X}$ is stable under $\smboxplus$, and it therefore suffices to use the evident symmetric monoidal functor
  \[
    \B\orthogonal^\oplus_{X_\derham} \to \classicaltruncation\myuline{\FSch}^{\Liouv,\boxplus}
  \]
  mapping a quadratic bundle $(Q,q) \to S$ over $S = \Spec A \in \Aff_{X_\derham}$ to the formal neighborhood $\widehat Q_0$ of its zero section, equipped with the Liouville section:
  \begin{equation}\label{liouvillesectionactionquad}
    \newcommand{\downsymbol}[1]{\ar[phantom, start anchor=center, end anchor=center]{d}[sloped]{#1}}
    \begin{tikzcd}[row sep=small]
      q \downsymbol{\in} \ar[mapsto]{r} & 0 \downsymbol{\in} & 0 \ar[mapsto]{l} \downsymbol{\in}
      \\
      \Symcompleted_A(Q^\dual) \ar[r] \downsymbol{\simeq} & \Sym_A^{\leq 1}(Q^\dual)  \downsymbol{\simeq} & A \ar{l} \downsymbol{=}
      \\
      \Gamma(\widehat Q_0, \structuresheaf_{\widehat Q_0}) \ar[r] & \derhamcomplex_{S/\widehat Q_0}^{\leq 1} & \derhamcomplex_{S/S}^{\leq 1}\rlap{,} \ar{l}
    \end{tikzcd}
  \end{equation}
  which is easily checked to be monoidal, as formal completion preserves fiber products.
\end{proof}

\subsection{Action of quadratic bundles on the Darboux stack}
\label{section-Action}

The goal of this section is to construct a natural action of quadratic bundle on the Darboux stack:

\begin{theorem}\label{actionondarbouxcharts}
  Let $X$ be a derived Deligne--Mumford stack with a $(-1)$\nobreakdashes-shifted exact symplectic form $\lambda$.
  The monoid  stack $\stackQuadnabla_X$ of \cref{stackquaddef}
  acts on $\Darbstack_X^{\lambda}$, in a way such that
  the action by $(\calQ,q) \in \stackQuadnabla_X(S)$ sends $S \to^i \formalU \to^f \affineline{}$ (forgetting the rest of the Liouville section) to
  \[
    S \lra^{(i,0)} \formalU \times_{S_\deRham} \widehat\formalQ_0 \lra^{f \smboxplus q} \affineline{}.
  \]
\end{theorem}

\begin{proof}
  Consider the morphism $\B\orthogonal^{\oplus}_{X_\derham} \to \stackfactor^{\Liouv,\boxplus}_{p_X}$ of \cref{orthogonaltoliouv}.
  Restricting to the small étale site of $X$ via $(-)_\derham^X$, we get a nullhomotopic composite morphism
  \begin{multline*}
    \stackQuadnablaoplus_X \subset \B\orthogonal^\oplus_{X_\derham} \circ (-)^X_\derham \to \stackfactor^{\Liouv,\boxplus}_{p_X} \circ (-)^X_\derham \to \stackforms^{2,\exact}_\relative(p_X((-)_\derham^X), -1)^+ \\\simeq \stackforms^{2,\exact}(-,-1)^+_{|X_\et^\daff}.
  \end{multline*}
  Recall from \cref{fixedexactformstack} the notation
  \[
    \stackfactor^{\Liouv,\lambda}_{p_X} \circ (-)^X_\derham \coloneqq \fiber\left(\stackfactor^{\Liouv}_{p_X} \circ (-)^X_\derham \to \stackforms^{2,\exact}(-,-1)_{|X_\et^\daff} \ni \lambda\right).
  \]
  We get this way an action of $\stackQuadnablaoplus_X$ on $\stackfactor^{\Liouv,\lambda}_{p_X|X_\et^\daff} \coloneqq \stackfactor^{\Liouv,\lambda}_{p_X} \circ (-)^X_\derham$ using the symmetric monoidal structure of \cref{constructionLiouville4marsmonoidal}
  \[
    \begin{tikzcd}
      \stackQuadnablaoplus_X \times \stackfactor^{\Liouv,\lambda}_{p_X|X_\et^\daff}\ar{r}  &  \stackfactor^{\Liouv,0}_{p_X|X_\et^\daff} \times \stackfactor^{\Liouv,\lambda}_{p_X|X_\et^\daff}\ar{r}{\boxplus} & \stackfactor^{\Liouv, \lambda}_{p_X|X_\et^\daff}.
    \end{tikzcd}
  \]
  Recall now from \cref{Liouvillesubfunctor} the definition of the Liouville substack $\Liouvstack_X^{\lambda}$ as a substack of $\stackfactor^{\Liouv,\lambda}_{p_X|X_\et^\daff}$. Let us show that it is preserved by the action. Unpacking the construction, the action of $(\calQ,q) \in \stackQuadnabla_X(S)$ is (forgetting the rest of the Liouville section)
  \[
     \left( S \to^i \formalU \to^f \affineline{} \right) \mapsto \left(S \lra^{(i,0)} \formalU \times_{S_\deRham} \widehat\formalQ_0 \lra^{f \smboxplus q} \affineline{}\right).
  \]
  In particular, since $\formalU \times_{S_\deRham} \widehat\formalQ_0$ is a smooth formal scheme and $(i,0)$ is a reduced equivalence, we see that $\Liouvstack^{\lambda}_X$ is preserved by the action.

  To conclude, it remains to show that the Darboux substack of \cref{definitiondarbouxstackasstacks} is also preserved by the action, or in other words, that the isotropic fibration structure of $S \to \formalU \times_{S_\deRham} \widehat\formalQ_0$ is non-degenerate if that of $S \to \formalU$ is. Let us assume the latter.
  Since $S \to \widehat\formalQ_0$ factors through $S_\deRham$ by construction, the morphism $S \to \formalU \times_{S_\deRham} \widehat\formalQ_0$ factors through $\formalU$.
  The isotropic distribution associated to our action thus yields a null-homotopy of the composite morphism
  \[
    \tangent_{S/\formalU \times_{S_\deRham} \widehat \formalQ_0} \simeq \tangent_{S/\formalU} \oplus Q[-1] \to \tangent_{S/\formalU} \to \tangent_S \simeq \cotangent_S[-1] \to \cotangent_{S/\formalU}[-1] \to \cotangent_{S/\formalU}[-1] \oplus Q^\dual.
  \]
  Since the action adds up Liouville structures, it adds up the isotropic structure and thus the associated boundary operators as well.
  In particular, by \cref{liouvillesectionactionquad}, the induced boundary operator $\cofiber(\tangent_{S/\formalU} \oplus Q[-1] \to \tangent_S) \simeq i^*\tangent_{\formalU} \oplus Q \simeq \cotangent_{S/\formalU}[-1] \oplus Q \to \cotangent_{S/\formalU}[-1] \oplus Q^\dual$ is identified with $\Id_{\cotangent_{S/\formalU}[-1]} \oplus \hessian(q)$.
  It is an equivalence as $q$ is assumed to be non-degenerate, and the action indeed preserves $\Darbstack_X^{\lambda}$.
\end{proof}

\begin{observation}\label{re-interpretationofactionpreservinglagrangians}
  \Cref{liouvtoLGmonoidal} implies that the action of quadratic bundles is compatible with the convolution $\smboxplus$ of $\LG$\nobreakdashes-pairs.
  In particular, \cref{actionondarbouxcharts} establishes the pattern observed in \cref{ambiguitypresentationdcrit}, namely, the exact symplectic identification of $\dCrit(\formalU,f)$ with $\dCrit(\formalU \smboxplus \formalQ, f \smboxplus q)$.
\end{observation}

\subsection{Transitivity of the action}
\label{section-transitivity}

The goal of this section is to prove the transitivity of the action of \cref{actionondarbouxcharts}. Since $\stackQuadnabla_X$ is a monoid, the notion of transitivity only makes sense in the group-completed sense:

\begin{theorem}[Transitivity]\label{thm-transitivity}
  Let $X$ be derived Deligne--Mumford stack equipped with a $(-1)$\nobreakdashes-shifted exact symplectic form $\lambda$.
  Then, the action of $\stackQuadnabla_X$ on $\Darbstack^{\lambda}_X$ is transitive.

  More precisely: for any $S\in X_\et^\daff$ with Darboux charts $\formalU$, $\formalV\in \Darbstack_X^\lambda(S)$, then, étale locally on $S$ there exists quadratic bundles $\formalQ_1, \formalQ_2\in \stackQuadnabla_X$ and equivalences  $\formalU\smboxplus \formalQ_1 \simeq \formalV\smboxplus \formalQ_2$ in $\Darbstack_X^\lambda$.

  In particular, the étale sheaf $\homotopysheaf_0(\Darbstack^{\lambda}_X/\stackQuadnabla_X)$ is trivial.
\end{theorem}

\begin{observation}\label{localformformorphisms1}
 Recall from \cref{definitiondarbouxstackasstacks} that we defined $\Darbstack_X^\lambda$ as a stack in \igroupoids.
 There is however a natural notion of morphisms by taking the corresponding full substack (in \icategories) of $\Liouvstack_X^\lambda$. Similarly, the stack in groupoids $\stackQuadnabla_X$ has a natural categorical enhancement, and the action of \cref{actionondarbouxcharts} is obviously compatible.
 In particular, for any $S \in X_\et^\daff$ and any quadratic bundle $(\formalQ,q) \in \stackQuadnabla_X(S)$, the canonical morphism $0 \to (\formalQ,q)$ of quadratic bundles induces a morphism
 \[
   \formalU = \formalU \smboxplus 0 \to \formalU \smboxplus \formalQ = \formalU \times_{S_\deRham} \widehat\formalQ_0 \in \Liouvstack_X^\lambda(S)
 \]
 for any $\formalU \in \Liouvstack_X^\lambda(S)$ (and a fortiori in $\Darbstack_X^\lambda(S)$).
\end{observation}

One of the main ingredients in proving \cref{thm-transitivity} is the following converse:

\begin{proposition}\label{cor:map=action}
  Let $S \in X_\et^\daff$ and $\rho \colon \formalU \to \formalV$ in $\Liouvstack_X^\lambda(S)$.
  If both $\formalU$ and $\formalV$ belong to $\Darbstack_X^\lambda(S)$, then $\rho$ is étale locally of the form $\formalU \to \formalU \smboxplus \formalQ$ as in \cref{localformformorphisms1} where $\formalQ$ is the normal bundle of $\rho$.
\end{proposition}
\begin{proof}
  Fix the notations
  \[
    \begin{tikzcd}[row sep=0]
      & \formalU\ar{dd}{\rho} \ar{dr}{f} \\
      S \ar{ur}{i} \ar{dr}[swap]{j} && \affineline{\basefield}\rlap{.} \\
      & \formalV \ar{ur}[swap]{g}
    \end{tikzcd}
  \]
  \begin{enumerate}[leftmargin=0pt,itemindent=5em,label={$\bullet$ Step \arabic*.}, ref={Step \arabic*}]
    \item\label{stepproofclosedimmersion}
          Consider the underlying map of isotropic distributions. As in \cref{diagrammapinducedatlevelofdistributionsbyuniversalproperty}, it fits in a sequence of pushout squares in $\Perf(S)$
          \[
            \begin{tikzcd}[row sep=scriptsize, column sep=scriptsize]
              \arrow[rounded corners, "\ell_\formalU", to path={|- ([yshift=1em]\tikztotarget.north)[near end]\tikztonodes -- (\tikztotarget)}]{rr} \tangent_{S/\formalU} \tikzcart \ar{r}{\tangent_{S/\rho}} \ar{d} &
              \tangent_{S/\formalV} \tikzcart \ar{r}{\ell_\formalV} \ar{d} &
              \tangent_S\ar{d}\ar[rounded corners, to path={-|(\tikztotarget)\tikztonodes}]{dr}{\ell_\formalU^\vee[-1]\circ \omega_S} &
              \\
              0 \ar{r} &
              i^\ast \conormalUV[-1] \ar{r} \ar{d} \tikzcart &
              i^\ast \tangent_\formalU \ar{d}[swap]{\tangentmap{\rho}} \ar{r}{\alpha_{\eta_\formalU}} \ar[phantom]{dr}{\scriptstyle(\sigma)} &
              \cotangent_{S/\formalU}[-1]
              \\
              & 0 \ar{r} &
              j^\ast \tangent_{\formalV} \ar{r}[swap]{\alpha_{\eta_\formalV}} &
              \cotangent_{S/\formalV}[-1] \ar{u}[swap]{(\tangent_{S/\rho})^\vee[-1]}
            \end{tikzcd}
          \]
          where the $2$\nobreakdashes-cell $(\sigma)$ stems from the homotopy $\eta_\formalU \sim \Sym^2((\tangent_{S/\rho})^\vee[-1])\circ \eta_\formalV $ (\ie the isotropic structure on the morphism $\rho$) and $\normal_{\formalU/\formalV} \coloneqq \cofiber(\tangentmap{\rho} \colon \tangent_\formalU \to \rho^\ast \tangent_\formalV)$.
          Since $\eta_\formalU$ is Lagrangian,  $\alpha_{\eta_\formalU}$ is an equivalence.
          Therefore, we find $\alpha_{\eta_\formalU}^{-1}\circ (\tangent_{S/\rho})^\vee[-1]\circ \alpha_{\eta_\formalV}\circ \tangentmap{\rho} \sim \Id$ so that
          $r_\rho \coloneqq \alpha_{\eta_\formalU}^{-1}\circ (\tangent_{S/\rho})^\vee[-1]\circ \alpha_{\eta_\formalV}$ is a retract of $\tangentmap{\rho}$.
          Using this retract we find a compatible direct sum decomposition $j^\ast \tangent_{\formalV}\simeq i^\ast \tangent_\formalU\oplus i^\ast \normal_{\formalU/\formalV}$
          and $i^\ast \normal_{\formalU/\formalV}$ is thus projective.
          It follows that $\normal_{\formalU/\formalV}$ is projective over $\formalU$. In particular it is concentrated in degree $0$ and $\rho$ is a closed immersion.

    \item(Special tangential splitting):
          Denote by $\conormalUV$ the conormal bundle of $\rho$ and consider the exact sequence $0 \to \conormalUV \to \rho^* \cotangent_\formalV \to \cotangent_\formalU \to 0$.
          As $\formalU$ is smooth, it always admits a $\formalU$\nobreakdashes-linear splitting. Choose such a splitting $r \colon \rho^* \cotangent_\formalV \to \conormalUV$.
          We in particular get a commutative diagram (to avoid cluttering the picture, we omit the pushforward functors)
          \[
            \begin{tikzcd}
              \conormalUV \ar{r} \ar{d} & \rho^* \cotangent_\formalV \ar{r} \ar{d} \ar[bend right=1.3em]{l}[swap]{r} & \cotangent_\formalU \ar{d} \\
              i^* \conormalUV \ar{r} & j^* \cotangent_\formalV \ar{r} \ar[bend right=1.3em]{l}[swap,pos=.6]{i^*(r)} & i^* \cotangent_\formalU\rlap{.}
            \end{tikzcd}
          \]
          The section $\rho^\ast dg$ of $\rho^* \cotangent_\formalV$ has a nullhomotopic image in $j^* \cotangent_\formalV$ (recall that $S \simeq \dCrit(g)$), so that $r(\rho^\ast dg)$ canonically lifts to a section $a_r$
          \[
            a_r \colon \basefield \to \derhamcomplex_{S/\formalU}^{\geq 1} \otimes_{\structuresheaf_\formalU} \conormalUV \simeq \fiber(\structuresheaf_\formalU\to \structuresheaf_S)\otimes_{\structuresheaf_\formalU} \conormalUV \simeq \fiber\left( \conormalUV \to i^* \conormalUV \right).
          \]
          Using the explicit model for $\structuresheaf_{S}$ of \cref{explicitmodulefunctionsdCrit}, the complex $\derhamcomplex_{S/\formalU}^{\geq 1}$ can be represented as $[\cdots \to \Lambda^2\tangent_\formalU \to \tangent_\formalU]$, and we thus get a surjective morphism
          \[
            \tangent_\formalU \otimes \conormalUV \twoheadrightarrow \cohomology^0(\derhamcomplex_{S/\formalU}^{\geq 1} \otimes \conormalUV).
          \]
          Denote by $b_r$ a lift of $a_r$ and by $r_0$ the twist of $r$ by $-b_r \colon \cotangent_\formalU \to \conormalUV$.
          By construction, we get a new splitting $r_0$ such that $a_{r_0}$ vanishes\footnote{The arguments of \ref{stepproofclosedimmersion} provide a structural splitting of $i^\ast \conormalUV \to j^\ast \cotangent_\formalV \to i^\ast \cotangent_U$. We suspect that this splitting coincides with $i^*(r_0)$.}.
    \item
          It follows from \ref{stepproofclosedimmersion} that $\formalV$ identifies with the formal neighborhood of $\formalU$ in $\formalV$.
          As both $\formalU$ and $\formalV$ are smooth, a standard deformation theory argument allows us to extend the tangential splitting $r_0$ to a geometric splitting: there exists an isomorphism $\alpha \colon \formalV \simeq \widehat \vectorbundle(\normal_{\formalU/\formalV})$ such that $\rho$ corresponds to the $0$\nobreakdashes-section and such that the splitting $\tangentmap{\alpha}^\dual \colon \rho^* \cotangent_\formalV \simeq \cotangent_\formalU \oplus \conormalUV$ coincides with $r_0$.

          We denote by $\pi \colon \formalV \simeq \widehat \vectorbundle(\normal_{\formalU/\formalV}) \to \formalU$ the projection to complete the diagram (where $f\pi \neq g$ in general)
          \[
            \begin{tikzcd}
              \mathllap{\widehat \vectorbundle(\normal_{\formalU/\formalV}) \simeq {}} \formalV \ar[start anchor={[xshift=-5pt,yshift=5pt]south east}]{dr}{g} \ar{d}{\pi} \\
              \formalU \ar[bend left]{u}{s_0 = \rho} \ar{r}{f} & \bbA^1\rlap{.}
            \end{tikzcd}
          \]
          Consider the commutative diagram
          \[
            \begin{tikzcd}
              & \conormalUV \ar{d} \ar{r} & i^* \conormalUV \ar{d}
              \\
              \structuresheaf_{\formalV} \ar{r} \ar{d} & \derhamcomplex^{\leq 1}_{\formalU/\formalV} \ar{r} \ar{d} & \derhamcomplex_{S/\formalV}^{\leq 1} \ar{d} & \derhamcomplex_{S}^{\leq 1} \ar{l} \ar{dl} \mathrlap{{} \ni \lambda_S}
              \\
              \structuresheaf_\formalU \ar{r}{\sim} & \derhamcomplex^{\leq 1}_{\formalU/\formalU} \ar{r} & \derhamcomplex_{S/\formalU}^{\leq 1}\rlap{.}
            \end{tikzcd}
          \]
          By construction, the function $h \coloneqq g - f \pi$ is mapped to $0$ in $\structuresheaf_\formalU$ and it therefore induces a section of $dh \colon \basefield \to \conormalUV$.
          Moreover, the Liouville structures on $\formalU$ and $\formalV$, as well as their compatibility with $\rho$ induce a nullhomotopy of $dh$ in $i^* \conormalUV$.
          Unraveling the construction, the corresponding lift to
          \[
            \derhamcomplex_{S/\formalU}^{\geq 1} \otimes_{\structuresheaf_\formalU} \conormalUV \simeq \fiber\left( \conormalUV \to i^* \conormalUV \right)
          \]
          is nothing but $a_{r_0} \sim 0$ (we use here that $r_0(dh) = r_0(dg)$).

          It follows that the Liouville structure on $\formalV$ is homotopic to the sum of the restriction of the Liouville structure of $\formalU$ along $\pi$ with the section of its Liouville complex given by
          \[
            \newcommand{\downsymbol}[1]{\ar[phantom, start anchor=center, end anchor=center]{d}[sloped]{#1}}
            \begin{tikzcd}[row sep=small]
              h \downsymbol{\in} \ar[mapsto]{r} & 0 \downsymbol{\in} \ar[mapsto]{r} & 0 \downsymbol{\in} & 0 \ar[mapsto]{l} \downsymbol{\in}
              \\
              \Symcompleted_\formalU(\conormalUV) \ar[r] \downsymbol{\simeq} & \Sym_\formalU^{\leq 1}(\conormalUV) \downsymbol{\simeq} \ar{r} & \derhamcomplex_{S/\formalV}^{\leq 1} \downsymbol{=} & \derhamcomplex_S^{\leq 1} \ar{l} \downsymbol{=}
              \\
              \Gamma(\formalV, \structuresheaf_{\formalV}) \ar[r] & \derhamcomplex_{\formalU/\formalV}^{\leq 1} \ar[r] & \derhamcomplex_{S/\formalV}^{\leq 1} & \derhamcomplex_S^{\leq 1}\rlap{.} \ar{l}
            \end{tikzcd}
          \]
          In virtue of \cref{liouvillesectionactionquad}, it is now enough to ensure that $h \colon \widehat \vectorbundle(\normal_{\formalU/\formalV}) \to \affineline{}$ is isomorphic to a quadratic form.
          This will be a consequence of the formal Morse lemma.
    \item(The function $h$ is relatively Morse)
          Consider the relative derived critical locus $\dCrit_{\formalU}(\formalV, h)$ (see \cref{definitionrelativedcrit}). Its relative tangent over $\formalU$ is the $2$\nobreakdashes-term complex
          \[
            \begin{tikzcd}
              \structuresheaf_{\dCrit_U(\formalV, h)} \otimes \normal_{\formalU/\formalV} \ar{r}{\hessian_h} &
              \structuresheaf_{\dCrit_U(\formalV, h)} \otimes \normal_{\formalU/\formalV}^\vee
            \end{tikzcd}
          \]
          where $\hessian_h$ is the Hessian of $h$.
          Because $\dCrit(\formalU, f) \simeq S \simeq \dCrit(\formalV,g)$, we have a commutative diagram where rows and columns are fiber sequences:
          \[
            \begin{tikzcd}[row sep=scriptsize]
              \tangent_S \ar[equal]{d} \ar{r} & i^* \tangent_\formalU \ar{d} \ar{r}{\hessian_f} &
              i^* \cotangent_\formalU \ar{d} \\
              \tangent_S \ar{r} \ar{d} & j^* \tangent_\formalV \ar{r}[swap]{\hessian_g} \ar{d} &
              j^* \cotangent_\formalV \ar{d} \\
              0 \ar{r} & i^* \normal_{\formalU/\formalV} \ar{r}[swap]{i^* \hessian_h} & i^* \conormalUV.
            \end{tikzcd}
          \]
          In particular, $i^* \hessian_h$ is an equivalence thus so is $\hessian_h$ by Nakayama.
          It follows that $\dCrit_\formalU(\formalV, h) \simeq \formalU$ and $h$ is thus relatively Morse (see \cref{definitionrelativelymorse}).
    \item(Conclusion)
          We now apply the formal Morse lemma (see \cref{relativeformalmorselemma}) and get a (non-linear) automorphism of $\widehat \vectorbundle(\normal_{\formalU/\formalV})$ identifying $h$ with a non-degenerate quadratic form $q$.
          Up to further étale localization, we can assume the corresponding orthogonal bundle is trivial: we get an isomorphism of quadratic bundles $\normal_{\formalU/\formalV} \simeq \formalU \times (\bbA^n, \sum x_i^2)$.
          Denoting by $\formalQ \coloneqq S_\deRham \times (\bbA^n, \sum x_i^2)$, we have $\formalV \simeq \formalU \smboxplus \formalQ$ in $\Darbstack_X^{\lambda}(S)$.
  \end{enumerate}
\end{proof}

\begin{proof}[Proof of \cref{thm-transitivity}]
  Since $\formalV$ is formally smooth and $\formalU$ is a formal thickening of $S$, \cite[\href{https://stacks.math.columbia.edu/tag/02GZ}{Tag 02GZ}]{stacks-project} guarantees the existence of a morphism of formal schemes $\theta \colon \formalU \to V$ which commutes with the inclusion of $S$:
  \[
    \begin{tikzcd}[row sep=0pt]
      & \formalU\ar{dd}{\theta} \\
      S \ar{ur}{i} \ar{dr}[swap]{j} \\
      & \formalV\rlap{.}
    \end{tikzcd}
  \]
  Étale-localizing on $S$, we can and do assume that $\cotangent_{\formalU}$ is trivial, $\cotangent_{\formalU}\simeq \structuresheaf_{\formalU}^{\oplus_{\dim \formalU}}\simeq \rho^\ast \structuresheaf_{\formalV}^{\oplus_{\dim \formalU}}$.
  In this case, we consider the object of $\Quad^\nabla(S)$ given by $\formalQ \coloneqq (\structuresheaf_{S_\derham}^{\oplus_{\dim \formalU}})^\dual \oplus \structuresheaf_{S_{\derham}}^{\oplus_{\dim \formalU}}$ equipped with the standard non-degenerate quadratic form
  \[
    \begin{pmatrix} 0 & \Id  \\ \Id & 0 \end{pmatrix}.
  \]
  It comes with an equivalence $Q \coloneqq \formalQ_{|S} \simeq i^\ast \cotangent_\formalU \oplus i^\ast \tangent_\formalU$.
  The action of \cref{actionondarbouxcharts} equips $\formalV \smboxplus \formalQ$ with a structure of Darboux chart and \cref{localformformorphisms1} renders the inclusion $s \colon \formalV \to \formalV \smboxplus \formalQ$ a morphism of Darboux charts.

  Now, the section $df$ of $\cotangent_\formalU$ combined with the zero section of $\tangent_\formalU$ gives a section
  \[
    \begin{tikzcd}[row sep=scriptsize]
      & \formalV \smboxplus \formalQ \tikzcart \ar[r] \ar{d} \tikzcart & \formalQ \ar{d}
      \\
      \formalU \ar{ur}{\mathllap{\rho \coloneqq {}}(df,0)} \ar{r}[swap]{\theta} & \formalV \ar[r] & \formalV_{\derham}\mathrlap{{} = S_{\derham}.}
    \end{tikzcd}
  \]
  and the null-homotopy $df_{|_S} \sim 0$ coming from the Darboux structure on $\formalU$ provides a 2-cell rendering the commutativity of the diagram
  \[
    \begin{tikzcd}
      S \ar{r}{j} \ar{d}[swap]{i} & \formalV \ar{d}{s}
      \\
      \formalU \ar{r}[swap]{\rho} \ar{ur}{\theta} & \formalV \smboxplus \formalQ\rlap{.}
    \end{tikzcd}
  \]
  The map $\rho$ thus obtained is not a priori a map of Darboux charts.
  The rest of the proof consists in changing the Liouville structure in $\formalV \smboxplus \formalQ$ in order to make $\rho$ a map of Darboux charts, while preserving the Liouville compatibility of $s$.
  Since the exact structure on $S$ is fixed, the failure of $\rho$ to be compatible with the Liouville structures is the failure to be compatible with the isotropic structures.
  Denote by $\eta_{\formalU}$, $\eta_{\formalV}$ and $\eta_{\formalV \smboxplus \formalQ}$ the respective isotropic structures:
  \[
    \begin{tikzcd}
      &&& \derhamcomplex_{S/\formalU}^{\geq 2}[1]
      \\
      \basefield \ar{r}[swap, near end]{\omega_S} \ar[bend right, ""{name=ETA}]{rr}[swap]{0} \ar[bend left=12pt, ""{name=ETAU}]{rrru}{0} & \derhamcomplex_S^{\geq 2}[1] \ar[r, ""{name=ETAU2,near end}]
      \ar[phantom, from=ETA, "{\phantom{\scriptstyle\eta} \rotatebox[origin=c]{90}{$\sim$} \scriptstyle \eta_{\formalV\smboxplus\formalQ}}"]
      \ar[phantom, from=ETAU, to=ETAU2, "{\phantom{\scriptstyle\eta_\formalU} \rotatebox[origin=c]{110}{$\sim$} \scriptstyle \eta_\formalU}"]
      & \derhamcomplex_{S/\formalV\smboxplus \formalQ}^{\geq 2}[1] \ar{ur}{ \rho^\ast} \ar{dr}[swap]{ s^\ast}
      \\
      &&& \derhamcomplex_{S/\formalV}^{\geq 2}[1]
    \end{tikzcd}
  \]
  with $\eta_\formalV = s^\ast \circ \eta_{\formalV\smboxplus\formalQ}$ being the horizontal composition of $\eta_{\formalV\smboxplus\formalQ}$ with $s^\ast$ (since $s$ is a morphism of Darboux charts).
  We get two null-homotopies $\eta_\formalU$ and $\rho^\ast \circ \eta_{\formalV\smboxplus\formalQ}$ of the composite $\basefield \to \derhamcomplex_{S/\formalU}^{\geq 2}[1]$.
  Their difference is a morphism $\delta = (\eta_\formalU - \rho^\ast \circ \eta_{\formalV\smboxplus\formalQ}) \colon \basefield \to \derhamcomplex_{S/\formalU}^{\geq 2}$, \ie, an element $\delta \in \operatorname{H^0}\left(\derhamcomplex_{S/\formalU}^{\geq 2}\right)$.
  Our main computation is the following:

  \begin{lemma}\label{lem-rhos0onto}
    There exists $\Delta \in \cohomology^0( \derhamcomplex_{S/\formalV\smboxplus \formalQ}^{\geq 2})$ such that
    \begin{enumerate}[label={\textrm{\upshape{(\alph*)}}}]
      \item The isotropic structure $\eta_{\formalV\smboxplus\formalQ} + \Delta$ on $\formalV\boxplus\formalQ$ is Lagrangian,
      \item $s^\ast \circ \Delta = 0$ in $\cohomology^0( \derhamcomplex_{S/\formalV}^{\geq 2})$,
      \item $\rho^\ast \circ \Delta = \delta$.
    \end{enumerate}
  \end{lemma}

  Let us first see how this lemma implies \cref{thm-transitivity}.
  Setting $\eta \coloneqq \eta_{\formalV\smboxplus\formalQ} + \Delta$ as a new Lagrangian structure on $\formalV\smboxplus\formalQ$, we get:
  \[
    \rho^\ast \circ \eta = \rho^\ast \circ \eta_{\formalV\smboxplus\formalQ} + \delta = \eta_\formalU \hspace{2em}\text{and}\hspace{2em}  s^\ast \circ \eta =  s^\ast \circ \eta_{\formalV\smboxplus\formalQ} + 0 = \eta_\formalV.
  \]
  In particular, we find a common upper bound $(\formalU, \eta_\formalU) \to (\formalV\smboxplus \formalQ, \eta) \leftarrow (\formalV, \eta_\formalV)$ in the \icategory of Darboux charts of $S$. We can now conclude using \cref{cor:map=action}.
\end{proof}

Before proving \cref{lem-rhos0onto} we will need a preliminary step related to the statement of the lemma:

\begin{construction}
  \label{twistinglagrangians}
  Let $(\rmE,\ell,\eta)$ be an isotropic distribution on $S=\Spec(A)$.
  The 2-cell $\eta$ is a null-homotopy of the composition \cref{restrictionofunderlyingformtofoliation} and the choice of any other symmetric bivector $\tau \colon \structuresheaf_S\to \Sym^2_A(\rmE^\vee[-1])$ can be understood as a self-null-homotopy of the zero morphism $0 \colon \structuresheaf_S\to \Sym^2_A(\rmE^\vee[-1])[1]$.
  Therefore, the sum $\eta+\tau$ (corresponding to the composition of null-homotopies) is a new isotropic structure on the same distribution $\ell \colon \rmE \to \tangent_S$.
  In this case the induced maps $\alpha_\eta,\alpha_{\eta+\tau} \colon  \rmN_\ell\to \rmE^\vee[-1]$ differ by the composition
  \begin{equation}\label{equationtwistinglagrangianbyhomotopy}
    \alpha_{\eta+\tau}=\alpha_\eta+\underline{\tau}\circ \partial_\ell,
    \hspace{2cm}
    \begin{tikzcd} \rmN_\ell \ar{r}{\partial_\ell} & \rmE[1] \ar{r}{\underline{\tau}} & \rmE^\vee[-1] \end{tikzcd}
  \end{equation}
  where $\underline{\tau}$ is the self-dual map induced by $\tau$, and $\partial_\ell$ the boundary operator of \cref{normalofdistribution}.

  Assume that $\eta$ is Lagrangian and let $\tau_\eta$ be the symmetric bivector given by the image of $\tau$ via the map $\Sym^2(\alpha_\eta^{-1})\colon \Sym^2(\rmE^\dual[-1])\to \Sym^2(\rmN_\ell)$.
  Then $\eta+\tau$ is Lagrangian (\ie $\alpha_{\eta+\tau}$ is an equivalence) if and only if
  \begin{equation}\label{testisequivalence}
    \alpha_\eta^{-1}\circ \alpha_{\eta+\tau}=\alpha_\eta^{-1}\circ(\alpha_\eta+\underline{\tau}\circ \partial_\ell)\underset{\cref{equationtwistinglagrangianbyhomotopy}}{=}\Id_{\rmN_\ell} + \underline{\tau}_\eta \circ \rmH_\eta
  \end{equation}
  is an equivalence, with $\rmH_\eta$ as in \cref{const:quadfromlagdist-lagrangian}. Notice that, since both $\rmH_\eta$ and $\underline{\tau}_\eta$ are self-dual, so is their composition.
\end{construction}

And a second result on connectivity estimates:

\begin{lemma}\label{lemmaconnectivityestimate}
  Let $\formalU \in \Darbstack_X^\lambda(S)$.
  Then the dg-module $\derhamcomplex_{S/\formalU}^{\geq p}$ is connective, \ie, it belongs to $\dgMod_\basefield^{\leq 0}$\footnote{in cohomological convention}.
\end{lemma}
\begin{proof}
  By definition (\cref{relativeclosedpforms}), $\derhamcomplex_{S/\formalU}^{\geq p}=|\DR(S/\formalU)|^{\geq p}[-2p]$ and because of the Lagrangian condition, as a graded module, we have
  \[\DR(S/\formalU)\simeq \bigoplus_{n\geq 0} \Sym^n_S(i^\ast \tangent_\formalU[2])=\bigoplus_{n\geq 0} \left(\Sym_S^n(i^\ast\tangent_\formalU)\right)[2n]\]
  We remark that $\left(\Sym_S^n(i^\ast\tangent_\formalU)\right)[2n][-2n]=\Sym_S^n(i^\ast\tangent_\formalU)$ is connective.  We fall in the following general situation: let $\rmE=\bigoplus_{n\geq 0} E_n \in \mixedgradedmodules_\basefield$ as in \cref{remindermixedgraded}. Suppose that for all weights $n\geq 0$, $E_n[-2n]\in \dgMod_\basefield^{\leq 0}$. Then for any $p\geq 0$, the totalization of \cref{reminderrealizations}, given by $|\rmE|^{\geq p}[-2p]$, is connective. Indeed, by construction, $|E|^{\geq p}[-2p]$ is the (product) totalization
  \[
    \totalization \left(
    \begin{tikzcd}[column sep=1.2em]
        E^p[-2p]\ar{r}{\epsilon}& E^{p+1}[-2(p+1)][1]\ar{r}{\epsilon}&\cdots\ar{r}{\epsilon}& E^{p+n}[-2(p+n)][n]\ar{r}{\epsilon}&\cdots
      \end{tikzcd}\right)\rlap{.}
  \]
\end{proof}

\begin{proof}[Proof of \cref{lem-rhos0onto}:]
  We divide the proof in two parts:
  \begin{enumerate}[leftmargin=0pt,itemindent=5em,label={$\bullet$ Step \arabic*.}, ref={Step \arabic*}]
    \item \label{firststeptransitivity} We first show that the morphism of mixed graded cdga's
          \[
            \chi \coloneqq (\rho^\ast, s^\ast) \colon \DR(S/\formalV\smboxplus\formalQ) \to \DR(S/\formalU) \oplus \DR(S/\formalV)
          \]
          admits a section $\sigma$ as a map of graded dg-modules (not compatible with the mixed structure). Indeed, by \cref{constructionrelativederham}, as map of graded modules, $\chi$ is induced by the map in weight $1$
          \[
            \chi^{(1)} \coloneqq \begin{pmatrix} \tangent_{S/s}^\vee[1] & \tangent_{S/\rho}^\vee[1] \end{pmatrix} \colon
            \begin{tikzcd}[cramped]
              \cotangent_{S/\formalV\smboxplus\formalQ}[1]\ar{r} & \cotangent_{S/\formalV}[1]\oplus \cotangent_{S/\formalU}[1]
            \end{tikzcd}
          \]
          where $\tangent_{S/s}$ and $\tangent_{S/\rho}$ are as in \cref{diagrammapinducedatlevelofdistributionsbyuniversalproperty}. Using the Lagrangian condition independently for $\formalV\oplus \formalQ$, $\formalU$, and $\formalV$, $\chi^{(1)}[-2]$ can be written as a block matrix
          \[
            \begin{pmatrix} A & B \\ C & D \end{pmatrix}\colon
            \begin{tikzcd}
              j^\ast \tangent_\formalV \oplus Q \ar{r} & j^\ast \tangent_\formalV \oplus i^\ast \tangent_\formalU
            \end{tikzcd}
          \]
          Now, since by construction $s \colon \formalV\to \formalV \smboxplus \formalQ$ is the $0$\nobreakdashes-section, we have $A = \Id$ and $B = 0$.
          It remains to describe the map $C$ and $D$. Recall that by construction we have $Q \simeq i^\ast\cotangent_\formalU \oplus i^\ast\tangent_\formalU$. We split $D$ in a block matrix $D = \begin{pmatrix} D_\cotangent & D_\tangent\end{pmatrix} \colon Q \simeq i^\ast\cotangent_\formalU \oplus i^\ast\tangent_\formalU \to i^\ast \tangent_\formalU$.
          Now, by compatibility with the boundary operators, and using \cref{const:quadfromlagdist-lagrangian}, we find a commutative diagram
          \begin{equation}\label{mapofdistributionsnotfoliations}
            \begin{tikzcd}[ampersand replacement=\&]
              i^\ast \tangent_\formalU \arrow[rounded corners, to path={|- ([yshift=1em]\tikztotarget.north)[near end]\tikztonodes -- (\tikztotarget)}]{rr}{\hessian_f} \ar{r}{\partial_{\ell_\formalU}}\ar{d}{\tangentmap{\rho}}
              \&
              \tangent_{S/\formalU}[1] \ar{r}{\alpha_{\eta_\formalU}^\dual}[swap]{\sim} \ar{d}{\tangent_{S/\rho}[1]}
              \&
              i^\ast \cotangent_\formalU \ar{d}{\left(\begin{smallmatrix} C^\dual \\ D_\cotangent^\dual \\ D_\tangent^\dual \end{smallmatrix}\right) }
              \\
              j^\ast\tangent_\formalV\oplus i^\ast\cotangent_\formalU \oplus i^\ast\tangent_\formalU \ar{r}[swap]{\partial_{\ell_{\formalV\smboxplus\formalQ}}} \arrow[swap, rounded corners, to path={|- ([yshift=-1em]\tikztotarget.south)[near end]\tikztonodes -- (\tikztotarget)}]{rr}{\mathllap{\hessian_{g\smboxplus q}={}}\begin{pmatrix} \hessian_g&0&0\\ 0&0&\Id\\ 0&\Id&0\end{pmatrix}}
              \&
              \tangent_{S/\formalV\smboxplus\formalQ}[1] \ar{r}{\sim}[swap]{\alpha^\dual_{\eta_{\formalV\smboxplus\formalQ}}}
              \&
              j^\ast\cotangent_\formalV\oplus i^\ast\tangent_\formalU \oplus i^\ast\cotangent_\formalU
            \end{tikzcd}
          \end{equation}
          By definition of $\rho$, the tangent map $\tangentmap{\rho}$ is the block matrix $\tangentmap{\rho} = \begin{pmatrix} \tangentmap{\theta} & \hessian_f & 0\end{pmatrix}$,
          where $\tangentmap{\theta}$ the tangent morphism of $\theta\colon\formalU\to\formalV$.
          Finally, the map of fiberwise distributions induced by $\theta$, $\tangent_{S/\theta}\colon \tangent_{S/\formalU}\to \tangent_{S/\formalV}$, allows us to complete diagram \cref{mapofdistributionsnotfoliations}
          \[
            \begin{tikzcd}[ampersand replacement=\&,row sep=3em, column sep=4em]
              i^* \tangent_\formalU \ar{r}{\hessian_f} \ar{d}[swap]{\left(\begin{smallmatrix} \tangentmap{\theta} \\ \hessian_f \\ 0  \end{smallmatrix}\right)}
              \ar[dd, "\tangentmap{\theta}", rounded corners, to path={ -| ([xshift=-4em]\tikztotarget.west)[near end,swap]\tikztonodes
                    -- (\tikztotarget)}]
              \&
              i^* \cotangent_\formalU \ar{d}{\left(\begin{smallmatrix} C^\dual \\ D_\cotangent^\dual \\ D_\tangent^\dual \end{smallmatrix}\right)}
              \ar[dd, "{\tau \coloneqq \alpha^\dual_{\eta_{\formalV}}\circ \tangent_{S/\theta}[1]\circ (\alpha_{\eta_\formalU}^\dual)^{-1}}", rounded corners, to path={ -| ([xshift=4em]\tikztotarget.east)[near end]\tikztonodes
                  -- (\tikztotarget)}]
              \\
              j^* \tangent_\formalV\oplus i^\ast\cotangent_\formalU \oplus i^\ast\tangent_\formalU \ar{r}{\left(\begin{smallmatrix} \hessian_g & 0 & 0 \\ 0& 0 & \Id \\ 0 & \Id & 0 \end{smallmatrix}\right)} \ar{d}[swap]{\left(\begin{smallmatrix} \Id & 0 & 0 \end{smallmatrix}\right)}
              \&
              j^* \cotangent_\formalV\oplus i^\ast\tangent_\formalU \oplus i^\ast\cotangent_\formalU
              \ar{d}{\left(\begin{smallmatrix} \Id & 0 & 0\end{smallmatrix}\right)}
              \\
              j^* \tangent_\formalV \ar{r}[swap]{\hessian_g}
              \&
              j^* \cotangent_\formalV
            \end{tikzcd}
          \]
          where the outer square is a pullback since the fibers of the horizontal arrows are constant equal to $\tangent_S$ as in \cref{lagrangianconditiondCritHessian}. Since the lower square is obviously a pullback so must the upper square be.
          We also have an obviously Cartesian square
          \[
            \begin{tikzcd}[ampersand replacement=\&,row sep=3em, column sep=4em]
              i^* \tangent_\formalU \ar{r}{\hessian_f} \ar{d}[swap]{\left(\begin{smallmatrix} \tangentmap{\theta} \\ \hessian_f \\ 0  \end{smallmatrix}\right)}
              \&
              i^* \cotangent_\formalU \ar{d}{\left( \begin{smallmatrix} \tau \\ 0 \\ \Id \end{smallmatrix} \right)}
              \\
              j^* \tangent_\formalV\oplus i^\ast\cotangent_\formalU \oplus i^\ast\tangent_\formalU \ar{r}{\left(\begin{smallmatrix} \hessian_g & 0 & 0 \\ 0& 0 & \Id \\ 0 & \Id & 0 \end{smallmatrix}\right)}
              \&
              j^* \cotangent_\formalV\oplus i^\ast\tangent_\formalU \oplus i^\ast\cotangent_\formalU\rlap{.}
            \end{tikzcd}
          \]
          This implies $C = \tau^\dual$, $D_\tangent = \Id$ and $D_\cotangent = 0$.
          In conclusion, we have
          \[
            \chi^{(1)}[2] =
            \begin{pmatrix}
              \Id        & 0 & 0   \\
              \tau^\dual & 0 & \Id
            \end{pmatrix} \colon j^* \tangent_\formalV \oplus i^\ast\cotangent_\formalU \oplus i^\ast\tangent_\formalU \to j^* \tangent_\formalV \oplus i^\ast \tangent_\formalU.
          \]
          It admits a section $\sigma$ given by the matrix
          \[
            \sigma \coloneqq
            \begin{pmatrix}
              \Id         & 0   \\
              0           & 0   \\
              -\tau^\dual & \Id
            \end{pmatrix} \colon j^* \tangent_\formalV \oplus i^\ast \tangent_\formalU \to j^* \tangent_\formalV \oplus i^\ast\cotangent_\formalU \oplus i^\ast\tangent_\formalU.
          \]
          This implies that, forgetting the mixed structure, the morphism of graded complexes $\chi = \Sym(\chi^{(1)})$ admits a section.

    \item  Since the realizations of \cref{reminderrealizations} preserves fibers and direct sums (being a right adjoint), $\chi$ induces morphisms
          \[
            \chi_p\colon\derhamcomplex_{S/\formalV\smboxplus\formalQ}^{\geq p}\to \derhamcomplex_{S/\formalU}^{\geq p}\oplus \derhamcomplex_{S/\formalV}^{\geq p}.
          \]
          By \cref{lemmaconnectivityestimate} both source and target of $\chi_p$ are connective. We now observe that $\fiber(\chi_p) \simeq |\fiber(\chi)|^{\geq p}[-2p]$ remains connective.
          Indeed, the fiber of $\chi$ is a mixed graded module $E$ given in weight $n\geq 0$ by
          \[
            E^n=\fiber\left(\Sym^n_S(j^\ast \tangent_\formalV\oplus Q) \to^{\chi} \Sym^n_S(j^\ast \tangent_\formalV)\oplus \Sym^n_S(j^\ast \tangent_\formalU)\right)[2n].
          \]
          Using the section of $\chi$ of \ref{firststeptransitivity}, we see that $E^n[-2n]$ is connective and we can conclude as in the proof of \cref{lemmaconnectivityestimate}. We obtain a morphism between the fiber sequences \cref{gradedpiecestruncatedderhamfiltration} and thus a morphism of (vertical) exact sequences
          \begin{equation}\label{diagramsurjectivearrows}
            \begin{tikzcd}[row sep=0.5cm, column sep=6em]
              \ar[two heads]{r}{\chi_3} \cohomology^0(\derhamcomplex_{S/\formalV\smboxplus\formalQ}^{\geq 3}) \ar{d}
              &
              \cohomology^0(\derhamcomplex_{S/\formalU}^{\geq 3})\oplus \cohomology^0(\derhamcomplex_{S/\formalV}^{\geq 3})\ar{d}
              \\
              \ar[two heads]{r}{\chi_2} \cohomology^0(\derhamcomplex_{S/\formalV\smboxplus\formalQ}^{\geq 2}) \ar[two heads]{d}
              &
              \cohomology^0(\derhamcomplex_{S/\formalU}^{\geq 2})\oplus \cohomology^0(\derhamcomplex_{S/\formalV}^{\geq 2}) \ar[two heads]{d}
              \\
              \ar[two heads]{r}[swap]{\chi^{(2)} = \Sym^2(\chi^{(1)})} \cohomology^0(\Sym^2_S(j^\ast \tangent_\formalV \oplus Q))
              &
              \cohomology^0(\Sym^2_S(i^\ast \tangent_\formalU ))\oplus \cohomology^0(\Sym^2_S(j^\ast \tangent_\formalV)) \ar[bend right=1em, start anchor={[yshift=.6em]west}, end anchor={[yshift=.6em]east}]{l}[swap]{\sigma^{(2)} \coloneqq \Sym^2(\sigma)}
            \end{tikzcd}
          \end{equation}
          where the double-headed arrows are surjections and $\sigma^{(2)} \coloneqq \Sym^2(\sigma)$ is the section of $\chi^{(2)}$ obtained in \ref{firststeptransitivity}.
    \item
          Consider now the element
          \[
            (\delta,0) \in \operatorname{H^0}\left( \derhamcomplex_{S/\formalU}^{\geq 2}\right)\oplus \cohomology^0(\derhamcomplex_{S/\formalV}^{\geq 2})
          \]
          and denote by $(\tau_{\eta_\formalU}, 0)$ its image in $\cohomology^0(\Sym^2_S(i^\ast \tangent_\formalU ))\oplus \cohomology^0(\Sym^2_S(j^\ast \tangent_\formalV))$.
          The element $\tau_{\eta_\formalU}$ is the bivector induced by $\delta$ on the underlying distribution of $\formalU$ (\cref{distributionassociatedtofibration}), as in \cref{twistinglagrangians}.
          A straightforward computation shows that the element $\sigma^{(2)}(\tau_{\eta_\formalU},0)$ coincides with the image of $\tau_{\eta_\formalU}$ via the canonical inclusion:
          \[
            \cohomology^0(\Sym^2_S(i^\ast \tangent_\formalU ))\subseteq \cohomology^0(\Sym^2_S(j^\ast \tangent_\formalV \oplus i^\ast \cotangent_\formalU \oplus i^\ast \tangent_\formalU)) = \cohomology^0(\Sym^2_S(j^\ast \tangent_\formalV \oplus Q)),
          \]
          \ie,
          \begin{equation}\label{formulaforsectionexplicit}
            \sigma^{(2)}(\tau_{\eta_\formalU},0) = \sigma
            \begin{pmatrix}
              \tau_{\eta_\formalU} & 0 \\ 0 & 0
            \end{pmatrix}
            \sigma^\vee
            =
            \begin{pmatrix}
              0 & 0 & 0 \\ 0 & \tau_{\eta_\formalU} & 0 \\ 0 & 0 & 0
            \end{pmatrix}.
          \end{equation}
          Now, denote by $\Delta_1$ a lift of $\sigma^{(2)}(\tau_{\eta_\formalU},0)$ to $\cohomology^0(\derhamcomplex_{S/\formalV\smboxplus\formalQ}^{\geq 2})$ in \cref{diagramsurjectivearrows}. By construction $(\delta,0) - \chi_2(\Delta_1)$ lifts to some class $a \in   \cohomology^0(\derhamcomplex_{S/\formalU}^{\geq 3})\oplus \cohomology^0(\derhamcomplex_{S/\formalV}^{\geq 3})$, which in turn lifts to $d \in \cohomology^0(\derhamcomplex_{S/\formalV\smboxplus\formalQ}^{\geq 3})$. We denote by $\Delta_2$ the image of $d$ in $\cohomology^0(\derhamcomplex_{S/\formalV\smboxplus\formalQ}^{\geq 2})$, and set $\Delta = \Delta_1 + \Delta_2$.
          By design, we have $\chi_2(\Delta) = (\delta, 0)$.
          This implies assertions (b) and (c) from the lemma.

          It remains to show (a) claiming that the isotropic structure $\eta_{\formalV\smboxplus\formalQ}+\Delta$ is Lagrangian.
          Since this a property of the underlying isotropic distribution, we find ourselves in the context of \cref{twistinglagrangians} and it suffices to show that the map \cref{testisequivalence} is an equivalence.
          By design, the image of $\Delta$ in $\cohomology^0(\Sym^2_S(j^\ast \tangent_\formalV \oplus Q))$ coincides with the image of $\Delta_1$ and thus recovers $\sigma^{(2)}(\tau_\eta, 0)$ described by the matrix of \cref{formulaforsectionexplicit}.
          Therefore, we are reduce to showing that the map
          \[
            \Id + \underline{\sigma}^{(2)}(\tau_{\eta_\formalU}, 0) \circ \rmH_{\eta_{\formalV\smboxplus\formalQ}}
          \]
          is an equivalence. It follows from the computation (relying on \cref{bilinearformisthehessian}):
          \[
            \underline{\sigma}^{(2)}(\tau_{\eta_\formalU}, 0) \circ \rmH_{\eta_{\formalV\smboxplus\formalQ}}=
            \begin{pmatrix} 0 & 0 & 0 \\ 0 & \tau_\delta & 0 \\ 0 & 0 & 0 \end{pmatrix}
            \begin{pmatrix}  \hessian_g & 0& 0 \\ 0 & 0 & \Id \\ 0 & \Id & 0 \end{pmatrix}
            =
            \begin{pmatrix} 0 & 0 & 0 \\ 0 & 0 & \tau_\delta \\ 0 & 0 & 0 \end{pmatrix}.
          \]
  \end{enumerate}
  This concludes the proof of \cref{lem-rhos0onto} and thus of \cref{thm-transitivity}.
\end{proof}

\section{\texorpdfstring{$\affineline{}$}{A¹}-Isotopies of Darboux Charts}
\label{sectionisotopies}

In \cref{thm-transitivity} above we saw that the quotient  $\Darbstack^{\lambda}_X/\stackQuadnabla_X$ is connected. Unfortunately, as \cref{ambiguitymorphisms} shows,  this quotient stack is not contractible. \Cref{isotopiessolveambiguity} suggests the next step: we need to identify isotopic maps of Darboux charts.
In  \cref{subsectionclassicalA1localization} we review the classical formulas of $\affineline{}$\nobreakdashes-localization as introduced in \cite{voevodsky-morel} and in  \cref{subsectionisotopies} we explain our new general procedure of isotopic localization. In \cref{darbouxstackuptoisitopysection}  we construct the new stack of Darboux charts  and quadratic bundles up to isotopies and finally in  \cref{sectionproofcontractibility} we prove  \cref{theoremcontractibility}.

\subsection{\texorpdfstring{Reminders on the classical $\affineline{}$}{A¹}-localization}
\label{subsectionclassicalA1localization}

Before introducing the new moduli of Darboux charts where isotopic maps are identified (\cref{darbouxstackuptoisitopysection} below), we need some preliminaries, starting with a well-known notion of $\affineline{}$\nobreakdashes-invariance.
Throughout this section we fix $Y$ a derived stack.
\begin{definition}\label{localA1MV}
  Let $G \in \PSh(\dAff_Y, \inftygpd)\simeq \PSh(\dAff_\basefield, \inftygpd)_Y$\footnote{cf. \cref{descriptionsofslicecategoriesofstacks}-(i).}.
  Recall from \cite{voevodsky-morel} that $G$ is said to be $\affineline{Y}$\nobreakdashes-local if the canonical morphism induced by composition with $\affineline{Y} \to Y$
  \[
    G \simeq \stackMap_Y(Y, G) \to \stackMap_Y(\affineline{Y},G)
  \]
  is an equivalence of presheaves. This is by definition equivalent to asking that for any $S \in \dAff_Y$, the morphism $G(S)\to G(S\times \affineline{\basefield})$ is an equivalence. Let $\ipshA{Y}\colon\PSh_{\affineline{Y}}(\dAff_Y, \inftygpd) \subseteq \PSh(\dAff_Y, \inftygpd)$ denote the inclusion of the full subcategory of $\affineline{Y}$\nobreakdashes-local presheaves.
  By the results of \cite{voevodsky-morel} (see also \cite{MR3281141}) this inclusion admits a left adjoint $\LpshA{Y}$. We denote its unit by $\eta_{\affineline{}} \colon \Id \to \LpshA{Y}$
\end{definition}

\begin{reminder}\label{reminderexplicitmodelusualA1localization}
  As explained in \cite{voevodsky-morel}, the functor $\LpshA{Y}$ of \cref{localA1MV} admits an explicit model using the \emph{algebraic simplices}
  \[
  \simplex^n_\basefield \coloneqq \Spec\left(\quot{\basefield[x_0,\cdots x_n]}{(x_0+\cdots + x_n=1)}\right).
  \]
  They fit in a cosimplicial diagram $\simplex^\bullet_\basefield \colon \Delta \to \dAff_\basefield$.
  For $G \in \PSh(\dAff_Y, \inftygpd)$, we have an equivalence in $\PSh(\dAff_Y, \inftygpd)$
  \[
   \ipshA{Y}\circ \LpshA{Y}(G) \simeq \colim_{[n]\in \Delta^\op} \stackMap_Y(\simplex^n_\basefield \times Y, G)\rlap{.}
  \]
  See \cite[\S 4.3]{Antieau2017} for details.
\end{reminder}

\begin{notation}\label{notationforrestrictionpresheavescategories}
  Let $\calC$ and $\calD$ be \icategories and  $G \colon \calC^\op \to \calD$ an \ifunctor. Given $S \in \calC$, we consider the restriction $G_{|_S}$ defined by the composition $(\calC_{/S})^\op \to \calC^\op \to \calD$.
\end{notation}

\begin{observation}\label{naiveA1localizationcommuteswithrestriction}
  Given $G \in \PSh(\dAff_Y, \inftygpd)$ and $S \in \dAff_Y$ we have
  \begin{equation}\label{restrictionformulanaive}
    \stackMap_Y(\affineline{Y}, G)_{|_S} \simeq \stackMap_S(\affineline{S}, G_{|_S})
  \end{equation}
  In particular, the restriction $(-)_{|_S} \colon \PSh(\dAff_Y, \inftygpd) \to \PSh(\dAff_S, \inftygpd)$ sends $\affineline{Y}$\nobreakdashes-local objects to $\affineline{S}$\nobreakdashes-local objects.
  Finally, we observe that $(-)_{|_S}$ commutes with the $\affineline{}$\nobreakdashes-localization formulas:
  \begin{multline*}
    \left(\ipshA{Y}\circ \LpshA{Y}(G)\right)_{|_S} \simeq \colim_{[n]\in \Delta} \stackMap_Y(\simplex^n_\basefield \times Y, G)_{|_S}
    \\
    \simeq \colim_{[n] \in \Delta} \stackMap_S(\simplex^n_\basefield \times S, G_{|_S}) \simeq \ipshA{S}\circ \LpshA{S}(G_{|_S})\rlap{.}
  \end{multline*}
\end{observation}

\begin{observation}\label{intersectionoflocalizations}
  Let $\calC$ be a presentable \icategory and let
  \[
    \begin{tikzcd}
      \calC_0\ar[hook]{r}[swap]{i_0} & \calC \ar[bend right=1.8em]{l}[swap]{\rmL_0}&   \calC_1\ar[hook]{r}[swap]{i_1} & \calC \ar[bend right=1.8em]{l}[swap]{\rmL_1}
    \end{tikzcd}
  \]
  be accessible reflexive localizations\footnote{See \cite[5.4.2.5, 5.5.4.2]{lurie-htt}} and consider the intersection
  \[
    \begin{tikzcd}
      \calC_0\cap \calC_1 \ar[hook]{d}{i_1'} \ar[hook]{r}{i_0'}
      & \calC_1 \ar[hook]{d}{i_1} \\
      \calC_0 \ar[hook]{r}[swap]{i_0} & \calC
    \end{tikzcd}
  \]
  By \cite[5.5.4.17, 5.5.4.18]{lurie-htt} each inclusion  admits a left adjoint $\rmL_0' \colon \calC_0 \to \calC_0 \cap \calC_1$, $\rmL_1' \colon \calC_1 \to \calC_0 \cap \calC_1$ and therefore the square of left adjoints commutes. In particular, this exhibits the inclusion $i \coloneqq i_1 \circ i_0' = i_0\circ i_1' \colon \calC_0 \cap \calC_1 \hookrightarrow \calC$ as an accessible reflexive localization with left adjoint $\rmL_1' \circ \rmL_0 = \rmL_0' \circ \rmL_1$.
\end{observation}

\begin{construction}\label{classicalA1localizationofsheaves}
  Consider the pullback
  \[
    \begin{tikzcd}[column sep=small]
      \Sh_{\affineline{Y}}(\dAff_Y, \inftygpd) \tikzcart \ar[hook]{r}{\ishA{Y}} \ar[hook]{d}& \Sh(\dAff_Y, \inftycats) \ar[hook]{d}{\ish{Y}} \\
      \PSh_{\affineline{Y}}(\dAff_Y, \inftygpd) \ar[hook]{r}{\ipshA{Y}} &  \PSh(\dAff_Y, \inftygpd)
    \end{tikzcd}
  \]
  where $\Sh(\dAff_Y, \inftygpd)$ denotes the \itopos of étale hypersheaves.
  We denote by $\LshA{Y}$ the left adjoint to the inclusion $\ishA{Y}\colon \Sh_{\affineline{Y}}(\dAff_Y, \inftygpd) \subseteq \Sh(\dAff_Y, \inftygpd)$, given by \cref{intersectionoflocalizations}.
\end{construction}

We will need the following construction to describe the functor $\LshA{Y}$ explicitly:

\begin{construction}\label{endofunctortransfiniteinduction}
  Let $\calC$, $\calC_0$ and $\calC_1$ be as in \cref{intersectionoflocalizations} and assume without loss of generality that $\calC_0$ and $\calC_1$ are both $\kappa$\nobreakdashes-accessible reflexive localizations for $\kappa$ a regular cardinal.
  Denote by $\eta_0 \colon \Id \to i_0 \circ \rmL_0$ and $\eta_1 \colon \Id \to i_1 \circ \rmL_1$ the unit natural transformations and consider the endofunctor $\Phi \coloneqq i_1 \circ \rmL_1 \circ i_0 \circ \rmL_0 \colon \calC \to \calC$ together with the natural transformation
  \[
    \begin{tikzcd}[column sep=large]
      \nu \colon \Id_{\calC} \ar{r}{i_1\circ \rmL_1 (\eta_0) \circ \eta_1} & \Phi
    \end{tikzcd}
  \]
  Let $\Seq[\kappa]$ be the poset of ordinal numbers $\gamma<\kappa$. Since $\kappa$ is regular, $\Seq[\kappa]$  is $\kappa$\nobreakdashes-small and $\kappa$\nobreakdashes-filtered. As in \cite[page 76]{voevodsky-morel} and since $\calC$ admits all $\kappa$\nobreakdashes-small colimits, we consider the $\Seq[\kappa]$\nobreakdashes-diagram of endofunctors of $\calC$ defined by transfinite induction as follows:
  \begin{itemize}
    \item If $\gamma=0$, $\Phi^0=\Id_\calC$;
    \item If $\gamma = \gamma' + 1$, we set $\Phi^{\gamma} \coloneqq \Phi \circ \Phi^{\gamma'}$;
    \item If $\gamma$ is a limit ordinal, we set $\Phi^{\gamma} \coloneqq \colim_{\gamma'<\gamma} \Phi^{\gamma'}$ with transition maps given by the iterations of $\nu$.
  \end{itemize}
  Finally, we define an endofunctor $\Phi^{\kappa}$ of $\calC$ as the colimit of this $\Seq[\kappa]$\nobreakdashes-diagram
  \begin{equation}\label{colimititerationformula}
    \Phi^{\kappa} \coloneqq \colim_{\gamma<\kappa} \Phi^{\gamma}
  \end{equation}
  together with the natural transformation $\mu \colon \Id_{\calC} \to \Phi^{\kappa}$.
\end{construction}

\begin{lemma}\label{formulaforiteratedlocalization}
  Let $\calC$, $\calC_0$ and $\calC_1$ be as in \cref{endofunctortransfiniteinduction}. Then the composition $i \circ \rmL \colon \calC \to \calC$ coincides with the endofunctor $\Phi^{\kappa} \colon \calC\to \calC$ of \cref{endofunctortransfiniteinduction}.
  In particular, the composition $i_0' \circ \rmL_0' \colon \calC_1\to \calC_1$ coincides with the composition
  \begin{equation}\label{formulaiteratedfinal}
    \rmL_1 \circ \Phi^{\kappa} \circ i_1 =(\rmL_1\circ i_0\circ \rmL_0\circ i_1)^{\kappa} \colon \calC_1\to \calC_1.
  \end{equation}
\end{lemma}
\begin{proof}
  A sketch of proof can be found in \cite[Lem.\,4.42.]{MR3281141} in the particular case of Nisnevich and $\affineline{}$\nobreakdashes-local sheaves. We give here a context-independent complete argument since we need it for hypersheaves on the big étale site of derived affine schemes. Let $X \in \calC$. By cofinality, the canonical map
  \[
    \Phi^{\kappa}(\nu_X) \colon \Phi^{\kappa}(X)\simeq \Phi^{\kappa}(\Phi(X))
  \]
  is an equivalence.
  Since $i_1$ is fully faithful and preserves $\kappa$\nobreakdashes-filtered colimits, we conclude that $\Phi^{\kappa}(X)$ belongs to the essential image of $i_1$.
  The commutativity of the diagram of natural transformations
  \begin{equation}\label{diagramnaturaltransformationsiteration}
    \begin{tikzcd}[column sep=large]
      \ar{dr}{\nu}\Id \ar{d}{\eta_0} \ar{r}{\eta_1}& i_1\circ \rmL_1\ar{d}{i_1\circ \rmL_1(\eta_0)}\\
      i_0 \circ \rmL_0 \ar{r}{\eta_1 \circ \Id_{i_0 \circ \rmL_0}} & i_1 \circ \rmL_1 \circ i_0 \circ \rmL_0 \mathrlap{{} \eqqcolon \Phi}
    \end{tikzcd}
  \end{equation}
  implies that by cofinality we can also write
  \[
    \Phi^{\kappa}(X) \simeq \Psi^{\kappa}(i_0 \circ \rmL_0(X))
  \]
  where $\Psi \coloneqq i_0\circ \rmL_0 \circ i_1\circ \rmL_1$ and $\Psi^{\kappa}$ is defined as in \cref{endofunctortransfiniteinduction}. Since $i_0$ is fully faithful and preserves $\kappa$\nobreakdashes-filtered colimits, it follows that $\Phi^{\kappa}(X)$ belongs to the essential image of $i_0$.
  We conclude that $\Phi^{\kappa}(\calC)\subseteq \calC_0\cap \calC_1$. In particular, if $\eta \colon \Id\to i \circ \rmL$ denotes the unit of the adjunction, the induced natural transformation
  \[
    \eta \circ \Id_{\Phi^{\kappa}} \colon \Phi^{\kappa} = \Id \circ \Phi^{\kappa} \to i \circ \rmL \circ \Phi^{\kappa}
  \]
  is an equivalence.

  Now, assume $X \in \calC_0 \cap \calC_1$.
  In particular, when seen as an object of $\calC$, both maps $\eta_0 \colon X \to i_0 \circ \rmL_0(X)$ and $\eta_1 \colon X \to i_1 \circ \rmL_1(X)$ are equivalences.
  In this case all the maps in \cref{diagramnaturaltransformationsiteration} are equivalences and so is $\nu_X \colon X\to \Phi(X)$ and therefore, so is $\mu_X \colon X\to \Phi^{\kappa}(X)$.
  This shows $\calC_0\cap \calC_1\subseteq \Phi^{\kappa}(\calC)$.
  We get $\Phi^{\kappa}(\calC) = \calC_0\cap \calC_1$.
  In particular, both natural transformations
  \[
    \mu \circ \Id_{i \circ \rmL} \colon i \circ \rmL \to \Phi^{\kappa} \circ i \circ \rmL
    \hspace{1em}\text{and}\hspace{.5em}
    \mu \circ \Id_{\Phi^{\kappa}} \colon \Phi^{\kappa}\to \Phi^{\kappa} \circ \Phi^{\kappa}
  \]
  are equivalences.
  In summary, for every $X\in \calC$ the natural transformations $\mu$ and $\eta$ (with its universal property exhibiting $\rmL$ as a left adjoint) provides commutative diagrams
  \[
    \begin{tikzcd}
      X \ar{r}{\eta_X} \ar{d}{\mu_X} & i \circ \rmL(X) \ar{dd}[near end]{i\circ\rmL(\mu_X)} \ar{dr}[near start]{\mu_{i\circ \rmL(X)}}[swap,sloped]{\sim} &
    \\
      \Phi^{\kappa}(X) \ar[crossing over]{rr}[near start]{\Phi^{\kappa}(\eta_X)} \ar{dr}[swap]{\eta_{\Phi^{\kappa}(X)}}[sloped]{\sim} && \Phi^{\kappa}(i\circ \rmL(X))
    \\
      & i \circ \rmL(\Phi^{\kappa}(X))\ar[dashed]{ur} &
    \end{tikzcd}
  \]
  In particular, we find that the following conditions are equivalent:
  \begin{enumerate}[label=(\alph*)]
    \item the map $\mu_X \colon X\to \Phi^{\kappa}(X)$ is an $\rmL$\nobreakdashes-equivalence, \ie $i\circ \rmL(\mu_X)$ is an equivalence;
    \item the dashed map is an equivalence;
    \item $\eta_X$ is a $\Phi^{\kappa}$\nobreakdashes-equivalence, \ie $\Phi^{\kappa}(\eta_X)$ is an equivalence.
  \end{enumerate}
  We conclude the proof by showing that indeed $\mu_X$ is an $\rmL$\nobreakdashes-equivalence.
  By construction, it suffices to show that $\nu_X \colon X \to \Phi(X)$ is an $\rmL$\nobreakdashes-equivalence. This is in turn tantamount to proving that for every object $Y\in \calC_0\cap \calC_1$, the composition
  \[
    \Map_{\calC}(\Phi(X), i_1\circ i_0'(Y))\to \Map_{\calC}(X,i_1\circ i_0'(Y))
  \]
  is an equivalence. We compute
  \begin{align*}
    \Map_{\calC}(\Phi(X)&, i_1 \circ i_0'(Y))
    \simeq \Map_{\calC}(i_1 \circ \rmL_1 \circ i_0 \circ \rmL_0(X), i_1 \circ i_0'(Y))
    \\& \simeq \Map_{\calC_1}(\rmL_1 \circ i_0 \circ  \rmL_0(X), i_0'(Y))
    \simeq \Map_{\calC}(i_0 \circ \rmL_0(X), i_1\circ i_0'(Y))
    \\& \simeq \Map_{\calC}(i_0 \circ  \rmL_0(X), i_0 \circ i_1'(Y))
    \simeq \Map_{\calC_0}(\rmL_0(X), i_1'(Y))
    \\& \simeq \Map_{\calC}(X,  i_0 \circ i_1'(Y)).
  \end{align*}
  The assertion concerning the formula \cref{formulaiteratedfinal} is now a direct computation using the fact both $i_1$ and $\rmL_1$ preserve $\kappa$\nobreakdashes-filtered colimits by assumption.
\end{proof}

\begin{corollary}\label{preservesproductsifbothdo}
  In the context of \cref{formulaforiteratedlocalization}, if both $\rmL_0$ and $\rmL_1$ preserve finite products, so do $\rmL_0'$ and  $\rmL$.
\end{corollary}
\begin{proof}
  The formula for $\rmL_{0}'$ established in \cref{formulaforiteratedlocalization} in terms of the $\kappa$\nobreakdashes-filtered colimit \cref{colimititerationformula}, implies (as $\kappa$\nobreakdashes-filtered colimits preserve $\kappa$\nobreakdashes-small limits \cite[5.3.3.3]{lurie-htt} and in particular, finite products) that $\rmL^{'}_{0}$ preserves finite products.
  Therefore, so does the composition $\rmL=\rmL_0'\circ \rmL_1$.
\end{proof}

\begin{observation}\label{lemmaappliedtoourcase}
  The localization
  \[
  \begin{tikzcd}
    \Sh(\dAff_Y, \inftygpd)\ar[hookrightarrow]{r}[swap]{\ish{Y}}&\ar[bend right=10]{l}[swap]{\Lsh{Y}}\PSh(\dAff_Y, \inftygpd)
    \end{tikzcd}
  \]
  is $\kappa$\nobreakdashes-accessible for some regular cardinal $\kappa$. At the same time, the reflexive localization
  \[
    \begin{tikzcd}
      \PSh_{\affineline{Y}}(\dAff_Y, \inftygpd) \ar[hookrightarrow]{r}[swap]{\ipshA{Y}}&\ar[bend right=10]{l}[swap]{\LpshA{Y}} \PSh(\dAff_Y, \inftygpd)
    \end{tikzcd}
  \]
  of \cref{localA1MV} is $\omega$\nobreakdashes-accessible since $\affineline{Y}$\nobreakdashes-local presheaves are stable under filtered colimits.
  In particular it is also $\kappa$\nobreakdashes-accessible for any regular cardinal $\kappa \geq \omega$. Therefore, \cref{formulaforiteratedlocalization} with $\rmL_0=\LpshA{Y}$, $i_0=\ipshA{Y}$, $\rmL_1=\Lsh{Y}$, $i_1=\ish{Y}$, implies that
  \[
  \ishA{Y}\circ \LshA{Y}\simeq \PhiA{Y}^\kappa: \Sh(\dAff_Y, \inftygpd)\to \Sh(\dAff_Y, \inftygpd)
  \]
  where $\PhiA{Y} \coloneqq \Lsh{Y}\circ \ipshA{Y}\circ \LpshA{Y}\circ \ish{Y}$.
\end{observation}

\begin{lemma}\label{A1localizationforsheavespreservesproducts}
  Both localizations $\LpshA{Y}$ and $\LshA{Y}$ (of \cref{classicalA1localizationofsheaves}) preserve finite products.
\end{lemma}
\begin{proof}
  The results of \cite[5.5.8.4, 5.5.8.11]{lurie-htt} combined with the explicit description of $\LpshA{Y}$ of \cref{reminderexplicitmodelusualA1localization} implies that $\LpshA{Y}$ preserves finite products.
  Since $\Sh(\dAff_Y,\inftygpd)$ is an hypercomplete \itopos, the sheafification functor also preserves products \cite[6.2.1.1]{lurie-htt}. Finally, we apply \cref{preservesproductsifbothdo}.
\end{proof}

\subsection{\texorpdfstring{$\affineline{}$}{A¹}-isotopic localization}
\label{subsectionisotopies}
For our main results we will only need to impose $\affineline{}$\nobreakdashes-invariance on the mapping spaces of categorical sheaves (such as the Darboux stack). For this purpose we will need some categorical preliminaries:

\subsubsection{Categorical suspension}

\begin{construction}[{Categorical suspension -- \cf \cite[\href{https://kerodon.net/tag/01LJ}{01LJ}]{kerodon}}]\label{constructioncategoricalsuspension}
  Let $E \in \inftygpd$ viewed as an object in $\inftycats$ via the inclusion $\inftygpd\subseteq \inftycats$. We define the \emph{categorical suspension} $\catsuspension E$ of $E$
  as the pushout in $\inftycats$
  \begin{equation}\label{pushoutcategoricalsuspension}
    \begin{tikzcd}
      E\times \partial \Delta^1 \ar[hook]{d} \ar{r} \tikzcocart & \partial \Delta^1\ar{d}\\
      E\times \Delta^1 \ar{r} & \catsuspension E
    \end{tikzcd}
  \end{equation}
  where $E \times \partial \Delta^1 \to \partial \Delta^1$ is the canonical projection. The \icategory $\catsuspension E$ contains two objects $0$ and $1$ and $\Map_{\catsuspension E}(0,1)\simeq E$.
  This construction defines an \ifunctor $\catsuspension (-) \colon \inftygpd \to \inftycats$ and since $\catsuspension \ast=\Delta^1$, it factors by $\catsuspension (-) \colon \inftygpd \to (\inftycats)_{\partial\Delta^1/-/\Delta^1}$.
\end{construction}

\begin{remark}
  The functor $\catsuspension (-)\colon \inftygpd\to (\inftycats)_{\partial\Delta^1/-/\Delta^1}$ preserves finite products:
  \[
    \catsuspension (E_1\times E_2)\simeq \catsuspension E_1\pullback{\Delta^1}\catsuspension E_2.
  \]
\end{remark}

\begin{observation}\label{sigmalocalisff}
  As explained in \cite[\href{https://kerodon.net/tag/01LL}{01LL}]{kerodon}, for any \icategory $\calD$ with objects $x$ and $y$, the pushout \cref{pushoutcategoricalsuspension} induces a fiber product in $\inftycats$:
  \[
    \begin{tikzcd}
      \Fun(E,\Map_{\calD}(x, y)) \tikzcart \ar{r}\ar{d}& \tikzcart  \Fun(\catsuspension E, \calD)\ar{d}\ar{r}& \Fun(E\times \Delta^1,\calD)\ar{d} \\
      \Delta^0\ar{r}{(x,y)}& \calD\times \calD \ar{r} & \Fun(E,\calD\times \calD)
    \end{tikzcd}
  \]
  where the central vertical map is given by composition with $\partial\Delta^1\to \catsuspension E$. In particular, the constant diagram functor $\calD\to \Fun(E,\calD)$ is fully faithful if and only if the functor $\Fun(\Delta^1, \calD)\to \Fun(\catsuspension E,\calD)$ induced by $\catsuspension E\to \Delta^1$ is an equivalence.
\end{observation}

\subsubsection{Locally $\affineline{}$-local presheaves of $\infty$-categories}

\begin{notation}\label{categoricalsuspensionpresheaf}
  Let $\calC \in \inftycats$. We (abusively) denote by $\catsuspension$ the \ifunctor
  \[
   \begin{tikzcd}
    \catsuspension \colon \PSh(\calC, \inftygpd) \ar{r}{\catsuspension \circ -} & \PSh(\calC, \inftycats).
   \end{tikzcd}
  \]
  In particular, for any $C \in \calC$, we have $\catsuspension C \in \PSh(\calC, \inftycats)$, where $C$ is identified with its image by the Yoneda functor.
\end{notation}

\begin{observation}\label{internalhomcategoricalpresheaves}
  Let $\calC \in \inftycats$ with finite products and $F \in \PSh(\calC, \inftycats)$.
  Then by the enriched Yoneda lemma \cite[\S 6.2.7]{hinich2020yoneda} and the fact that the constant diagram functor is left adjoint to the global sections functor, we find for any $I \in \inftycats$ and $C \in \calC$, formulas for the internal-hom in $\PSh(\calC, \inftycats)$: for any $S \in \calC$
 \[
 \stackMap(I, F)(S) = \Fun(I, F(S)) \hspace{1em}\text{and}\hspace{.5em} \stackMap(C, F)(S) =  F(C \times S)
 \]
 (where the constant diagram and Yoneda functors are implicit).
\end{observation}

\begin{definition}\label{definitionlocallyA1local}
  Let $Y$ be a derived stack.
  Let $F \in \PSh(\dAff_Y, \inftycats)$.
  We say that $F$ is \emph{locally $\affineline{Y}$\nobreakdashes-local} if the natural transformation $F \to \stackMap(\affineline{Y}, F)$ (induced by composition with $\affineline{Y} \to \ast$) is fully faithful when evaluated on each object of $\dAff_Y$.
  We denote by
  \[
  \ipshcatSA{Y}\colon  \PSh_{\catsuspension \affineline{Y}}(\dAff_Y, \inftycats) \subseteq \PSh(\dAff_Y, \inftycats)
  \]
  the inclusion of the full subcategory spanned by locally $\affineline{Y}$\nobreakdashes-local presheaves.
\end{definition}

\begin{proposition}\label{localizationlocallyA1local}
  A categorical presheaf $F \in \PSh(\dAff_Y, \inftycats)$ is locally $\affineline{Y}$\nobreakdashes-local if and only if $F$ is local with respect to the canonical map $\catsuspension \affineline{Y} \to \catsuspension \ast=\Delta^1$ in the sense of \cite[5.5.4.1]{lurie-htt}, \ie the induced map
  \[
    F^{\Delta^1}=\stackMap(\catsuspension \ast,F)\to \stackMap(\catsuspension \affineline{Y}, F)
  \]
  is an equivalence.
  The inclusion $\ipshcatSA{Y}$ thus admits a left adjoint $\LpshcatSA{Y}$ presenting locally $\affineline{Y}$\nobreakdashes-local presheaves as a presentable localization obtained by inverting the class of maps $\{ S\times \catsuspension \affineline{Y}\to S\times \Delta^1\}_{S\in \dAff_Y}$ \footnote{remark that the class of maps in indexed by S in dAffY but is not itself in dAffY}.
\end{proposition}
\begin{proof}
  According to \cref{sigmalocalisff,internalhomcategoricalpresheaves}, it follows from \cref{definitionlocallyA1local} that $F$ is locally $\affineline{Y}$\nobreakdashes-local if and only if the diagram
  \[
    \begin{tikzcd}
      \stackMap(\Delta^1, F) \ar{d} \ar{r} & \stackMap(\Delta^1 \times \affineline{Y}, F) \ar{d} \\
      F \times F \ar{r} & \stackMap(\affineline{Y},F) \times \stackMap(\affineline{Y},F)
    \end{tikzcd}
  \]
  is Cartesian. Since $\catsuspension \affineline{\basefield} \simeq \Delta^1 \times \affineline{\basefield} \amalg_{\partial \Delta^1 \times \affineline{\basefield}} \partial \Delta^1$, the second half of the statement follows by applying \cite[5.5.4.15, 5.5.4.20]{lurie-htt} to the class $\{S \times \catsuspension \affineline{\basefield} \to S \times \catsuspension \ast = S \times \Delta^1\}_{S \in \dAff_Y}$.
\end{proof}

\begin{construction}\label{mappingpresheaf}
  Let $\calC$ be a category with finite limits and $F\in \PSh(\calC, \inftycats)$. Let $S \in \calC$ and let $x_S, y_S \in F(S)$ be two sections. We consider the presheaf of \igroupoids $\Mappresheaf_{F_{|_S}}(x_S,y_S) \in \PSh(\calC/S, \inftygpd)$ given by the fiber product
  \[
    \begin{tikzcd}
      \Mappresheaf_{F_{|_S}}(x_S,y_S) \ar{r} \ar{d} \tikzcart & \stackMap(\Delta^1,F)_{|_S}\mathrlap{{}\simeq \stackMap(\Delta^1, F_{|_S})} \ar{d}
      \\
      \ast \ar{r}{(x,y)} & F_{|_S} \times F_{|_S}\rlap{.}
    \end{tikzcd}
  \]
  By construction, for any $S' \in \calC/S$ we have $\Mappresheaf_{F_{|_S}}(x_S,y_S)(S') \simeq \Map_{F(S')}(x_{|{S'}}, y_{|{S'}})$.
\end{construction}

\begin{observation}\label{Sigmalocalandlocal}
  Let $F \in \PSh(\dAff_Y, \inftycats)$.
  It follows from \cref{definitionlocallyA1local} and the restriction formula \cref{restrictionformulanaive} that $F$ is \emph{locally $\affineline{Y}$\nobreakdashes-local}, if and only if for every $S\in \dAff_Y$ and for every $x_S, y_S\in F(S)$, the induced map of presheaves on $\dAff_S$ of \cref{mappingpresheaf}
  \[
    \Mappresheaf_{F_{|_S}}(x_S,y_S) \to \stackMap_S(\affineline{S}, \Mappresheaf_{F_{|_S}}(x_S,y_S))
  \]
  is an equivalence, \ie the presheaf $\Mappresheaf_{F_{|_S}}(x_S,y_S)$ is $\affineline{S}$\nobreakdashes-local in the sense of \cref{localA1MV}.
\end{observation}

\subsubsection{Locally $\affineline{}$-local presheaves of $\infty$-groupoids}

\begin{construction}\label{localizationpresheavesofgroupoids}
  The inclusion $\inftygpd\subseteq \inftycats$ admits a left adjoint sending $\calC$ to the localization $\calC[W_{all}^{-1}]$ inverting all $1$\nobreakdashes-morphisms in $\calC$ \cite[4.1.7.1]{lurie-ha}. As $\ipshcatSA{Y}:\PSh_{\catsuspension \affineline{Y}}(\dAff_Y, \inftycats)\subseteq \PSh(\dAff_Y, \inftycats)$ is a presentable reflective localization (cf.\cref{localizationlocallyA1local}), it follows from \cite[5.5.4.17]{lurie-htt} that the full subcategory of $\PSh(\dAff_Y, \inftygpd)$ given by intersection
  \[
    \PSh_{\catsuspension \affineline{Y}}(\dAff_Y, \inftygpd) \coloneqq \PSh_{\catsuspension \affineline{Y}}(\dAff_Y, \inftycats)\cap \PSh(\dAff_Y, \inftygpd)
  \]
  is again a presentable reflexive localization.
  We denote by $\LpshSA{Y}$ the left adjoint to the inclusion $\ipshSA{Y}\colon\PSh_{\catsuspension \affineline{Y}}(\dAff_Y, \inftygpd)\subseteq \PSh(\dAff_Y, \inftygpd)$.
\end{construction}

\begin{observation}\label{A1localimpliesSigmaA1local}
  Let $F\in \PSh(\dAff_Y, \inftygpd)$. Then, if $F$ is $\affineline{Y}$\nobreakdashes-local in the sense of \cref{localA1MV} then it is locally $\affineline{Y}$\nobreakdashes-local \cref{definitionlocallyA1local}.
\end{observation}

In the spirit of \cref{reminderexplicitmodelusualA1localization}, we will establish an explicit formula for the functor $\LpshSA{Y}$. We start by collecting an elementary fact:

\begin{lemma}\label{localizationpizero}
  For any $F \in \PSh(\dAff_Y, \inftygpd)$, the canonical morphism $\pi_0 \circ F \to \pi_0 \circ \ipshSA{Y}\circ \LpshSA{Y}(F)$ is an isomorphism (of presheaves of sets).
\end{lemma}
\begin{proof}
  Consider the inclusion $\Sets \subseteq \inftygpd$ and the intersection
  \[
    \begin{tikzcd}[column sep=small]
      \PSh_{\catsuspension \affineline{Y}}(\dAff_Y, \Sets) \ar[hook]{r}{u} \ar[hook]{d}{v}& \PSh_{\catsuspension \affineline{Y}}(\dAff_Y, \inftygpd) \ar[hook]{d}{\ipshSA{Y}}
      \\
      \PSh(\dAff_Y, \Sets) \ar[hook]{r}{z} & \PSh(\dAff_Y, \inftygpd)
    \end{tikzcd}
  \]
  It follows from \cref{definitionlocallyA1local} that any presheaf of sets is automatically \emph{locally $\affineline{Y}$\nobreakdashes-local}, so that $v$ is an equality.
  The above diagram consists of right adjoints and therefore their left adjoints commute as well:
  \[
    \begin{tikzcd}[column sep=small]
      \PSh_{\catsuspension \affineline{Y}}(\dAff_Y, \Sets) & \ar{l}[swap]{\widetilde{\pi_0}} \PSh_{\catsuspension \affineline{Y}}(\dAff_Y, \inftygpd)
      \\
      \PSh(\dAff_Y, \Sets)\ar[equals]{u} & \ar{l}[swap]{\pi_0} \ar{u}[swap]{\LpshSA{Y}} \PSh(\dAff_Y, \inftygpd)\rlap{.}
    \end{tikzcd}
  \]
  We then compute using that $\ipshSA{Y}$ is fully faithful:
  \[
    \pi_0 \circ \ipshSA{Y} \circ \LpshSA{Y} = \widetilde{\pi_0}\circ \LpshSA{Y} \circ \ipshSA{Y} \circ \LpshSA{Y} = \widetilde{\pi_0} \circ \LpshSA{Y}= \pi_0. \qedhere
  \]
\end{proof}

\begin{construction}\label{notationconstantobjects}
  Let $E \colon \Delta^\op \to \inftygpd$ be a simplicial object. We denote by $E^{\constant} \colon \Delta^\op \to \inftygpd$ the new simplicial object obtained from $E$ by replacing $E_n$ by the essential image of $E_0 \to E_n$ (induced by the unique map $[n] \to [0]$ in $\Delta$).
\end{construction}

\begin{observation}\label{takingconstantisconstantonpizero}
  Let $E$ be as in \cref{notationconstantobjects}. Any morphism $[0] \to [n]$ induces a retract of $E_0 \to E_n$.
  In particular, $\pi_0(E_0)$ injects into $\pi_0(E_n)$, thus giving a bijection $\pi_0(E_0) \simeq \pi_0(E_n^\constant)$.
  As a consequence, the functor
  \[
    \begin{tikzcd}[column sep=small]
      \Delta^\op \ar{r}{E^{\constant}} & \inftygpd \ar{r}{\pi_0} & \Sets
    \end{tikzcd}
  \]
  is the constant functor.
\end{observation}

\begin{construction}\label{singularcomplexconstant}
  Consider once more the cosimplicial object $\simplex^\bullet_\basefield \colon \Delta \to \dAff_\basefield$ of algebraic simplices (see \cref{reminderexplicitmodelusualA1localization}).
  Given $F \in \PSh(\dAff_Y, \inftygpd)$, we consider the simplicial object
  \[
    \stackMap^{\constant}_Y(\simplex_\basefield^\bullet \times Y,F) \coloneqq (-)^\constant \circ \stackMap_Y(\simplex_\basefield^\bullet \times Y,F) \colon \Delta^\op \to \PSh(\dAff_Y, \inftygpd).
  \]
  Finally, we consider the endofunctor $\rmT_{\catsuspension\affineline{Y}}$ of $\PSh(\dAff_Y, \inftygpd)$, defined by the colimit formula
  \[
    F \mapsto \rmT_{\catsuspension\affineline{Y}}(F) \coloneqq \colim \stackMap^{\constant}_Y(\simplex_\basefield^\bullet \times Y,F)
  \]
  together with the natural transformation $\zeta \colon \Id\to \rmT_{\catsuspension\affineline{Y}}$ induced by the identification $\simplex^0_\basefield=\Spec(\basefield)$.
\end{construction}

\begin{observation}\label{remarkTlocal}
Notice that $F$ is \emph{locally $\affineline{}$\nobreakdashes-local} (\cref{definitionlocallyA1local}) precisely when the first degeneracy map given by the pullback $F\to \stackMap^{\constant}_Y(\affineline{Y},F)$ is an equivalence.
\end{observation}

Our goal is to show that the functor $\rmT_{\catsuspension\affineline{Y}}$ of \cref{singularcomplexconstant} coincides with the endofunctor $\ipshSA{Y}\circ \LpshSA{Y}$ of \cref{localizationpresheavesofgroupoids}. We start by showing that $\rmT_{\catsuspension\affineline{Y}}$  also defines an accessible reflexive localization of $\PSh(\dAff_Y, \inftygpd)$. For that, we will need two key results:

\begin{lemma}\label{remarkOmegapreservessiftedcolimits}
  Let $E_\bullet \colon \Delta^\op \to \inftygpd$ and suppose that the functor $\pi_0(E_\bullet)$ is constant.
  In particular, since $\pi_0$ commutes with colimits, we have
  \[
    \pi_0(\colim E_\bullet) \simeq \pi_0 E_0 \simeq \cdots \simeq \pi_0 E_n.
  \]
  Let $e_0 \in \pi_0 E_0$ and denote by $e_n$ its image in $\pi_0 E_n$ and $e \in \pi_0(\colim E_\bullet)$.
  Then
  \[
    \Omega_{e}\left(\colim E_\bullet\right) \simeq \colim_{[n] \in \Delta^\op} \Omega_{e_n} E_n.
  \]
\end{lemma}
\begin{proof}
  Indeed, since $\pi_0 E_\bullet$ is constant, the decomposition in connected components $E_0 = \coprod_{v \in \pi_0 E_0} E_0^v$ induces a decomposition of the simplicial object $E_\bullet = \coprod_{v \in \pi_0 E_0} E^v_\bullet$ as a coproduct of simplicial objects.
  Now, the choice of the point $v$ in $E_0$ determines, via the unique map $E_0 \to E_n$ a lift of $E^v_\bullet$ to \emph{pointed} connected spaces.
  Finally, we evoke the fact $\Omega$ commutes with sifted colimits on pointed connected spaces \cite[5.2.6.18]{lurie-ha}.
\end{proof}

\begin{observation}\label{restrictioncommuteswithA1localizationonpresheaves}
  Let $F \in \PSh(\dAff_Y, \inftycats)$ and $S \in \dAff_Y$. The restriction $F \mapsto F_{|_S}$ (see \cref{notationforrestrictionpresheavescategories}) sends locally $\affineline{Y}$\nobreakdashes-local presheaves to locally $\affineline{S}$\nobreakdashes-local presheaves.
  Moreover, as in \cref{naiveA1localizationcommuteswithrestriction}, we have from the definitions:
  \[
    \left(\rmT_{\catsuspension\affineline{Y}}(F)\right)_{|_S} \simeq \rmT_{\catsuspension\affineline{S}}(F_{|_S}).
  \]
\end{observation}

\begin{lemmalist}{Let $F\in \PSh(\dAff_Y, \inftygpd)$.}
    \item \label{mappingspacesofTareA1localizations-pi0} The map $\pi_0(\zeta_F) \colon \pi_0(F)\to \pi_0(\rmT_{\catsuspension\affineline{Y}}(F))$ is an isomorphism of presheaves.
    \item \label{keyformulaforA1localizationsrelation} Let $S \in \dAff_Y$ and $x_S \in F(S)$. Then, we have a canonical equivalence of presheaves on $\dAff_S$
    \[
      \Omega_{x_S} \rmT_{\catsuspension\affineline{S}}(F_{|_S}) \simeq \ipshA{S}\circ \LpshA{S}\left(\Omega_{x_S} (F_{|_S})\right)
    \]
    with $\ipshA{S} \circ \LpshA{S}$ as in \cref{localA1MV} and \cref{naiveA1localizationcommuteswithrestriction}.
\end{lemmalist}
\begin{proof}
  Since $\pi_0(F)$ is determined objectwise, (a) is a direct consequence of \cref{takingconstantisconstantonpizero} and the fact $\pi_0$ commutes with colimits.
  Let us now address (b): since this is a statement about presheaves, it amounts to proving the statement when evaluated on each object $S'\in \dAff_S$.
  Let $x_{S'} \coloneqq (x_S)_{|S'}$.
  By \cref{singularcomplexconstant}, we have
  \[
    \left(\Omega_{x_S} \rmT_{\catsuspension\affineline{S}}(F_{|_S})\right)(S') = \Omega_{x_{S'}} \left(\colim F(\simplex^\bullet_\basefield \times S')^{\constant}\right).
  \]
  By \cref{takingconstantisconstantonpizero}, the simplicial set $\pi_0F(\simplex^\bullet_\basefield \times S')^{\constant}$ is constant with value $\pi_0(F(S'))$.
  \cref{remarkOmegapreservessiftedcolimits} thus implies
  \[
    \Omega_{x_{S'}} \left(\colim F(\simplex^\bullet_\basefield \times S')^{\constant}\right) \simeq \colim \Omega_{x_{S'}} F(\simplex^\bullet_\basefield \times S')^{\constant}.
  \]
  Finally, by definition of $(-)^\constant$, we get
  \begin{multline*}
    \left(\Omega_{x_S} \rmT_{\catsuspension\affineline{S}}(F_{|_S})\right)(S') \simeq \colim \Omega_{x_{S'}} F(\simplex^\bullet_\basefield \times S')^{\constant} \simeq \colim \Omega_{x_{S'}} F(\simplex^\bullet_\basefield \times S')
    \\
    \simeq \ipshA{S}\LpshA{S}\left(\Omega_{x_S} (F_{|_S})\right)(S')
  \end{multline*}
where the last equivalence follows from \cref{reminderexplicitmodelusualA1localization}. Since this chain of equivalences is natural in $S'$ we have obtained the claimed natural equivalence of presheaves on $\dAff_S$.
\end{proof}

\begin{proposition}\label{explicitformulaTislocalization}
The endofunctor $\rmT_{\catsuspension\affineline{Y}}$ together with the natural transformation $\zeta$, defines an accessible reflexive localization of $\PSh(\dAff_Y, \inftygpd)$.
\end{proposition}
\begin{proof}
  We use the criteria of \cite[5.2.7.4-(3)]{lurie-htt}: it suffices to check that for any $F$, the two maps $
  \rmT_{\catsuspension\affineline{Y}}(\zeta_F),\, \zeta_{\rmT_{\catsuspension\affineline{Y}}(F)}$ fitting in the commutative diagram
  \[
    \begin{tikzcd}
      F\ar{r}{\zeta_F}\ar{d}{\zeta_F}& \rmT_{\catsuspension\affineline{Y}}(F)\ar{d}{\zeta_{\rmT_{\catsuspension\affineline{Y}}(F)}}\\
      \rmT_{\catsuspension\affineline{Y}}(F)\ar{r}{  \rmT_{\catsuspension\affineline{Y}}(\zeta_F)}&\rmT_{\catsuspension\affineline{Y}}(\rmT_{\catsuspension\affineline{Y}}(F))
    \end{tikzcd}
  \]
  are equivalences in $\PSh(\dAff_Y, \inftygpd)$.
  Since equivalences of presheaves are determined pointwise, it suffices to check that when evaluated on $S\in \dAff_Y$, both maps:
  \begin{enumerate}[label={(\roman*)}]
    \item \label{substepenumpi0} induce equivalences on $\pi_0$;
    \item \label{substepenumOmega} for every connected component, the induced map of loop spaces is an equivalence.
  \end{enumerate}
  Assertion \ref{substepenumpi0} follows from \cref{mappingspacesofTareA1localizations-pi0}.
  To show \ref{substepenumOmega}, let $x_S\in F(S)$ and consider the induced maps of loop prestacks on $\dAff_S$
  \[
    \Omega_{x_S} \rmT_{\catsuspension\affineline{S}}(\zeta_F) \hspace{1em}\text{and}\hspace{1em} \Omega_{x_S}\zeta_{\rmT_{\catsuspension\affineline{S}}(F_{|_S})}:\Omega_{x_S} \rmT_{\catsuspension\affineline{S}}(F_{|_S})\to \Omega_{x_S} \rmT_{\catsuspension\affineline{S}}\left(\rmT_{\catsuspension\affineline{S}}(F_{|_S})\right).
  \]
  Via \cref{keyformulaforA1localizationsrelation}, these maps can be identified with the two natural maps associated to the localization $\LpshA{S}$ of \cref{localA1MV}:
  \[
    \ipshA{S} \LpshA{S}(\eta_{\affineline{}}) \hspace{1em}\text{and}\hspace{1em} (\eta_{\affineline{}})_{ \ipshA{S} \LpshA{S}}: \ipshA{S} \LpshA{S}\left(\Omega_{x_S} (F_{|_S})\right)\to  \ipshA{S} \LpshA{S}\left( \ipshA{S} \LpshA{S}\left(\Omega_{x_S} (F_{|_S})\right)\right).
  \]
  Since we know that $\LpshA{S}$ is a localization, we have by \cite[5.2.7.4-(3)]{lurie-htt} that these two maps are equivalences. This concludes the proof.
\end{proof}

Finally, we compare the two localizations:

\begin{proposition}\label{comparisonoftwonaivesigmaA1localizations}
  The two localizations $\rmT_{\catsuspension\affineline{Y}}$ and $\LpshSA{Y}$ of $\PSh(\dAff_Y, \inftygpd)$ coincide.
\end{proposition}
\begin{proof}
  By \cite[5.2.7.8, 5.2.2.2]{lurie-htt}, it suffices to check that the two classes of local objects coincide.
  First we show that for any $F$, $\rmT_{\catsuspension\affineline{Y}}(F)$ is locally $\affineline{Y}$\nobreakdashes-local: indeed, following \cref{Sigmalocalandlocal}, it suffices to check that for every object $S \in \dAff_Y$, and $x_S \in F(S)$ the induced map on loops
  \[
    \Omega_{x_S} \rmT_{\catsuspension\affineline{S}}(F_{|_S})\to   \Omega_{x_S}\stackMap_S(\affineline{S}, \rmT_{\catsuspension\affineline{S}}(F_{|_S}))=\stackMap_S(\affineline{S},\Omega_{x_S} \rmT_{\catsuspension\affineline{S}}(F_{|_S}))
  \]
  is an equivalence of presheaves of \igroupoids on $\dAff_S$.
  This then follows from \cref{keyformulaforA1localizationsrelation}.

  Conversely, if $F$ is locally $\affineline{}$\nobreakdashes-local then the map $\zeta_{F} \colon F \to \rmT_{\catsuspension\affineline{Y}}( F)$ is an equivalence by \cref{remarkTlocal}.
\end{proof}

\begin{observation}\label{A1locallylocallocalizationpreservesproducts}
  As in the proof of \cref{A1localizationforsheavespreservesproducts}, the explicit formula $\ipshSA{Y}\circ \LpshSA{Y}\simeq \rmT_{\catsuspension\affineline{Y}}$, combined with \cite[5.5.8.4, 5.5.8.11]{lurie-htt},  implies that $\LpshSA{Y}$ preserves finite products.
\end{observation}

\subsubsection{Locally $\affineline{}$-local sheaves}

We now adapt the discussion of the previous paragraph to sheaves. First we discuss sheaves of categories.

\begin{construction}\label{localizationofsheaves1}
  Consider the pullback
  \begin{equation}\label{localizationofsheaves}
    \begin{tikzcd}[column sep=small]
      \Sh_{\catsuspension{\affineline{Y}}}(\dAff_Y, \inftycats)\tikzcart \ar[hook]{r}\ar[hook]{d}{\ishcatSA{Y}}& \PSh_{\catsuspension{\affineline{Y}}}(\dAff_Y, \inftycats)\ar[hook]{d}{\ipshcatSA{Y}}\\
      \Sh(\dAff_Y, \inftycats)\ar[hook]{r}{\ishcat{Y}}&  \PSh(\dAff_Y, \inftycats)\rlap{.}
    \end{tikzcd}
  \end{equation}
  As in \cref{intersectionoflocalizations}, each inclusion functor in \cref{localizationofsheaves} admits a left adjoint and the square of left adjoints commutes. We denote by $\LshcatSA{Y}$ the left adjoint to the inclusion $\ishcatSA{Y}\colon\Sh_{\catsuspension{\affineline{Y}}}(\dAff_Y, \inftycats)\subseteq \Sh(\dAff_Y, \inftycats)$.
\end{construction}

We now discuss sheaves of \igroupoids:

\begin{construction}\label{A1localizationofsheavesofgroupoids}
  As in \cref{intersectionoflocalizations}, the localization $\LpshSA{Y}$ of \cref{localizationpresheavesofgroupoids} induces a localization at the level of sheaves of \igroupoids
  \[
   \LshSA{Y} \colon  \Sh(\dAff_Y, \inftygpd)\to \Sh_{\catsuspension{\affineline{Y}}}(\dAff_Y, \inftygpd)
  \]
  with $\ishSA{Y}\circ \LshSA{Y}$ given by the iteration \cref{formulaiteratedfinal} of the endofunctor
  \begin{equation}\label{formulaforendofunctorlocalizationSigmaA1beforeiteration}
    \PhiSA{Y} \coloneqq \Lsh{Y}\circ \ipshSA{Y}\circ \LpshSA{Y}\circ \ish{Y}:\Sh(\dAff_Y, \inftygpd)\to \Sh(\dAff_Y, \inftygpd).
  \end{equation}
  Moreover, since the diagram of inclusions
  \[
    \begin{tikzcd}[column sep=-3em,row sep=small]
      \PSh(\dAff_Y, \inftygpd) \ar[hook]{rr} && \PSh(\dAff_Y, \inftycats) \ar[hook, from=dd]  &
      \\
      & \Sh(\dAff_Y, \inftygpd) \ar[hook,crossing over]{rr} \ar[hook]{ul} && \ar[hook]{ul}  \Sh(\dAff_Y, \inftycats)
      \\
      \ar[hook]{uu} \PSh_{\catsuspension{\affineline{Y}}}(\dAff_Y, \inftygpd) \ar[hook]{rr} && \PSh_{\catsuspension{\affineline{Y}}}(\dAff_Y, \inftycats) &
      \\
      & \ar[hook,crossing over]{uu} \Sh_{\catsuspension{\affineline{Y}}}(\dAff_Y, \inftygpd) \ar[hook]{rr} \ar[hook]{ul} && \ar[hook]{uu} \ar[hook]{ul} \Sh_{\catsuspension{\affineline{Y}}}(\dAff_Y, \inftycats)
    \end{tikzcd}
  \]
  commutes, so does the diagram of left adjoints.
\end{construction}

\begin{observation}\label{A1localizationpreservespizeroforsheaves}
  The computation of \cref{localizationpizero} remains valid for sheaves if we replace the functor $\pi_0$ by its sheaf-theoretic counterpart $\homotopysheaf_0 \colon \Sh(\dAff_Y, \inftygpd)\to \Sh(\dAff_Y, \Sets)$ and $\LpshSA{Y}$ by the functor $\LshSA{Y}$ of \cref{A1localizationofsheavesofgroupoids}, \ie if
  \[
    \begin{tikzcd}
      \Sh_{\catsuspension{\affineline{Y}}}(\dAff_Y, \inftygpd) \ar[hook]{r}{\ishSA{Y}} & \Sh(\dAff_Y, \inftygpd)
    \end{tikzcd}
  \]
  denotes the inclusion, then we have canonical equivalences
  \begin{equation}\label{formulaA1localizationpreservespizerosheaf}
    \homotopysheaf_0 \circ \Phi_{\catsuspension{\affineline{Y}}} \simeq \homotopysheaf_0 \hspace{1cm}\text{ and by induction } \hspace{1cm}
    \homotopysheaf_0 \circ \ishSA{Y} \circ \LshSA{Y} \simeq \homotopysheaf_0.
  \end{equation}
\end{observation}

\begin{lemma}\label{suspensedA1localizationpreservesfiniteproducts}
  The functors $\PhiSA{Y}$ and  $\LshSA{Y}$ of \cref{A1localizationofsheavesofgroupoids} preserve finite products.
\end{lemma}
\begin{proof}
  We argue as in the proof of \cref{A1localizationforsheavespreservesproducts}: by \cref{A1locallylocallocalizationpreservesproducts},  $\ipshSA{Y}\circ \LpshSA{Y}=\rmT_{\catsuspension\affineline{Y}}$ preserves finite products. By \cite[6.2.1.1]{lurie-htt}, the same holds for the hyper-sheafification functor.
  We conclude using the filtered colimit formula of \cref{endofunctortransfiniteinduction} and \cref{formulaforiteratedlocalization}.
\end{proof}

Finally, we establish the sheaf-version of \cref{keyformulaforA1localizationsrelation}:

\begin{lemma}\label{lemmamappingspacesA1isotopiclocalizationsareA1localizations}
  Let $F\in \Sh(\dAff_Y, \inftygpd)$, $S\in \dAff_Y$ and $x_S\in F(S)$.
  Then, we have canonical equivalences of hypersheaves on $\dAff_S$
  \[
    \Omega_{x_S} \PhiSA{S}(F_{|_S})\simeq \PhiA{S}\left(\Omega_{x_S} (F_{|_S})\right)
  \]
  and
  \[
    \Omega_{x_S} \ishSA{S}\LshSA{S}(F_{|_S})\simeq \ishA{S}\LshA{S}\left(\Omega_{x_S} (F_{|_S})\right).
  \]
\end{lemma}

\begin{proof}
  We have, for any pointed sheaf $x \colon \ast \to F$ over $\dAff_S$:
  \begin{equation}\label{transfinitestepsuccessor}
    \Omega_x \PhiSA{S}(F_{|_S}) \coloneqq \Omega_x \Lsh{S} \ipshSA{S}\LpshSA{S}\ish{S}(F_{|_S})\simeq \Lsh{S}\Omega_x \ipshSA{S}\LpshSA{S}\ish{S}(F_{|_S})
  \end{equation}
  where second  equivalence follows the sheafification functor is a left exact localization. Using \cref{comparisonoftwonaivesigmaA1localizations}, the RHS becomes equivalent to $\Lsh{S}\Omega_x \rmT_{\catsuspension\affineline{S}}\ish{S}(F_{|_S})$ which by \cref{keyformulaforA1localizationsrelation} is equivalent to $\Lsh{S}\ishA{S}\LshA{S}\ish{S}(\Omega_x F_{|_S}) \eqqcolon \PhiA{S}(\Omega_x F_{|_S})$. This establishes the first formula.

We now prove the second formula: fix an ordinal $\gamma$ such that both localizations are $\gamma$\nobreakdashes-accessible, and a pointed sheaf $x \colon \ast \to F_{|_S}$ over $\dAff_S$.
  If $\gamma = \gamma' + 1$ and if $\Omega_x (\PhiSA{S})^{\gamma'}(F_{|_S}) \simeq (\PhiA{S})^{\gamma'} \Omega_x(F_{|_S})$, then applying \cref{transfinitestepsuccessor} to the pointed sheaf $(\PhiSA{S})^{\gamma'}(F_{|_S})$ yields $\Omega_x (\PhiSA{S})^{\gamma}(F_{|_S}) \simeq (\PhiA{S})^{\gamma} \Omega_x(F_{|_S})$.
  If $\gamma$ is a limit ordinal, we have, as the poset $\{\gamma'<\gamma\}$ is filtered:
  \[
    \Omega_x (\PhiSA{S})^{\gamma}(F_{|_S}) \coloneqq \Omega_x \colim_{\gamma'<\gamma} (\PhiSA{S})^{\gamma'}(F_{|_S}) \simeq \colim_{\gamma'<\gamma} \Omega_x (\PhiSA{S})^{\gamma'}(F_{|_S}).
  \]
  In particular, assuming $\Omega_x (\PhiSA{S})^{\gamma'}(F_{|_S}) \simeq (\PhiA{S})^{\gamma'} \Omega_x(F_{|_S})$ for any $\gamma'<\gamma$, we find $\Omega_x (\PhiSA{S})^{\gamma}(F_{|_S}) \simeq (\PhiA{S})^{\gamma} \Omega_x(F_{|_S})$.

We conclude by transfinite induction, using \cref{formulaforiteratedlocalization}.
\end{proof}

\subsection{Darboux charts and quadratic bundles up to \texorpdfstring{$\affineline{}$}{A¹}-Isotopies}
\label{darbouxstackuptoisitopysection}

In this section we use \cref{subsectionisotopies} to construct a new version of the Darboux stack where we contract all $\affineline{}$\nobreakdashes-isotopies in the mapping sheaves.
Let us start with a construction implementing restriction to the étale site:

\begin{construction}\label{compatibilityofpizerosheavesunderrestriction}
  Let $X$ be a derived Deligne-Mumford stack.
  Consider the functors
  \begin{align*}
  (-)^X_\derham\colon X_\et^\daff &\to \dSt_{X_\derham} \supset \dAff_{X_\derham} \\
   S &\mapsto S_\derham.
  \end{align*}
(\cf \cref{notationforrestrictiontoXet})  together with the associated restriction functors, yielding a zig-zag
  \[
   \PSh(X_\et^\daff, \inftygpd) \from \PSh(\dSt_{X_\derham}, \inftygpd) \to \PSh(\dAff_{X_\derham}, \inftygpd)
  \]
  preserving hypersheaves (with the canonical topology on $\dSt_{X_\derham}$) and $\pi_0$\nobreakdashes-presheaves.
  By \cite[Prop.\,7.1]{zbMATH07387425}, those restriction functors commute with colimits, hypersheafification and truncations.
  Since $\Sh(\dSt_{X_\derham}, \inftygpd) \simeq \Sh(\dAff_{X_\derham}, \inftygpd)$ \footnote{cf. \cref{descriptionsofslicecategoriesofstacks}}, they induce a restriction at the level of sheaves
  \[
    \hspace{2em}
    \begin{tikzcd}[column sep=tiny,row sep=-.2em]
      \mathllap{ (-)_{|X_\et^\daff} \colon{}} \Sh(\dAff_{X_\derham}, \inftygpd) \ar[phantom]{r}{\simeq} & \Sh(\dSt_{X_\derham}, \inftygpd) \ar{r} & \Sh(X_\et^\daff, \inftygpd)
    \\
      F \ar[mapsto]{rr} && \left(S \mapsto \Map_{X_\derham}(S_\derham, F) \right)
    \end{tikzcd}
  \]
  commuting with colimits and with taking $\pi_0$\nobreakdashes-sheaves.
\end{construction}

\begin{constructions}{\label{constructionofquotientA1stacksgeneral}Let $X$ be a derived Deligne--Mumford stack and let $F\in \Sh(\dAff_{X_\derham}, \inftygpd)$.
  Let $\homotopysheaf\subseteq \homotopysheaf_0(F){|{X_\et^\daff}}$ be a subsheaf of the restriction where $\homotopysheaf_0$ denotes the $\pi_0$\nobreakdashes-sheaves (unambiguous by \cref{compatibilityofpizerosheavesunderrestriction}) and consider $F_{\homotopysheaf}$ defined as the fiber product in $\Sh(X_\et^\daff,\inftygpd)$:
  \[
  \begin{tikzcd}[ampersand replacement=\&]
    F_{\homotopysheaf} \ar{d} \ar{r} \tikzcart \& F_{|{X_\et^\daff}} \ar{d}
    \\
    \homotopysheaf \ar{r} \& \homotopysheaf_0(F)_{|{X_\et^\daff}}\rlap{.}
  \end{tikzcd}
  \]
  We construct:}
  \item \label{A1quotientgeneral} the \emph{$\affineline{}$\nobreakdashes-isotopic quotient} of $F_{\homotopysheaf}$, denoted $F_{\homotopysheaf}^{\affineline{}}$, as the sheaf of \igroupoids on $X_\et^\daff$ given by the fiber product
  \[
    \begin{tikzcd}[column sep=scriptsize]
      F_{\homotopysheaf}^{\affineline{}} \ar[hook]{rr} \ar{d} \tikzcart
      &&
      \left( \PhiSA{X_\derham} \left(F\right)\right)_{|{X_\et^\daff}} \ar{d}
      \\
      \mathllap{\homotopysheaf = {}} \homotopysheaf_0\left(F_{\homotopysheaf}\right) \ar[hook]{r}
      &
      \homotopysheaf_0\left(F\right)_{|{X_\et^\daff}} \ar{r}{\cref{formulaA1localizationpreservespizerosheaf}}[swap]{\sim}
      &
      \homotopysheaf_0\left(\PhiSA{X_\derham} \left(F\right)\right)_{|{X_\et^\daff}} \rlap{.}
    \end{tikzcd}
  \]
  By construction, it comes with a canonical morphism $F_{\homotopysheaf} \to F_{\homotopysheaf}^{\affineline{}}$ which induces an isomorphism on $\homotopysheaf_0$.

  \item \label{A1localizationsigmageneral} the \emph{$\catsuspension\affineline{}$\nobreakdashes-localization} of $F_{\homotopysheaf}$ , denoted $F_{\homotopysheaf, \catsuspension\affineline{}}$, as the sheaf of \igroupoids on $X_\et^\daff$ given by the fiber product
  \[
    \begin{tikzcd}[column sep=scriptsize]
      F_{\homotopysheaf, \catsuspension\affineline{}} \ar[hook]{rr} \ar{d} \tikzcart
      &&
      \left( \ishSA{X_\derham}\LshSA{X_\derham}  \left(F\right)\right)_{|{X_\et^\daff}} \ar{d}
      \\
      \mathllap{\homotopysheaf = {}} \homotopysheaf_0\left(F_{\homotopysheaf}\right) \ar[hook]{r}
      &
      \homotopysheaf_0\left(F\right)_{|{X_\et^\daff}} \ar{r}{\cref{formulaA1localizationpreservespizerosheaf}}[swap]{\sim}
      &
      \homotopysheaf_0\left(\ishSA{X_\derham}\LshSA{X_\derham}  \left(F\right)\right)_{|{X_\et^\daff}} \rlap{.}
    \end{tikzcd}
  \]
  It comes with canonical morphisms
  \begin{equation}\label{definitionA1darbouxstackCartesianmap}
    \begin{tikzcd}F_{\homotopysheaf}\ar{r}& F_{\homotopysheaf}^{\affineline{}}\ar{r}& F_{\homotopysheaf, \catsuspension\affineline{}}\end{tikzcd}
  \end{equation}
  in $\Sh(X_\et^\daff, \inftygpd)$, which induce isomorphisms on $\homotopysheaf_0$.
\end{constructions}

\begin{observation}\label{observationmapsareA1localizations}
In the context of \cref{constructionofquotientA1stacksgeneral}, let $\sigma$ be a section of $F_{\homotopysheaf}$. Then \cref{lemmamappingspacesA1isotopiclocalizationsareA1localizations} provides equivalences
  \begin{multline}
    \label{automorphismsinthedarbouxstacksareA1localizations}
    \Autpresheaf_{F_{\homotopysheaf}^{\affineline{}}}(\sigma) \coloneqq \Omega_\sigma F_{\homotopysheaf}^{\affineline{}} \simeq \Omega_\sigma \left( \PhiSA{X_\derham} (F)\right)_{|{X_\et^\daff}}
    \\
    \simeq \left(\PhiA{X_\derham} \Omega_\sigma \, F\right)_{|{X_\et^\daff}}\simeq \left(\PhiA{X_\derham} \Autpresheaf_{ F}(\sigma)\right)_{|{X_\et^\daff}}
  \end{multline}
  \begin{equation}
    \label{automorphismsinthedarbouxstacksareA1localizationsla}
    \llap{and }\Autpresheaf_{F_{\homotopysheaf, \catsuspension\affineline{}}}(\sigma) \coloneqq \Omega_\sigma F_{\homotopysheaf, \catsuspension\affineline{}}
    \simeq \cdots \simeq \left(\ishA{X_\derham}\LshA{X_\derham} \Autpresheaf_{ F}(\sigma)\right)_{|{X_\et^\daff}}.
  \end{equation}
\end{observation}

We can finally introduce the isotopic quotient of the Darboux stack:

\begin{definitions}{Let $X$ be a Deligne-Mumford equipped with a $(-1)$\nobreakdashes-shifted exact symplectic form $\lambda$.
  Consider the stack $F \coloneqq \stackfactor^{\Liouv,\lambda,\simeq}_{p_X} \colon \dAff_{X_\derham}^\op \subset \dSt_{X_\derham}^\op \to \inftygpd$ obtained by extracting the maximal \igroupoid of the stack defined by the fiber product \cref{definitionLiouvillestackalphafiberproductconstruction} and define $\homotopysheaf \coloneqq \homotopysheaf_0(\Darbstack_X)$, together with its inclusion $\homotopysheaf \subseteq \homotopysheaf_0\left(F\right)_{|{X_\et^\daff}}$.
  We define: }
 \item \label{Liouvilleuptoisotopy} $\DA \coloneqq F_{\homotopysheaf}^{\affineline{}}$ as in \cref{A1quotientgeneral};
\item \label{LiouvilleuptoisotopySigmalocalization} $\Darbstack_{X, \catsuspension\affineline{}}^{\lambda} \coloneqq F_{\homotopysheaf, \catsuspension\affineline{}}$ as in \cref{A1localizationsigmageneral}.
\end{definitions}

Now we address versions of the $\affineline{}$\nobreakdashes-isotopic quotient and $\affineline{}$\nobreakdashes-localization of the stack of quadratic bundles with flat connection of \cref{stackquaddef}:

\begin{definitions}{Let $X$ be a derived Deligne-Mumford stack.
      Consider the stack $F=\B\orthogonal_{X_\derham} \colon \dSt_{X_\derham}^\op \to \inftygpd$ mapping $S \to X_\derham$ to $\Map(S, \B\orthogonal)$.
      By construction, the stack $\stackQuadnabla_X$ (see \cref{stackquaddef}) is a substack of the restriction of $\B\orthogonal_{X_\derham}$ along $(-)_\derham \colon X_\et^\daff \to \dSt_{X_\derham}$. Setting $\homotopysheaf \coloneqq \homotopysheaf_0(\stackQuadnabla_X)$ together with the inclusion $\homotopysheaf \subset \homotopysheaf_0\left(\B\orthogonal_{X_\derham}\right)_{|{X_\et^\daff}}$, we define:}

  \item \label{A1contractionofQuad} $\QA \coloneqq F_{\homotopysheaf}^{\affineline{}}$ as in \cref{A1quotientgeneral};

  \item \label{A1contractionofQuadSigmalocalization} $  \stackQuadnabla_{X, \isotopies} \coloneqq F_{\homotopysheaf, \catsuspension\affineline{}}$  as in \cref{A1localizationsigmageneral}.
\end{definitions}

\begin{remark}
  Contrary to $ \stackQuadnabla_{X, \isotopies}$ and $\Darbstack^{\lambda}_{X,\isotopies}$, the stacks $\QA$ and $\DA$ are not themselves locally $\affineline{}$\nobreakdashes-local in the sense of \cref{definitionlocallyA1local}. Nevertheless, they are sufficient to identify isotopic families of automorphisms as in \cref{isotopiessolveambiguity}.
  We will use them to glue local invariants of singularities in the next section and in future work.
\end{remark}

Finally, we observe that

\begin{proposition}\label{actiondescendstoisotopyquotient}
  The action of \cref{actionondarbouxcharts} descends to the $\affineline{}$\nobreakdashes-isotopic quotients and localizations
  \begin{gather*}
    \QA \times \DA \to \DA
  \\
    \stackQuadnabla_{X,\isotopies} \times \Darbstack^{\lambda}_{X,\isotopies} \to \Darbstack^{\lambda}_{X,\isotopies}.
  \end{gather*}
\end{proposition}
\begin{proof}
  Using \cref{suspensedA1localizationpreservesfiniteproducts}, the monoid structure of $\B\orthogonal$ of \cref{definitionBO} induces a commutative monoid structure on  both $\QA$ and $\stackQuadnabla_{X,\isotopies}$ and the canonical maps \cref{definitionA1darbouxstackCartesianmap}
  \[
    \stackQuadnabla_X\to \QA\to \stackQuadnabla_{X, \isotopies}
  \]
  are  maps of commutative monoids.
  Moreover, the action of $\B\orthogonal_{X_\derham}$ on $\stackfactor^{\Liouv,\lambda,\simeq}_{p_X}$ of \cref{actionondarbouxcharts} induces an action between their $\affineline{}$\nobreakdashes-quotients and $\Sigma \affineline{}$\nobreakdashes-localizations.
  Finally, the isomorphism \cref{formulaA1localizationpreservespizerosheaf} guarantees that this action preserves the Darboux substacks $\DA$ and $\Darbstack^{\lambda}_{X,\isotopies}$ and is compatible with canonical maps \cref{definitionA1darbouxstackCartesianmap}.
\end{proof}

\subsection{The contractibility theorem}
\label{sectioncontractibilitytheoremdetailedformulation}

Let us now give a more precise formulation of \cref{theoremcontractibility}:

\begin{theorem}\label{contractibilitytheoremdetailedversion}
  Let $X$ be a $(-1)$-shifted exact symplectic derived Deligne--Mumford stack (with form $\lambda$). Then the quotient stack is contractible
  \[
   \quotA \simeq \ast_X.
  \]
\end{theorem}

The proof of \cref{contractibilitytheoremdetailedversion}  will last until the end of \cref{sectionisotopies}.
We start by reducing the it to the technically more accessible \cref{prop:A1-automorphisms}.
This will require the technical notion of minimal models (also called minimal critical charts) at the neighborhood of a point.
The notion was introduced by Brav, Bussi and Joyce in \cite[Def. 2.13]{MR3904157}.

\begin{definition}[Brav-Bussi-Joyce]\label{definitionminimalmodel}
  Let $X$ be an exact $(-1)$\nobreakdashes-symplectic derived scheme (with form $\lambda$) and $x \in X$ be geometric point. A \emph{minimal model} of $X$ at $x$ is an affine Zariski open neighborhood $S=\Spec(A) \subset X$ of $x$ equipped with a critical chart $(\formalU ,f) \in \Darbstack^\lambda_X(S)$ such that the Hessian of $f$ at $x$ vanishes.
\end{definition}

\begin{example}
The critical chart of \cref{ambiguitymorphisms} is minimal at $x=0$.
\end{example}

\begin{theorem}[{\cite[Thm.\,4.1]{MR3904157}}]
  Every geometric point $x \in X$ admits a minimal model.
\end{theorem}

\begin{proposition}\label{prop:A1-automorphisms}
  Let $x \in X$ be a geometric point and $(S=\Spec(A), \formalU_0)$ a minimal model at $x$.
  Let $(Q,q)$ be a (non-degenerate) quadratic vector space over $\basefield$ and let $\formalQ$ be the induced trivial quadratic bundle over $S_{\derham}$. Set $\formalU \coloneqq \formalU_0 \smboxplus \formalQ$ and
  denote by $\stackAut(\formalQ)$ the stack (over $S$) of automorphisms of $\formalQ$ in $\stackQuadnabla_X$, and by $\stackAut(\formalU)$ the stack of automorphisms of $\formalU$ in $\Darbstack^{\lambda}_X$.
  Then, the induced morphism
  \[
    \stackAut(\formalQ) \to \stackAut_\Darb(\formalU)
  \]
  is a $\PhiA{}$\nobreakdashes-equivalence at $x \in S$, \ie
  $
    \PhiA{S}(\stackAut(\formalQ))_{x} \simeq \PhiA{S}(\stackAut_\Darb(\formalU))_{x}.
  $
\end{proposition}

We will delay the proof of \cref{prop:A1-automorphisms} until the next section.
We focus for now on proving it implies \cref{contractibilitytheoremdetailedversion}. This will in turn require a bit of preliminaries on group completions.

\subsubsection{Isotopies and group completion}
We will need some understanding of the group completion of the monoid objects $\QA$  of \cref{A1contractionofQuad}.

\begin{notation}
 Let $M$ be a monoid object in stacks. We denote by $M^+ \simeq \Omega \B M$ its group-completion.
 If $E$ is acted upon by $M$, we denote by $E^+$ the base change
 \[
  E^+ \coloneqq E \otimes_M M^+.
 \]
\end{notation}

\begin{lemma}\label{inversehyperbolic}
  Denote by $Q_1$ the trivial quadratic bundle with fiber $\basefield$ and form $x^2$.
  There are natural equivalences (of étale hypercomplete stacks over $X_\et^\daff$)
  \begin{align*}
     \left(\QA\right)^+ &\simeq \colim \left[
    \begin{tikzcd}[cramped,ampersand replacement=\&]
      \QA\ar{r}{- \oplus Q_1} \& \QA\ar{r}{- \oplus Q_1} \& \cdots
    \end{tikzcd} \right]\rlap{\hspace{1em}and}\\
     \left(\DA\right)^+ &\simeq \colim \left[
     \begin{tikzcd}[cramped,ampersand replacement=\&]
       \DA\ar{r}{- \oplus Q_1} \& \DA \ar{r}{- \oplus Q_1} \& \cdots
     \end{tikzcd}
     \right].
  \end{align*}
\end{lemma}
\begin{proof}
\newcommand{\localquad}{\QA}
\newcommand{\localquadplus}{\left(\QA\right)^+}
  The second equivalence follows from the first.  Denote by $C$ the colimit in stacks
  \[
    C \coloneqq \colim \left[
    \begin{tikzcd}[cramped]
      \localquad \ar{r}{- \oplus Q_1} & \localquad \ar{r}{- \oplus Q_1} & \cdots
    \end{tikzcd}
    \right].
  \]
  By cofinality, we can replace $Q_1$ in the colimit diagram by $Q_2 \coloneqq Q_1 \oplus Q_1$.
  Choosing a square root of $-1 \in \basefield$ yields an isomorphism $Q_2 \simeq E \oplus E^\dual$, where $E$ is trivial of rank $1$, and $E \oplus E^\dual$ is equipped with its canonical form.
  Consider the $1$\nobreakdashes-parameter family of endomorphisms of $E^3$ given by
  \[
    f_t \coloneqq
    \begin{pmatrix}
      1-t^3 & t & 0 \\
      0 & 1 - t^3 & t \\
      t(3-3t^3 +t^6) & 0 & 1-t^3.
    \end{pmatrix}
  \]
  It satisfies $\det(f_t) = 1$ and $f_0 = \Id$, while $f_1$ is the matrix of a $3$\nobreakdashes-cycle permuting the coordinates.
  The associated $\affineline{}$\nobreakdashes-family of automorphisms of $Q_2^{\oplus 3}$ interpolates between the identity and the $3$\nobreakdashes-cycle permutation.
  By \cref{automorphismsinthedarbouxstacksareA1localizationsla}, the $3$\nobreakdashes-cycle is homotopic to the identity in $\localquad$. By applying \cite[Thm 2.14]{MR3281141} pointwise\footnote{See also the survey \cite[Prop 6]{nikolaus2017group}.}, we get that the colimit \emph{in prestacks}
  \[
    C_\PSh \coloneqq \colim \left[
    \begin{tikzcd}[cramped]
      \localquad \ar{r}{- \oplus Q_1} & \localquad \ar{r}{- \oplus Q_1} & \cdots
    \end{tikzcd}
    \right] \in \PSh(X_\et^\daff, \inftygpd)
  \]
  carries a symmetric monoidal structure such that
  \begin{enumerate}
    \item the canonical morphism $\localquad \to C_\PSh$ is monoidal and
    \item the image of $Q_1$ is tensor-invertible as a section of $C_\PSh$.
  \end{enumerate}
  The monoidal prestack $C_\PSh$ is moreover universal for those properties.
  Since stackification preserves finite products, the same holds for the stackification $C$ of $C_\PSh$.

  Finally, the image of the quadratic bundle $Q_1$ in $\localquadplus$ is invertible, we get a canonical monoidal morphism
  \[
    \Theta \colon C \to \localquadplus.
  \]
  Moreover, as any section of $\localquad$ is locally trivial by definition (see also \cref{formulaA1localizationpreservespizerosheaf}), it is locally a power of $Q_1$.
  It follows that every section of $C$ is locally tensor-invertible, and thus globally tensor-invertible.
  In particular, $C$ satisfies the universal property of the group-completion, so $\Theta$ is an equivalence.
  See also \cite[Lem.\,3.8(1)]{zbMATH07623470} for a similar result.
\end{proof}

\begin{observation}\label{PropositionidentificationGrothendieckWittgroup}
  This lemma implies in particular that for any point $x \in X$, we have $\pi_0\left((\QA)^{+}_x , 0\right) \simeq \integers$.
  This could be obtained directly, as $\pi_0\left((\QA)^{+}_x , 0\right) $ identifies with the Grothendieck--Witt group of $\basefield$ and is thus isomorphic to $\integers$ via the rank morphism.
\end{observation}

\subsubsection{Deducing contractibility}

\begin{proof}[\textbf{Proof of \cref{contractibilitytheoremdetailedversion} assuming \cref{prop:A1-automorphisms}}]

\newcommand{\quotientA}{\quotA}
  Since we are working with hypersheaves, to prove the quotient stack $\quotientA$ is contractible, it suffices to prove that its stalks at geometric points are contractible.
  We therefore fix $x \in X$ a geometric point.
  First, \cref{thm-transitivity} implies that the stalk $(\quotientA)_x$ is connected.
  We can thus compute its higher homotopy groups at any given point.
  Let $\formalU_0$ be a minimal chart at $x$ (up to shrinking $X$, we may assume $\formalU_0$ is defined over $X$).
  We will prove
  \[
    \pi_n\left(\left(\quotientA\right)_x, [\formalU_0]\right) \simeq \homotopysheaf_n\left(\quotientA, [\formalU_0]\right)_x \simeq 0,
  \]
  where we write $\homotopysheaf_n$ for the homotopy sheaves. The Cartesian square
  \[
    \begin{tikzcd}
     \left( \QA\right)^{+} \ar{r}{\formalU_0 \oplus -} \ar{d} \tikzcart & \left(\DA\right)^{+} \ar{d} \\
      * \ar{r}[swap]{[\formalU_0]} & \quotientA
    \end{tikzcd}
  \]
  induces a long exact sequence of abelian sheaves (or of pointed sheaves for $n - 1 = 0$)
  \[
    \begin{tikzcd}[cramped,column sep=small, row sep=small]
      \homotopysheaf_n\left( \left( \QA\right)^{+} , 0\right) \ar{r}{\formalU_0 \oplus -} &
      \homotopysheaf_n\left(\left(\DA\right)^{+} , \formalU_0\right) \ar{r}
      \arrow[d, phantom, ""{coordinate, name=Z}] &
      \homotopysheaf_n\left(\quotientA, [\formalU_0]\right)
      \arrow[dll, rounded corners, overlay, to path={ -- ([xshift=2ex]\tikztostart.east) |- (Z) -| ([xshift=-2ex]\tikztotarget.west) -- (\tikztotarget)}]
      \\
      \homotopysheaf_{n-1}\left( \left( \QA\right)^{+} , 0\right) \ar{r} &
      \homotopysheaf_{n-1}\left(\left(\DA\right)^{+} , \formalU_0\right) \ar{r} &
      \homotopysheaf_{n-1}\left(\quotientA, [\formalU_0]\right)\rlap{.}
    \end{tikzcd}
  \]
 Since taking the stalks at $x$ is given by a filtered colimit, it commutes with taking loops and therefore, we see that it is enough to prove the map
  \[
    \alpha_n \coloneqq \formalU_0 \oplus - \colon \pi_n\left( \left( \QA\right)^{+}_x , 0\right) \to \pi_n\left(\left(\DA\right)^{+}_x , \formalU_0\right)
  \]
  is an isomorphism (of groups for $n \ge 1$, of pointed sets for $n = 0$).

 We already know (by \cref{thm-transitivity}) that $\alpha_0$ is surjective.
  Moreover, by \cref{PropositionidentificationGrothendieckWittgroup}, we have $\pi_0\left( \left( \QA\right)^{+}_x , 0\right)\simeq \integers$ via the rank morphism.
  The dimension provides an invariant of Darboux charts, which induces a map
  \[
    \dim \colon \pi_0\left(\left(\DA\right)^{+}_x, \formalU_0\right) \to \integers.
  \]
  The composite map
  \[
    \begin{tikzcd}
      \integers \simeq \pi_0\left( \left( \QA\right)^{+}_x , 0\right) \ar{r}{\alpha_0} & \pi_0\left(\left(\DA\right)^{+}_x, \formalU_0\right) \ar{r}{\dim} & \integers
    \end{tikzcd}
  \]
  identifies with the shift by $\dim \formalU_0$ and is thus injective. It follows that $\alpha_0$ is a bijection.

 We now focus on $\alpha_n$ for $n \geq 1$. For $p \geq 0$, denote by $Q_p$ the trivial quadratic bundle with fiber $\basefield^p$ and form $\sum x_i^2$, and by $\formalU_p \coloneqq \formalU_0 \oplus Q_p$.
  By \cref{inversehyperbolic}, we can identify $\alpha_n$ as the colimit morphism $\alpha_n \simeq \colim_p \alpha_n^p$:
  \[
    \begin{tikzcd}[row sep=small]
      \pi_n\left(\left(\QA\right)_{x} , 0 \right) \ar{d} \ar{r}{\alpha_n^0} &   \pi_n\left(  \left(\DA\right)_x  , \formalU_0 \right) \ar{d}\\
      \vdots \ar{d}& \vdots \ar{d}\\
      \pi_n\left(\left(\QA\right)_{x} , Q_p \right) \ar{d}{- \oplus Q_1} \ar{r}{\alpha_n^p = \formalU_0 \oplus -} &   \pi_n\left(  \left(\DA\right)_x , \formalU_p \right) \ar{d}{- \oplus Q_1} \\
      \pi_n\left(\left(\QA\right)_{x} , Q_{p+1} \right) \ar{r}{\alpha_n^{p+1}} \ar{d} &   \pi_n\left(\left(\DA\right)_x  , \formalU_{p+1} \right) \ar{d} \\
      \vdots\ar{d}&\vdots \ar{d}\\
      \pi_n\left( \left( \QA\right)^{+}_x, 0\right) \ar{r}{\alpha_n}&\pi_n\left(  \left(\DA\right)^{+}_x , \formalU_0\right)\rlap{.}
    \end{tikzcd}
  \]
  Now, using the assumption that $\formalU_0$ is minimal, by \cref{prop:A1-automorphisms}, the morphism $\stackAut(Q_p) \to \stackAut(\formalU_p)$ induces an equivalence at the stalk at  $x$ after applying  $\PhiA{}$
  \[
    \PhiA{}\left(\stackAut(Q_p)\right)_{x} \to^\sim  \PhiA{}\left(\stackAut(\formalU_p)\right)_{x}.
  \]
  But thanks to \cref{automorphismsinthedarbouxstacksareA1localizations}, this means that we have an equivalence at the level of loop stacks
 \[
 \left(\Omega_{Q_p} \QA\right)_x\simeq \left(\Omega_{\formalU_p}\DA\right)_x
 \]
 In particular, the morphism $\alpha_n^p$ is an isomorphism for any $n$ and any $p$:
  \[
    \begin{tikzcd}[column sep=.3em]
      \pi_n\left(\left(\QA\right)_{x}, Q_p\right) \ar[phantom]{r}{\simeq} \ar{d}[swap]{\alpha_n^p} & \pi_{n-1}\left(  \PhiA{}\left(\stackAut(Q_p)\right)_{x}, \Id \right) \ar{d}[sloped]{\sim}
      \\
      \pi_n\left(\left(\DA\right)_x, \formalU_p \right) \ar[phantom]{r}{\simeq} & \pi_{n-1}\left(\PhiA{}\left(\stackAut(\formalU_p)\right)_{x}., \Id \right).
    \end{tikzcd}
  \]
  It follows that $\alpha_n$ is itself an isomorphism as well. The result follows.
\end{proof}

We conclude this section with a version of  \cref{contractibilitytheoremdetailedversion} for $\catsuspension\affineline{}$\nobreakdashes-localizations. This will not be needed in this paper but we keep it for future reference:

\begin{proposition}\label{contractibilitytheoremdetailedversion2}
  Let $X$ be a derived Deligne--Mumford stack equipped with a $(-1)$\nobreakdashes-shifted exact symplectic form $\lambda$. Then the quotient stack
  \begin{equation}
  \label{quotientstackSigmalocalizationcontractible}\quot{\left(\Darbstack^{\lambda}_{X,\isotopies}\right)}{\left(\stackQuadnabla_{X,\isotopies}\right)} \simeq \constantsheaf \ast_X.
  \end{equation}
  is contractible.
\end{proposition}
\begin{proof}
  The proof runs essentially the same steps as in the proof of \cref{contractibilitytheoremdetailedversion}, assuming \cref{prop:A1-automorphisms}. First, since the maps \cref{definitionA1darbouxstackCartesianmap}  are isomorphism on the sheaves $\homotopysheaf_0$, the transitivity of \cref{thm-transitivity} implies that the quotient (hyper)stack of \cref{quotientstackSigmalocalizationcontractible} is connected. We now fix $\formalU_0$ a minimal model at $x$. All it remains is to show that the map of loop stacks induced by the action
  \[
\Omega_{0}  \left(\stackQuadnabla_{X,\isotopies}\right)^+\to \Omega_{\formalU_0}\left(\Darbstack^{\lambda}_{X,\isotopies}\right)^+
  \]
  is an equivalence after taking stalks as $x$. But now the proof of \cref{inversehyperbolic} combined with \cref{automorphismsinthedarbouxstacksareA1localizationsla}, implies that the group completions are also given by the filtered colimit formulas
  \begin{align*}
    \left(\stackQuadnabla_{X,\isotopies}\right)^+ &\simeq \colim \left[
    \begin{tikzcd}[cramped,ampersand replacement=\&]
      \stackQuadnabla_{X,\isotopies}\ar{r}{- \oplus Q_1} \& \stackQuadnabla_{X,\isotopies}\ar{r}{- \oplus Q_1} \& \cdots
    \end{tikzcd} \right]\rlap{\hspace{1em}and}\\
    \left(\Darbstack^{\lambda}_{X,\isotopies}\right)^+ &\simeq \colim \left[
    \begin{tikzcd}[cramped,ampersand replacement=\&]
      \Darbstack^{\lambda}_{X,\isotopies}\ar{r}{- \oplus Q_1} \& \Darbstack^{\lambda}_{X,\isotopies} \ar{r}{- \oplus Q_1} \& \cdots
    \end{tikzcd}
    \right].
  \end{align*}
  and therefore we are reduce to showing that the action map
  \[
    \Omega_{Q_p} \stackQuadnabla_{X,\isotopies} \to \Omega_{\formalU_p}\Darbstack^{\lambda}_{X,\isotopies}
  \]
  is an equivalence at the stalk at $x$. We conclude combining \cref{prop:A1-automorphisms}, \cref{A1contractionofQuadSigmalocalization}, \cref{automorphismsinthedarbouxstacksareA1localizationsla} and \cref{lemmaappliedtoourcase}.
\end{proof}

\subsection{Linearization and proof of the contractibility theorem}
\label{sectionproofcontractibility}

The goal of this section is to prove \cref{prop:A1-automorphisms}.
The key idea is to prove that the stack of endomorphisms of a Darboux chart is (locally) an affine space.
Recall that such an endomorphism $\phi \colon \formalU \to \formalU$ amounts to a commutative diagram
\[
  \begin{tikzcd}[row sep=0]
    & \formalU \ar{dr}{f} \ar{dd}{\phi} \\
    X \ar{ru}{i} \ar{dr}[swap]{i} && \affineline{} \\
    & \formalU \ar{ur}[swap]{f}
  \end{tikzcd}
\]
and some compatibilities of the symplectic data. It is relatively straightforward to interpolate morphisms under $X$ (see \cref{lem:endunderXlinear} below).
The main difficulty comes from the necessity to preserve the function $f$.
This is dealt with by a linearization argument explained below.

\subsubsection{Connections, exponential isomorphisms}
  We fix $\formalU$ a smooth formal affine scheme.
  Let $\widehat \Delta$ denote the formal neighborhood of its diagonal and let $\widehat \tangentbundle_\formalU$ denote the formal neighborhood of the zero vector field.

\begin{definition}
  Let $\nabla$ be a (torsion-free) connection on $\formalU$.
  We denote by
  \[
    \exp_\nabla \colon \widehat \tangentbundle_\formalU \to^\sim \widehat \Delta
  \]
  the exponential isomorphism associated to $\nabla$ (see \eg \cite{Kapranov1999}) fitting in a commutative diagram
  \[
    \begin{tikzcd}
      \formalU \ar{r}{\Delta} \ar{d}[swap]{s_0} & \widehat \Delta \ar{d}{p_1} \\ \widehat \tangentbundle_\formalU \ar{r} \ar{ur}[swap,near start]{\exp_\nabla}[sloped]{\sim} & \formalU,
    \end{tikzcd}
  \]
  Recall that it is determined by the formula\footnote{The usual $\dfrac{1}{n!}$ factor only appears after the symmetrization map $\Sym^n \Omega_\formalU \to \Omega_\formalU^{\otimes n}$.}
  \[
   \begin{tikzcd}[row sep=-.5em, column sep=small]
    \structuresheaf_\formalU \ar{r}{p_2^*} & \structuresheaf_{\widehat \Delta} \ar{r}{\exp_\nabla^*} & \Symcompleted \Omega_\formalU
    \\
    f \ar[mapsto]{rr} && \sum d^n_\nabla f
   \end{tikzcd}
  \]
  where $d^n_\nabla \colon \structuresheaf_\formalU \to \Sym^n \Omega_\formalU$ is the $n$\nobreakdashes-th derivative operator associated to $\nabla$.
  Note that $d^0_\nabla f = f$ and $d^1_\nabla f = df$ by construction.
  We denote by $\log_\nabla \colon \widehat \Delta \to \widehat \tangentbundle_\formalU$ the inverse of $\exp_\nabla$.
\end{definition}

\begin{observation}
  Considered over $\formalU$, the formal neighborhood $\widehat \Delta$ classifies jets of endomorphisms of $\formalU$, seen as their graph.
  We will sometimes write $\widehat \stackEnd(\formalU)$ for $\widehat \Delta$ in order to emphasize this point.
  In the presence of a connection $\nabla$, we thus have
  \[
    \exp_\nabla \colon \widehat \tangentbundle_\formalU \simeq \widehat \Delta = \widehat \stackEnd(\formalU).
  \]
\end{observation}
The next lemma extends this equivalence to the stack of endomorphism fixing a reduced equivalent closed derived subscheme.

\begin{notations}{}
  \item Let $S \to \formalU$ be a morphism. We denote by $\stackEnd_{S/}(\formalU)$ the stack over the affine site of $\formalU$ classifying endomorphisms of $\formalU$ fixing $S$.
  \item Let $E$ be a perfect complex over $\formalU$. We denote by $\vectorbundle_\formalU(E)$ the associated linear stack over $\formalU$.
  In particular, we have $\tangentbundle_{\formalU} \coloneqq \vectorbundle_\formalU(\tangent_\formalU)$.
\end{notations}

\begin{lemma}\label{lem:endunderXlinear}
  Let $\formalU$ be a smooth formal scheme and let $i \colon S \to \formalU$ be an lci derived closed subscheme.
  Denote by $I_S$ the homotopy fiber of $\structuresheaf_\formalU \to i_* \structuresheaf_S$ and fix a connection $\nabla$ on $\formalU$.
  If $i$ is a reduced equivalence, then there is an equivalence
  \[
    \exp_\nabla^S \colon \vectorbundle_\formalU\left(\tangent_\formalU \otimes I_S\right) \simeq \stackEnd_{S/}(\formalU) :\, \log_\nabla^S.
  \]
  Moreover, the natural morphisms $\vectorbundle_\formalU(\tangent \otimes I_S) \to \tangentbundle_\formalU$ and $\stackEnd_{S/}(\formalU) \to \stackEnd(\formalU)$ factor through $\widehat \tangentbundle_\formalU$ and $\widehat \stackEnd(\formalU)$ respectively, and the following diagram commutes:
  \[
    \begin{tikzcd}[row sep=tiny]
      &
      \vectorbundle_\formalU\left(\tangent_\formalU \otimes I_S\right) \ar{dd}{\exp_\nabla^S}[swap,sloped]{\sim} \ar{r} &
      \widehat \tangentbundle_\formalU \ar{dd}{\exp_\nabla}[swap,sloped]{\sim}
    \\
      \formalU \ar{ur}{s_0} \ar{dr}[swap]{\{\Id\}}
    \\
      &
      \stackEnd_{S/}(\formalU) \ar{r} &
      \widehat \stackEnd(\formalU)\rlap{.}
    \end{tikzcd}
  \]
\end{lemma}

\begin{sublemma}\label{mappingstacktoformalneigbourhood}
  Consider $U \from^i S \to V$ be morphisms of derived stacks. Denote by $\widehat V$ the formal neighborhood of $S$ in $V$.
  If $i \colon S \to U$ is a reduced equivalence, then the morphism
  \[
    \theta \colon \stackMap_{S/}(U, \widehat V) \to \stackMap_{S/}(U, V)
  \]
  is an equivalence.
\end{sublemma}

\begin{proof}
  By definition, $\widehat V \coloneqq V \times_{V_\deRham} S_\deRham$, so the map $\theta$ fits into a Cartesian square:
  \[
    \begin{tikzcd}
      \stackMap_{S/}(U, \widehat V) \ar{r}{\theta} \ar{d} \tikzcart & \stackMap_{S/}(U, V) \ar{d} \\
      \stackMap_{S/}(U, S_\deRham) \ar{r} & \stackMap_{S/}(U, V_\deRham)\rlap{.}
    \end{tikzcd}
  \]
  However, since $i \colon S \to U$ is a reduced equivalence, both the stacks $\stackMap_{S/}(U, S_\deRham)$ and $\stackMap_{S/}(U, V_\deRham)$ are contractible, and $\theta$ is therefore an equivalence.
\end{proof}

\begin{proof}[Proof of \cref{lem:endunderXlinear}]
  By definition, the stack $\stackEnd_{S/}(\formalU)$ is the relative mapping stack
  \[
    \stackEnd_{S/}(\formalU) \coloneqq \stackMap_{S/-/\formalU}(\formalU, \formalU \times \formalU)
  \]
  where $\formalU \times \formalU$ is seen as a stack over $\formalU$ via the first projection (so an endomorphism is represented by its graph) and the map $S \to \formalU \times \formalU$ is $i$ followed by the diagonal immersion.

  Since $i \colon S \to \formalU$ is a reduced equivalence, \cref{mappingstacktoformalneigbourhood} implies
  \[
    \stackEnd_{S/}(\formalU) \coloneqq \stackMap_{S/-/\formalU}(\formalU, \formalU \times \formalU) \simeq \stackMap_{S/-/\formalU}\left(\formalU, \widehat \Delta\right).
  \]
  In particular, the morphism
  \[
    \stackEnd_{S/}(\formalU) = \stackMap_{S/-/\formalU}(\formalU, \formalU \times \formalU) \to \stackMap_{\formalU}(\formalU, \formalU \times \formalU) = \stackEnd(\formalU)
  \]
  naturally factors through $\widehat \stackEnd(\formalU) = \widehat \Delta$.

  Fix now a connection $\nabla$ on $\formalU$.
  Its exponential isomorphism and \cref{mappingstacktoformalneigbourhood} yield
  \[
    \stackEnd_{S/}(\formalU) \simeq \stackMap_{S/-/\formalU}\left(\formalU, \widehat \Delta\right) \simeq \stackMap_{S/-/\formalU}\left(\formalU, \widehat \tangentbundle_\formalU\right) \simeq \stackMap_{S/-/\formalU}(\formalU, \tangentbundle_\formalU).
  \]
  In turn, the stack $\stackMap_{S/-/\formalU}(\formalU, \tangentbundle_\formalU)$ classifies vector fields nullhomotopic on $S$, so morphisms $\Omega_\formalU \to I_S$.
  We deduce the equivalence $\stackEnd_{S/}(\formalU) \simeq \vectorbundle_\formalU(\tangent_\formalU \otimes I_S)$, using that $I_S$ is perfect (recall that we assumed the map $i$ to be lci so that $\structuresheaf_S$ is perfect over $\structuresheaf_\formalU$).
\end{proof}

\subsubsection{Taylor expansion and function-preserving endomorphisms}

We will now be interested in endomorphisms preserving a given function.
We thus fix a function $f \colon \formalU \to \affineline{}$.

Define on $\widehat \stackEnd(\formalU) = \widehat \Delta$ the function $F = (-f, f) = p_2^* f - p_1^* f$.
\begin{lemma}
  The zero locus of $F \colon \widehat \stackEnd(\formalU) \to \affineline{}$ is equivalent to the stack of jets of function-preserving endomorphisms
  \[
    \widehat \stackEnd_{/\affineline{}}(\formalU) \simeq \zerolocus(F).
  \]
\end{lemma}

\begin{proof}
  It boils down to the equivalence
  \[
    \widehat \stackEnd_{/\affineline{}}(\formalU) = \stackMap_\formalU(\formalU, \widehat \Delta) \times_{\stackMap_\formalU(\formalU, \affineline{\formalU})} \{ 0 \} \simeq \widehat \Delta \times_{\affineline{\formalU}} \{F\} = \zerolocus(F).\qedhere
  \]
\end{proof}

In the presence of a connection $\nabla$, the restriction of $p_2^* f$ along $\exp_\nabla$ is called the (global) $\nabla$\nobreakdashes-Taylor expansion of $f$.
It is given by the series
\[
  \taylorexpansion_\nabla(f) \coloneqq \exp_\nabla^* p_2^* f = \sum_{n \geq 0} d^n_\nabla f \in \Symcompleted \Omega_\formalU.
\]
Since $d_\nabla^0 f = f$, the function $F$ thus corresponds to the series $\taylorexpansion_\nabla^{\geq1}(f) \coloneqq \exp_\nabla^*(F) = \sum_{n \geq 1} d^n_\nabla f$ and is in particular not linear.
Its zero locus, classifying endomorphisms preserving the function, is therefore not linear a priori.

Assuming further we have an lci closed derived subscheme $i \colon S \to \formalU$ such that $i$ is a reduced equivalence, we can use \cref{lem:endunderXlinear} to study function-preserving and $S$\nobreakdashes-fixing endomorphisms.
In the specific situation where $S$ is the derived critical locus of $f$, we can go further and consider isotropic morphisms.

\subsubsection{Isotropic morphisms}

Fix $\formalU$ a smooth formal scheme and $f \colon \formalU \to \affineline{}$ a function.
We will now be interested in linearizing the stack of endomorphisms of $\formalU$ seen as a critical chart of $\dCrit(f)$ -- \ie of endomorphisms fixing $S \coloneqq \dCrit(f)$ and preserving the isotropic structure.

Assume that the inclusion $i \colon S \coloneqq \dCrit(f) \to \formalU$ is a reduced equivalence.
Denote by $\cdots \to I_S^{(n)} \to I_S^{(n-1)} \to \cdots \to I_S^{(0)} = \structuresheaf_\formalU$ the relative Hodge filtration defined by $I_S^{(n)} \coloneqq \derhamcomplex^{\geq n}_{S/\formalU}$.
Recall that we have $I_S^{(1)} = I_S$ as well as a fiber sequence for each integer $n$:
\[
  I^{(n+1)}_S \to I_S^{(n)} \to i_* i^* \Sym^n(\tangent_\formalU).
\]
By construction, the image of the symplectic form $\omega_S \colon \basefield \to \derhamcomplex^{\geq 2}_{S}[1]$ of $S \coloneqq \dCrit(f)$ in $I_S{(2)}[1] = \derhamcomplex^{\geq 2}_{S/\formalU}[1]$ is nullhomotopic:
\[
  \begin{tikzcd}
    \basefield \ar{r} \ar[bend right]{rr}[yshift=1.5em,xshift=-.2em,rotate=-90]{\displaystyle \sim}[yshift=.3em,xshift=-.7em]{\eta}[swap]{0} & \derhamcomplex^{\geq 2}_{S}[1] \ar{r} & \derhamcomplex^{\geq 2}_{S/\formalU}[1] \mathrlap{{} = I^{(2)}_S[1].}
  \end{tikzcd}
\]

\begin{definition}
  Fix an endomorphism $\phi \colon \formalU \to \formalU$ fixing $S \coloneqq \dCrit(f)$.
  The isotropy defect of $\phi$ is the difference $\phi^*\eta - \eta \in I^{(2)}_S[1][-1] = I^{(2)}_S$.We denote by $\delta_\eta$ the morphism
  \[
    \delta_\eta \colon \stackEnd_{S/}(\formalU) \to \vectorbundle_\formalU(I^{(2)}_S)
  \]
  mapping an endomorphism $\phi$ to its isotropy defect.
  By construction, the zero-locus of $\delta_\eta$ is the stack of endomorphisms of $\formalU$ fixing $S$ and preserving the isotropic structure $\eta$: \ie isotropic morphisms.
\end{definition}

\begin{observations}{}
  \item By \cref{cor:map=action}, an isotropic endomorphism is in fact formally étale, and is thus invertible.
  \item An isotropic morphism in particular preserves the function $f$, as the following diagram commutes
  \[
    \begin{tikzcd}
      \stackEnd_{S/}(\formalU) \ar{d}{\delta_\eta} \ar{r} & \widehat \stackEnd(\formalU) \ar{d}{F}
    \\
      \vectorbundle_\formalU(I^{(2)}_S) \ar{r}{I_S^{(2)} \to \structuresheaf_\formalU} & \affineline{\formalU} \mathrlap{{} = \vectorbundle_\formalU(\structuresheaf_\formalU).}
    \end{tikzcd}
  \]
\end{observations}

\begin{lemma}\label{lem-describedeltaeta}
  Fix $\nabla$ a torsion-free connection on $\formalU$. The composite morphism
  \[
    \begin{tikzcd}
      \vectorbundle_\formalU(\tangent_\formalU \otimes I_S) \ar{r}{\exp_\nabla^S}[swap]{\sim} & \stackEnd_{S/}(\formalU) \ar{r}{\delta_\eta} & \vectorbundle_\formalU(I^{(2)}_S)
    \end{tikzcd}
  \]
  maps $\varphi \colon \Omega_\formalU \to I_S$ to $\sum_{n \geq 1} K_n(\varphi)$ where
  \[
  \begin{array}{ll}
    \begin{tikzcd}[ampersand replacement=\&,column sep={5em,between origins}]
      K_n(\varphi) \colon \structuresheaf_\formalU \ar{r}{d_\nabla^nf} \&[1em] \Sym^n \Omega_\formalU \ar{r}{\Sym^n \varphi} \&[2em] \Sym^n I_S \ar{r} \& I_S^{(n)} \ar{r} \& I_S^{(2)}
    \end{tikzcd}
    &
    \text{for } n \geq 2 \text{ and}
    \\
    \begin{tikzcd}[ampersand replacement=\&,column sep={5em,between origins}]
    K_1(\varphi) \colon \structuresheaf_\formalU \ar{r}{df} \&[1em] I_S \Omega_\formalU \ar{r}{\varphi} \&[2em] I_S^{(2)}.
    \end{tikzcd}
    &
  \end{array}
  \]
\end{lemma}
\begin{proof}
  By definition of $\exp_\nabla^S$, its image of $\varphi \colon \Omega_\formalU \to I_S$ is the endomorphism $e_\varphi \colon \structuresheaf_\formalU \to \structuresheaf_\formalU$ given by $g \mapsto \sum_{n \geq 0} u_n \varphi^{(n)}(d_\nabla^n g)$ (where $\varphi^{(n)} \colon \Sym^n \Omega_\formalU \to I_S^{(n)}$ and $u_n \colon I_S^{(n)} \to \structuresheaf_\formalU$).
  Computing the isotropy defect from this formula yields the announced result.
\end{proof}

\subsubsection{Linearization of the isotropy defect}

\begin{proposition}\label{linearizetaylorexpansion}
  Let $\formalU_0$ be a smooth formal scheme over $\basefield$ and $f_0 \colon \formalU_0 \to \affineline{}$ a function.
  Let $(Q,q)$ a $\basefield$\nobreakdashes-vector space equipped with a non-degenerate quadratic form $q$.
  Denote by $\formalU$ the formal neighborhood of the zero-section in $\formalU_0 \times Q$ and $f \colon \formalU \to \affineline{}$ the function $f_0 + q$.
  Fix a connection $\nabla_0$ on $\formalU_0$ and extend it (using the trivial connection of $Q$) to a connection $\nabla$ on $\formalU$.
  Let $x \colon \Spec \basefield \to S \subset \formalU_0 \subset \formalU$ be a $\basefield$\nobreakdashes-point.
  Assume that
  \begin{enumerate}[label={\textrm{\upshape{(\alph*)}}}]
    \item The canonical inclusion $i_0 \colon S \coloneqq \dCrit(f_0) \to \formalU_0$ is a reduced equivalence;
    \item The Hessian of $f$ vanishes at $x \in S \subset \formalU_0$ (so that $\formalU_0$ is minimal at $x$, as a critical chart of $S$).
  \end{enumerate}
  Then the form $df \colon \tangent_\formalU \to \structuresheaf_\formalU$ factors via $I_S \coloneqq \fiber(\structuresheaf_\formalU \to i_* \structuresheaf_S)$ and there exists a (non-linear) endomorphism $L \colon \vectorbundle_\formalU(\tangent_\formalU \otimes \tangent_\formalU) \to \vectorbundle_\formalU(\tangent_\formalU \otimes \tangent_\formalU)$ such that
  \begin{enumerate}[label={\textit{(\roman*)}},ref={\textit{(\roman*)}}]
    \item\label{linearizetaylorexpansion-commutes}
      The following diagrams commute
      \[
        \begin{tikzcd}
          \formalU \ar{r}{s_0} \ar{d}[swap]{s_0} & \vectorbundle_\formalU(\tangent_\formalU \otimes \tangent_\formalU) \ar{dl}{L} \ar{d} \\
          \vectorbundle_\formalU(\tangent_\formalU \otimes \tangent_\formalU) \ar{r} & \formalU
        \end{tikzcd}\hspace{1cm}
        \begin{tikzcd}
          \vectorbundle_\formalU(\tangent_\formalU \otimes \tangent_\formalU) \ar{r}{\Id \otimes df} \ar{d}{L} & \vectorbundle_\formalU(\tangent_\formalU \otimes I_S) \ar{d}{\delta_\eta} \\
          \vectorbundle_\formalU(\tangent_\formalU \otimes \tangent_\formalU) \ar{r}[swap]{df \otimes df} & \vectorbundle_\formalU(I_S^{(2)})\rlap{.}
        \end{tikzcd}
      \]
    \item\label{linearizetaylorexpansion-fix}
      The subspace $\vectorbundle_\basefield(\tangent_{\formalU_0,x} \otimes Q)$ of $\vectorbundle_\basefield(\tangent_{\formalU,x} \otimes \tangent_{\formalU,x})$ is fixed by $L_{|x}$.
    \item\label{linearizetaylorexpansion-etale}
      Let $E$ be the fiber product $\tangent_\formalU \otimes \tangent_\formalU \times_{x_* x^*(\tangent_\formalU \otimes \tangent_\formalU)} x_* (\tangent_{\formalU_0,x} \otimes Q)$ and let $Y \coloneqq \vectorbundle_\formalU(E)$.
      By \ref{linearizetaylorexpansion-fix}, it is stable under $L$.
      Moreover, the induced morphism $L \colon Y \to Y$ is formally étale in the neighborhood of its fiber over $x$.
    \item\label{linearizetaylorexpansion-artin}
      There exists an étale neighborhood $T$ of $x \in \formalU$ such that $L_{|T} \colon Y \times_\formalU T \to Y \times_\formalU T$ is invertible.
  \end{enumerate}
\end{proposition}

\begin{proof}
  Consider the formula
  \begin{multline*}
    L(\omega \otimes \lambda) = \omega \otimes \lambda \\ + \sum_{n \geq 2} \sum (d^n_\nabla f_{(1)} \otimes df) \otimes^{\mathrm s} \cdots \otimes^{\mathrm s} (d^n_\nabla f_{(n-2)} \otimes df) \otimes^{\mathrm s} (d^n_\nabla f_{(n-1)} \otimes \omega) \otimes^{\mathrm s} (d^n_\nabla f_{(n)} \otimes \lambda),
  \end{multline*}
  where $d^n_\nabla f_{(j)}$ denote Sweedler components for $d^n_\nabla f$.
  Since $df$ has values in the pro-nilpotent radical of $\structuresheaf_\formalU$, the above series defines a morphism $\Omega_\formalU \otimes \Omega_\formalU \to \Sym_\formalU(\Omega_\formalU \otimes \Omega_\formalU)$, leading to the announced (non-linear) endomorphism of $\vectorbundle_\formalU(\tangent_\formalU \otimes \tangent_\formalU)$.
  Assertion \ref{linearizetaylorexpansion-commutes} follows directly from the definition and \cref{lem-describedeltaeta}.

  To prove \ref{linearizetaylorexpansion-fix}, we first observe that $df$ vanishes at $x$, so that every term corresponding to $n \geq 3$ in the above formula vanishes.
  The term for $n = 2$ restricted to $\vectorbundle_k(\tangent_{\formalU_0,x} \otimes Q)$ also vanishes since the Hessian $d^2_{\nabla_0}(f_0)$ vanishes at $x$ by assumption.

  We now focus on \ref{linearizetaylorexpansion-etale}.
  Let $\delta$ denote the relative de Rham differential on $\vectorbundle_\formalU(\tangent_\formalU \otimes \tangent_\formalU)$ over $\formalU$.
  Since $df$ vanishes at $x$, we find on $\vectorbundle_\basefield(\tangent_{\formalU_0,x} \otimes Q)$:
  \begin{align}
    D_L^\dual(\delta(\omega \otimes \tau)) = {}& \delta(\omega \otimes \tau) \tag{A} \\
    &+ \sum \delta(d^2_\nabla f_{(1)} \otimes \omega) \otimes^{\mathrm{s}} (d^2_{\nabla} f_{(2)})_{|\tangent_{\formalU_0,x}} \otimes \tau_{|Q} \tag{B} \\
    &+ \sum (d^2_{\nabla} f_{(1)})_{|\tangent_{\formalU_0,x}} \otimes \omega_{|Q} \otimes^{\mathrm{s}} \delta(d^2_\nabla f_{(2)} \otimes \tau) \tag{C} \\
    &+ \sum \delta(d^3_\nabla f_{(1)} \otimes df) \otimes^{\mathrm{s}} (d^3_{\nabla} f_{(2)})_{|\tangent_{\formalU_0,x}} \otimes \omega_{|Q} \otimes^{\mathrm{s}} (d^3_{\nabla} f_{(3)})_{|\tangent_{\formalU_0,x}} \otimes \tau_{|Q}. \tag{D}
  \end{align}
  With $df$ vanishing at $x$ and $\delta$ being $\structuresheaf_\formalU$\nobreakdashes-linear, we deduce that the term (D) vanishes.

  We split $d^2_\nabla f \in \Sym^2_\formalU(\Omega_\formalU)$ into the sum $d^2_\nabla f = d^2_{\nabla_0} f_0 + d^2 q$, where $d^2_{\nabla_0} f_0 \in \Sym_{\formalU_0}^2(\Omega_{\formalU_0})$ and $d^2 q \in \Sym^2_k(Q^\dual)$.
  As $d^2_{\nabla_0} f_0$ vanishes at $x$, we can choose a decomposition $\sum_i a_i \otimes b_i$ of $d^2_{\nabla_0} f_0 \in \Sym^2_{\formalU_0}(\Omega_{\formalU_0}) \subset \Omega_{\formalU_0}^{\otimes 2}$ where each $b_i$ vanishes at $x$.
  Again, $\delta$ being $\structuresheaf_\formalU$\nobreakdashes-linear allows us to conclude that both (B) and (C) vanish as well.
  It follows that the cotangent morphism of $L$ relatively to $\formalU$ is the identity at every point of $\vectorbundle_\basefield(\tangent_{\formalU_0,x} \otimes Q)$, so that $L$ is indeed formally étale.

  At last, assertion \ref{linearizetaylorexpansion-artin} follows from Artin's algebraization theorem (\cf \cite{MichaelArtin921}).
  Let $S^{(2)}$ denote the second infinitesimal neighborhood of $\Crit(f)$ in $\formalU$ and $Y^{(2)} \coloneqq Y \times_\formalU S^{(2)}$.
  Denote by $Y_x$ the fiber of $Y$ at $x \in S^{(2)} \subset \formalU$. By definition, we have $Y_x = \vectorbundle_\basefield(\tangent_{\formalU_0,x} \otimes Q)$.
  Since $L \colon Y \to Y$ is the identity on $Y_x$ and is étale in its neighborhood, it restricts to an automorphism of the formal neighborhood $\widehat Y_x$ of $Y_x$ in $Y^{(2)}$.
  Consider the moduli problem $\calF$ (over $\formalU$) of inverses of $L \colon Y^{(2)} \to Y^{(2)}$.
  Since $\widehat Y_x$ identifies with the restriction of $Y^{(2)}$ on the formal neighborhood $\widehat x$ of $x$ in $S^{(2)}$, we get a formally universal section $\sigma \colon \widehat x \to \calF$.
  As $\calF$ is locally of finite presentation over $S^{(2)}$ and $S^{(2)}$ is affine and excellent, Artin algebraization (over $S^{(2)}$) implies there exists an affine scheme $T^{(2)}$ and a factorization $\widehat x \to T^{(2)} \to \calF$ of $\sigma$, where $\widehat x \to T^{(2)}$ is étale.
  It follows that $T^{(2)} \to \calF \to S^{(2)}$ is étale in the neighborhood of $x \in T^{(2)}$.
  Shrinking $T^{(2)}$ if necessary, we have an étale neighborhood $T^{(2)} \to S^{(2)}$ of $x$ such that $L|_{T^{(2)}} \colon Y \times_\formalU T^{(2)} \to Y \times_\formalU T^{(2)}$ is invertible.
  Denote by $T \to \formalU$ the unique étale neighborhood of $x$ such that $T \times_{\formalU} S^{(2)} \simeq T^{(2)}$ (recall that $S^{(2)}_\red = S_\red = \formalU_\red$).
  The fact that $L|_{T^{(2)}}$ is invertible then implies that $L|_T \colon Y \times_\formalU T \to Y \times_\formalU T$ is invertible as well.
\end{proof}

\subsubsection{More on minimal models}

\begin{lemma}\label{AutUtoAutQ}
  Let $S$ be a $(-1)$\nobreakdashes-symplectic derived scheme, $x \in S$ a geometric point and $\formalU_0$ be a minimal model of $S$ at $x$.
  Let $(Q,q)$ be a non-degenerate quadratic $\basefield$\nobreakdashes-vector space and denote by $\calQ$ the trivial quadratic bundle with connection on $\formalU$, with fiber $Q$.
  Let $\formalU \coloneqq \formalU_0 \oplus \calQ \coloneq \formalU_0 \smboxplus \widehat Q_0$ (where $\widehat Q_0$ is the formal neighborhood of $0$ in $Q$).

 For any automorphism $\phi$ of $\formalU$ (as a section of $\Darbstack_X^\lambda(S)$), the tangent automorphism
  \[
    \tangentmap{\phi,x} \colon \tangent_{\formalU,x} = \tangent_{\formalU_0,x} \oplus Q \to \tangent_{\formalU_0,x} \oplus Q = \tangent_{\formalU,x}
  \]
  decomposes as a matrix
  \[
    \tangentmap{\phi,x} =
    \begin{pmatrix}
      \Id_{\tangent_{\formalU_0,x}} & B \\ 0 & \tangentmap{\phi,x}^Q
    \end{pmatrix}
  \]
  where $\tangentmap{\phi,x}^Q \colon Q \to Q$ is an orthogonal automorphism.
  The assignment $\phi \mapsto \tangentmap{\phi,x}^Q$ is functorial and defines a map of group stacks $\theta \colon \stackAut_\Darb(\formalU) \to \delta^x_{\Aut(Q)}$, where $\delta^x_{\Aut(Q)}$ denotes the skyscraper sheaf with stalk $\Aut(Q) \simeq \orthogonal(\dim Q, \basefield)$ at $x$.
\end{lemma}
\begin{proof}
  Fix an automorphism $\phi$.
  It induces an isomorphism of exact sequences
  \[
    \begin{tikzcd}
      \cohomology^0 \tangent_{S,x} \ar{r} \ar[equals]{d} & \tangent_{\formalU,x} \ar{r}{\hessian_{f,x}} \ar{d}{\tangentmap{\phi,x}}[sloped,swap]{\sim} & \cotangent_{\formalU,x} \ar{d}{(\tangentmap{\phi,x}^\dual)^{-1}}[sloped,swap]{\sim} \\
      \cohomology^0 \tangent_{S,x} \ar{r} & \tangent_{\formalU,x} \ar{r}[swap]{\hessian_{f,x}} & \cotangent_{\formalU,x}.
    \end{tikzcd}
  \]
  Using the decomposition $\tangent_{\formalU,x} = \tangent_{\formalU_0,x} \oplus Q$ and the fact that $\formalU_0$ is minimal at $x$, we see that the morphism $\cohomology^0 \tangent_{S,x} \to \tangent_{\formalU,x}$ identifies with the inclusion of $\tangent_{\formalU_0,x}$.
  Write the above morphisms as block matrices:
  \[
    \tangentmap{\phi,x} = \begin{pmatrix} A & B \\ C & D \end{pmatrix} \hspace{1cm}\text{and} \hspace{1cm} \hessian_{f,x} = \begin{pmatrix} \hessian_{f_0,x} & 0 \\ 0 & \hessian_q \end{pmatrix} = \begin{pmatrix} 0 & 0 \\ 0 & \hessian_q \end{pmatrix},
  \]
  we find $A = \Id_{\tangent_{\formalU_0,x}}$ and $C = 0$. We can then compute:
  \begin{multline*}
    \begin{pmatrix} 0 & 0 \\ 0 & \hessian_q \end{pmatrix}
    = \hessian_{f,x}
    = \tangentmap{\phi,x}^\dual \hessian_{f,x} \tangentmap{\phi,x}
    = \begin{pmatrix} \Id & 0 \\ B^\dual & D^\dual \end{pmatrix} \begin{pmatrix} 0 & 0 \\ 0 & \hessian_q \end{pmatrix} \begin{pmatrix} \Id & B \\ 0 & D \end{pmatrix}
    \\
    = \begin{pmatrix} 0 & 0 \\ 0 & D^\dual \hessian_q D \end{pmatrix}.
  \end{multline*}
  It follows that $D$ is an orthogonal automorphism of $Q$.
  This defines the announced morphism of group stacks $\theta \colon \stackAut_\Darb(\formalU) \to \delta^x_{\Aut(Q)}$.
\end{proof}

\begin{lemma}\label{thetaA1sectionatx}
  In the situation of \cref{AutUtoAutQ}, for any integer $n$, the composition
  \[
    \begin{tikzcd}
      \stackMap(\bbA^n, \stackAut(\calQ)) \ar{r}{\formalU_0 \oplus -} & \stackMap(\bbA^n, \stackAut_\Darb(\formalU)) \ar{r}{\theta} & \stackMap(\bbA^n, \delta^x_{\Aut(Q)})
    \end{tikzcd}
  \]
  in $\St_\formalU$ induces an isomorphism between the stalks at $x$.
  In particular since colimits commute with stalks, the composition
  \[
    \begin{tikzcd}
      \PhiA{}\left(\stackAut(\calQ)\right)_{x} \ar{r}{\formalU_0 \oplus -} &  \PhiA{}\left(\stackAut_\Darb(\formalU)\right)_{x} \ar{r}{\theta} &  \PhiA{}\left(\delta^x_{\Aut(Q)}\right)_{x}
    \end{tikzcd}
  \]
 is an equivalence.
\end{lemma}
\begin{proof}
  The first half of the statement follows from the Riemann--Hilbert correspondence (see \cref{constantquadstack}), while the second is deduced from \cref{reminderexplicitmodelusualA1localization}.
\end{proof}

\subsubsection{End of the proof of \cref{prop:A1-automorphisms}}
From \cref{thetaA1sectionatx}, it suffices to prove that the morphism
\begin{equation}\label{thetalocallyA1eq}
  \PhiA{}\left(\stackAut_\Darb(\formalU)\right)_{x} \lra  \PhiA{}\left(\delta^x_{\Aut(Q)}\right)_{x}.
\end{equation}
induced by $\theta \colon \stackAut_\Darb(\formalU) \to \delta^x_{\Aut(Q)}$ is an equivalence.
Denote by $\calF$ the fiber of $\theta$, by $\stackEnd_{S/}^\sharp(\formalU)$ the fiber product
\[
\stackEnd_{S/}^\sharp(\formalU) \coloneqq \stackEnd_{S/}(\formalU) \times_{\vectorbundle_\formalU(\tangent_\formalU \otimes I_S)} \vectorbundle_\formalU(\tangent_\formalU \otimes \tangent_\formalU),
\]
and by $\stackAut^\sharp(\formalU)$ and $\calF^\sharp$ its intersection with $\stackAut_\Darb(\formalU)$ and $\calF$ respectively.

Consider the following diagram over $\formalU$, where each square is Cartesian\footnote{Note that the bottom pentagon is not a Cartesian square.} and where each down-right arrow of the shape $Z \to M$ represents the inclusion of the zero-locus of the morphism down to $\vectorbundle_\formalU(I_S^{(2)})$ restricted to $M$:
\[
  \begin{tikzcd}[column sep={3.9em,between origins},row sep={3.7em,between origins}]
    & \calF^\sharp \ar{dd} \ar{dl} \ar[dashed]{rr}
    && Z_{Y} \ar{rr}{L} \ar{dr} \ar{dd}
    && {Z(df \otimes df)} \ar{dr}
    &
  \\
    \calF \ar{dd}
    &&&& Y \ar{dd} \ar{rr}{L}
    && Y \ar{dd}
  \\[-2em]
    & \stackAut^\sharp(\formalU) \ar{dl} \ar{rr}{\sim}
    && Z_\tangent \ar{dl} \ar{dr}
  \\
   \stackAut_\Darb(\formalU) \ar{rr}{\sim} \ar{dr}
   && Z_I \ar{dr}
   && \vectorbundle_\formalU (\tangent_\formalU \otimes \tangent_\formalU) \ar{dl} \ar{rr}[swap]{L}
   && \vectorbundle_\formalU (\tangent_\formalU \otimes \tangent_\formalU) \ar{dd}[swap]{df \otimes df}
  \\
    & \stackEnd_{S/}(\formalU) \ar{rr}{\sim}[swap]{\exp_\nabla^S} \ar[drrrrr,to path={|- (\tikztotarget) [near end]\tikztonodes}, rounded corners, "\delta_\eta"]
    && \vectorbundle_\formalU (\tangent_\formalU \otimes I_S)
  \\[-1.4em]
    &&&&&& \vectorbundle_\formalU(I_S^{(2)})\rlap{.}
  \end{tikzcd}
\]
\begin{lemma}\label{technicallemmaaboutetaleneigh}
  Assume we have an étale morphism $z = (z_1, \dots, z_n) \colon \formalU_0 \to \affinespace^n$ (this always exists up to replacing $\formalU_0$ with an étale neighborhood of $x$) and set $\nabla_0$ to be the pullback of the trivial connection on $\affinespace^n$.

 The composite morphism $\Aut^\sharp(\formalU) \to \vectorbundle_\formalU(\tangent_\formalU \otimes \tangent_\formalU) \to x_* \vectorbundle_k(\tangent_{\formalU,x} \otimes \tangent_{\formalU,x})$ associates to an automorphism $\phi \colon \formalU \to \formalU$ the bivector corresponding to the morphism
  \[
    \begin{tikzcd}
      \Omega_{\formalU,x} \ar{r}{\tangentmap\phi^\dual - \Id} & \Omega_{\formalU,x} \simeq \Omega_{\formalU_0,x} \oplus Q \ar{r} & Q \ar{r}{\sim}[swap]{q} & Q^\dual \ar{r} & \tangent_{\formalU,x}.
    \end{tikzcd}
  \]
\end{lemma}

\begin{proof}
  Denote by $(z_{n+1}, \dots z_{n+p})$ coordinates of $Q = k^p$ such that $q = \sum_{i > n} z_i^2$.
  For $\phi \in \stackAut^\sharp(\formalU)$, the image in $\vectorbundle_\formalU(\tangent_\formalU \otimes I_S)$ classifies the morphism $\Omega_\formalU \to I_S$ given by $dz_i \mapsto \phi(z_i) - z_i$.
  Its restriction to the critical locus induces a morphism
  \[
    \Omega_\formalU \to i_*i^* I_S \simeq i_*i^* \Sym(\tangent_\formalU[1])[-1] \to i_* i^* \tangent_\formalU
  \]
  which naturally identifies with the projection $\Omega_\formalU \to I_S \to I_S / I_S^{(2)} = i_* \conormal_{S/\formalU} \simeq i_* i^* \tangent_\formalU$.
  It follows that:
  \begin{enumerate}
   \item The image of $\phi \in \stackAut^\sharp(\formalU)$ to a section of $i^*(\tangent_\formalU \otimes \tangent_\formalU)$ only depends on the image of $\phi$ in $\stackAut_\Darb(\formalU)$ and
   \item it moreover identifies with the morphism $\tangentmap\phi^\dual - \Id \colon \Omega_\formalU \to i_*\conormal_{S/\formalU} \simeq i_*i^* \tangent_\formalU$.
  \end{enumerate}
  Unwinding the equivalence $\conormal_{S/\formalU} \simeq i^* \tangent_\formalU$ and the shape of $\tangentmap\phi$ at $x$, we get the announced result.
\end{proof}

\Cref{technicallemmaaboutetaleneigh} together with the definitions of $\calF$ and $Y$ imply that the composition $\calF^\sharp \to \vectorbundle_\formalU(\tangent_\formalU \otimes \tangent_\formalU)$ factors through $Y$ (and thus through $Z_Y$), and that the obtained square
\[
  \begin{tikzcd}
    \calF^\sharp \ar{r} \ar{d} & Z_Y \ar{d} \\ \stackAut^\sharp(\formalU) \ar{r}{\sim} & Z_\tangent
  \end{tikzcd}
\]
is Cartesian.
In particular, using \cref{linearizetaylorexpansion}\ref{linearizetaylorexpansion-artin}, we get an equivalence (up to replacing $\formalU$ by an étale neighborhood of $x$)
\[
 \calF^\sharp \simeq Z_Y \simeq Z(df \otimes df) = \vectorbundle_\formalU(\fiber(df \otimes df \colon \tangent_\formalU \otimes \tangent_\formalU \to I_S^{(2)})).
\]
The right hand side being strongly $\affineline{}$\nobreakdashes-contractible (over $\formalU$), in particular we have
\[
  \PhiA{\formalU}(\calF^\sharp)\simeq \PhiA{\formalU}(\vectorbundle_\formalU(\fiber(df \otimes df \colon \tangent_\formalU \otimes \tangent_\formalU \to I_S^{(2)})))\simeq \PhiA{\formalU}(\formalU)\simeq \ast
\]
in the category of stacks over $\formalU$.

\begin{sublemma}\label{lemmaA1localizationsnaive}
  Let $\rmE \to \rmE'$ be a map of perfect complexes over $\formalU$ whose fiber is connective and let $\phi \colon \vectorbundle_\formalU(\rmE) \to \vectorbundle_\formalU(\rmE')$ be the induced map of derived vector bundles.
  Then any pullback of $\phi$ is an $\PhiA{\formalU}$\nobreakdashes-equivalence
  \begin{proof} We use the notations in \cref{localA1MV} and \cref{lemmaappliedtoourcase}. Fix a stack $Y$ over $\formalU$ and a morphism $Y \to \vectorbundle_\formalU(\rmE')$. Denote by $\rmN_\bullet$ the nerve of $\phi$. By assumption, the map $\ish{\formalU}(\phi)$ is an epimorphism of presheaves over $\formalU$, and therefore so is the pullback. We obtain
  \[
    \ish{\formalU}(Y)\simeq \colim \, \left(\ish{\formalU}Y \times_{\ish{\formalU}\vectorbundle_\formalU(\rmE')} \ish{\formalU}\rmN_\bullet\right) \in \PSh(\dAff_\formalU, \inftygpd).
  \]
  Applying the left adjoint $\LpshA{\formalU}$ we get
  \[
    \LpshA{\formalU} \ish{\formalU}(Y)
    \simeq \colim \, \LpshA{\formalU}\left(\ish{\formalU}Y \times_{\ish{\formalU}\vectorbundle_\formalU(\rmE')} \ish{\formalU}\rmN_\bullet\right) \in \PSh_{\affineline{\formalU}}(\dAff_\formalU, \inftygpd).
  \]
  Now, each boundary operator $\rmN_{n+1} \to \rmN_n$ is linear and admits a section, and is thus a strong $\affineline{\formalU}$\nobreakdashes-homotopy equivalence. In particular the induced morphism $Y \times_{\vectorbundle_\formalU(\rmE')} \rmN_{n+1} \to Y \times_{\vectorbundle_\formalU(\rmE')} \rmN_{n}$  is also a strong $\affineline{\formalU}$\nobreakdashes-homotopy equivalence in presheaves and therefore becomes an equivalence after applying  $\LpshA{\formalU}$. This colimit diagram is thus constant and we find
  \[
    \LpshA{\formalU}\ish{\formalU}(Y) \simeq \colim \, \LpshA{\formalU}\left(\ish{\formalU}Y \times_{\ish{\formalU}\vectorbundle_\formalU(\rmE')} \ish{\formalU}\rmN_\bullet\right) \simeq  \LpshA{\formalU}\ish{\formalU}\left(Y \times_{\vectorbundle_\formalU(\rmE')} \vectorbundle_{\formalU}(\rmE)\right).
  \]
  Applying $\Lsh{\formalU}\ipshA{\formalU}$, we get
  $
    \PhiA{\formalU}(Y) \simeq \PhiA{\formalU}\left(Y \times_{\vectorbundle_\formalU(\rmE')} \vectorbundle_{\formalU}(\rmE)\right).
  $
  \end{proof}
\end{sublemma}

\cref{lemmaA1localizationsnaive} implies the morphism $\calF^\sharp \to \calF$ is an $\PhiA{\formalU}$\nobreakdashes-equivalence.
We get that $\PhiA{\formalU}\left(\calF\right)_{x}$ is contractible.
Since the group morphism \cref{thetalocallyA1eq} admits a product-preserving section (\cf \cref{thetaA1sectionatx}), we find
\[
 \PhiA{\formalU}\left(\stackAut_\Darb(\formalU)\right)_{x} \simeq   \PhiA{\formalU}\left(\calF\right)_{x} \times   \PhiA{\formalU}\left(\delta^x_{\Aut(Q)}\right)_{x} \simeq   \PhiA{\formalU}\left(\delta^x_{\Aut(Q)}\right)_{x}.
\]
\cref{prop:A1-automorphisms} is thus proven, and so is \cref{theoremcontractibility} (a.k.a. \cref{contractibilitytheoremdetailedversion}). \qed

\section{Application: Behrend's function and vanishing cycles}
\label{sectionapplications}
We now discuss the applications of \cref{theoremcontractibility}.
We will construct 
examples of invariants that can be glued on $(-1)$\nobreakdashes-symplectic schemes.

\subsection{Critical invariants}\label{sectioncriticalinvariants}
Many singularity invariants are defined for algebraic LG pairs. We formalize here their properties, and show that in some cases, those invariants only depend on the formal neighborhood of the critical locus.
As a corollary, we will obtain that those invariants (\eg Milnor number, perverse sheaves or mixed Hodge modules of vanishing cycles) are also defined for locally algebraizable formal LG pairs -- see \cref{algebraicinvariants1cat} below.

\begin{notation}
  Denote by $\LGet_\basefield$ the site whose objects are \emph{smooth} affine algebraic LG-pairs $(U,f)$ over $\basefield$ such that $\Crit(f)_\red \subset f^{-1}(0)$ and whose morphisms are (potential-preserving) étale morphisms between them.
  Taking the critical locus yields a functor $\Crit \colon \LGet_\basefield \to \Affet_\basefield$ to the site of affine schemes with étale morphisms between them.
\end{notation}

\begin{definition}\label{criticalinvariant}
  Let $\calF$ be a stack\footnote{Recall that in this paper, all stacks satisfy hyperdescent.} $\left(\Affet_\basefield \right)^\op \to \inftygpd$.
  An \emph{$\calF$\nobreakdashes-valued algebraic critical invariant} $I$ is a morphism in $\Sh(\LGet_\basefield, \inftygpd)$:
  \[
  I \colon \constantsheaf *_{\LGet_\basefield} \to \calF \circ \Crit
  \]
  where $\constantsheaf *_{\LGet_\basefield}$ is the final constant stack.
  We denote by $\criticalinvariants_\calF$ the $\infty$\nobreakdashes-groupoid of $\calF$\nobreakdashes-valued algebraic critical invariants.
\end{definition}

\begin{observation}\label{remarkprecisedescriptionofcriticalinvariants}
  An $\calF$\nobreakdashes-valued algebraic critical invariant $I$ can equivalently be seen as a section
  \[
    \begin{tikzcd}
      &\int \calF\ar{d}{\pi_{\calF}}\\
      \LGet_\basefield\ar[dashed]{ur}{I}\ar{r}{\Crit}& \Affet_\basefield
    \end{tikzcd}
  \]
  where $\pi_{\calF}$ is the right fibration classifying $\calF$ and such that $I$ sends all morphisms in $\LGet_\basefield$ to $\pi_{\calF}$-Cartesian edges.
  Informally, the invariant $I$ maps a smooth affine algebraic LG-pair $(U, f)$ to a section $I_{U,f}$ of $\calF(\Crit(f))$, in a way compatible with étale maps.
  Note in particular that if $V \to U$ is an étale neighborhood of $\Crit(f)$ in $U$, then $I_{U,f}$ and $I_{V, f_{|V}}$ are equivalent.
\end{observation}

\begin{examples}{\label{examplesofalgebraiccriticalinvariants}Examples of algebraic critical invariants in the sense of \cref{criticalinvariant} include:}
  \item The Milnor number: Here $\calF$ is the sheaf of constructible $\integers$\nobreakdashes-valued functions.
  \item \label{vanishingcyclesarealgebraiccriticalinvariants} Vanishing cycles: To any smooth scheme $U$ equipped with a function $f \colon U \to \affineline{}$ with sole critical value $0$, we can associate a perverse sheaf on the étale site of $U$ (with coefficients in $A = \integers$ or $\rationals$) of \emph{vanishing cycles} $\PJoyce_{U,f} \coloneqq \phi_f(\constantsheaf{A}_U[\dim U])$.
  It is in fact supported on the critical locus $\Crit(f)$ and we can see it as an object in $\Perv(\Crit(f))$.
  Moreover, the perverse sheaf of vanishing cycles has a structure of mixed Hodge module, and may be seen as an object in $\MHM(\Crit(f))$.

  Let $\mathcal{F}=\stackPerv$ (respectively $\mathcal{F}=\MHMstack$) denote the stack $ (\Affet_\basefield)^\op \to \gpd$ classifying perverse sheaves (respectively graded polarized mixed Hodge modules \cite{MorihikoSaito911})\footnote{Since both these stacks are 1-truncated, by \cite[6.5.2.9]{lurie-htt} they satisfy hyperdescent.}.
  Following \cref{remarkprecisedescriptionofcriticalinvariants}, the assignment sending an algebraic $\LG$\nobreakdashes-pair $(U,f)$ to the perverse sheaf $\PJoyce_{U,f}\in \calF(\Crit(U,f))$ defines an algebraic critical invariant $I_{\PJoyce}\in \criticalinvariants_\calF $ in the sense of \cref{criticalinvariant}.
  We refer to \cite[Thm.\,2.11(iii)]{MR3353002} for the required functorialities.

  \item \label{exmotivicvanishingcycles}Motivic vanishing cycles: Here $\calF$ is the $\infty$\nobreakdashes-stack of Voevodsky motives (see \eg \cite{JosephAyoub916}).
  \item Matrix factorizations: Here $\calF$ is the $\infty$\nobreakdashes-stack of $2$\nobreakdashes-periodic dg-categories.
\end{examples}

The first two examples above will be studied in \cref{examplemilnornatural,naturaltransformationperversesection}.
The example of matrix factorizations will be addressed in \cite{hennionholsteinrobaloII}.

\begin{definitions}{}
  \item A locally algebraizable smooth formal $\LG$\nobreakdashes-pair is a pair $(\formalU, f)$ with $\formalU$ a smooth formal scheme and $f \colon \formalU \to \affineline{}$ such that
  \begin{enumerate}[nosep, itemsep=.3\baselineskip]
    \item $\Crit(f)$ is a scheme,
    \item $\Crit(f)_\red \subset f^{-1}(0)$, \ie $0$ is the only possible critical value and
    \item such that étale-locally on $\formalU$, there exists a smooth scheme $U$ with a function $g \colon U \to \affineline{}$, and an isomorphism between $\formalU$ and the formal neighborhood of a closed subscheme of $U$ identifying $f$ with the restriction of $g$.
  \end{enumerate}
  \item\label{entirelycritical} A locally algebraizable smooth formal $\LG$\nobreakdashes-pair $(\formalU, f)$ is \emph{entirely critical} if the closed immersion $\Crit(f) \hookrightarrow \formalU$ is a reduced equivalence.
  \item Let $\LGfet_\basefield$ be the category of entirely critical locally algebraizable smooth formal $\LG$\nobreakdashes-pairs, with étale adic and function-preserving morphisms between them. By construction, it comes with a critical locus functor that we denote by
  \[
    \Critf \colon \LGfet_\basefield \to \Affet_\basefield.
  \]
\end{definitions}

The following proposition justifies that for the purposes of this work we restrict attention to locally algebraizable $\LG$-pairs:

\begin{proposition}\label{darbouxchartsarelocallyalgebraizable}
  Let $X$ be $(-1)$\nobreakdashes-shifted exact symplectic Deligne-Mumford stack. Then for every $S\in X_\et^\daff$, every object $(\formalU,f)\in \Darbstack_X$ is étale locally algebraizable.
\end{proposition}
\begin{proof}
  By definition of the Darboux stack, we have $S\simeq \dCrit(\formalU,f)$ and $S_\red\simeq \formalU_\red$ so that $S$ and
  $\formalU$ have the same geometric points. Let $x\in\formalU$ be a geometric point.
  By the Darboux lemma of \cite{MR3904157}, there exists an étale neighborhood $S'$ of $x$ in $S$, a smooth affine scheme $V$, together with a function $g$ and a derived symplectomorphism $\dCrit(V,g) \simeq S'$.
  Moreover, by \cref{dCritonlydependsformal}, we can replace $V$ by the formal completion $\formalV \coloneqq \widehat{V}$. Now, by the transitivity of the action of quadratic forms of \cref{thm-transitivity}, we can relate the induced étale cover $\formalU'\to \formalU$ to $\formalV$, namely, there exists quadratic bundles $\formalQ_1$ and $\formalQ_2$ on $S'$ and equivalences of Darboux charts $\formalU'\smboxplus\formalQ_1 \simeq \formalV\smboxplus \formalQ_2$. Since $\formalV$ is algebraizable, so is $\formalV \smboxplus \formalQ_2$ and therefore, $\formalU'\smboxplus\formalQ_1$.
  Finally, by \cref{localformformorphisms1}, the inclusion along the zero section $\formalU' \hookrightarrow \formalU' \smboxplus \formalQ_1$ defines a morphism of Darboux charts and exhibits $\formalU'$ as algebraizable.
\end{proof}

\begin{lemmalist}{}
  \item The formal completion at the critical locus defines a functor $\formalizeLG \colon \LGet_\basefield \to \LGfet_\basefield$.
  \item Let $\psi \colon \formalU \to \formalV \to \affineline{}$ be an étale adic morphism of formal LG-pairs. If $\formalV \to \affineline{}$ is (locally) algebraizable, then $\formalU \to \affineline{}$ admits a (local) algebraization such that the morphism $\psi$ lifts to an algebraic morphism.
  \item Étale, adic and function-preserving covers generate a Grothendieck topology on $\LGfet_\basefield$.
  \item The functor $\formalizeLG$ is continuous as a map of sites.
\end{lemmalist}

\begin{proof}
  The functoriality of (a) is obvious.
  We focus on (b). Up to replacing $\formalV$ by an étale covering and $\formalU$ by the induced étale covering, we can and do assume that $\formalV$ is globally algebraizable.
  Let thus $(V, f)$ be an LG-pair with $\formalV \simeq \formalizeLG(V,f)$.
  By Artin's algebraization (see \cite[Thm.\,1.1]{JarodAlper839}), there exists a commutative diagram
  \[
    \begin{tikzcd}
      \formalU \ar{r}{\psi}[swap]{\et.} \ar[hook]{d} & \formalV \ar[hook]{d} \\ U \ar{r}{\et.} & V
    \end{tikzcd}
  \]
  which identifies $\formalU$ with $\formalizeLG(U, g|_U)$.

  To prove (c), we first observe that (b) implies that fiber products of locally algebraizable formal LG pairs along étale and adic morphisms are still locally algebraizable.
  Étale, adic and function-preserving covers are clearly stable under base change and composition, and we get indeed a topology.

 The functor $\formalizeLG$ maps étale morphisms to étale and adic morphisms, and preserves surjections and fiber products. It thus indeed defines a continuous morphism of sites and (d) holds.
\end{proof}

We can now extend \cref{criticalinvariant} to the formal setting:
\begin{definition}\label{definitionofformalcriticalinvariant}
  Let $\calF$ be a stack $\left(\Affet_\basefield\right)^\op \to \inftygpd$.
  An \emph{$\calF$\nobreakdashes-valued formal critical invariant} $I$ is a morphism of stacks  $I \colon *_{\LGfet_\basefield} \to \calF \circ \Critf$.
  We denote by $\criticalinvariantsf_\calF$ the $\infty$\nobreakdashes-groupoid of $\calF$\nobreakdashes-valued formal critical invariants.
  
  The formal completion functor $\formalizeLG \colon \LGet_\basefield \to \LGfet_\basefield$ induces an obvious forgetful functor $\criticalinvariantsf_\calF \to \criticalinvariants_\calF$.
\end{definition}

\begin{lemma}\label{LGoverAffetCartesian}
  The functor $\Critf \colon \LGfet_\basefield \to \Affet_\basefield$ is a Cartesian fibration in groupoids.
  It classifies a stack $\stackLG{}{} \colon (\Affet_\basefield)^\op \to \mathsf{Gpd}$ such that for $S \in \Affet_\basefield$, the $1$\nobreakdashes-groupoid $\stackLG{S}{} \coloneqq \stackLG{}{}(S)$ classifies triples $(\formalU, f, \alpha)$ where $(\formalU, f) \in \LGfet_\basefield$ and $\alpha$ is an isomorphism $\alpha \colon \Critf(f) \simeq S$ satisfying:
  \begin{enumerate}[nosep,itemsep=0.3\baselineskip]
   \item $\Critf(f) \to \formalU$ is a reduced equivalence and
   \item $\Critf(f)_\red \subset f^{-1}(0)$.
  \end{enumerate}
\end{lemma}
\begin{proof}
  As objects of $\LGfet_\basefield$ are assumed to be entirely critical (\cf \cref{entirelycritical}), the fact that $\Critf$ is a Cartesian fibration follows directly from \cref{equivalenceetalesites}.
  Moreover, any morphism in a fiber $\Critf^{-1}(S)$ is both formally étale and a reduced equivalence. It is thus invertible by \cref{etaleformalcompletionisnoop}.
  The descent condition is straightforward to check.
\end{proof}

\begin{observation}\label{explicitdescriptionofformalcriticalinvariant}
  As in \cref{remarkprecisedescriptionofcriticalinvariants}, the data of an $\calF$-valued formal critical invariant amounts to a section
  \[
  \begin{tikzcd}
    &\int \calF\ar{d}{\pi_{\calF}}\\
    \LGfet_\basefield \ar[dashed]{ur}{I} \ar{r}{\Critf} & \Affet_\basefield \rlap{.}
  \end{tikzcd}
  \]
  In particular, it amounts to a morphism of stacks $\stackLG{}{} \to \calF$.
\end{observation}

\begin{construction}\label{formalinvariantgivesdarbouxinvariant}
  Let $X$ be $(-1)$\nobreakdashes-shifted symplectic Deligne-Mumford stack and denote by $\lambda$ its canonical exact structure of \cref{canonicalexactstructure}.
  Recall the map of stacks on $X_\et^\daff$ of \cref{naturaltransformationfromLiouvilletoLG2} $\Darbstack^{\lambda}_X \to \stackLGder{X}{}$.
  By \cref{darbouxchartsarelocallyalgebraizable,remarkcriticalvalues}, we now know that this map factors through the full substack $\stackLGder{X}{\localg} \subseteq \stackLGder{X}{}$ spanned by locally algebraizable $\LG$-pairs with only $0$ as critical value.
  This latter stack comes with an obvious morphism $\stackLGder{X}{\localg} \to (\stackLG{}{})_{|{X_\et^\daff}}$ forgetting the derived structures involved in the LHS, where $(\stackLG{}{})_{|{X_\et^\daff}}$ denotes the restriction of $\stackLG{}{}$ along the truncation functor $\classicaltruncation \colon X_\et^\daff \to \Affet$.
  All in all, we get a natural morphism
  \[
    \Darbstack^\lambda_X \to (\stackLG{}{})_{|{X_\et^\daff}}.
  \]
  Let now $I$ be an $\calF$-valued formal critical invariant.
  Composing the above with (the restriction of) the stack morphism of \cref{explicitdescriptionofformalcriticalinvariant} yields a natural morphism
  \begin{equation}\label{mapfromdarbouxtopervofstacks}
    \Darbstack^\lambda_X \to (\stackLG{}{})_{|{X_\et}} \to^I \calF_{|X_\et^\daff}.
  \end{equation}
  In particular, any formal critical invariant induces an invariant out of the Darboux stack.
\end{construction}

From the above discussion, formal critical invariants will play a crucial role in using \cref{theoremcontractibility}.
However, most examples of critical invariants are a priori only defined algebraically.
We will thus be interested in lifting some invariants from $\criticalinvariants_\calF$ to $\criticalinvariantsf_\calF$, which is the content of the following section.

\subsection{From algebraic critical invariants to formal critical invariants}\label{sectionalgebraictoformulainvariants}
\subsubsection{The case of discrete critical invariants}

We focus in this subsection on the case of critical invariants with values in a set-valued sheaf.

\begin{lemma}\label{lem-injectivepi0}
  Let $(U,f)$ and $(V,g)$ be algebraic LG pairs and let $\phi \colon \formalizeLG(U,f) \simeq \formalizeLG(V,g)$ be an isomorphism of formal LG pairs.
  Then, there exist a smooth scheme $W$ and étale morphisms $U \from W \to V$ over $\affineline{}$ such that
  \begin{enumerate}
    \item the induced morphisms $\Crit(f) \from \Crit(f|_W) \to \Crit(g)$ are isomorphisms;
    \item the following triangle commutes
    \[
    \begin{tikzcd}[row sep=tiny]
      & \Crit(f|_W) \ar{dl} \ar{dr} & \\ \Crit(f) \ar{rr}[swap]{\phi} && \Crit(g).
    \end{tikzcd}
    \]
  \end{enumerate}
\end{lemma}

\begin{warning}\label{warningapproximation}
  In this statement, the étale morphisms $U \from W \to V$ induce isomorphisms $\formalizeLG(U,f) \simeq \formalizeLG(W, f|_W) \simeq \formalizeLG(V,g)$. However, the composition does not need to be $\phi$, and it only coincides with $\phi$ at the level of critical loci.
  One could improve on the above statement and require that this composition agrees with $\phi$ up to a fixed order $N \in \bbN$, but the scheme $W$ then depends on $N$.
\end{warning}

\begin{proof}
  This is a consequence of Artin's approximation, using arguments due to Artin \cite[Cor (2.6)]{MichaelArtin886}. We apply \cite[Prop.\,2.12]{JarodAlper839} with $N = 2$, $S = \affineline{}$, $\calX = U \times_{\affineline{}} V$, $\calZ = \formalizeLG(U,f)$, $\calZ_0 = \Crit(f)$ and $\eta \colon \calZ = \formalizeLG(U,f) \to U \times_{\affineline{}} V = \calX$ induced by $\phi$ and the canonical inclusions.
  We get $\xi \colon W \to U \times_{\affineline{}} V$ and a closed immersion $\Crit(f) \to W$ such that $\Crit(f) \to W \to U \times_{\affineline{}} V$ coincides with the restriction of $\eta$.
  Furthermore, the $2^\mathrm{nd}$ infinitesimal neighborhood $W_2$ of $\Crit(f)$ in $W$ is also identified with that of $\Crit(f)$ in $U$ and with that of $\Crit(g)$ in $V$. It follows that the composite map $W \to U \times_{\affineline{}} V \to U$ (\resp $V$) are indeed étale along $\Crit(f)$.
  Up to shrinking $W$, we may assume the maps $U \from W \to V$ are étale everywhere.
\end{proof}

\begin{lemma}\label{sheavesLG}
  Consider the commutative diagram
  \[
    \begin{tikzcd}[row sep=tiny]
      \LGet_\basefield \ar{rr}{\Crit} \ar{dr}[swap]{\formalizeLG} && \Affet_\basefield \\
      & \LGfet_\basefield \rlap{.} \ar{ur}[swap]{\Critf}
    \end{tikzcd}
  \]
  For each one of those sites $\calS$, we denote by $*_\calS$ the final constant sheaf (of sets).
  The canonical morphism $\Crit^{-1}(*_{\LGet_\basefield}) \to \Critf^{-1}(*_{\LGfet_\basefield}) \in \Sh(\Affet_\basefield, \Sets)$ is an isomorphism.
\end{lemma}

\begin{proof}
  Denote by $\calD$ the full subcategory of $\LGfet_\basefield$ consisting of globally algebraizable LG pairs, endowed with the induced Grothendieck topology.
  As $\calD$ topologically generates $\LGfet_\basefield$, those sites are equivalent and their categories of sheaves are equivalent (\cf \cite[Exp. III Thm 4.1]{MR0354652}).
  Since $\formalizeLG$ factors through $\calD$, we can replace $\LGfet_\basefield$ by $\calD$.

  It now suffices to prove that the analogous statement for presheaf pullbacks and study the morphism $\tau \colon \Crit^\dagger(*_{\LGet_\basefield}) \to \Critf^\dagger(*_{\calD})$.
  We fix $S \in \Affet_\basefield$ and compute the $S$\nobreakdashes-sections of those presheaves.
  Denote by $\calC \coloneqq \LGet_\basefield$ and $\calE \coloneqq \Affet_\basefield$. Let $\calC_S \coloneqq \calC \times_\calE \undercat{S}{\calE}$ and $\calD_S \coloneqq \calD \times_\calE \undercat{S}{\calE}$. We get
  \[
    \Crit^\dagger(*_{\calC})(S) = \colim_{\calC_S} * \in \Sets \hspace{2em} \text{and} \hspace{2em} \Critf^\dagger(*_{\calD})(S) = \colim_{\calD_S} * \in \Sets.
  \]
  The canonical morphism $\tau$ is thus simply the functor $\formalizeLG_S \colon \calC_S \to \calD_S$ on connected components.
  By definition, every formal LG pair of $\calD$ is algebraizable, so that $\formalizeLG_S$ is essentially surjective. Moreover, \cref{lem-injectivepi0} implies that $\formalizeLG_S$ is injective at the level of connected components. We get
  \[
    \Crit^\dagger(*_{\calC})(S) \simeq \pi_0(\calC_S) \simeq \pi_0(\calD_S) \simeq \Critf^\dagger(*_{\calD})(S).\qedhere
  \]
\end{proof}

\begin{warning}\label{warningonlysheaves}
  The statement analogous to \cref{sheavesLG} for stacks (either $1$\nobreakdashes- or higher) does not hold. Attempting to prove such a statement leads to studying the higher homotopy groups of the categories $\calC_S$ and $\calD_S$ above. It is true that $\pi_1(\calC_S)$ injects into $\pi_1(\calD_S)$ (at any base point), but the morphism is not surjective in general: some automorphisms of formal LG pairs are transcendental in nature and cannot be algebraized. This is also related to \cref{warningapproximation}.
  We will circumvent this issue in the next subsection.
\end{warning}

\begin{corollary}\label{formalcriticalinvariant}
  Let $\calF$ be a (set-valued) sheaf on $\Affet_\basefield$. Any $\calF$\nobreakdashes-valued critical invariant $I$ induces a natural formal invariant
  \[
    I^\mathsf{f} \colon \ast_{\LGfet_\basefield} \to \calF \circ \Critf.
  \]
  Unpacking the definitions, for any locally algebraizable formal LG pair $(\formalU, f)$, we get a section $I^\mathsf{f}_{\formalU, f} \in \calF(\Critf(\formalU, f))$, in a way compatible with restrictions along étale morphisms.
\end{corollary}
\begin{proof}
  The critical invariant $I \colon \ast_{\LGet_\basefield} \to \calF \circ \Crit$ corresponds, by adjunction, to a morphism $\Crit^{-1}(\ast_{\LGet_\basefield}) \to \calF$.
  Using \cref{sheavesLG}, we get $\Critf^{-1}(\ast_{\LGfet_\basefield}) \simeq \Crit^{-1}(\ast_{\LGet_\basefield}) \to \calF$, which in turn amounts to the announced formal critical invariant be adjunction.
\end{proof}

\begin{definition}
  Let $\calF$ be a (set-valued) sheaf on $\Affet_\basefield$ and $I$ an $\calF$\nobreakdashes-valued critical invariant.
  We say that $I$ is \emph{stable} if for any $(U, f) \in \LGet_\basefield$ and any quadratic bundle $(\pi \colon Q \to U,q)$ over $U$, we have
  \[
    I_{Q,q + f \pi} = I_{U,f} \in \calF(\Crit(f)).
  \]
\end{definition}

\begin{corollary}
  Let $X$ be a $(-1)$\nobreakdashes-shifted symplectic derived Deligne--Mumford stack and let $\calF$ be a sheaf (of sets) over $\Affet_\basefield$.
  Any $\calF$\nobreakdashes-valued critical invariant $I$ induces a morphism of stacks over the small étale site of $X$
  \[
    \begin{tikzcd}
      \Darbstack_X^\lambda \ar{r}{I} & \calF_{|{X_\et^\daff}}
    \end{tikzcd}
  \]
  which coincides with $I$ on any (local) algebraization of a Darboux chart, where $\calF_{|{X_\et^\daff}}$ is the restriction of $\calF$ to the small étale site of (the truncation of) $X$.

 If moreover $I$ is stable, then it induces a global section $I_X \in \calF(X)$ uniquely determined by the property that its restriction to an algebraic critical chart $\dCrit(f)$ is $I_{U,f}$, where $f \colon U \to \affineline{}$.
\end{corollary}

\begin{proof}
  This first half of the statement follows directly from \cref{formalcriticalinvariant,formalinvariantgivesdarbouxinvariant}. If furthermore $I$ is stable, then the map $I \colon \Darbstack_X^\lambda \to \calF_{|{X_\et^\daff}}$ factors canonically as
  \[
    \Darbstack_X^\lambda \to \truncationleq{0}\Darbstack_X^\lambda \to \quot{\truncationleq{0}\Darbstack_X^\lambda}{\truncationleq{0}\stackQuadnabla_X} \simeq *_X \to \calF_{|{X_\et^\daff}}.\qedhere
  \]
\end{proof}

\begin{example}[Behrend's function]\label{examplemilnornatural}
  Consider $\calF$ the sheaf classifying $\integers$\nobreakdashes-valued constructible functions.
  For any LG pair $(U,f)$, denote by $\nu_{U,f}$ the function $\nu_{U,f} = (-1)^{\dim U - 1} \chi(\phi_f)$ on $\Crit(f)$.
  It defines a stable $\myuline \integers$\nobreakdashes-valued critical invariants, and therefore a function $\nu_X \in \cohomology_\et^0(X, \integers)$ for any $(-1)$\nobreakdashes-symplectic scheme $X$.
  This function coincides with Behrend's function \cite{MR2600874} as the functions agrees locally.
\end{example}

\subsubsection{Isotopy-invariance.}
Similarly to \cite{MR3353002}, we will also be interested in gluing sheaves of vanishing cycles or mixed Hodge modules.

We will thus need a categorified version of \cref{sheavesLG}.
As explained in \cref{warningonlysheaves}, some extra care is required as some formal automorphisms of LG pairs may not be algebraic.
They are however isotopic to algebraic automorphisms by \cref{theoremcontractibility} (see also \cite[Prop.\,3.4]{MR3353002}), and working up to $\affineline{}$\nobreakdashes-isotopy is enough.

\begin{definitions}{\label{definitionofalternativeA1localization}%
  Fix an affine scheme $T$.}
  \item Denote by $\Affet_T$ the category of affine schemes over $T$, with étale morphisms between them.
  \item For any scheme $U$ (\resp formal scheme $\formalU$) smooth over $T$ equipped with a function $f \colon U \to \affineline{\basefield}$ (\resp $f \colon \formalU \to \affineline{}$), we denote by $\Critrelative{T}(f)$ (\resp $\Critfrelative{T}(f)$) the relative critical locus of $f$:
  \[
    \Critrelative{T}(f) \coloneqq \{ df = 0 \in \Omega_{U/T} \}.
  \]
  \item We denote by $\Affet_{\affineline{T}}$ the category of affine schemes $U$ over $\affineline{T}$ smooth over $T$, with étale morphisms (over $\affineline{T}$) between them.
  Consider the pullback
  \[
    \begin{tikzcd}[column sep=large,row sep=scriptsize]
      \calE_T \ar{r} \ar{d} \tikzcart & \Affet_{\affineline{T}} \ar{d}{\Critrelative{T}} \\
      \Affet_\basefield \ar{r}[swap]{U \mapsto U \times T} & \Affet_T\rlap{.}
    \end{tikzcd}
  \]
  We denote by $\LGetcst_{T} \subseteq \calE_T$ the essential image of $\LGet_\basefield \subset \Affet_{\affineline{\basefield}} \simeq \calE_\basefield$ under the pullback functor $\calE_\basefield \to \calE_T$.
  \item We define similarly $\LGfetcst_{T}$ as the essential image of $\LGfet_\basefield$ in the analogous pullback $\calE'_T$:
  \[
    \begin{tikzcd}[column sep=large,row sep=scriptsize]
      \calE_T' \ar{r} \ar{d} \tikzcart & \FSch^\et_{\affineline{T}} \ar{d}{\Critfrelative{T}} \\
      \Affet_\basefield \ar{r}[swap]{U \mapsto U \times T} & \FSch^\et_T\rlap{.}
    \end{tikzcd}
  \]
\end{definitions}

\begin{construction}\label{constructionA1naivelocalizationofLG}
  Recall the algebraic simplices $[n] \mapsto \simplex^n \coloneqq \simplex^n_\basefield$ (\cf \cref{reminderexplicitmodelusualA1localization}). They induce simplicial diagrams of $1$-categories $\Delta^\op \to \CAT \subset \inftycats$
  \[
    [n] \mapsto \LGetcst_{\simplex^n}
    \hspace{.7em}\text{and}\hspace{.8em}
    [n] \mapsto \LGfetcst_{\simplex^n}.
  \]
  We denote by $\LGetA_\basefield$ and $\LGfetA_\basefield$ their respective colimits in $\inftycats$:
  \[
    \LGetA_\basefield \coloneqq \colim_{[n] \in \Delta} \LGetcst_{\simplex^n},  \hspace{0.5cm}  \LGfetA_\basefield \coloneqq \colim_{[n] \in \Delta} \LGfetcst_{\simplex^n}.
  \]
  By construction, they come with natural functors
  \[
    \begin{tikzcd}
      \LGet_\basefield \ar{r}{\locA} \ar{dr}[swap]{\formalizeLG} & \LGetA_\basefield \ar{rr}{\CritA} \ar{dr}[near end]{\formalizeLG^{\affineline{}}} && \Affet_\basefield \\
      & \LGfet_\basefield \ar{r}{\locA} & \LGfetA_\basefield\rlap{.} \ar{ur}[swap]{\CritfA}
    \end{tikzcd}
  \]
\end{construction}

\begin{observation}\label{mapsareisotopyclasses}
  Since the simplicial diagram $\LGetcst_{\simplex^\bullet}$ is constant on objects, by \cite[Lem.\,1.3.4.10]{lurie-ha}, $\LGetA_\basefield$ can informally be described as the \icategory whose object are those of $\LGet_\basefield$ and arrows are $\affineline{}$\nobreakdashes-isotopy classes of morphisms in $\LGet_\basefield$, where the isotopies are assumed to be trivial on the critical loci.
  Its higher categorical structure is however difficult to describe.
\end{observation}

\begin{definition}\label{definitionofisotopyinvariantcriticalinvariant}
  Let $\calF$ be a (higher) stack over $\Affet_\basefield$. An \emph{isotopy-invariant} $\calF$\nobreakdashes-valued algebraic critical invariant is a natural transformation
  \[
    I \colon *_{\LGetA_\basefield} \to \calF \circ \CritA.
  \]
  We denote by $\criticalinvariantsA_\calF$ their $\infty$\nobreakdashes-groupoid. Similarly, a formal isotopy-invariant critical invariant is a natural transformation $*_{\LGfetA_\basefield} \to \calF \circ \CritfA$ and $\criticalinvariants^{\mathsf{f},\affineline{}}_\calF$ denotes their $\infty$\nobreakdashes-groupoid.
\end{definition}

\begin{observation}\label{compositionsunderlyingcriticalinvariants}
  Composition with the functors $\locA$ of \cref{definitionofalternativeA1localization}
  \[
    \theta \colon \criticalinvariantsA_\calF \to \criticalinvariants_\calF \hspace{2cm} \theta_\mathsf{f}\colon \criticalinvariants^{\mathsf{f},\affineline{}}_\calF \to \criticalinvariants^{\mathsf{f}}_\calF
  \]
  yields underlying $\calF$\nobreakdashes-valued algebraic (\resp formal) critical invariant in the sense of \cref{criticalinvariant}.
\end{observation}

\begin{lemma}\label{LGoverAffetCartesianT}
  For any $T \in \Aff_\basefield$, the functor $p \colon \LGetcst_{T} \to \Affet_\basefield$ induced by $\Critrelative{T}$ is a Cartesian fibration in $1$-groupoids.
  The associated functor satisfies étale hyperdescent.
\end{lemma}
\begin{proof}
  This is similar to \cref{LGoverAffetCartesian}.
  Any object of $\LGetcst_{T}$ is of the form $\formalU \times T$ with a function $f \colon \formalU \to \affineline{}$, such that the closed immersion $\Critf(f) \times T \simeq \Critfrelative{T}(f_T) \to \formalU \times T$ is a reduced equivalence, where $f_T = f \times \Id_T \colon \formalU \times T \to \affineline{T}$.
  Moreover, a morphism $(\formalU \times T, f) \to (\formalV \times T, g)$ of $\LGetcst_{T}$ is an étale map $\phi \colon \formalU \times T \to \formalV \times T$ over $\affineline{T}$ such that the induced morphism
  \[
    \Critf(f) \times T \simeq \Critfrelative{T}(f_T) \to \Critfrelative{T}(g_T) \simeq \Critf(g) \times T
  \]
  is of the form $p(\phi) \times \Id_T$.
  It then follows from \cref{equivalenceetalesites} that every such $\phi$ is $p$-Cartesian and that $p$ is a Cartesian fibration in $1$-groupoids.
  The descent condition is straightforwardly verified.
\end{proof}
\begin{notation}
  For any $[n] \in \Delta$, we denote by $\stackLG{}{n}$ the stack described by the Cartesian fibration of \cref{LGoverAffetCartesianT} applied to $T = \simplex^n$.
\end{notation}

\begin{observation}\label{explicitdescriptionofdataofisotopycriticalinvariantfromcolimit}
  As in \cref{explicitdescriptionofformalcriticalinvariant}, the data of an \emph{isotopy-invariant} $\calF$-valued formal critical invariants correspond to the data of section
  \[
    \begin{tikzcd}
      & \int \calF\ar{d}{\pi_{\calF}} \\
      \LGfetA_\basefield\ar[dashed]{ur}{I}\ar{r}{\CritfA} & \Affet_\basefield\rlap{.}
    \end{tikzcd}
  \]
  But since $\LGfetA_\basefield$ is defined as a colimit in $\inftycats$, the data of $I$ amounts to a compatible system of functors $I_n$
  \[
    \begin{tikzcd}
      \mathllap{\int \stackLG{}{n} = {}} \LGfetcst_{\simplex^n} \ar[dashed]{r}{I_n} & \int \calF
    \end{tikzcd}
  \]
  over $\Affet_\basefield$.
  It thus amounts to a compatible system of morphisms of stacks $I_n \colon \stackLG{}{n} \to \calF$, and therefore to a morphism
  \[
    \begin{tikzcd}
      \colim \stackLG{}{\bullet} \ar{r}{I} & \calF \in \Sh(\Affet_\basefield, \inftygpd).
    \end{tikzcd}
  \]
\end{observation}

This observation allows us to relate $\affineline{}$-invariant critical invariants to the $\affineline{}$-isotopic localization of \cref{subsectionisotopies}.

\begin{lemma}\label{keyrelationbetweenLGstacks}
  Let $X$ be $(-1)$\nobreakdashes-shifted symplectic Deligne-Mumford stack and denote by $\lambda$ its canonical exact structure of \cref{canonicalexactstructure}.
  Every isotopy invariant formal critical invariant $I$ induces a map
  $\DA  \to \calF_{|{X_\et^\daff}}$ in $ \Sh(X_\et^\daff,\inftygpd)$.
\end{lemma}
\begin{proof} 
	\newcommand{\lgA}{\stackLGder{X}{\localg,\affineline{}}} 
  Recall the stack $\stackLGder{X}{\localg}$ on $X_\et^\daff$ of \cref{formalinvariantgivesdarbouxinvariant}, defined as a full substack of the restriction of the stack $\calG \coloneqq \stackfactor^{\LG}_{p_X}$ over $X_\derham$ to $X_\et^\daff$.
  Setting $\homotopysheaf \coloneqq \homotopysheaf_0(\stackLGder{X}{\localg}) \subset \homotopysheaf_0(\calG_{|X_\et^\daff})$, we introduce a new stack on $X_\et^\daff$:
  \[
    \lgA \coloneqq \calG^{\affineline{}}_{\homotopysheaf}
  \]
  as in \cref{A1quotientgeneral}. Alternatively, using the explicit formula for the functor
  \[
    \PhiSA{X_\derham}(\calG) \simeq \Lsh{X_\derham}\left(\colim \stackMap^{\constant}_{X_\derham}\left(\simplex^\bullet \times X_\derham,\ish{X_\derham}(\calG)\right)\right)
  \]
  obtained combining \cref{formulaforendofunctorlocalizationSigmaA1beforeiteration}, with \cref{singularcomplexconstant,keyformulaforA1localizationsrelation}, the stack $\lgA$ can also be described as a colimit of the simplicial diagram: indeed setting (using the notation from \cref{constructionofquotientA1stacksgeneral})
  \[
    \stackLGder{X}{\localg,n} \coloneqq \stackMap^{\constant}_{X_\derham}\left(\simplex^n \times X_\derham,\ish{X_\derham}(\calG)\right)_{\homotopysheaf}
  \]
  for each $[n] \in \Delta^\op$, we obtain
  \[
     \lgA \simeq \colim_{[n] \in \Delta^\op} \stackLGder{X}{\localg,n} \in \Sh(X_\et^\daff, \inftygpd).
  \]
  We further observe the existence of a natural forgetful map of simplicial objects in stacks on $X_\et^\daff$
  \[
    \stackLGder{X}{\localg,\bullet} \to (\stackLG{}{\bullet})_{|{X_\et^\daff}},
  \]
  which allows using \cref{explicitdescriptionofdataofisotopycriticalinvariantfromcolimit} to induce, from every isotopy invariant formal critical invariant $I$, a morphism in $\Sh(X_\et^\daff,\inftygpd)$
  \[
    \begin{tikzcd}
      \lgA \simeq \colim \stackLGder{X}{\localg,\bullet} \ar{r} & \colim (\stackLG{}{\bullet})_{|{X_\et^\daff}} \ar{r}{I} & \calF_{|{X_\et^\daff}}.
    \end{tikzcd}
  \]
  Finally, by the functoriality of $\PhiSA{}$ and by \cref{A1quotientgeneral} and \cref{Liouvilleuptoisotopy}, the map \cref{naturaltransformationfromLiouvilletoLG} inducing \cref{naturaltransformationfromLiouvilletoLG2} and \cref{mapfromdarbouxtopervofstacks} descends to a morphism
  \[
    \begin{tikzcd}
      \DA \ar{r} & \lgA \ar{r} & \calF_{|{X_\et^\daff}} \in \Sh(X_\et^\daff,\inftygpd).
    \end{tikzcd}\qedhere
  \]
\end{proof}

\subsubsection{The case of critical invariants with values in 1-groupoids.}
For the rest of this section, we will restrict ourselves to $\calF$\nobreakdashes-valued invariants for $\calF$ a $1$\nobreakdashes-stack.

\begin{lemma}\label{a1invariantimage}
  Let $\calF$ be a $1$\nobreakdashes-stack over $\Affet_\basefield$.
  The forgetful functor $\theta \colon \criticalinvariantsA_\calF \to \criticalinvariants_\calF$ is fully faithful.
  Its essential image consists of critical invariants $I$ with the following additional property:
  \begin{itemize}[leftmargin=2em]
    \item[$(*)$] For any $\LG$ pairs $(U,f)$ and $(V,g)$ and any $\affineline{}$\nobreakdashes-family of function-preserving étale morphisms $(\phi_t \colon U \to V)_{t \in \affineline{}}$, if $\phi_t|_{\Crit(U,f)} \colon \Crit(U,f) \to \Crit(V,g)$ is independent of $t \in \affineline{}$, then the induced isomorphism
          \[
            I(\phi_t) \colon I_{U,f} \simeq \left(\phi_t|_{\Crit(U,f)}\right)^* I_{V,g}
          \]
          is independent of $t$.
  \end{itemize}
\end{lemma}

\begin{proof}
  As $\calF$ is a $1$\nobreakdashes-stack, the notion of isotopy-invariant critical invariants can be expressed in terms of the $1$\nobreakdashes-truncation of the \icategory $\LGetA_\basefield$. At the level of $1$\nobreakdashes-truncations, the functor $\locA$ identifies with the functor being the identity on objects and sending a morphism to its isotopy class.
  The result follows.
\end{proof}

\begin{proposition}\label{algebraicinvariants1cat}
  Let $\calF$ be a $1$\nobreakdashes-stack. Composition with $\formalizeLG^{\affineline{}}$ induces an equivalence (of $1$\nobreakdashes-groupoids) $\criticalinvariants^{\mathsf{f},\affineline{}}_\calF \to^\sim \criticalinvariantsA_\calF$.
\end{proposition}

\begin{lemma}\label{etaleequalizerisotopy}
  Consider a diagram of smooth affine schemes
  \[
    \begin{tikzcd}
      U \ar[shift right=2pt]{r}[swap]{a} \ar[shift left=2pt]{r}{b} & V \ar{r}{g} & \affineline{}
    \end{tikzcd}
  \]
  where both $a$ and $b$ are étale and $g \circ a = g \circ b$. Set $f \coloneqq g \circ a = g \circ b$ and assume further that the two morphisms $\Crit(f) \to V$ induced by $a$ and $b$ coincide.

  Assume there exists an isotopy $\varphi \colon \formalizeLG(U, f) \times \affineline{} \to V$ such that $\varphi_0 = a|_{\formalizeLG(U, f)}$, $\varphi_1 = b|_{\formalizeLG(U, f)}$ and such that $\varphi_t$ is étale, function-preserving and the identity on $\Crit(f)$ for every $t$,
  then there exists an étale neighborhood $W$ of $\Crit(f)$ in $U$ and an isotopy $\psi \colon W \times \affineline{} \to V$ such that $\psi_0 = a|_W$, $\psi_1 = b|_W$ and $\psi_t$ is étale, function-preserving and the identity on $\Crit(f)$ for every $t$.
\end{lemma}

\begin{proof}
  Denote by $E$ the stack of function-preserving étale isotopies between $a$ and $b$ fixing $\Crit(f)$:
  \[
    \begin{tikzcd}[column sep={between origins, 5.5em}]
      E \ar{rr} \tikzcart \ar{d} && \Map_{\Crit(f) \times \affineline{}/-/U}^\et\left(U \times \affineline{}, V \times_{\affineline{}} U\right) \ar{d}{\ev_0,\ev_1} \\
      U \ar{rr}{\{(a,b)\}} & {} & \Map_{\Crit(f)/-/U}^\et\left(U, V \times_{\affineline{}} U\right)^2\rlap{,}
    \end{tikzcd}
  \]
  where the superscript $(-)^\et$ signifies we restrict to the substacks whose points have underlying morphisms $U \to V$ that are étale.
  It is locally of finite presentation and we can apply \cite[Prop.\,2.12]{JarodAlper839} to $N = 2$, $S = \affineline{}$, $\calX = E$, $\calZ = \formalizeLG(U,f)$, $\calZ_0 = \Crit(f)$ and $\eta \colon \calZ = \formalizeLG(U,f) \to E = \calX$ the morphism classifying the isotopy $\varphi$.

  We get $W \to E$ and a factorization of $\calZ_2 \to W$ of the restriction of $\eta$ on the second infinitesimal neighborhood of $\calZ_0$ in $\calZ$.
  This implies that the composite morphism $W \to E \to U$ is étale on $\calZ_0 = \Crit(f)$.
  Up to shrinking $W$, we can assume $W \to U$ is étale everywhere. We get the announced isotopy.
\end{proof}

\begin{proof}[Proof of \cref{algebraicinvariants1cat}]
  We follow the same first steps as for \cref{sheavesLG}.
  Writing $\calD^{\affineline{}}$ for the full subcategory of $\LGfetA_\basefield$ spanned by algebraizable formal LG pairs and $\calC^{\affineline{}} \coloneqq \LGetA_\basefield$, we are reduced to studying, for any $S \in \Affet_\basefield$, the functor
  \[
    \begin{tikzcd}
      \formalizeLG^{\affineline{}}_S \colon \calC^{\affineline{}}_S \coloneqq \calC^{\affineline{}} \times_{\Affet} \undercat{S}{\Affet} \ar{r} &
      \calD^{\affineline{}} \times_{\Affet} \undercat{S}{\Affet} \eqqcolon \calD^{\affineline{}}_S.
    \end{tikzcd}
  \]
  The proof of \cref{sheavesLG} implies that this functor induces a bijection on the sets of isomorphism classes. Because $\calF$ is $1$\nobreakdashes-truncated, it suffices to prove the above functor $\formalizeLG^{\affineline{}}_S$ induces an isomorphism on $\pi_1$ (where $\pi_1$ of a category is the fundamental group of its geometric realization).

 Let us borrow the notations $\calC$, $\calC_S$, $\calD$ and $\calD_S$ from the proof of \cref{sheavesLG}.
  We fix an object $\tau$ in $\calC_S^{\affineline{}}$ and focus on $\pi_1(\calC_S^{\affineline{}}, \tau) \to \pi_1(\calD_S^{\affineline{}}, \formalizeLG_S^{\affineline{}}(\tau))$.
  The object $\tau$ lifts to an object $(S \to V \to^g \affineline{}) \in \calC_S$ which we also denote by $\tau$.
  By \cref{mapsareisotopyclasses}, the induced morphism $\pi \colon \pi_1(\calC_S, \tau) \to \pi_1(\calC_S^{\affineline{}}, \tau)$ is surjective (any zigzag can be lifted) and we have a commutative diagram
  \[
    \begin{tikzcd}
      \pi_1(\calC_S, \tau) \ar[two heads]{d}{\pi} \ar{r} & \pi_1(\calD_S, \formalizeLG_S(\tau)) \ar[two heads]{d}
      \\
      \pi_1(\calC_S^{\affineline{}}, \tau) \ar{r} & \pi_1(\calD_S^{\affineline{}}, \formalizeLG_S^{\affineline{}}(\tau)).
    \end{tikzcd}
  \]
  A class in $\pi_1(\calC_S, \tau)$ (\resp in $\pi_1(\calD_S, \formalizeLG_S(\tau))$) can be represented as a zig-zag of étale maps over $\affineline{}$ and under $S$, starting and ending at $V$ (\resp at $\formalizeLG(V,g)$).
  Since étale pullbacks exist in $\calC$ (\resp in $\calD$), any class in $\pi_1(\calC_S, \tau)$ (\resp in $\pi_1(\calD_S, \formalizeLG_S(\tau))$) can be represented as a single span
  \[
    \begin{tikzcd}[column sep = small]
      U \ar[shift right=2pt]{r}[swap]{a} \ar[shift left=2pt]{r}{b} & V
    \end{tikzcd}
    \hspace{2em} \left( \text{\resp }
    \begin{tikzcd}[column sep = small]
      \formalU \ar[shift right=2pt]{r}[swap]{a} \ar[shift left=2pt]{r}{b} & \formalizeLG(V,g)
    \end{tikzcd}\right)
  \]
  with $g \circ a = g \circ b$ and commuting under $S$.
  Such a loop is $\affineline{}$\nobreakdashes-contractible (so give rise to the trivial class in $\pi_1(\calC_S^{\affineline{}}, \tau)$ -- \resp in $\pi_1(\calD_S^{\affineline{}}, \formalizeLG_S(\tau))$) if and only if $a$ and $b$ are isotopic on some étale neighborhood of the image of $S$ in $U$ (resp.\
  $\formalU$).
  \cref{etaleequalizerisotopy} then implies that $\pi_1(\calC_S^{\affineline{}}, \tau) \to \pi_1(\calD_S^{\affineline{}}, \formalizeLG_S(\tau))$ is injective.

  It is also locally surjective by \cite[Prop.\,3.4]{MR3353002} (every formal automorphism is locally isotopic to an algebraic one).
  The result follows by adjunction as in \cref{sheavesLG} and \cref{formalcriticalinvariant}.
\end{proof}

\subsection{Vanishing cycles}\label{sectionrecoverJoyce}
We now apply the results of \cref{sectionalgebraictoformulainvariants} to the case of vanishing cycles (\cref{vanishingcyclesarealgebraiccriticalinvariants})

\begin{proposition}\label{MHMisformal}
   Let $\mathcal{F}=\stackPerv$ (respectively $\mathcal{F}=\MHMstack$). Then, the algebraic critical invariant $I_{\PJoyce}\in \criticalinvariants_\calF$ of \cref{vanishingcyclesarealgebraiccriticalinvariants} is in the full subcategory $\criticalinvariantsA_{\mathcal{F}}\subseteq \criticalinvariants_\calF$ (cf. \cref{definitionofisotopyinvariantcriticalinvariant}).
\end{proposition}
\begin{proof}
  Vanishing cycles define algebraic critical invariants (cf. \cref{vanishingcyclesarealgebraiccriticalinvariants}).  To show they are isotopy invariant we use \cref{a1invariantimage}. We first note that the (shifted) pull-back $\Perv(X) \to \Perv(X \times \affineline{})$ is fully faithful (see \eg \cite{PortaHaine}), thus $\affineline{}$\nobreakdashes-families of endomorphisms of perverse sheaves are constant.
  As the forgetful functor $\MHM \to \Perv$ is faithful the result also follows for $\MHMstack$. 
\end{proof}

\begin{remarks}{}
  \item\label{remvanishingformal} The fact that vanishing cycles are a formal invariant has been proven for perverse sheaves by Berkovich in \cite{MR1262943, Berkovich906}.
  \item In particular, \cref{MHMisformal} implies that vanishing cycles in mixed Hodge modules are canonically defined for any function on a locally algebraizable smooth formal scheme.
  \item
 \cref{MHMisformal}  may also be recovered from \cite[Thm.\,4.2 and \S2.10]{MR3353002}. In fact, there it is shown that vanishing cycles only depend on the third order thickening of the critical locus, using a third order version of Artin's approximation in the special case of $\LG$\nobreakdashes-pairs.
\end{remarks}

\begin{construction}\label{constructionofperversesheavesasmapfromdarbouxstack}
Using \cref{MHMisformal}, we can lift the invariant of vanishing cycles along the equivalence of \cref{algebraicinvariants1cat}
\[
  \criticalinvariantsA_{\stackPerv} \simeq \criticalinvariants^{\mathsf{f},\affineline{}}_{\stackPerv}
\]
and via \cref{keyrelationbetweenLGstacks} obtain a morphism of stacks on the small étale site of a $(-1)$\nobreakdashes-symplectic Deligne--Mumford stack $X$
\begin{equation}\label{mapofstacksdarbouxtopervesesheaves}
  \begin{tikzcd}
    \DA \ar{r}{\PJoyce} & \stackPerv_X.
  \end{tikzcd}
\end{equation}
\end{construction}

\begin{remark}
  \cref{MHMisformal} and the uniqueness of the lift to formal schemes offered by the \cref{algebraicinvariants1cat} imply that the vanishing cycles invariants encoded by the map of stacks \cref{mapofstacksdarbouxtopervesesheaves} coincide with the ones of \cite{MR3353002}.
\end{remark}

Next we discuss how the map of stacks of \cref{constructionofperversesheavesasmapfromdarbouxstack} behaves with respect to Thom-Sebastiani.
Recall that $\PJoyce_{U,f} \coloneqq \phi_f(\constantsheaf{A}_U[\dim U])$ for $A = \integers$ or $\rationals$.

\begin{lemma}[Thom--Sebastiani isomorphisms and Knörrer periodicity]\label{thomsebastianiforperversesheaves}
  Let $(\formalU,f)$ be a locally algebraizable formal Landau--Ginzburg pair.
  \begin{enumerate}[label=\textrm{\upshape{(\alph*)}}, ref=\textrm{\upshape{(\alph*)}}]
    \item\label{enumTS} If $(\formalV,g)$ is another locally algebraizable formal Landau-Ginzburg pair, then there is a natural isomorphism
          \[
            \PJoyce_{\formalU \times \formalV, f \smboxplus g} \simeq \PJoyce_{\formalU,f} \boxtimes_{\constantsheaf{A}} \PJoyce_{\formalV,g}
          \]
          of perverse sheaves over $\Crit(f) \times \Crit(g) = \Crit(f \smboxplus g)$.
          It is furthermore associative when taking the product of more than two LG pairs.
    \item\label{enumTSrelquad} Let $(\pi \colon Q \to \formalU, q)$ a quadratic bundle over $\formalU$. There is a natural isomorphism (of perverse sheaves over $\Crit(\formalU, f) = \Crit(Q, f\pi + q)$)
          \[
            \PJoyce_{Q, f \pi + q} \simeq \PJoyce_{\formalU, f} \otimes \PJoyce_{Q, q}[-\dim \formalU]
          \]
          where we write $\PJoyce_{Q,f\pi+q}$ for the perverse sheaf associated to the total space of $Q$.
          This equivalence is furthermore compatible with direct sums of quadratic bundles.
    \item\label{enumknorrer} Let $(Q, q)$ be a quadratic bundle over $\formalU$. The datum of a maximal isotropic sub-bundle $\Lambda \subset Q$ induces a natural isomorphism
          \[
            \PJoyce_{Q,q} \simeq \PJoyce_{\formalU, 0} \simeq \constantsheaf{A}_\formalU[\dim \formalU].
          \]
  \end{enumerate}
\end{lemma}
\begin{proof}
  We start by proving \ref{enumTS}.
  As the formal LG pairs are locally algebraizable, \cite{massey:thomsebastiani} shows there are local Thom--Sebastiani isomorphisms $\tau_\mathrm{loc} \colon \PJoyce_{\formalU \times \formalV, f \smboxplus g} \simeq \PJoyce_{\formalU,f} \boxtimes_A \PJoyce_{\formalV,g}$.
  By \cref{lem-injectivepi0}, any two local algebraization are étale locally isomorphic, and thus the local sections defined above are independent of the choice of algebraization and glue to a global isomorphism
  \[
    \tau \colon \PJoyce_{\formalU \times \formalV, f \smboxplus g} \simeq \PJoyce_{\formalU,f} \boxtimes_{\constantsheaf{A} }\PJoyce_{\formalV,g}.
  \]
  We deduce the associativity from the analogous property of the Thom--Sebastiani isomorphisms associated to regular functions.

 Concerning \ref{enumTSrelquad} and \ref{enumknorrer},
  we may establish independence of the algebraization as for the case \ref{enumTS}.
  This reduces us to the case where $\formalU$ is globally algebraizable, which may be deduced from \cite[Thm 5.4]{MR3353002}.
  To compare our setting with that of \cite{MR3353002}, we observe that the perverse sheaf $\PJoyce_{Q, q}$ corresponds to $\myuline A[\dim \formalU] \otimes_{\quot{\integers}{2}} \det Q^\dual$.
  This can be seen by comparing the canonical bundle of the Milnor fiber with $\det Q^\dual$.

  For the sake of completeness we will now give alternative proofs for the globally algebraic cases of \ref{enumTSrelquad} and \ref{enumknorrer}. Firstly, \ref{enumknorrer} follows immediately from the fact that $\PJoyce_{Q, q}$ corresponds to $\myuline A[\dim \formalU] \otimes_{\quot{\integers}{2}} \det Q^\dual$
  as a maximal isotropic subbundle of $Q$ trivializes $\det Q$.

  The following proof of \ref{enumTSrelquad} was suggested to us by Nick Rozenblyum. Let $(U, f)$ be an LG pair and $(Q, q)$ a quadratic bundle over $U$. The idea is to construct a canonical map $ \PJoyce_{U,f} \otimes \PJoyce_{Q,q}[-\dim U] \to \PJoyce_{Q, f \pi + q}$ and show it is locally an isomorphism. Write $d_U = \dim U$ and $d_Q = \dim Q$.
  We will work with sheaves on $Q$. We note that $\pi^*\PJoyce_{U,f} \cong \pi^*\phi_{f} \constantsheaf A_U[d_U] \cong \PJoyce_{Q,f\pi}[d_U-d_Q]
  $ where the second isomorphism is smooth base change for vanishing cycles which follows from \cite[Expos\'e XIII, Equation 2.1.7.2]{MR0354657}.

  Now we construct the desired map $\Theta \colon \PJoyce_{Q,f\pi} \otimes \PJoyce_{Q,q}[-d_Q] \to \PJoyce_{Q, f \pi + q}$ and show it is locally a quasi-isomorphism. To do this we will apply Thom--Sebastiani to the product of $Q$ with itself. We consider the following commutative diagram:
  \[
    \begin{tikzcd}
      \Crit(f)
      \ar[hook]{r}{i}
      \ar[hook]{dr}[swap]{\iota}
      &
      Z(f \pi)\cap Z(q)
      \ar[hook]{r}{\Delta^{\times}}
      \ar[hook]{d}
      &
      Z(f\pi) \times Z(q)
      \ar[hook]{d}{k} & \Crit(f\pi) \times U \ar[hook']{dl}{j} \ar[hook']{l} \\
      & Z(f\pi + q) \ar[hook]{r}{\Delta^{\boxplus}} \ar[hook]{d}
      & Z(f\pi \boxplus q) \ar[hook]{d}\\
      & Q \ar[hook]{r}{\Delta} & Q \times Q\\
    \end{tikzcd}
  \]
  and we consider $\Crit(f) = \Crit(f\pi) \cap U \subset Q$. We now consider the following composition of maps of constructible sheaves over $\Crit(f\pi) \times \Crit(q) = \Crit(f\pi) \times U$ (because of the shifts the objects on the left are not perverse sheaves):
  \[
    \PJoyce_{Q, f\pi} \boxtimes \PJoyce_{Q, q}[-d_Q]
    \to_{\tau}^\sim
    j^*\PJoyce_{Q \times Q, f\pi \boxplus q}[-d_Q]
    \to_\lambda
    j^*\phi_{f\pi \boxplus q}(\Delta_*\constantsheaf{A}_Q[d_Q])
    \to_\beta^\sim
    j^* \Delta^{\boxplus}_* \PJoyce_{U,f\pi+q}.
  \]
  Here $\tau$ is the Thom--Sebastiani isomorphism {\textrm{\upshape{(a)}}}, $\lambda$ is induced  by diagonal adjunction $(\constantsheaf A_Q[d_Q] \boxtimes \constantsheaf  A_Q[d_Q])[-d_Q] = \constantsheaf A_{Q \times Q}[d_Q] \to \Delta_* \constantsheaf A_Q[d_Q]$ and $\beta$ is proper base change for vanishing cycles, \eg \cite[4.2.11]{Dimca2004}. Then we obtain $\Theta$ by restricting $\beta \circ \lambda \circ \tau$ to the diagonal via $\Delta^*$. First observe that   $\PJoyce_{Q,f\pi} \otimes \PJoyce_{Q,q}[-d_Q]$ and $ \PJoyce_{Q, f \pi + q}$ are supported on $\Crit(f) \subset U$ so we may compare after restricting along $i$ respectively $\iota$.

  Writing $\psi$ for $(\Delta^{\times}\circ i)^*(\beta \circ \lambda \circ \tau)$ we obtain
  \[
    \begin{tikzcd}[column sep=scriptsize]
      i^*(\PJoyce_{Q,f\pi} \otimes \PJoyce_{Q,q}[-d_Q])
      \ar[equals]{d}
      \ar{rr}{\Theta} &  & \iota^* \PJoyce_{Q, f\pi + q} \\
      (\Delta ^{\times} \circ i)^*(\PJoyce_{Q,f\pi} \boxtimes \PJoyce_{Q,q}[-d_Q])
      \ar{r}{\psi}
      &  (\Delta^{\times} \circ i)^* k^* \Delta^{\boxplus}_* \PJoyce_{Q, f\pi+q}\ar{r}{\cong}
      & \iota^* (\Delta^{\boxplus})^* \Delta^{\boxplus}_* \PJoyce_{Q, f\pi+q}
      \ar{u}[sloped]{\sim}
    \end{tikzcd}
  \]
  Here the right vertical arrow is an isomorphism as $\Delta^{\boxplus}$ is a closed immersion.
  It remains to show that $\Theta$ is an isomorphism by checking locally, \ie if $Q = U \times Q_0$ with $U$ contractible and $q = q_0 \circ \pi_Q$ where the projections are $\pi_U$ and $\pi_Q$.

  All the constituent parts of $\Theta$ except for $\lambda$ are isomorphisms. Let us consider its shift
  \[
    \phi_{f \pi \boxplus q}\constantsheaf A_{U \times Q_0 \times U \times Q_0} \to \phi_{f\pi \boxplus q} (\Delta_*\constantsheaf A_{U \times Q_0}).
  \]
  Using the description $\phi_g(\mathcal F) \cong R\Gamma_{\Re(g) \geq 0}(\mathcal F)|_{Z(g)}$ from \cite[Exercise VIII.13]{Kashiwara1990}
  it suffices to compare
  $R\Gamma_{Z} A_{Q\times Q}$
  and $R\Gamma_{Z} \Delta_*\constantsheaf A_Q
    = \Delta_* R\Gamma_{\Delta^{-1}(Z)}(\constantsheaf A_Q)$  where $Z = \{f \pi_U \pi_1 + q_0 \pi_Q \pi_2 \geq 0\}$ and correspondingly $\Delta^{-1}(Z) = \{f \pi_U + q \geq 0\}$.
  However, the pairs $(U \times Q_0 \times U \times Q_0, Z)$ and $(U \times Q_0, \Delta^{-1}(Z))$ are homotopy equivalent as the former is a product of the latter with $U \times Q_0$.
  This completes the proof.
\end{proof}

\begin{observation}\label{perversetorsion}
  As a consequence, the choice of the maximal isotropic sub-bundle $(1,i)$ in $Q \oplus Q$ puts $\PJoyce_{Q, q}[-d_U]$ in the $2$\nobreakdashes-torsion Picard group of complexes of constructible sheaves (for $\otimes$), which identifies with $\mu_2$\nobreakdashes-bundles. Under this identification, the morphism $\PJoyce_{Q,q}$ corresponds to $\det Q[d_U]$.
\end{observation}

\begin{corollary}[of \cref{thomsebastianiforperversesheaves}]\label{corollaryisJoyce}
The morphism
$\begin{tikzcd}[column sep=small,cramped]
    \DA \ar{r}{\PJoyce} & \stackPerv_X
  \end{tikzcd}$
  \cref{mapofstacksdarbouxtopervesesheaves} descends to a morphism
  \[
    \begin{tikzcd}
      \constantsheaf *_X \simeq \displaystyle \quot{\DA}{\QA} \ar{r}{\PJoyce} & \displaystyle \quot{\stackPerv_X}{\B\mu_2}.
    \end{tikzcd}
  \]
  Thus defining a global section of $\quot{\stackPerv_X}{\B\mu_2}$, \ie a twisted perverse sheaf on $X$.
\end{corollary}
\begin{proof}
  \cref{thomsebastianiforperversesheaves}\ref{enumTSrelquad} (together with the observations that vanishing cycles are a formal invariant and formal completion computes with pullbacks)
  shows that $\PJoyce$ is compatible with the action of $\QA$ on $\DA$ from \cref{actiondescendstoisotopyquotient}.
  
  \cref{perversetorsion} then shows that $\QA$ is sent to 2-torsion in perverse sheaves.
\end{proof}

As an immediate consequence of \cref{corollaryisJoyce} we obtain a slight improvement of the main theorem of \cite{MR3353002}:

\begin{theorem}\label{naturaltransformationperversesection}
  Let $X$ be a $(-1)$\nobreakdashes-shifted symplectic derived Deligne--Mumford stack.
  \begin{enumerate}
    \item There is a canonical twisted perverse sheaf $\PJoyce_X \in \quot{\stackPerv_X}{\B\mu_2}(X)$
    \item The obstruction to untwisting $\PJoyce_X$, classified by a class $\beta \in \cohomology^2(X, \mu_2)$, corresponds to the composite morphism
          \[
            \constantsheaf *_X \to \quot{\stackPerv_X}{\B\mu_2} \to \B\B\mu_2 = \Kspace(\mu_2, 2).
          \]
          It corresponds to the $\mu_2$\nobreakdashes-gerbe of square roots of the canonical bundle $\canonicalbundle_X = \det \cotangent_X$.
    \item In particular, we recover \cite[Thm.\,6.9]{MR3353002}: fixing a square root $\canonicalbundle_X^{\quot 1 2}$ defines a lifting
    \[
      \begin{tikzcd}
        & \stackPerv_X \ar{d} \ar{r} \tikzcart  & \constantsheaf *_X \ar{d}
      \\
        *_X \ar[dashed]{ur}{}\ar{r} & \quot{\stackPerv_X}{\B\mu_2} \ar{r} & \B\B\mu_2 \mathrlap{{} = \Kspace(\mu_2, 2).}
      \end{tikzcd}
    \]
    defining a perverse sheaf $\PJoyce_X$ on $X$ such that for any critical chart
    \[
      X \leftarrow S \coloneqq \dCrit(f) \hookrightarrow U \to^f \affineline{},
    \]
    if $\canonicalbundle_U|_S \simeq \canonicalbundle_X^{\quot 1 2} |_S$ as a square root of $\canonicalbundle_X|_S$ then $\PJoyce_X|_S \simeq \PJoyce_{U,f}$.
  \end{enumerate}
\end{theorem}

\begin{remarks}{}
	\item These arguments (together with \cref{MHMisformal}) apply mutatis mutandis to the mixed Hodge modules of vanishing cycles, thus leading to a global mixed Hodge module $\PJoyce^{\MHM}_X \in \MHM(X)$ provided a square root $\canonicalbundle_X^{\quot 1 2}$.
  \item Our proof of \cref{naturaltransformationperversesection} is not entirely independent of \cite{MR3353002}: in order to show that each $\affineline{}$\nobreakdashes-invariant critical invariants induces a formal invariant we use \cite[Prop.\,3.4]{MR3353002} (see the proof of \cref{algebraicinvariants1cat}).
  This is not really necessary to get the result for perverse sheaves of vanishing cycles (see \cref{remvanishingformal})
  The extension to mixed Hodge modules however relies on some results of \cite{MR3353002}.
  \item The reader may be tempted to apply \cref{theoremcontractibility} to glue motivic vanishing cycles (as in \cref{exmotivicvanishingcycles}).
  Such an application is surely possible, but not as straightforward as one would hope.
  Indeed, one would need either to extend the work of Ayoub--Gallauer--Vezzani \cite{JosephAyoub925} so we could define motivic vanishing cycles of \emph{non-adic} functions on formal schemes, or to generalize \cref{algebraicinvariants1cat} for higher categorical $\calF$'s.
\end{remarks}

\subsection{Universal orientation data and Grothendieck--Witt groups}\label{sectionuniversalorientation}

In this section, we define universal orientation data and describe their ties to Grothendieck--Witt groups.
The results are not used in the rest of the article, and the uninterested reader can safely ignore this section.

\begin{definition}[Universal orientation obstruction]\label{defuniversalorientation}
  Let $X$ be a $(-1)$\nobreakdashes-shifted exact symplectic derived Deligne--Mumford stack (with form $\lambda$).
  The \emph{universal orientation obstruction} of $X$ is the morphism
  \[
    \universalobstruction_X \colon \constantsheaf \ast_X \simeq  \quotA \lra \B \left( \QA \right).
  \]
  Universal orientation data are given by a nullhomotopy of the universal orientation obstruction.
\end{definition}

\begin{remarks}{}
\item In particular, if $X$ is equipped with universal orientation data, any invariant out of $\DA$ (satisfying Thom--Sebastiani) glues into a global invariant of $X$.
\item\label{remarkorientationondCrit} Whenever $X$ is presented as a \emph{global} derived critical locus, the given presentation $X\simeq \dCrit(\formalU, f)$ also defines a universal orientation data in the sense of \cref{defuniversalorientation}.
\end{remarks}

The next proposition, in conjunction with \cite[Thm.\,1.2]{SchlichtingTripathi}, tells us universal orientation data live in a semi-topological version of  Grothendieck--Witt spectra.
It would be interesting to have a better understanding of those objects and their ties to the very deep subject of Grothendieck--Witt groups and spectra.

\begin{proposition}
  Let $\orthogonal_{\infty,X_\derham}$ be the group stack $\dAff_{X_\derham}^\op \to \groups(\inftygpd)$ defined by the colimit
  \[
    \orthogonal_{\infty,X_\derham} \coloneqq \left( \colim \left[\orthogonal_1 \to \orthogonal_2 \to \cdots\right] \right)_{|\dAff_{X_\derham}}.
  \]
  Then, there is an equivalence
  \[
    \B \left(\QA\right) \simeq \B \constantsheaf \integers \times \B \B \left(\PhiA{}(\orthogonal_{\infty,X_\derham})\right)_{|{X_\et^\daff}}
  \]
  (with $\PhiA{}$ as in \cref{lemmaappliedtoourcase}).
\end{proposition}
\begin{proof}
  The stack $\B\left( \QA \right)$ is connected, and the question thus reduces to studying its loop stack, which by \cref{inversehyperbolic} is equivalent to
  \[
  \left(\QA\right)^+ \simeq \colim \left[
  \begin{tikzcd}[cramped]
    \QA \ar{r}{- \oplus Q_1} & \QA \ar{r}{- \oplus Q_1} & \cdots
  \end{tikzcd}
  \right].
  \]
  By definition of $\QA$ as a fiber product (\cf \cref{A1contractionofQuad}), and using the notations of \cref{A1localizationofsheavesofgroupoids} we get
  \[
  \begin{tikzcd}
    \left(\QA\right)^+ \ar{r} \ar{d} \tikzcart
    &
    \left( \PhiA{} \left(\colim \left[\B\orthogonal_{X_\derham} \to \cdots\right]\right)\right)_{|{X_\et^\daff}} \ar{d}
    \\
    \mathllap{\constantsheaf \integers \simeq {}} \homotopysheaf_0\left(\colim \left[\stackQuadnabla_X \to \cdots \right] \right)  \ar[hook]{r}
    &
    \homotopysheaf_0\left(\PhiA{} \left(\colim \left[\B\orthogonal_{X_\derham} \to \cdots\right]\right)\right)_{|{X_\et^\daff}} \rlap{.}
  \end{tikzcd}
  \]
  The colimit $\colim \left[\B\orthogonal_{X_\derham} \to \cdots\right]$ is equivalent to $\constantsheaf \integers \times \B \orthogonal_{\infty,X_\derham}$.
  The result then follows from \cref{lemmamappingspacesA1isotopiclocalizationsareA1localizations,suspensedA1localizationpreservesfiniteproducts}.
\end{proof}

\appendix
\numberwithin{equation}{section}

\section{Relative formal Morse lemma} 

\begin{definitions}{Let $\pi \colon \formalV \to \formalU$ be a smooth morphism of smooth formal schemes and let $h \colon \formalV \to \affineline{}$ be a function.}
  \item \label{definitionrelativedcrit} The \emph{relative derived critical locus} of $h$ is the derived intersection
  \[
    \begin{tikzcd}
      \dCrit_\formalU(\formalV, h) \ar{r} \ar{d} \tikzcart & \formalV \ar{d}{dh} \\
      \formalV \ar{r}{0} & \cotangentstack \formalV/\formalU \mathrlap{{} \coloneqq \Spec_{\formalV}\left( \Sym_\formalV \tangent_{\formalV/\formalU} \right).}
    \end{tikzcd}
  \]
  \item \label{definitionrelativelymorse} We say that $h$ is \emph{relatively Morse} if $\dCrit_\formalU(\formalV,h) \to \formalV \to \formalU$ is étale. In particular, if $\pi$ is a reduced equivalence, this is equivalent to $\dCrit_\formalU(\formalV,h) \simeq \formalU$.
\end{definitions}

The following is an adaptation of the argument in  \cite[page 513]{zbMATH06362646}:

\begin{lemma}[Relative formal Morse lemma]\label{relativeformalmorselemma}
  Let $\pi \colon \formalV \to \formalU$ be a smooth morphism of smooth formal schemes and let $h \colon \formalV \to \affineline{}$. Assume that $\pi$ is a reduced equivalence, that $h|{\Crit(h)} = 0$ and that $h$ is relatively Morse.
  Then there exists an isomorphism $\theta \colon \formalV \simeq \widehat \vectorbundle_\formalU(\normal)$ where $\normal$ is the normal bundle of $\formalU \simeq \dCrit_\formalU(\formalV, h) \to \formalV$, such that
  \begin{enumerate}[label={(\alph*)}]
    \item The following diagram commutes
          \[
            \begin{tikzcd}[row sep=0]
              \dCrit_\formalU(\formalV, h) \ar{r} \ar{dd}[sloped]{\sim} & \formalV \ar{dr} \ar{dd}[sloped]{\sim}[swap]{\theta}
              \\
              & & \formalU.
              \\
              \formalU \ar{r}{s_0} & \widehat \vectorbundle_\formalU(\normal) \ar{ur}
            \end{tikzcd}
          \]
    \item The function $q \coloneqq h \circ \theta^{-1} \colon \widehat \vectorbundle_\formalU(\normal) \to \affineline{}$ is a non-degenerate quadratic form.
  \end{enumerate}
\end{lemma}

\begin{proof}
  Using smoothness we first choose an arbitrary isomorphism $\formalV \simeq \widehat \vectorbundle_\formalU(\normal)$ as in Step 3 of \cref{cor:map=action} to reduce the problem to $\formalV = \widehat \vectorbundle_\formalU(\normal)$.
  The map $\normal \to \normal^\dual$ induced by the Hessian of $h$ defines a non-degenerate quadratic form $q$ on $\normal$.

  Write $h = \sum_{k \geq 0} h_k \in \Symcompleted(\normal^\dual) \simeq \structuresheaf_\formalV$, where $h_k \in \Sym^k(\normal^\dual)$ for all $k$.
  By assumption, we have $h_0 = 0$ and, as $h$ is relatively Morse, $h_1 = 0$, and $h_2$ is the non-degenerate quadratic form $q$.

  The key step is to construct a family of adic automorphisms $\tau_k$ of $\Symcompleted(\normal^\dual)$ mapping $h = q + \sum_{k \geq 3} h_k$ to $q + \Symcompleted^{\geq k}(\normal^\dual)$.
  We do so inductively. We set $\tau_3$ to be the identity. Fix $k \geq 3$ and assume $\tau_k$ is constructed.
  Let $f \in \Sym^{k}(\normal^\dual)$ such that $\tau_k(h) \in q + f + \Symcompleted^{\geq k+1}(\normal^\dual)$ and denote by $\lambda$ its image under the inclusion
  \[
    \begin{tikzcd}[row sep=0]
      \Sym^k(\normal^\dual) \hookrightarrow \normal^\dual \otimes \Sym^{k-1}(\normal^\dual) \ar{r}{q^\flat \otimes \Id}[swap]{\sim} & \normal \otimes \Sym^{k-1}(\normal^\dual)\hspace{8em}\\
      & \mathrlap{\hspace{-4em}\simeq \Hom(\normal^\dual, \Sym^{k-1}\normal^\dual).}
    \end{tikzcd}
  \]
  Let $\phi_k$ be the filtered automorphism of $\Symcompleted(\normal^\dual)$ determined by
  \[
    \phi_k = \Id - \lambda \colon \normal^\dual \to \normal^\dual \oplus \Sym^{k-1}(\normal^\dual) \subset \Symcompleted(\normal^\dual).
  \]
  By construction, we have $\phi_k(q + f) \in q + \Symcompleted^{\geq k+1}(\normal^\dual)$.
  To verify this claim we may consider an étale neighborhood where $q = \sum_i y_i^2$ for suitable coordinates, then $\sum_i(y_i-\lambda(y_i))^2+f \in q + \Symcompleted^{\geq k+1}(\normal^\dual)$ as $\sum_i 2 y_i \lambda(y_i) = f$ by definition of $\lambda$.
  We can thus set $\tau_{k+1} = \phi_k \circ \tau_k$.

  Now each $\tau_k$ induces an automorphism $\bar \tau_k$ of $\Sym(\normal^\dual)/\Sym^{\geq k-1}(\normal^\dual)$ such that $\bar \tau_k(h) = q$. Since $\phi_k$ induces the identity on $\Sym(\normal^\dual)/\Sym^{\geq k-1}(\normal^\dual)$, this system of automorphisms is compatible with the projections, and thus induces an automorphism
  \[
    \tau \colon
    \begin{tikzcd}[cramped,column sep=large]
      \displaystyle \lim_k \quot{\Sym(\normal^\dual)}{\Sym^{\geq k}(\normal^\dual)} \ar{r}{\lim \bar \tau_{k+1}} & \displaystyle \lim_k \quot{\Sym(\normal^\dual)}{\Sym^{\geq k}(\normal^\dual)}.
    \end{tikzcd}
  \]
  We have by construction $\tau(h) = q$. Setting $\theta = \tau^{-1}$ completes the proof.
\end{proof}

\printbibliography

\end{document}